\documentclass[a4paper, 11pt,reqno]{amsart}
\usepackage{amsmath}
\usepackage{amssymb}
\usepackage{amsfonts}
\usepackage{amscd}
\usepackage{amsthm}
\usepackage{mathtools}
\usepackage{mathrsfs}
\usepackage{multirow,bigdelim}
\usepackage{MnSymbol}
\usepackage[bb=boondox,bbscaled=.95,cal=boondoxo]{mathalfa}
\usepackage[numbers]{natbib}
\usepackage{xcolor}
\usepackage{url}
\usepackage{hyperref}
\usepackage{tikz}
\usetikzlibrary{cd}
\usetikzlibrary{decorations.pathmorphing}
\usetikzlibrary{patterns}
\newcommand{\bfA}{{\rm\bf A}}

\newcommand{\CC}{{\rm\bf C}}
\newcommand{\RR}{{\rm\bf R}}
\newcommand{\QQ}{{\rm\bf Q}}
\newcommand{\ZZ}{{\rm\bf Z}}
\newcommand{\NN}{{\rm\bf N}}
\newcommand{\GG}{{\rm\bf G}}
\newcommand{\FF}{{\rm\bf F}}

\newcommand{\PP}{{\rm\bf P}}

\newcommand{\OO}{\mathcal{O}}

\DeclareMathOperator{\Spec}{\mathrm{Spec}}
\DeclareMathOperator{\diag}{\mathrm{diag}}

\DeclareMathOperator{\GL}{\mathrm{GL}}
\DeclareMathOperator{\SL}{\mathrm{SL}}

\DeclareMathOperator{\Oo}{\mathrm{O}}
\DeclareMathOperator{\U}{\mathrm{U}}
\DeclareMathOperator{\SO}{\mathrm{SO}}

\DeclareMathOperator{\Sp}{\mathrm{Sp}}

\DeclareMathOperator{\Mod}{Mod}

\DeclareMathOperator{\Sch}{Sch}

\DeclareMathOperator{\End}{\mathrm{End}}
\DeclareMathOperator{\Hom}{\mathrm{Hom}}
\DeclareMathOperator{\iEnd}{\mathcal E\textit{nd}\,}
\DeclareMathOperator{\iHom}{\mathcal H\!\textit{om}\,}
\DeclareMathOperator{\Gr}{\mathrm{Gr}}

\DeclareMathOperator{\Der}{{\mathcal D\!\textit{e}\textit{r}}}
\DeclareMathOperator{\D}{{\mathcal D}}
\DeclareMathOperator{\Pp}{{\mathcal P}}

\DeclareMathOperator{\Image}{\mathrm{Im}}
\DeclareMathOperator{\op}{\mathrm{op}}
\DeclareMathOperator{\rank}{\mathrm{rank}}

\DeclareMathOperator{\Gal}{\mathrm {Gal}}
\DeclareMathOperator{\ind}{\mathrm{ind}}
\DeclareMathOperator{\res}{\mathrm{Res}}

\DeclareMathOperator{\act}{\mathrm{act}}
\DeclareMathOperator{\id}{\mathrm{id}}
\DeclareMathOperator{\Lie}{\mathrm{Lie}}

\DeclareMathOperator{\Frac}{\mathrm{Frac}}

\DeclareMathOperator{\pr}{\mathrm{pr}}

\DeclareMathOperator{\rtype}{\mathrm{strtype}}
\DeclareMathOperator{\type}{\mathrm{type}}

\DeclareMathOperator{\SPS}{\mathrm{SPS}}

\DeclareMathOperator{\Proj}{\mathrm{Proj}}
\DeclareMathOperator{\spl}{\mathrm{spl}}

\newcommand{\liea}{{\mathfrak {a}}}

\newcommand{\liec}{{\mathfrak {c}}}

\newcommand{\lieg}{{\mathfrak {g}}}
\newcommand{\lieh}{{\mathfrak {h}}}
\newcommand{\liek}{{\mathfrak {k}}}
\newcommand{\liel}{{\mathfrak {l}}}

\newcommand{\lies}{{\mathfrak {s}}}
\newcommand{\liet}{{\mathfrak {t}}}

\newcommand{\liep}{{\mathfrak {p}}}
\newcommand{\lieq}{{\mathfrak {q}}}

\newcommand{\lieu}{{\mathfrak {u}}}

\newcommand{\lieso}{{\mathfrak {so}}}

\newcommand{\Aa}{{\mathcal {A}}}
\newcommand{\Bb}{{\mathcal {B}}}
\newcommand{\Cc}{{\mathcal {C}}}
\newcommand{\Ee}{{\mathcal {E}}}
\newcommand{\Ff}{{\mathcal {F}}}
\newcommand{\Gg}{{\mathcal {G}}}
\newcommand{\Ii}{\mathcal{I}}
\newcommand{\Ll}{{\mathcal {L}}}
\newcommand{\Mm}{{\mathcal {M}}}
\newcommand{\Nn}{{\mathcal {N}}}
\newcommand{\Rr}{{\mathcal {R}}}
\newcommand{\Vv}{{\mathcal {V}}}
\DeclareMathOperator{\Bun}{\mathrm{Bun}}
\DeclareMathOperator{\can}{\mathrm{can}}

\DeclareMathOperator{\cl}{\mathrm{cl}}
\DeclareMathOperator{\DR}{DR}

\DeclareMathOperator{\ev}{\mathrm{ev}}
\DeclareMathOperator{\Ker}{\mathrm{Ker}}
\DeclareMathOperator{\Map}{Map}
\DeclareMathOperator{\cmod}{-\mathrm{mod}}
\DeclareMathOperator{\PA}{PA}

\DeclareMathOperator{\qc}{qc}

\DeclareMathOperator{\ord}{\mathrm{ord}}
\DeclareMathOperator{\Sym}{Sym}
\DeclareMathOperator{\TDO}{TDO}

\DeclareMathOperator{\Tor}{Tor}
\DeclareMathOperator{\Ss}{{\mathscr S}}

\newcommand{\lieS}{{\mathfrak {S}}}

\newcommand{\xcong}[1]{%
	\mathrel{\tikz[baseline=0pt] {
			\node[above] at (0,1.2ex) (a) {\(\scriptstyle #1\)};
			\draw[preaction={
				transform canvas={yshift=-.5ex},
				draw,
				decorate,
				decoration={lineto}},
			preaction={
				transform canvas={yshift=-1ex},
				draw,
				decorate,
				decoration={lineto}}]
			(a.south west) .. controls +(.25,.15) and +(-.25,-.15) .. (a.south east);
}}}

\theoremstyle{plain}
\newtheorem{theorem}{Theorem}[subsection]
\newtheorem{lemma}[theorem]{Lemma}
\newtheorem{corollary}[theorem]{Corollary}
\newtheorem{proposition}[theorem]{Proposition}

\newtheorem{condition}[theorem]{Condition}
\newtheorem{problem}[theorem]{Problem}
\newtheorem{claim}[theorem]{Claim}
\newtheorem{variant}[theorem]{Variant}
\newtheorem{strategy}[theorem]{Strategy}
\newtheorem{prop-defn}[theorem]{Proposition-Definition}
\newtheorem{definition}[theorem]{Definition}
\newtheorem{construction}[theorem]{Construction}
\theoremstyle{remark}
\newtheorem{example}[theorem]{Example}
\newtheorem{remark}[theorem]{Remark}

\bibpunct{[}{]}{,}{n}{}{,}

\numberwithin{equation}{subsection}

\begin{document}
	
	\title{Families of twisted $\D$-modules\\and\\arithmetic models of Harish-Chandra modules}
	
	\author{Takuma Hayashi, Fabian Januszewski}
	\address{Department of Pure and Applied Mathematics, Graduate School of Information Science and Technology, Osaka University, 1-5 Yamadaoka, Suita, Osaka 565-0871, Japan}
	\email{hayashi-t@ist.osaka-u.ac.jp}
	\address{Institut f\"ur Mathematik, Fakult\"at EIM, Paderborn University, Warburger Str.\ 100, 33098 Paderborn, Germany}
	\email{fabian.januszewski@math.uni-paderborn.de}
	
	
	\begin{abstract}
		We develop a theory of tdos and twisted $\mathcal D$-modules over general base schemes with a focus on functorial asp
		ects. In particular, we introduce a flat base change functor and establish its compatibility with globalization and direct image functors. We also study forms of closed $K$-orbits of $\theta$-stable parabolic subgroups in the total flag variety. We apply these two developed theories to give a geometric construction of half-integral models of cohomologically induced modules. With a view towards arithmetic applications, we further demonstrate desirable properties of the constructed half-integral models, such as projectivity over the base and torsion-free relative Lie algebra cohomology.
	\end{abstract}
	
	\maketitle
	
	{\small
		\tableofcontents
	}
	
	\section*{Introduction}
	
	\subsection{Arithmetic background}
	
	Langlands' representation-theoretic perspective on the theory of automorphic forms led to fundamental conjectures relating automorphic representations, $L$-functions and motives.
	
	While automorphic representations in the Langlands Programme naturally live over $\CC$, rational and integral structures arising from the cohomology of arithmetic groups and coherent cohomology of Shimura varieties play a central role in arithmetic applications of automorphic forms.
	
	For example, the notion of a $p$-adic automorphic form is an inherently integral concept and is intimately linked to integral structures on the cohomology of arithmetic groups and to integral models of Shimura varieties, as well as to $p$-adic Galois representations.
	
	The study of the integrality properties of elliptic modular forms via their Fourier expansion at $\infty$ goes back to Shimura in the 1950ies \cite[Th\'eor\`eme 3]{shimura1959} and found its analog in the study of integral structures on the non-archimedean part $\Pi_f$ of an \lq{}algebraic\rq{} automorphic representation $\Pi=\Pi_\infty\otimes\Pi_f$.
	
	\smallskip
	In the case of automorphic representations $\Pi$ of $\GL_2/\QQ$ associated to an elliptic cusp form, not much is lost by focusing on $\Pi_f$, since the minimal $K$-type in $\Pi_\infty$ decomposes canonically into two lines, giving rise to two essentially canonical embeddings $\iota_\pm\colon\Pi_f\to\Pi$. Moreover, the $L$-functions associated to elliptic cusp forms are represented by good test vectors supported in these images of $\iota_\pm$, as is the contribution of $\Pi_f$ to the (cuspidal) cohomology of the arithmetic groups $\Gamma\subseteq\GL_2(\ZZ)$.
	
	Such properties are implicitly exploited for example in \cite{eichler1957,shimura1959,manin1972,zagier1977}. The number of results implicitly or explicitly relying on the above-mentioned properties is too large to be listed in its entirety.
	
	\smallskip
	For more general connected reductive groups $G/\QQ$ the analogous statements are no more true: Minimal $K$-types are higher-dimensional, in general we do not know whether we find good test vectors supported in the minimal $K$-types for a given integral representation of an associated automorphic $L$-function, and cuspidal cohomological representations $\Pi$ of $G$ contribute to cohomology in the so-called \lq{}cuspidal range\rq{}, i.\,e.\ no more in a single cohomological degree (cf.\ \cite{clozel1990} for a discussion of the case $G=\res_{F/\QQ}\GL_n$). A conjecture of Venkatesh addresses the issue of higher multiplicities of $\Pi_f$ in various cohomological degrees \cite{venkatesh2017,venkatesh2018ICM,venkatesh2019}.
	
	From a conceptual point of view, restricting attention to rational and integral structures on $\Pi_f$ is not satisfactory either: Any embedding $\Pi_f\to\Pi$ depends on the choice of an element in $\Pi_\infty$, and a fortiori a rational or integral structure on $\Pi_f$ does not globalize to $\Pi$ in a canonical way.
	
	\smallskip
	To address these complications, a better understanding of $\Pi_\infty$ from an arithmetic point of view is required. This motivates the study of rational and integral structures on the Harish-Chandra module $\Pi_\infty^{(K)}$ of $K$-finite vectors in $\Pi_\infty$, which in our present applications is always a cohomologically induced module (cf.\ Section \ref{sec:kitchen} below). For example, suitable arithmetic forms, i.\,e.\ integral models of $\Pi_\infty^{(K)}$ naturally allow us to globalize integral structures on $\Pi_f$. The resulting global integral structures will be unique up to scalars in $\CC^\times$.
	
	Therefore different normalizations of such global integral structures will give rise to automorphically defined periods, which in certain situations can be expected to be related to motivic periods as defined by Deligne in \cite{deligne1979} for example.
	
	Together with an appropriate well-behaved integral notion of $(\lieg,K)$-cohomology, a global integral structure on $\Pi$ will give rise to an integral structure on the $(\lieg,K)$-cohomology of $\Pi$. This in turn then induces an integral structure on a subspace of the isotypic component of $\Pi_f$ inside the cohomology of arithmetic groups in every single degree (cf.\ \cite{januszewskisheaves}).
	
	In that sense, global integral structures allow to transcend individual cohomological degrees and are therefore expected to be related to Venkatesh's Conjecture \cite{venkatesh2017,venkatesh2018ICM,venkatesh2019} on the action of the derived Hecke algebra on such isotypic components (cf.\ \cite{januszewskirationality,januszewskisheaves}).
	
	\smallskip
	Before studying such rational and integral structures, we need to establish their existence. Rational structures have been constructed previously by Harder, Harris and the second author \cite{harderraghuram,harris2013,harris2013erratum,januszewskirationality}, each using different methods: Harder used a $\ZZ[1/2]$-integral analog of real parabolic induction for $\GL_n/\QQ$, Harris sketched an approach via a $\QQ$-rational Beilinson--Bernstein theory, and Januszewski used a $\QQ$-rational analog of the Zuckerman functor. Together with the refinements laid out in \cite{januszewskisheaves}, we now have a satisfactory theory of rational models for $(\lieg,K)$-modules.
	
	\smallskip
	The construction of integral structures on infinite-dimensional modules is a non-trivial and more delicate matter. While rational structures on absolutely irreducible modules are unique up to homotheties, this is no more true for integral models. For example, any rational model is also an integral model, and torsion may be present as well. Therefore a mere construction of an integral model is only interesting if it can be shown that the constructed integral model satisfies additional properties, such as being projective over the base for example.
	
	Projectivity over the base excludes the rational models as integral models and also implies the absence of torsion. Additionally, we want integral $(\lieg,K)$-cohomology of integral models to be non-trivial and projective as well.
	
	\smallskip
	Our current understanding suggests that the integral analogs of cohomological induction are not suitable for the construction of integral models, since it appears difficult to establish the desired properties. Therefore we give another construction of integral models here, based on twisted $\D$-modules, which will allow us to construct suitable integral models and prove the above-mentioned desired properties.
	
	\smallskip
	In order to execute this construction, we first need to develop a suitable theory of twisted $\D$-modules which serves our needs.
	
	\smallskip
	Then in order to apply such a theory of twisted $\D$-modules to the construction of integral models of cohomologically induced $(\lieg,K)$-modules, we need to identify appropriate generalizations of partial flag varieties and $K$-orbits therein: The natural solution is given by moduli spaces of certain parabolic subgroups. To get integral models of characters of Levi subgroups or of the coefficient representations, we also need line bundles on the partial flag schemes (integral models of partial flag varieties), which was studied in \cite{hayashilinebdl}. 
	
	\smallskip
	With all these ingredients at hand, we can then construct integral models as outlined as global sections of twisted $\D$-modules over integer rings and prove that they are projective over the base. As mentioned above, we also show that their relative Lie algebra cohomology is well behaved. For the latter, it is crucial that we can control the corresponding integral structures on minimal $K$-types via our theory of twisted $\mathcal D$-modules over schemes.
	
	\smallskip
	We will discuss our approach to twisted $\D$-modules in the next section and will return the study of orbits in partial flag schemes and to current and potential future applications in number theory and representation theory before outlining the structure of the paper at the end of the introduction.

	\subsection{Twisted $\D$-modules}
	
	\subsubsection{Background}
	Classical $\D$-modules play a central role in representation theory of complex semisimple Lie algebras and real reductive Lie groups in the framework of of Beilinson--Bernstein's localization theory (\cite{beilinsonbernstein1981,beilinsonbernstein1993}). The untwisted $\D$-modules, i.\,e.\ modules over the classical sheaf $\D$ of differential operators, give rise to modules with trivial infinitesimal character.
	
	\smallskip
	To approach a broader class of representations, i.\,e.\ to allow for more general infinitesmal characters, Beilinson--Bernstein introduced twisted differential operators: Tdos (sheaves of twisted differential operators) on smooth complex algebraic varieties $X$ are sheaves of rings with homomorphisms from the structure sheaves $\OO_X$ of $X$ which are \'etale locally isomorphic to $\D_X$ (\cite[section 1]{beilinsonbernstein1993}).
	Subsequently, Kashiwara and Beilinson--Bernstein gave an axiomatic formalism of tdos on smooth algebraic varieties over fields $F$ of characteristic zero without reference to $\D_X$ and developed their general theory in \cite{kashiwara1989,beilinsonbernstein1993}. 
	For their classification and functoriality properties, the notions of Picard algebroids and $\Mm^\bullet$-torsors are introduced in \cite{beilinsonbernstein1993}, where $\Mm^\bullet$ is a cochain complex concentrated in positive degrees. In fact, the groupoid of tdos on $X$ is equivalent to that of $\Omega^{\geq 1}_{X/F}$-torsors, where $\Omega^{\geq 1}_{X/F}$ is the stupidly truncated de Rham complex of $X$ over $F$ at degree 0 (\cite[Lemma 2.1.4]{beilinsonbernstein1993}).
	
	\smallskip
	Then we can simply define a twisted $\D$-module as a module over a tdo. A remarkable feature of the theory of twisted $\D$-modules are functoriality properties of their derived categories (the derived direct image functors $f_+$ and the derived pullback $Lf^\ast$). The basic references are \cite{boreletal, milicic, hottaetal2008, hechtetal, kashiwaratanisaki}. They are key tools to construct rational forms of Harish-Chandra modules. In fact, motivated by the study of periods, Harris proposed to construct and classify irreducible equivariant twisted $\D$-modules on the flag variety (\cite[Theorem 2.1.3 (3)]{harris2013}). In particular, he proposed to construct rational forms of discrete (fundamental) series representations by using the direct image functor from closed orbits (\cite[Theorem 2.1.3 (4)]{harris2013}). Since then, these were amended by \cite{harris2013erratum}, which concerns rationality of objects appearing in the construction (see also \cite{hayashilinebdl} and section \ref{sec:orbit}).
	
	\smallskip
	However, a proper discussion of the base change formalism was lacking in the literature. We need this in order to verify that Harris' Harish-Chandra modules are rational forms of the expectedly-corresponding Harish-Chandra modules over $\CC$.
	
	\smallskip
	As explained in the previous section, our aim in this paper is to construct arithmetic forms of cohomologically induced modules by developing a general theory of twisted $\D$-modules over schemes as a generalization of the classical theory over fields $F$ of characterisic zero, where from the scheme-theoretic point of view, it can be thought of as a theory over $\Spec F$. At the same time, we also aim to introduce the base change formalism of tdos and twisted $\D$-modules.
	
	\subsubsection{Towards a general theory of twisted $\D$-modules}
	Our first task is to introduce a suitable notion of tdos for general smooth morphisms $X\to S$. In this paper, we define tdos over $S$ as a generalization of those over fields $F$ of characteristic zero by simply replacing the base $F$ with $S$ in \cite[Definition 2.3.3]{kashiwara1989} and \cite[Definition 2.1.1]{beilinsonbernstein1993}. The classical examples of tdos are Grothendieck's sheaf of differential operators on the structure sheaf $\OO_X$ if every nonzero integer is invertible on $S$ and the sheaf of PD differential operators (\cite[\S 4]{berthelotogus}) if there exists a positive integer $n$ such that $n$ vanishes on $S$. We adopt it as the definition of tdos on smooth schemes $X$ over general base schemes $S$ in this paper for compatibility with the formalism of $(\lieg,K)$-modules over commutative rings in \cite{harderraghuram,hayashi2019,hayashi2018}. We remark that Grothendieck's sheaf of differential operators is not a tdo in our sense in general since it may not be generated by the structure and tangent sheaves.

	Similarly, we define Picard algebroids and $\Omega^{\geq 1}_{X/S}$-torsors by replacement of $F$ with $S$ in \cite[Definition 2.1.3, 2.1.5]{beilinsonbernstein1993} (cf.~\cite[(16.5.15)]{ega44}). Then the equivalence of the groupoids of tdos, Picard algebroids, and $\Omega^{\geq 1}_{X/S}$-torsors still holds. 
	
	We define the base change functor in terms of $\Omega^{\geq 1}_{X/S}$-torsors, and transfer it to Picard algebroids and tdos. This is meaningful because of its local description: If $X$ and $S$ are affine, we can see the equivalence of tdos on a smooth scheme and those for its coordinate ring which are defined in the same manner as in the sheaf setting, which we refer to as the ring setting. Then the base change functor along a morphism of affine schemes $S'=\Spec k'\to \Spec k=S$ is given by $k'\otimes_k(-)$ in terms of the ring setting.

	Thanks to the base change functor, we can think of tdos as follows: For instance, we put $S=\Spec\ZZ$. This is one-dimensional with a point for every prime ideal $(p)\subseteq\ZZ$. If $p=0$ (resp. $p>0$), the corresponding point is $\Spec\QQ$ (resp.~$\Spec\ZZ/p\ZZ$) as a scheme. Hence one can view $\Aa$ as a {\em compatible family} $\{\Aa_{(p)}\}$ of tdos over $\Spec\ZZ/p\ZZ$ ($p>0$) extending a given module over the generic fiber $\Spec\QQ$ ($p=0$) in a suitable algebro-geometric (i.\,e.\ sheaf-theoretic) sense. For instance, if $\Aa$ is the sheaf of `untwisted' differential operators, so are the fibers.
	
	Twisted $\D$-modules are just modules over tdos. Using the base change of tdos, we can define that of twisted $\D$-modules in the standard manner. Hence one can imagine twisted $\D$-modules as families of twisted $\D$-modules over fields.

	\subsubsection{Tdos and twisted $\D$-modules over schemes}
	
	In the previous section, we proposed the definitions of tdos and twisted $\D$-modules over schemes $S$, and explained how to define the base change functors.
	
	\smallskip
	Our next step is to show that the fundamental aspects of the theory of tdos and twisted $\D$-modules in the literature generalize to general base schemes by using unbounded derived categories. On first sight, this appears to be straightforward at the formal level.

	\smallskip
	At the constructive level, the description of the pullback of tdos in \cite[Lemma 2.2]{beilinsonbernstein1993} is, however, not available (think over fields of positive characteristics). This causes technical difficulties in the definition of the
	left action of the pullback $f^\cdot\Aa$ on the transfer bimodule $\Aa_{X\to Y}$ for a morphism $X\to Y$ of smooth $S$-schemes and a tdo $\Aa$ on $Y$. We circumvent this issue by considering the actions of Picard algebroids.
	
	\smallskip
	Another technical issue arises in the context of the functoriality of the derived direct image functors in $X$ because of the (potential) failure of the projection formula (cf.~\cite[Proof of Proposition 1.5.21]{hottaetal2008}). We prove the required version of the projection formula under a mild condition.
	
	\smallskip
	Similarly, we need a finiteness condition for a generalization of the classical fact that the direct image functor along an affine immersion of smooth varieties over fields of characteristic zero is $t$-exact (\cite[section A.3.3]{hechtetal}, \cite[Chapter IV, Corollary 6.5]{milicic}). That is, $H^p (i_+\Mm)$ vanishes for an affine immersion $i:Y\to X$ of smooth schemes over a base $S$, a tdo $\Aa$ on $X$, a left $i^\cdot\Aa$-module $\Mm$, and $p\neq 0$. This fact is technical but important in representation theory since this makes analysis of certain representations constructed by using the direct image functor $i_+$ easy (see \cite[Proposition 2.6]{hechtetal}). This is also used for the classification of simple holonomic $\D$-modules (see \cite[Theorem 3.4.2]{hechtetal}). This vanishing theorem is nice when we work with equivariant twisted $\D$-modules since it turns out to lift $i_+$ to a functor between the abelian categories of quasi-coherent equivariant twisted $\D$-modules; Otherwise, we do not have a lift since the abelian category of equivariant twisted $\D$-modules does not have good homological properties (cf.~section \ref{sec:future}). We prove this for a general base $S$ under a reasonable condition.
	
	\smallskip
	In this paper, we prove extensions for the most fundamental classical statements under suitable hypotheses. However, some of them fail for reasons other than finiteness. For example, we give counterexamples to Kashiwara's theorem for closed immersions over $\ZZ$ and $\ZZ/p\ZZ$ for a prime $p$.
	
	\smallskip
	An expert of $\D$-modules over $\CC$ will certainly miss certain statements and probably would have put another emphasis. For example, we have the six operations if we restrict ourselves to holonomic twisted $\D$-modules over $\CC$. We do not discuss holonomic modules in the present paper since they are beyond the scope of our applications.

	\smallskip
	We also study the base change functors. In particular, we prove the flat base change theorems for the direct image and the globalization functors. Here the globalization functor for $x:X\to S$ is the usual direct image functor $x_\ast$ (or its right derived functor) to the base scheme $S$, which generalizes the global section functor in the classical theory.
	
	\smallskip
	As a quick conclusion, we achieve these goals in the present work and give a properly \lq{}geometric\rq{} theory of twisted $\D$-modules in Grothendieck's sense, i.e., a well-behaved (flat) base change formalism. Since we are guided by the fundamental principles of Grothendieck's theory of schemes, we are optimistic that the resulting generality will have applications beyond the construction of integral models of Harish-Chandra modules.

	\subsubsection{Twisted $\D$-modules over Dedekind schemes}
	
	Let us explain the subsequent development because we make use of it in applications to cohomologically induced modules later: It is a fundamental question that the theory of twisted $\D$-modules over schemes proposed above, in particular, the direct image functor really produces new objects over more general bases than $\Spec F$ above.
	
	After we established the generalities of twisted $\D$-modules in the present paper, the first author proved:
	
	\begin{theorem}[{\cite[Theorem 1.1, the end of section 1]{hayashifil}}]\label{thm:fil}
		Let $K$ be a smooth affine group scheme over a Dedekind scheme $S$, $i:Y\hookrightarrow X$ be a $K$-equivariant closed immersion of smooth $K$-schemes over $S$, and $\Aa$ be a $K$-equivariant tdo on $X$. Write $x:X\to S$ for the structure morphism. Suppose that the following conditions are satisfied:
		\begin{enumerate}
			\renewcommand{\labelenumi}{(\roman{enumi})}
			\item $X$ is quasi-compact over $S$;
			\item $Y$ is proper over $S$.
		\end{enumerate}
		Let $\Mm$ be a $K$-equivariant left $i^\cdot\Aa$-module which is locally free of finite rank over $\OO_Y$. Then there exists a natural exhaustive $K$-invariant filtration $F_\bullet i_+\Mm$ on the direct image $i_+\Mm$ such that the $p$th associated graded sheaf
		$\Gr^p x_\ast i_+\Mm$
		to the induced $K$-invariant filtration $x_\ast F_\bullet i_+\Mm$ on $x_\ast i_+\Mm$ is a locally free of finite rank over $\OO_S$ for every integer $p$.
	\end{theorem}
	
	When $S=\Spec F$ for a field $F$ of characteristic zero, the filtration $F_\bullet i_+\Mm$ is given by the so-called normal degree (cf.~\cite[Appendix A.3.3]{hechtetal}, \cite[Chapter II, Proposition 7.7]{bien1990}, \cite[Section 3]{oshima2015}). For right modules, we generalize it to an arbitrary base $S$ by using the given filtration on $\Aa$ and the left action of $i^\cdot\Aa$ on $\Aa_{Y\to X}$. A similar construction works for left modules as well.
	
	\smallskip
	In our main application of the present paper, $S$ will be affine with coordinate ring, say, $k$ a Dedekind domain. We identify $x_\ast i_+\Mm$ with the $k$-module $\Gamma(X,i_+\Mm)$ of global sections to use the induced filtration on $\Gamma(X,i_+\Mm)$. We find in virtue of Theorem \ref{thm:fil} that $\Gamma(X,i_+\Mm)$ is projective as a $k$-module (\cite[Lemma 5.4, Corollary 1.3]{hayashifil}).

	\subsection{Cohomologically induced modules and their geometric construction}
	
	According to the work of Vogan--Zuckerman, irreducible unitary representations of a connected real semisimple Lie group of finite center with nonzero relative Lie algebra cohomology are realized as $A_{\lieq}(\lambda)$, i.\.e.\ the cohomological induction of a character (\cite[Theorem 5.6]{voganzuckerman1984}). They also computed the relative Lie algebra cohomology of $A_{\lieq}(\lambda)$ (\cite[Theorem 5.5]{voganzuckerman1984}). 
	
	\smallskip
	For applications to the relative Lie algebra cohomology, it is worthy of record that only the minimal $K$-type contributes to this cohomology. Similarly, it is folklore that irreducible essentially unitarizable representations of certain real reductive Lie groups with nonzero cohomology are realized as cohomologically induced modules (cf.~Theorem \ref{thm:cohomological->aq}). Their cohomology has a similar description (see Theorem \ref{thm:a_q->cohomological}). Facts on their irreducibility and unitarizability are collected in \cite{knappvogan}.

	\smallskip
	For arithmetic applications, we wish to have arithmetic forms of cohomologically induced modules. There are three general ways known to construct cohomologically induced modules over the complex numbers. 
	
	One is Langlands' construction using real parabolic induction (\cite[(11.200)]{knappvogan}). The second way is to use the cohomological induction after Zuckerman. One can construct their rational forms by using the cohomological induction functor over fields of characteristic zero introduced by the second author in \cite{januszewskirationality} (see \cite{hayashi2018} for the base change property over rings). We remark that we cannot work directly within $\RR$ since nontrivial parabolic subalgebras $\lieq$ appearing in this construction are not defined over $\RR$.
	
	Nonetheless the second author pointed out that some cohomologically induced modules are defined over the real numbers and in certain situations even over the rational numbers. Discrete series representations of $\GL_2(\RR)$ are typical examples. The second author produced smaller forms by the Galois descent. One can also construct integral forms by using the cohomological induction functor over commutative rings introduced by the first author in \cite{hayashi2019}. They are hard to analyze since the derived construction produces huge torsion.

	\begin{remark}[{\cite[Proposition 6.3]{voganzuckerman1984}}, {\cite[Theorem 11.225]{knappvogan}}]
		The equivalence of these two constructions over the complex numbers is known as the transfer theorem.
	\end{remark}
	
	\begin{remark}
		Harder used parabolic induction to construct half-integral models in the case of $G=\GL_n$ over $\QQ$ \cite{harderraghuram}. Despite the transfer theorem over $\CC$, Harder's models are different from the models constructed by the first author in \cite{hayashi2019}.
	\end{remark}
	
	The third way is to use twisted $\D$-modules which we execute in this paper. To simplify notation, we write $-\otimes_\RR\CC=(-)_\CC$. The corresponding classical construction over $\CC$ is summarized in
	
	\begin{theorem}[{\cite[Theorem 4.3, Appendix B]{hechtetal}}, {\cite[Theorem 5.4 and Corollary 5.5]{kitchen2012}}, {\cite[(1.1)]{oshima2013}}, see also Corollary \ref{cor:duality}]\label{thm:duality}
		Let $G$ be a connected real reductive algebraic group with a Cartan involution $\theta$. Let $\lieg$ be the Lie algebra of $G$. Let $K\subset G$ be the fixed point subgroup by $\theta$. Let $H\subset G$ be a fundamental Cartan subgroup. Let $Q'$ be a $\theta_\CC$-stable parabolic subgroup containing $H_\CC$. Let $u$ be the dimension of the unipotent radical of $\bar{Q}'\cap K_\CC$.
		
		Let $i:K_\CC/(\bar{Q}'\cap K_\CC)\hookrightarrow G_\CC/\bar{Q}'$ denote the $K_\CC$-orbit map at the base point. Let $\Aa'$ be a $G_{\CC}$-equivariant tdo on $G_\CC/\bar{Q}'$, and $\Mm'$ be a $K_{\CC}$-equivariant integrable left $i^\cdot \Aa'$-connection. Let $M'$ be the fiber of $\Mm'$ at the base point. 
		\begin{enumerate}
			\item There exists an isomorphism
			$H^\bullet(G_\CC/\bar{Q}',i_+\Mm')
			\cong\Ll_{u-\bullet}(M')$
			of $(\lieg_\CC,K_{\CC})$-modules, where $\Ll_\bullet$ is the (co)homological induction functor in \cite{knappvogan}.
			\item Assume that $\Mm'$ is a line bundle. Consider the filtration $F_\bullet i_+\Mm'$ by the normal degree (\cite[Chapter II, Proposition 7.7]{bien1990} or \cite[section 3]{hayashifil}, see also \cite[Section 3]{oshima2015} for the affine setting). If $M'$ satisfies the conditions of \cite[Proposition 10.24]{knappvogan}, the global section module $\Gamma(K_\CC/(\bar{Q}'\cap K_\CC), F_0i_+\Mm')$ is the subspace of minimal $K_\CC$-types of $\Ll_{u}(M')$.
		\end{enumerate}
	\end{theorem}
	
	For the relation with the right derived version $\mathcal{R}^\bullet$ we refer to \cite[Theorems 8.21]{knappvogan}. In virtue of (2) and Vogan--Zuckerman's proof of \cite[Theorem 5.5]{voganzuckerman1984}, we can describe the relative Lie algebra cohomology in geometric terms of the realization (1). That is, assume the condition of (2). Let $V$ be a $(\lieg_\CC,K_\CC)$-module which is irreducible as a $\lieg_\CC$-module. Under suitable conditions on $M'$ and $V$ (see Theorem \ref{thm:a_q->cohomological} for the full statement), we have
	\[H^\bullet(\lieg_\CC,K_\CC,\Ll_{u}(M')\otimes_\CC V)
	\cong \Hom_{K_\CC}(\wedge^\bullet \lieg_\CC/\liek_\CC,\Gamma(K_\CC/(\bar{Q}'\cap K_\CC), F_0i_+\Mm')\otimes_\CC V),
	\]
	where $\liek$ is the Lie algebra of $K$. Once we obtain integral forms of these objects, we can expect a similar description of their relative Lie algebra cohomology groups. We could simplify this over $\CC$ as in \cite[Theorem 5.5]{voganzuckerman1984} by \cite[Proposition 3.6]{voganzuckerman1984}. Another remarkable point is that no derived operations are involved for cohomologically induced modules in this construction.

	\smallskip
	Motivated by these facts, we construct models of cohomologically induced modules over smaller commutative rings by defining the geometric ingredients (partial flag varieties, line bundles on them, and orbits) over such rings and using the theory of twisted $\D$-modules over schemes.

	\smallskip
	In order to set up this construction, we rely on partial flag schemes which are introduced as natural generalizations of partial flag varieties over algebraically closed fields. In particular, they are forms of complex partial flag varieties if the base is a subring of $\CC$.
	
	\smallskip
	Partial flag schemes are moduli schemes of geometric conjugacy classes of parabolic subgroups. To be precise, for a reductive group scheme $G$ over an arbitrary base scheme $S$ (\cite[D\'efinition 2.7]{sga3-19}), the moduli scheme $\Pp_G$ of parabolic subgroups of $G$ was introduced by the Grothendieck school, which may be called the total flag scheme of $G$ (\cite[Section 3.2, Th\'eor\`eme 3.3]{sga3-26}).
	
	\smallskip
	If we write $t$ for the quotient map $\Pp_G\to \type G$ by the conjugate action of $G$ on $\Pp_G$, a partial flag scheme of $G$ is the fiber of $t$ at an $S$-point of $\type G$, say, $\Pp_{G,x}$ for $x\in (\type G)(S)$. This formalism is crucial for our construction: Our intergral models of $(\lieg,K)$-modules will arise as the globalization of a certain equivariant twisted $\D$-module on $\Pp_{G,x}$.
	
	\smallskip
	In \cite{hayashilinebdl}, the first author constructed equivariant line bundles on partial flag schemes by Galois descent. We need the descent argument because we do not have base points in general. 
	
	\smallskip
	The moduli scheme of Borel subgroups of a reductive group scheme $G$ (\cite[Corollaire 5.8.3]{sga3-22}) is a typical example of partial flag schemes (see \cite[Theorem B.6]{hayashikgb} for a more explicit example without base points). It is evident that the flag scheme of $G$ has $S$-points if and only if $G$ admits a Borel subgroup over $S$. If we do not have base points, the construction of associated bundles is not available.
	
	\subsection{Models of orbits and orbit decompositions}\label{sec:orbit}
	
	In the last section, we proposed to construct arithmetic forms of cohomologically induced modules by application of the theory of twisted $\D$-modules over schemes to forms of the geometric ingredients. We saw how to get forms of partial flag varieties and line bundled on them. It remains to discuss what forms of orbits are.
	
	\smallskip
	To explain our idea, let us start with the projective line with the canonical action of the orthogonal group $\Oo(2)$ as a toy model. Over the field $\CC$ of complex numbers, we have a single closed $\Oo(2,\CC)$-orbit in the complex projective line. This consists of $\sqrt{-1}$ and $-\sqrt{-1}$. Since the complex conjugation respects this closed subvariety, it has a real form which is a closed subvariety of the real projective line. It is even defined over the ring $\ZZ$ of integers. In fact, we have the reduced closed subscheme $Z$ of the integral projective line $\PP^1_{\Spec\ZZ}=\Proj\ZZ[x,y]$ defined by the homogeneous polynomial $x^2+y^2$. As a scheme, we have $Z\cong\Spec\ZZ[t]/(t^2+1)$. In terms of the moduli description in \cite[Chapter II, Theorem 7.1]{hartshorneag}, the closed immersion is given by $(\ZZ[t]/(t^2+1),1,t)$. As for the group action, consider the canonical action of the general linear group scheme $\GL_2$ on $\PP^1_{\Spec\ZZ}$. Define a subgroup scheme $\Oo(2)\subset\GL_2$ as a group copresheaf by
	\[\Oo(2,R)=\{g\in\GL_2(R):~g^T=g^{-1}\},\]
	where $g^T$ denotes the transpose of $g$. Put the action of $\Oo(2)$ on $\PP^1_{\Spec\ZZ}$ by restriction. Then $Z$ is $\Oo(2)$-invariant. We notice that $Z$ has no $\ZZ$-point and does not even have real points. For this reason, $Z$ cannot be realized as an $\Oo(2)$-orbit in $\PP^1_{\Spec\ZZ}$.
	
	Likewise, we consider the action of the special orthogonal group $\SO(2)$ on the projective line. Over the complex numbers, each of $\{\sqrt{-1}\}$ and $\{-\sqrt{-1}\}$ is a closed $\SO(2,\CC)$-orbit. Since the complex conjugation switches them, they are not defined over the real numbers but their union is so. At the end, we obtain the same real form as above.
	
	\smallskip
	In the representation-theoretic perspectives at the level of orbits in the philosophy of the Beilinson--Bernstein equivalence, these observations correspond to the fact that discrete series representations of the special linear Lie group $\SL_2(\RR)$ are not defined over the real numbers but of $\GL_2(\RR)$ are so. The same reasoning over $\Spec\ZZ[1/2]$ tells us that discrete series representations of $\GL_2(\RR)$ (with appropriate central actions) should be defined over $\ZZ[1/2]$. Here we invert $2$ in order to avoid the singularity. In this base case, one can check by hand that this holds true. Alternatively, at the level of Lie algebras, we can construct a $\ZZ[1/2]$-form by using the induction functor $\ind$ over $\ZZ[1/2,\sqrt{-1}]$ in \cite{hayashi2019} and taking the restriction $\res_{\ZZ[1/2,\sqrt{-1}]/\ZZ[1/2]}$ in \cite{hayashi2018}. We conclude that $Z\otimes_{\ZZ}\ZZ[1/2]$ is the scheme we want to consider in the cases of $\GL_2$ and $\SL_2$.

	\smallskip
	As first generalization, we need real forms of $K_\CC$-orbits of $\theta$-stable parabolic subgroups in $\Pp_{G_\CC}$ for a general connected real reductive algebraic group $G$ with an involution $\theta$ and the fixed point subgroup $K=G^\theta$. We obtain those via Galois descent by studying the action of complex conjugation on the set of these orbits. Let us assume for simplicity that $K$ is connected. In this situation, a combinatorial classification of $K_\CC$-orbits was obtained by Brion--Helminck (\cite[Theorem 9.1]{brionhelminck}).
	
	It is well-known that for any $\theta_\CC$-stable parabolic subgroup $Q'$, $Q'\cap K_\CC$ is a parabolic subgroup of $K_\CC$. In particular, the $K_\CC$-orbit containing $Q'$ is closed. A possible clear explanation for this fact is as follows: recall that the parabolic subgroups of $G_\CC$ are exactly the subgroups of the form $P_{G_\CC}(\mu)$, where $\mu$ are cocharacters of $G_\CC$ (for the notation $P_{G_\CC}(\mu)$ and this fact we refer to \cite[the paragraph below Proof of Lemma 2.1.4 and Proposition 2.2.9]{conradpseudo} respectively). If this is $\theta_\CC$-stable, we can arrange $\mu$ so that $\mu$ factors through $K_\CC$. Then we have $P_{G_\CC}(\mu)\cap K_\CC=P_{K_\CC}(\mu)$. This is a key observation in \cite[Proof of Theorem 9.1]{brionhelminck}. Let us also note that once we know $\mu$ factors through $K_\CC$, the present argument goes without touching $\theta$. Namely, for a general connected reductive subgroup $K'\subset G_\CC$ and a cocharacter $\mu$ of $L'$, $P_{G_\CC}(\mu)\cap L'$ is a parabolic subgroup of $L'$. To promote \cite[Theorem 9.1]{brionhelminck}, we shall call a parabolic subgroup of this form stable relative to $L'$.

	By virtue of this fact, as complex variety, the moduli variety of $\theta_\CC$-stable parabolic subgroups of $G_\CC$ is the disjoint union of its $K_\CC$-orbits. In particular, its fppf quotient by $K_\CC$ is represented by a disjoint union of copies of $\Spec\CC$ indexed by the orbits, which we denote by $\rtype(G_\CC,K_\CC)$. The action of complex conjugation on the set of orbits determines a Galois action on $\rtype(G_\CC,K_\CC)$. Hence we obtain a real form $\rtype(G,K)$ by Galois descent.
	
	Moreover, Galois descent also provides us with a quotient map $rt$ from the moduli variety of $\theta$-stable parabolic subgroups of $G$ to $\rtype(G,K)$. Then real forms of $\theta_\CC$-stable $K_\CC$-orbits arise as the fibers of $rt$ at real points of $\rtype(G,K)$. More generally, taking all the fibers we obtain the minimal $K$-invariant decomposition of the moduli variety of $\theta$-stable parabolic subgroups of $G$. The previous example of $\SL_2$ is recovered by restricting ourselves to Borel subgroups. For general perfect ground fields of characteristic not two, the same argument using Galois descent applies.
	
	\smallskip
	In the present paper, we work over general base schemes (with $1/2$ in the construction of models for symmetric pairs). Motivated by the previous discussion, we introduce the notion of stable parabolic subgroups for a reductive group scheme $G$ and a smooth closed subgroup scheme $K$: We say a parabolic subgroup of $G$ is stable if it is \'etale locally of the form $P_G(\mu)$ for a cocharacter $\mu$ of $K$.
	
	This notion of stable parabolic subgroups defines a moduli problem, giving rise to the moduli space of stable parabolic subgroups. These moduli spaces are the appropriate generalization of those of $\theta_\CC$-stable parabolic subgroups. Along these lines we discuss \'etale local $K$-orbit decomposition and forms of each local $K$-orbit. As a result, we obtain forms of closed $K_\CC$-orbits of $\theta_\CC$-stable parabolic subgroups in Theorem \ref{thm:duality}.
	
	Putting all these ingredients together, we obtain arithmetic forms of cohomologically induced modules along the lines outlined above.

	\subsection{Future directions}\label{sec:future}
	
	\begin{description}
		\item[Deformation and contraction families] During the preparation of the present manuscript, Bernstein--Higson--Subag worked out deformation families and elementary examples of contraction families for some forms of classical groups over the complex and real projective lines in \cite{bernsteinetala,bernsteinetalb}. Subsequently, Barbasch--Higson--Subag generalized it to construct a real algebraic family of algebraic groups which interpolates between mutually dual real forms of a given complex algebraic group in \cite{barbaschetal}. Tan--Yao--Yu studied in \cite{tanetal2016} deformations of $\mathcal D$-modules in the case of $\SL_2(\RR)$ without providing formal justifications for the formalism they apply. We hope that our work contributes to a unification of these different approaches. The notion of families of Harish-Chandra pairs introduced by Bernstein--Higson--Subag corresponds to modules over the complex affine or projective lines in our context respectively. In \cite{hayashicontraction}, the first author gave purely algebraic proofs for results of \cite{barbaschetal} to allow us to work over general commutative ground rings with $1/2$. With the theory we develop here, we may work over the affine or projective line over $\ZZ[1/2]$.
		\item[General $K$-orbits in the flag scheme]
		Let us go back to the setting of $\PP^1_{\Spec\ZZ}$ in section \ref{sec:orbit}. One can easily show that the complementary subspace $U$ to $Z$ in $\PP^1_{\Spec\ZZ}$ is $\Oo(2)$-invariant. Moreover, $\Oo(2,F)$ acts transitively on $U(F)$ for every algebraically closed field $F$. One think of the set-theoretic decomposition $\PP^1_\ZZ=\Spec\ZZ[\sqrt{-1}]\coprod U$ as an integral form of the $\Oo(2,\CC)$-orbit decomposition of the complex flag variety of $\GL_2$. In fact, this decomposition is uniform in fields in the sense that we have $\PP^1_\ZZ(F)=(\Spec\ZZ[\sqrt{-1}])(F)\coprod U(F)$ for every field $F$ (\cite[Theorem 2.1]{hayashikgb}). The first author improved this idea to establish a decomposition of the flag scheme of $\SL_3$ over $\ZZ\left[1/2\right]$ (\cite[Theorem 1.1]{hayashikgb}). More generally, it is interesting to study forms of the $K$-orbit decomposition of the flag variety for arithmetric structures of the dual of cohomologically induced modules and irreducible Harish-Chandra modules arising from general $K$-orbits. The first author is currently working on it (cf.\ \cite{hayashiuniform}).
		\item[Localization of the equivariant Zuckerman functor]
		Mili\v{c}i\'{c}--Pand\v{z}i\'{c} made a proposal to work with equivariant derived categories for a functorial approach to \cite[Theorem 4.3]{hechtetal} from perspectives of homological algebra, and established the localization theory of the equivariant Zuckerman functor in the full flag case (\cite{milicicpandzic1998}).
		
		If we work on general $K_\CC$-orbits in partial flag varieties, the duality theorem of \cite[Theorem 4.3]{hechtetal} does not work since $K_\CC$-orbits are not imbedded affinely in general. Following Mili\v{c}i\'{c}--Pand\v{z}i\'{c}'s approach, Kitchen proved an equivariant derived version of \cite[Theorem 4.3]{hechtetal} (\cite[Theorem 5.4, Corollary 5.2, see also Corollary 5.6]{kitchen2012}).
		
		If we aim to work on forms of $K_\CC$-orbits in partial flag varieties and attached (equivariant complexes of) representations, a possible fundamental approach is to work out the theory of equivariant derived categories over schemes and their flat base change theorems. The study of derived functorialities of twisted $\D$-modules in the present paper could serve as a solid foundation of such an approach.
		
		\item[Number-theoretic applications]
		Our models of cohomologically induced modules offer a pathway to the construction of half-integral analogs of the rational structures on automorphic representation obtained by the second author in \cite{januszewskirationality}. This is executed in \cite{januszewskisheaves}, where global $1/N$-integral models of cohomological automorphic representations are constructed.
		An application of global automorphic $1/N$-integral structures is a canonical $v$-integral normalization of the canonical period matrices defined in loc.\ cit.\ for $v\nmid N$ a place of the field of definition $\QQ(\Pi)$ of the automorphic representation $\Pi$ under considaration.
		
		The same ideas extend to coherent cohomology of Shimura varieties as outlined in \cite[section 8.5]{januszewskirationality}.
		
		In both cases, global integral models induce new integral structures in the cohomology of artihmetic groups (as well as in coherent cohomology) which transcend cohomological degrees. Therefore they are expected to fit into Venkatesh's conjectural the framework \cite{venkatesh2017,venkatesh2018ICM,venkatesh2019} on the action of the derived Hecke algebra on cohomology.
	\end{description}

	\subsection{Results and organization of this paper}\label{sec:organization}
	
	We give a more detailed overview of the contents and main results of this paper. Sections \ref{sec:tdo_and_pic_alg}--\ref{sec:operationswithDmodules} treat generalities on tdos and twisted $\D$-modules over schemes. We collect preliminary definitions and related elementary facts in Appendix \ref{appendix}. Section \ref{sec:derivation} gives preliminary arguments on differential calculus which will be used in sections \ref{sec:tdo_and_pic_alg}--\ref{sec:operationswithDmodules}. The other sections of Appendix will be explained below in relation to how they are invoked in sections \ref{sec:tdo_and_pic_alg}--\ref{sec:operationswithDmodules}.

	\smallskip
	Section \ref{sec:tdo_and_pic_alg} is devoted to generalities on tdos and Picard algebroids over schemes. We adopt the straightforward generalization of \cite[Definitions 2.1.1, 2.1.3]{beilinsonbernstein1993} as their definitions (sections \ref{sec:tdo}, \ref{sec:picalg} respectively). For convenience to the reader, we review filtrations of sheaves in section \ref{sec:filt} for section \ref{sec:tdo}. For local study in terms of commutative rings (e.g.~Variant \ref{var:comparisonf_+f^+;imm}), tdos and Picard algebras for smooth homomorphisms of commutative rings are introduced in parallel (the ring setting).
	
	We give basic definitions and operations of torsors in section \ref{sec:torsor} for study of Picard algebroids in sections \ref{sec:picalg} and \ref{sec:functorspicalgtdo}.
	
	In section \ref{sec:tdo_and_pic_alg}, we confirm that large parts of the theory of tdos in \cite{kashiwara1989,beilinsonbernstein1993} work over general base schemes. Therefore overlapped proofs will be skipped. For applications to $(\lieg,K)$-modules, we treat equivariant structures in section \ref{sec:equivtdo}. For this, the canonical equivariant structure on the (co)tangent sheaves is recorded in section \ref{sec:equivsheaf}.

	Towards the construction of the bimodule structure on $\Aa_{X\to Y}\coloneqq f^\ast\Aa$ for a morphism $f:X\to Y$ of smooth schemes over an arbitrary base $S$ and a tdo $\Aa$ on $Y$, we give its left action in terms of the Picard algebroid corresponding to $\Aa$ (Construction \ref{cons:ambientfE}, Theorem \ref{thm:fcdotispicalg/sch}).
	
	\smallskip
	A notable new addition to the literature and critical to our applications is the introduction of the base change functor in section \ref{sec:Sbasechange}. This is achieved by defining the corresponding functor to $\Omega^{\geq 1}_{X/S}$-torsors, where for a smooth morphism $X\to S$, $\Omega^{\geq 1}_{X/S}$ is the truncated de Rham complex at the zeroth degree. The base change functors for tdos and Picard algebroids are locally given by the usual base change $k'\otimes_k-$ for a homomorphism $k\to k'$ of base rings in terms of the ring setting (section \ref{sec:Sbasechange}). As a fundamental property, we see that tdos form a stack over the fpqc site of the base scheme (Theorem \ref{thm:tdofpqcstack}). The main result of section \ref{sec:tdo_and_pic_alg} is the base change property of pullbacks:

	\begin{theorem}[Theorem \ref{thm:basechangevspullbackfortdo}]\label{thm:bc_tdo}
		Let $S'\to S$ be a morphism of schemes. Let $f:X\to Y$ be a morphism of smooth $S$-schemes, and $f':X'\coloneqq X\times_S S'\to Y\times_S S'\reflectbox{$\coloneqq$} Y'$ be its base change. Let $s_X:X'\to X$ and $s_Y:Y'\to Y$ be the canonical projections.
		
		Let $\Aa$ be a tdo on $Y$. Then there is a natural isomorphism
		$s^\ast_X f^\cdot \Aa\cong (f')^\cdot s^\ast_Y\Aa$,
		where $s^\ast_X$ and $s^\ast_Y$ (resp.~$f^\cdot$ and $(f')^\cdot$) are the base change (resp.~pullback) functors.
	\end{theorem}
	
	This is a crucial ingredient for base change properties of twisted $\D$-modules and their direct image functors.
	
	\smallskip
	Section \ref{sec:Dmodules} contains mostly elementary but technical generalities of twisted $\mathcal D$-modules necessary for the derived functoriality statements in section \ref{sec:operationswithDmodules}. In section \ref{sec:twistedd-mod}, we give some elementary examples of twisted $\D$-modules highlighting different phenomena from the classical theory. Then we give categorical and homological study of twisted (quasi-coherent) $\D$-modules in sections \ref{sec:twistedd-mod} and \ref{sec:qcd-mod}. For a digression, we also give preliminary arguments on twisted coherent $\D$-modules on locally Noetherian smooth schemes in section \ref{sec:coh} for future applications to their characteristic schemes.

	An important result in section \ref{sec:Dmodules} is to show equivalence of the settings of sheaves and rings for affine schemes at the categorical level (Lemma \ref{lem:modulecataffinecase}). As a consequence, we see the equivalence of the definitions of quasi-coherent twisted $\D$-modules in \cite[(5.1.3)]{ega1} and \cite[Notation 1.4.1]{hottaetal2008} over general bases (Proposition \ref{prop:defnofqcoh}). This is helpful since on the one hand \cite[Notation 1.4.1]{hottaetal2008} is useful for applications of algebro-geometry arguments and on the other hand it is sometimes easier to verify that a given twisted $\D$-module is quasi-coherent in the sense of \cite[(5.1.3)]{ega1}.
	
	\smallskip
	The category of main interest to us in the next section is the full subcategory of the unbounded derived category of left $\Aa$-modules consisting of cohomologically quasi-coherent complexes for a tdo $\Aa$, which we denote by $D_{\qc}(\Aa)$ and $D(\Aa)$ respectively. As a technical fact, we show that a q-injective (resp.~q-flat) complex of a twisted $\D$-modules on a smooth scheme $X$ is so as a complex of $\OO_X$-modules, where $\OO_X$ is the structure sheaf of $X$ (Proposition \ref{prop:resolutions}). This is important for comparison of derived operations of twisted $\D$-modules and of $\OO_X$-modules mentioned below.
	
	\smallskip
	In section \ref{sec:operationswithDmodules}, we study functorial operations of twisted $\mathcal D$-modules: globalization $Rx_\ast$, localization $\Aa\otimes^L_{x^{-1}x_\ast\Aa} x^{-1}(-)$, base change $Ls^\ast_X$, base restriction $R(s_X)_\ast$, derived inverse image $Lf^\ast$, and derived direct image $f_+$. We investigate to which extent the classical theorems on (twisted) $\D$-modules remain valid in our setting. We also study the relation of the (flat) base change functor with other operations (flat base change theorems). For this, preliminary generalities on unbounded derived categories of sheaves of modules and their functoriality are collected in section \ref{sec:homologicalalgebra}.

	\smallskip
	We define the globalization functor $x_\ast$ as the direct image functor to the base scheme (section \ref{sec:localization}). To be precise, for a smooth scheme $X$ over a base $S$ with the structure morphism $x:X\to S$ and a tdo $\Aa$ on $X$, we define $x_\ast$ as a functor from the category of left $\Aa$-modules to that of $x_\ast\Aa$-modules in the standard manner. Then the localization is defined as its right adjoint. We see that the derived functor $Rx_\ast:D(\Aa)\to D(x_\ast\Aa)$ is a lift of
	$Rx_\ast:D(\OO_X)\to D(\OO_S)$ (\eqref{diag:compatibilityofglobalization}). The symbols $D(-)$ denotes the corresponding unbounded derived categories of modules. We see as a consequence that $Rx_\ast$ respects complexes which are quasi-coherent over the structure sheaves if $x$ is concentrated, i.e., quasi-compact and quasi-separated (Proposition \ref{prop:x_ast;qc_pres}). We also see that the localization functor respects cohomologically quasi-coherent complexes as well under this assumption (Proposition \ref{prop:x_ast;qc_pres}, Remark \ref{rmk:onDaffinei}).

	\smallskip
	We introduce the transfer bimodule in section \ref{sec:bimodules} to define the derived inverse image functor $Lf^\ast$ and the derived direct image functor $f_+$ for a morphism $f:X\to Y$ of smooth schemes over an arbitary base $S$ (sections \ref{sec:derivedinverse}, \ref{sec:deriveddirectimage}). We see that $Lf^\ast$ lifts the classical derived inverse image functor $Lf^\ast:D(\OO_Y)\to D(\OO_X)$. In particular, $Lf^\ast$ respects cohomologically quasi-coherent complexes, and it is natural in $X$ (Theorem \ref{thm:inverseimage}).
	
	\smallskip
	In principle, one could study functoriality and preservation of cohomologically quasi-coherent complexes in a similar way to the case of $S=\Spec F$, where $F$ is a field of characteristic zero. One could prove the vanishing of cohomology for affine immersions as well.
	However, we sometimes need some finiteness conditions on $Rf_\ast$ for technical reasons. We collect preliminary discussions on them in section \ref{sec:boundRf_ast}. For simple explanation, let us restrict ourselves to the following cases from Examples \ref{ex:cloimm_dim0}, \ref{ex:findim}, Corollary \ref{cor:loc_noe_fin_dim} if needed:
	\begin{enumerate}
		\item[(i)] $f$ is a closed immersion;
		\item[(ii)] $X$ is Noetherian of finite dimension;
		\item[(iii)] $f$ is quasi-compact, and $X$ is locally Noetherian of finite Krull dimension.
	\end{enumerate}
	We say for short that $f$ satisfies (LB) if one of these conditions holds. In particular, $f$ is concentrated in this case. The main results of section \ref{sec:deriveddirectimage} are:

	\begin{theorem}[Theorems \ref{thm:qcohpreservation}, \ref{thm:compositionlaw}]\label{thm:main_f_+}
		Let $S$ be an arbitrary scheme.
		\begin{enumerate}
			\item For a morphism $f:X\to Y$ of smooth $S$-schemes and a tdo $\Aa$ on $Y$, the functor $f_+:D(f^\cdot\Bb)\to D(\Aa)$ respects cohomologically quasi-coherent complexes.
			\item Let $i:Y\to X$ be an affine immersion of smooth $S$-schemes with (LB), and $\Aa$ be a tdo on $X$. Then we have $H^p(i_+\Mm)=0$ for every left $i^\cdot\Aa$-module $\Mm$ and $p\neq 0$.
			\item Let $f:X\to Y$ and $g:Y\to Z$ be morphisms of smooth $S$-schemes, and $\Aa$ be a tdo on $Z$. If $f$ satisfies (LB) then $(g\circ f)_+=g_+\circ f_+$ on $D_{\qc}(g^\cdot f^\cdot \Aa)$.
		\end{enumerate}
	\end{theorem}
	
	Though statements are given for right modules, they immediately imply the variants for left modules as stated above by definition of $f_+$ for left modules. We note that (2) is reduced to the ring setting after Zariski localization of $X$ and $S$. We need the hypothesis (LB) to verify compatibility of $i_+$ with its version in the ring setting (Theorem \ref{thm:comparisonf_+f^+}, Variant \ref{var:comparisonf_+f^+;imm}).
	
	\smallskip
	In contrast, Kashiwara's theorem fails partly or almost completely over general bases. Let $i:Y\to X$ be a closed immersion of smooth schemes over a scheme $S$, and $\Aa$ be a tdo on $X$. Then $i_+$ is still faithful (Corollary \ref{cor:kashiwara} (1)), but no longer full in general. A counterexample for $S=\Spec\ZZ/p\ZZ$ with $p$ prime is given in terms of the ring setting in Example \ref{ex:counterexofkashiwaraeq}. The basic reason for the failure is the vanishing of some positive integers. We prove that $i_+$ is fully faithful if we restrict ourselves to torsion-free sheaves over $\ZZ$ (Corollary \ref{cor:kashiwara} (2)). Moreover, the essential image of quasi-coherent $i^\cdot\Aa$-modules along $i_+$ consists of quasi-coherent $\Aa$-modules supported on $Y$ if $X$ is locally Noetherian and all nonzero integers are invertible on $X$ (Theorem \ref{thm:kashiwaraequiv}, Corollary \ref{cor:kashiwara_ess_img}). We give a counterexample to this fact for $S=\Spec\ZZ$ (Example \ref{ex:counterexofessimg}). The failure of Kashiwara's theorem leads to a failure of the generalization of the classical base change theorem of \cite[Theorem 1.7.3]{hottaetal2008}. In this paper, we call it the $X$-base change theorem to emphasize that it is a base change theorem in $X$; base change theorems in $S$ which we will explain below are sometimes called $S$-base change theorems. A counterexample to the $X$-base change theorem over schemes is given in Example \ref{ex:counterexofbcthm}. However, we still have:
	
	\begin{theorem}[Theorem \ref{thm:xbc}]
		Suppose that we are given a Cartesian diagram of $S$-schemes
		\[\begin{tikzcd}
			\tilde{X}\ar[r, "\tilde{p}"]\ar[d, "\tilde{f}"']
			&\tilde{Y}\ar[d, "f"]\\
			X\ar[r, "p"]&Y.
		\end{tikzcd}\]
		Assume that the following conditions are satisfied:
		\begin{enumerate}
			\renewcommand{\labelenumi}{(\roman{enumi})}
			\item All of $X$, $Y$, $\tilde{X}$, and $\tilde{Y}$ are smooth over $S$.
			\item $p$ is isomorphic to a projection morphism in the category of morphisms of $S$-schemes.
			\item $f$ is concentrated.
		\end{enumerate}
		Let $\Aa$ be a tdo on $Y$. Then there is an equivalence
		$p^\ast\circ f_+\simeq \tilde{f}_+\circ \tilde{p}^\ast$
		on $D_{\qc}(f^\cdot\Aa)$.
	\end{theorem}
	
	This is enough in our applications because we can then lift the functor $i_+$ to the equivariant setting for an affine immersion $i$ with (LB):
	
	\begin{corollary}[Theorem \ref{thm:inductionofequivtwistedDmod}]\label{cor:equiv}
		Let $G$ ba a smooth $S$-affine group scheme over an arbitrary base $S$, $i:Y\hookrightarrow X$ be a $G$-equivariant affine immersion of smooth $G$-schemes over $S$, and $\Aa$ be a $G$-equivariant tdo on $X$. Suppose that $i$ and $G\times_Si$ satisfy (LB). Then for a $G$-equivariant left $i^\cdot\Aa$-module $\Mm$, $i_+\Mm$ is naturally endowed with the structure of a $G$-equivariant left $\Aa$-module.
	\end{corollary}
	
	\smallskip
	As a new aspect of twisted $\D$-modules, we introduce the adjunction of base change and restriction functors which is an enhancement of the classical one: For a morphism $S'\to S$ of schemes and a tdo $\Aa$ on a smooth $S$-scheme $X$, we define an adjunction
	\[Ls^\ast_X:D(\Aa)\rightleftarrows D(s^\ast_X\Aa):R(s_X)_\ast,\]
	and prove that it is a lift of
	$Ls^\ast_X:D(\OO_X)\rightleftarrows D(\OO_{X'}):R(s_X)_\ast$
	(Proposition \ref{prop:derivedadjunction_basechange}, recall Theorem \ref{thm:bc_tdo} for our notations). As its immediate consequence, we find that twisted $\D$-modules form a stack over the fpqc site of the base $S$ (Theorem \ref{thm:fpqcDmoduledescent}).
	
	\smallskip
	Our main results on $Ls^\ast_X$ are flat base change theorems of the globalization and the derived direct image functor:
	
	\begin{theorem}[Theorem \ref{thm:flatbasechangeofglobalsectionI}, \ref{thm:qcohpreservation}]\label{thm:bc}
		Let $s:S'\to S$ be a flat morphism of schemes, and $\Bb$ be a tdo on a smooth $S$-scheme $Y$ with the structure morphism $y$. Define $Y'$ and $s_Y$ as in Theorem \ref{thm:bc_tdo}. Write $y'$ for the canonical projection $Y'\to S'$.
		\begin{enumerate}
			\item If $y$ is concentrated, we have a canonical equivalence
			$Ry'_\ast\circ s^\ast_Y \simeq s^\ast \circ Ry_\ast$ on $D_{\qc}(\Bb)$.
			\item Let $f:X\to Y$ be a concentrated morphism of smooth $S$-schemes, and $s_X:X\times_S S'\to X$ be the projection. Set $f'=f\times_S S'$. Then we have a canonical equivalence
			$s^\ast_Y\circ f_+\simeq f'_+\circ s^\ast_X$ on $D_{\qc}(f^\cdot\Bb)$.
		\end{enumerate}
	\end{theorem}
	
	Part (1) is an immediate consequence of the corresponding statement for $\OO_Y$-modules (\cite[Proposition (3.9.5)]{lipman2009}). Part (2) is verified again by reducing to the cases of closed immersions and projections (Theorems \ref{thm:Sbasechange1}, \ref{thm:Sbasechange2}). With the classical faithfully flat descent, we see that the $\Aa$-affinity (\cite[Definition 1.4.2]{hottaetal2008}, see also Definition \ref{def:Daffine}) is invariant under base change of fields (Corollary \ref{cor:A-affine_is_geometric}). The base restriction functor $(s_X)_\ast$ is used for its proof.
	
	\smallskip
	The important conclusion of the preceding arguments for our purpose is to obtain $(\lieg,K)$-modules over commutative rings in the sense of \cite{hayashi2019} via the direct image functor of twisted $\D$-modules and the global section functor:
	
	\begin{theorem}[Theorems \ref{thm:inductionofequivtwistedDmod}, \ref{thm:globalization}]
		Let $i:Y\hookrightarrow X$ be an affine immersion of smooth $k$-schemes over an arbitrary commutative ring $k$, equipped with actions of respective smooth affine group schemes $K,G$ over $k$. Assume that $i$ and $K\times_S i$ satisfy (LB). Suppose that we are given a homomorphism $K\to G$ for which $i$ intertwines the actions. Let $\Aa$ be a $G$-equivariant tdo on $X$. Then for a $K$-equivariant left $i^\cdot\Aa$-module $\Mm$, the $k$-module $\Gamma(X,i_+\Mm)$ is naturally equipped with the structure of a $(\lieg,K)$-module. Moreover, this construction commutes with flat base change functors.
	\end{theorem}
	
	In fact, the statement extends to schemes. A general formalism on $(\lieg,K)$-modules over schemes is collected at the beginning of section \ref{sec:(g,K)-module}.
	
	\smallskip
	In section \ref{sec:descentoforbit}, we study the moduli space of stable parabolic subgroups. Let $G$ be a reductive group scheme over an arbitrary base $S$, and $K$ be a smooth closed subgroup scheme. We introduce the notion of stable parabolic subgroups of $G$ for $(G,K)$ (see the last paragraph of section \ref{sec:orbit} or Definition \ref{defn:stable}). We verify in Proposition \ref{prop:stable_vs_theta_stable} that this is a generalization of $\theta$-stable parabolic subgroups of $G$ for an involution $\theta$ when $2$ is invertible on $S$, i.e., a parabolic subgroup of $G$ is $\theta$-stable if and only if it is stable with respect to $(G,K)$, where $K$ is an open and closed subgroup scheme of the fixed point subgroup scheme of $G$ by $\theta$. Let us call it the symmetric setting for convenience.
	
	\smallskip
	Let $\Pp_{G}^{K-\mathrm{st}}$ be the moduli space of stable parabolic subgroups of $G$. This is equipped with an action of $K$ by conjugation. We are interested in its fppf (\'etale) local $K$-orbit decomposition and a form of each local $K$-orbit. Let $\rtype(G,K)$ be the fppf (\'etale) quotient of $\Pp_{G}^{K-\mathrm{st}}$ by the conjugate action of $K$, and $rt:\Pp_{G}^{K-\mathrm{st}}\to \rtype(G,K)$ be the canonical quotient map. Our main result in this section is:
	\begin{theorem}[Theorem \ref{thm:orbit}]
		Let $G$ be a reductive group scheme, and $K$ be a smooth closed subgroup scheme with reductive unit component such that the component group scheme $\pi_0(K)$ is finite \'etale. 
		\begin{enumerate}
			\renewcommand{\labelenumi}{(\arabic{enumi})}
			\item The $S$-space $\Pp_{G}^{K\mathrm{-st}}$ is represented by a smooth closed subscheme of $\Pp_G$.
			\item The $S$-space $\rtype(G,K)$ is represented by a finite \'etale $S$-scheme.
			\item Stable parabolic subgroups $P$ and $P'$ are \'etale locally $K$-conjugate to each other if and only if $rt(P)=rt(P')$.
			\item The morphism $rt:\Pp^{K\mathrm{-st}}_G\to \rtype(G,K)$ is smooth projective and surjective.
		\end{enumerate}
	\end{theorem}
	
	We then obtain $S$-form of \'etale local $K$-orbits by taking the fiber $\Pp_{G,x}^{K\mathrm{-st}}$ of $rt$ at $x\in \rtype(G,K)(S)$. We define $gt:\rtype(G,K)\to \type G$ by descending the canonical embedding $\Pp_{G}^{K\mathrm{-st}}\hookrightarrow\Pp_G$. This is to specify the ``ambient'' partial flag scheme of $\Pp_{G}^{K\mathrm{-st}}$. In fact, we obtain a closed immersion $i_x:\Pp_{G,x}^{K\mathrm{-st}}\hookrightarrow\Pp_{G,gt(x)}$ (Theorem \ref{thm:orbit} (6)). This models a closed $K$-orbit of a partial flag scheme of $G$.
	
	\smallskip
	In section \ref{sec:rem}, we consider the symmetric setting with $K$ reductive. We give a classification of $K$-orbits of $\theta$-stable parabolic subgroups under reasonable splitting conditions in the form of \cite[Proof of Theorem 9.1]{brionhelminck} (Theorem \ref{thm:solutiontotheclassificationproblem}). As application, we compute $\rtype(G,K)$ and $gt$ for the standard $\ZZ[1/2]$-forms of classical Lie groups in \cite{hayashi2019} (section \ref{sec:halfintegralformofrelativetype}). For computation of descent, a few elementary generalities and examples are collected in section \ref{sec:absolutedescent}.
	
	\smallskip
	Section \ref{sec:arithmeticmodelsofaq(lambda)} is devoted to applications of the previously developed theories to construction of arithmetic models of cohomologically induced modules. Based on the preceding ideas, we introduce a general geometric construction of $(\lieg,K)$-modules over schemes in section \ref{sec:general_construction}. In sections \ref{sec:kitchen}, we collect general facts on cohomologically induced modules, based on the algebro-geometric language in \cite{adamstaibi2018}.
	
	\smallskip
	Then as a byproduct of the preceding arguments, we obtain one of the main results throughout this paper, i.\,e., we establish a general formalism to get forms of cohomologically induced modules from those of their geometric ingredients:
	
	\begin{theorem}[Corollary \ref{cor:sformofcohohomologicalinduction1}]\label{thm:main}
		Suppose that we are given a commutative diagram
		\[\begin{tikzcd}
			\Spec\CC\ar[r]\ar[d, "s"']&\Spec\RR\ar[d]\\
			\tilde{S}\ar[r]&S
		\end{tikzcd}\]
		of $\ZZ[1/2]$-schemes with $s$ flat. Let $G$ be a reductive group scheme, equipped with an involution $\theta$ with $\theta\times_S\Spec\RR$ Cartan in the sense of \cite{adamstaibi2018} (see also section \ref{sec:kitchen}). Let $K$ be an open and closed subgroup scheme of the fixed point subgroup scheme of $G$ by $\theta$. Assume that $\pi_0(K)\times_S \tilde{S}$ is finite \'etale over $\tilde{S}$. Let $x\in\rtype(G,K)(\tilde{S})$, and $\Aa$ be a $G\times_S \tilde{S}$-equivariant tdo on $\Pp_{G\times_S \tilde{S},gt(x)}$.
		
		Choose a maximal torus $T\subset K\times_S \Spec\RR$ and a complex $\theta\times_S \Spec\CC$-stable parabolic subgroup $Q'\subset G\times_S \Spec\CC$ containing $T\times_{\Spec\RR}\Spec\CC$ such that $rt(\bar{Q}')=x|_{\Spec\CC}$. Let $u$ be the dimension of the unipotent radical of $Q'\cap K\times_S \Spec\CC$. 		
		
		Let $\Mm$ be a $K\times_S\tilde{S}$-equivariant left $i^\cdot_{x}\Aa$-module which is locally free of finite rank over $\OO_{\Pp^{K\times_S \tilde{S}\mathrm{-st}}_{G\times_S\tilde{S},x}}$. Let $M'$ denote the geometric fiber of $\Mm$ at $\bar{Q}'\in \Pp^{K\times_S \tilde{S}\mathrm{-st}}_{G\times_S\tilde{S},x}(\CC)$. Then
		\[R^\bullet (p_{G\times_S \tilde{S},x})_\ast(i_{x})_+\Mm\]
		exhibits an $\tilde{S}$-form of $\Ll_{u-\bullet}(M')$.
	\end{theorem}

	Theorem \ref{thm:fil} guarantees that the resulting $\tilde{S}$-forms are new to rational forms of cohomologically induced modules from Theorem \ref{thm:fil}. For instance, if $\tilde{S}$ is the spectrum of a Dedekind domain $k'$ (of dimension one), the corresponding $k'$-module is projective.
	
	\smallskip
	In section \ref{sec:gl_n} we discuss the Weil restriction of $\GL_n$ over totally real or CM number fields to $\QQ$. This is motivated by automorphic representation theory. We make an attempt to give numerical data to Theorem \ref{thm:main} in these cases. We do not work with arbitrary totally imaginary number fields but only CM fields since we need well behaved forms of maximal compact subgroups (see Proposition \ref{prop:U(n)} and its paragraph below)
	
	For simplicity, we only consider the sheaves of differential operators on equivariant line bundles over the Weil restriction of partial flag schemes of $\GL_n$. Accordingly, we put the restriction of these line bundles for $\Mm$. This is enough for our applications to cuspidal cohomological automorphic representations of $\GL_n$.
	
	\smallskip
	Towards this, we determine the minimal forms of characters of base changes of the Weil restriction $\res_{E/\QQ} T_{\spl}$ from $E$ to $\QQ$ of a split torus $T_{\spl}$ for a number field $E$ in terms of those over its Galois hull $E^{\Gal}$ and their Galois twists (Theorem \ref{thm:minimal_form}). Then we obtain equivariant line bundles on partial flag schemes by using \cite[Proposition 2.3.7, Example 2.3.8]{hayashilinebdl}.
	
	Since we use fundamental Cartan subgroups for parameterization of cohomologically induced modules in \cite{voganzuckerman1984} (see also Corollary \ref{cor:irreducible_unitary}), we compare character groups of the Weil restriction of a certain concrete fundamental Cartan subgroup and of the split maximal torus of diagonal matrices in the totally real case (section \ref{sec:totally_real}). For CM number fields, we may and do use the same maximal torus since the torus of diagonal matrices is a fundamental Cartan subgroup in this case.
	
	Then we obtain forms of cohomologically induced modules over certain explicitly defined commutative ring (Theorems \ref{thm:gl_n_totally_real}, \ref{thm:AqlambdaGLnCM}). They are canonical in the sense that they are determined from parabolic types of $\GL_n$ and the characters $\lambda$ over $E^{\Gal}$ since we have only one closed homogeneous subscheme in each partial flag scheme attached to a stable parabolic subgroup (see Propositions \ref{prop:orbit_GL_n} and \ref{prop:numberoforbits1}, Corollary \ref{cor:numberoforbits}).
	
	To apply generalities on cohomologically induced modules in \cite{knappvogan}, fix an embedding $E^{\Gal}\hookrightarrow\RR$ (resp.~$E^{\Gal}\hookrightarrow \CC$) if $E$ is totally real (resp.~CM). We choose a positive system fixed by (a form of) the standard Cartan involution. We fix an embedding $E^{\Gal}\hookrightarrow\RR$ (resp.~$E^{\Gal}\hookrightarrow \CC$) if $E$ is totally real (resp.~CM) to give a numerical condition on $\lambda$ so that $\lambda$ is strictly dominant and that the corresponding cohomologically induced module is essentially unitary (Definitions \ref{def:A(V)_totally_real_even}, \ref{def:A(V)_totally_real_odd}, \ref{def:A(V)_CM}). In particular, our partial flag scheme is forced to be full by the regularity of $\lambda$. Hence these forms are determined only from such $\lambda$.
	
	\smallskip
	Finally, we discuss the cohomology groups of our models of cohomologically induced modules in section \ref{sec:(g,K)-cohomology}. Towards applications to automorphic representations, we introduce the notion of $(\lieg,\lies,K)$-cohomology, which generalizes the $(\lieg,K_\infty)$-cohomology in number theory. That is, we take the action of the center into consideration. For instance, if $G=\GL_2$ over $\RR$, $K_\infty$ is the Lie subgroup of $\GL_2(\RR)$ generated by $\Oo(2)$ and the unit component of the center of $\GL_2(\RR)$. We record how to compute it in terms of the classical relative Lie algebra cohomology (Proposition \ref{prop:relativeliealgcoh}).
	
	Then we collect results of the Vogan--Zuckerman theory for possibly disconnected real reductive Lie groups arising from real reductive algebraic groups for convenience to the reader (Theorems \ref{thm:a_q->cohomological}, \ref{thm:cohomological->aq}). Thanks to their results, we find that the forms in Definitions \ref{def:A(V)_totally_real_even}, \ref{def:A(V)_totally_real_odd}, \ref{def:A(V)_CM} are cohomological, and that they are determined by its coefficient through its lowest weight. That is, for each finite dimensional representation $V$ of lowest weight $-\lambda$ with $\lambda$ in Definitions \ref{def:A(V)_totally_real_even}, \ref{def:A(V)_totally_real_odd}, \ref{def:A(V)_CM}, there is only one irreducible essentially unitary cohomological representation whose coefficient $V$ (see Corollary \ref{cor:unique} for a general statement). We can construct its form without any more choices as in Definitions \ref{def:A(V)_totally_real_even}, \ref{def:A(V)_totally_real_odd}, \ref{def:A(V)_CM}.
	
	\smallskip
	Our final main result is to give an explicit formula for the relative Lie algebra cohomology of our models of cohomologically induced modules under mild conditions. In particular, we prove that they are finitely generated and projective. To save space, we denote $(-)_{R'}=-\otimes_R R'$ for a homomorphism $R\to R'$ of commutative rings.
	
	\begin{theorem}[Variant \ref{var:gsKcohomology}]
		Consider the setting of Theorem \ref{thm:main}. Assume that $S,\tilde{S}$ are affine with coordinate rings $k,\tilde{k}$ respectively. Write $H$ for the centralizer of $T$. Fix a Borel subgroup $B'\subset G_\CC$ with $T_\CC\subset B'\subset Q'$. Let $\rho$ be the half sum of roots of $B'$ with respect to $H_\CC$.
		
		Let $V$ be a $(\lieg_{\tilde{k}},K_{\tilde{k}})$-module which is finitely generated and projective as a $\tilde{k}$-module.
		Suppose that the following conditions are satisfied:
		\begin{enumerate}
			\renewcommand{\labelenumi}{(\roman{enumi})}
			\item $\tilde{k}$ is a Dedekind domain;
			\item the map $\tilde{k}\to \CC$ corresponding to $s$ is injective;
			\item $\Mm=\Ll$ is a line bundle on $\Pp^{K_{\tilde{k}}\mathrm{-st}}_{G_{\tilde{k}},x}$;
			\item $V_\CC$ is irreducible as a $\lieg_\CC$-module;
			\item There exists a one-dimensional $(\lieg_\CC,K_\CC)$-module $\tau$ such that $M\otimes_\CC\tau$ is unitary;
			\item Let $\lambda$ be the character of the Lie algebra $\lieh_\CC$ of $H_\CC$ for which $\lieh_\CC$ acts on $M$. Then the pairing of $\alpha^\vee$ and $\lambda+\rho$ does not belong to $(-\infty,-1]$ for any coroot $\alpha^\vee$ of the unipotent radical of $Q'$.
		\end{enumerate}
		Let $F_\bullet \Gamma(\Pp_{G_{\tilde{k}},gt(x)},(i_x)_+\Ll)$ be the filtration on $\Gamma(\Pp_{G_{\tilde{k}},gt(x)},(i_x)_+\Ll)$ in Theorem \ref{thm:fil}.
		
		Let $\mu$ be the lowest weight of $V_\CC$. Let $\lies$ be a $\tilde{k}$-submodule of the center of $\lieg_{\tilde{k}}$. Then we have
		\begin{flalign*}
			&H^\bullet(\lieg_{\tilde{k}},\lies_{\tilde{k}}, K_{\tilde{k}},\Gamma(\Pp_{G_{\tilde{k}},gt(x)},(i_x)_+\Ll)\otimes_{\tilde{k}} V)\\
			&\cong 
			\Hom_{K_{\tilde{k}}} (\wedge^{\bullet}\lieg_{\tilde{k}}/(\liek_{\tilde{k}}+\lies_{\tilde{k}}),
			F_0\Gamma(\Pp_{G_{\tilde{k}},gt(x)},(i_x)_+\Ll)\otimes_{\tilde{k}} V).
		\end{flalign*}
		Moreover, it vanishes unless $\lambda=-\mu$.
		In particular, the cohomology has no torsion.
	\end{theorem}
	
	This is verified as a consequence of Theorems \ref{thm:main}, \ref{thm:fil}, and the proof of \cite[Theorem 5.5]{voganzuckerman1984}. We assume (iii)-(vi) in order to apply results of unitarizability and computation of the cohomology over $\CC$ in \cite{knappvogan} and \cite{voganzuckerman1984} respectively. Assumption (i) is used for Theorem \ref{thm:fil}. We assume (ii) in order to derive the result from that over $\CC$.	
	
	\smallskip
	In our treatment we suppose a fluent understanding of modern scheme-theoretic language. Throughout the paper we make frequent use of results from Grothendieck--Dieudonn\'e's EGA I-IV \cite{ega1,ega2,ega31,ega42,ega44} and SGA 3 \cite{sga3-2,sga3-4,sga3-5,sga3-6b, sga3-11,sga3-12, sga3-14,sga3-19,sga3-22,sga3-24,sga3-26}. For the general theory of tdos, we content us mostly with \cite{kashiwara1989,beilinsonbernstein1993}. For the functorial study of twisted $\mathcal D$-modules, we basically follow the approaches of \cite{boreletal, milicic, hottaetal2008, hechtetal, kashiwaratanisaki}. For the derived operations, we also refer to \cite{lipman2009} in order to eliminate Noetherian and other finiteness assumptions.

	{\em Acknowledgments.} The first author is grateful to Yoshiki Oshima for valuable comments. He also thanks Hiroki Kato, Masatoshi Kitagawa, and Daichi Takeuchi for helpful discussions. The second author thanks Haruzo Hida for insightful discussions and Leticia Barchini, Michael Harris and Guido Kings for helpful remarks on preliminary drafts of this manuscript. T. H. was supported by JSPS KAKENHI Grant Numbers JP21J00023 and JP22KJ2045. F. J. was supported by the Deutsche Forschungsgemeinschaft (DFG, German Research Foundation) -- SFB-TRR 358/1 2023 -- 491392403.
	
	\section*{Notation and Convention}
	
	\subsection*{Elementary concepts}
	Let $\ZZ$ (resp.\ $\NN$, $\QQ$, $\RR$, $\CC$) be the set of integers (resp.\ nonnegative integers, rational numbers, real numbers, complex numbers). For a prime $p$, $\FF_p$ denotes the finite field of order $p$.
	
	For a positive integer $n$, let $\lieS_n$ denote the $n$-th symmetric group.
	
	For a topological space $X$, write $\ZZ_X$ for the constant sheaf on $X$ attached to $\ZZ$.
	
	For a commutative ring $R$ (resp.\ a sheaf $\Rr$ of commutative rings), let $R^\times$ (resp.\ $\Rr^\times$) denote the group (resp.\ the group sheaf) of units of $R$ (resp.\ $\Rr$).
	
	For a finite dimensional vector space $V$ over a field $F$, we denote its dimension by $\dim_F V$.

	For an integral domain $k$, we denote the field of fractions of $k$ by $\Frac(k)$.
	
	In this paper, every filtration of an abelian group or a sheaf of abelian groups is assumed to be positive and increasing. That is, $F_{-1}A=0$ and $F_iA\subset F_{i+1}A$ for every filtration $F_\bullet A$ on an object $A$ and $i\in\ZZ$. For a filtered object $F_\bullet A$, its associated graded object will be denoted by $\Gr^\bullet A$.
	
	\subsection*{Categories}
	
	For an object $X$ of a category, let $\id_X$ denote its identity map. We may omit $X$ if it is obvious by the context. For objects $X,Y$ of a category, we denote the projection morphisms $X\times Y\to X$ and $X\times Y\to Y$ by $\pr_1$ and $\pr_2$ respectively (if the product exists).
	
	For an additive functor $F$ of abelian categories, we denote the attached functor between the homotopy categories of complexes of their objects by the same symbol $F$.

	We try to explain our notations on modules and sheaves simultaneously here. To distinguish the setting of sheaves from that of rings, we will use calligraphy letters for sheaves in applications.
	
	Let $\Vv$ be a Grothendieck abelian closed symmetric monoidal category. If necessary, see \cite[section 1.3]{hayashi2019} for basic definitions concerning symmetric monoidal categories. We say that a monoidal category is abelian if it is an abelian category and the monoidal structure $-\otimes-$ is additive in each variable. For the definition of Grothendieck abelian categories, see \cite{tohoku}. The basic example will be the category of abelian groups (resp.~sheaves of abelian groups on a topological space $X$) with the monoidal structures $\otimes_{\ZZ}$ (resp.~$\otimes_{\ZZ_X}$) and the usual symmetry constraint. Our main example will be the category of modules over a commutative algebra object $R$ in $\Vv$ with $-\otimes_R-$.
	
	For a pair $X,Y$ of objects of $\Vv$, we denote its symmetry constraint $X\otimes Y\cong Y\otimes X$ by $C_{X,Y}$.
	
	We denote the monoidal dual by $(-)^\vee$. That is, for $X\in\Vv$, set $X^\vee=\Map(X,I)$, where $\Map(-,-)$ is the internal Hom, and $I$ is the unit object of $\Vv$.

	Let $A$ be an algebra object of $\Vv$. We denote the multiplication morphism $A\otimes A\to A$ of $A$ by $m_A$. We put a Lie bracket $\left[-,-\right]$ on $A$ by the commutator unless mentioned otherwise. Let $A^{\op}$ denote the opposite monoid object. For example, if $\Vv$ is the category of modules over a commutative ring $R$, $A^{\op}$ is the opposite $R$-algebra to a given $R$-algebra $A$.
	
	We write $\Mod(A)$ (resp.\ $\Mod_{\mathrm{r}}(A)$) for the category of left (resp.\ right) $A$-modules. For a pair $X,Y$ of objects of $\Mod(A)$ or $\Mod_{\mathrm{r}}(A)$, we denote the set of morphisms from $X$ to $Y$ by $\Hom_{A}(X,Y)$. If $X=Y$, write $\End_{A}(X)=\Hom_{A}(X,X)$.
	
	Let $K(A)$ and $K_{\mathrm{r}}(A)$ be the homotopy categories of complexes of left and right $A$-modules respectively. We denote the unbounded derived category of left (resp.\ right) $A$-modules by $D(A)$ (resp.\ $D_{\mathrm{r}}(A)$). We remark that these derived categories are essentially locally small since $\Mod(\Aa)$ and $\Mod_{\mathrm{r}}(\Aa)$ are Grothendieck abelian (see \cite[section 1.5]{tohoku} if necessary). For each integer $p$, let $D^{\geq p}(A)$ (resp.\ $D^{\leq p}(A)$) be the full subcategory of $D(A)$ consisting of complexes cohomologically concentrated in degree $\geq p$ (resp.\ $\leq p$). Set $D^+(A)=\cup_{p\in\ZZ} D^{\geq p}(A)$ and $D^-(A)=\cup_{p\in\ZZ} D^{\leq p}(A)$. Define $D^{\geq p}_{\mathrm{r}}(A)$, $D^{\leq p}_{\mathrm{r}}(A)$, $D^+_{\mathrm{r}}(A)$, and $D^-_{\mathrm{r}}(A)$ in a similar way. We will call isomorphisms in derived categories equivalences in this paper for perspectives of homotopical algebra.
	
	Let $R$ be a commutative algebra object of $\Vv$. For an $R$-module $M$, let $\Sym_R M$ denote its symmetric product. For a nonnegative integer $n$, write $\Sym^n_R M$ for the degree $n$ part of $\Sym_R M$.

	\subsection*{Derivation}

	Let $k\to R$ be a homomorphism of commutative rings. For an $R$-module $M$, let $\mathrm{Der}_{k}(R,M)$ denote the $R$-module of $k$-derivations of $M$. Write $\Theta_{R/k}=\mathrm{Der}_k(R,R)$. More notations on derivations and related concepts (e.g.~differentials) are collected in Appendix \ref{sec:derivation}.
	
	\subsection*{Topological spaces}
	
	For a continuous map $f:X\to Y$ of topological spaces and an open subset $U\subset Y$, let $f_U:f^{-1}(U)\to U$ be the restriction of $f$.
	
	\subsection*{Sheaves}

	Let $\Ff$ be a sheaf on a topological space $X$, and $U$ be an open subset of $X$. We write $\Ff(U)= \Gamma(U,\Ff)$. We will refer to an element $a\in\Ff(U)$ as a local section $a\in\Ff$ when we do not intend to emphasize $U$. 
	
	Let $X$ be a topological space, and $\Aa$ be a sheaf of rings on $X$. For left or right $\Aa$-modules $\Mm$ and $\mathcal N$, let $\iHom_{\Aa}(\Mm,\mathcal N)$ denote the sheaf of $\Aa$-linear homomorphisms from $\Mm$ to $\Nn$. For a closed immersion $Z\hookrightarrow X$ and a left or right $\Aa$-module $\Ff$ on $X$, let $\Gamma_Z\Ff\subset\Ff$ denote the subsheaf of local sections of $\Ff$ supported in $Z$ (under the identification of $Z$ with a closed subset of $X$). It is evident that $\Gamma_Z\Ff$ is an $\Aa$-submodule of $\Ff$.
	
	For a ringed space $X$, let $\OO_X$ denote the structure sheaf of $X$ unless specified otherwise. For a morphism
	$f:(X,\OO_X)\to (S,\OO_S)$
	of ringed spaces, let its structure homomorphism $f^{-1}\OO_S\to\OO_X$ by $f^\sharp$.
	
	\subsection*{Schemes}
	Let $(X,\OO_X)$ be a scheme, and $\Aa$ be a sheaf of rings on $X$, equipped with a homomorphism $\OO_X\to \Aa$. Then we denote the full subcategory of $\Mod(\Aa)$ (resp.\ $\Mod_{\mathrm{r}}(\Aa)$) consisting of left (resp.\ right) $\Aa$-modules which are quasi-coherent as left (resp.\ right) $\OO_X$-modules by $\Mod_{\qc}(\Aa)$ (resp.\ $\Mod_{{\mathrm{r}},\qc}(\Aa)$). Similarly, for $\bullet\in\{\emptyset,+,-\}$, write $D^\bullet_{\qc}(\Aa)\subset D^\bullet(\Aa)$ and $D^\bullet_{{\mathrm{r}},\qc}(\Aa)\subset D^\bullet_{\mathrm{r}}(\Aa)$ for the full subcategories of complexes whose cohomology degreewisely lies in $\Mod_{\qc}(\Aa)$ and $\Mod_{{\mathrm{r}},\qc}(\Aa)$ respectively. Although these notations may be confusing in general since one can define the notion of quasi-coherent $\Aa$-modules without $\OO_X$ in general (see Definition \ref{defn:qcoh}), we will see that these two notions agree at least in our applications (Proposition \ref{prop:defnofqcoh}, Lemma \ref{lem:qcohmodoverqcohalg}).
	
	For a nonnegative integer $n$ and a scheme $S$, let $\bfA^n_S$ and $\PP^n_S$ denote the affine and projective $n$-spaces over $S$ respectively.
	
	Let $S$ be a scheme. Write $\Sch_S$ for the category of schemes over $S$. We call a presheaf on $\Sch_S$ an $S$-space. If a test $S$-scheme is $T=\Spec R$ for some commutative ring $R$ then for an $S$-space $X$, we sometimes denote $X(T)=X(R)$. For a set $C$, let $C_S$ denote the attached constant sheaf on $\Sch_S$ in the big \'etale topology.
	For an $S$-scheme $X$ and a geometric point $\bar{s}$ of $S$, i.e., a morphism to $S$ from the spectrum of an algebraically closed field, let $X_{\bar{s}}$ be the pullback of $X$ by $\bar{s}$. We call $X_{\bar{s}}$ the geometric fiber of $X$ at $\bar{s}$.
	
	Let $S'\to S$ be a morphism of schemes, and $X'$ be an $S'$-space. Let $\res_{S'/S}X$ denote the Weil restriction of $X'$, i.e., the $S$-space defined by
	$(\res_{S'/S}X')(T)=X'(S'\times_S T)$
	for $S$-schemes $T$. If $S=\Spec k$ and $S'=\Spec k'$ for some commutative rings $k,k'$, we sometimes denote $\res_{S'/S}=\res_{k'/k}$. For an automorphism $\sigma$ of $S'$ over $S$, define an auto-functor
	$(-)^\sigma:\Sch_{S'}\to\Sch_{S'}$
	as the base change by $\sigma$.
	
	\subsection*{Group schemes}
	For a scheme $S$, let $\GG_{m,S}$ denote the multiplicative group scheme over $S$. 
	
	For a group scheme $G$ over $S$, write $X^\ast(G)$ for the character group of $G$, i.e., the group of homomorphisms $G\to \GG_{m,S}$. If we are given a split torus $H$ over a connected scheme $S$, then we denote the cocharacter group of $H$ by $X_\ast(H)$.
	
	For an $S$-affine smooth group scheme over $S$, we denote its Lie algebra (see \cite[sections 2-4]{sga3-2}) by the corresponding German letter unless specified otherwise. For example, the Lie algebra of an $S$-affine smooth group scheme $G$ over $S$ will be expressed as $\lieg$.
	
	For a scheme $X$ over $S$, equipped with an action of a group scheme $H$ over $S$, $H\backslash X$ will denote the fppf quotient sheaf (see \cite{sga3-4} for the general formalism).
	
	For a smooth group scheme $G$ over a scheme $S$, we denote its unit component in the sense of \cite[D\'efinition 3.1]{sga3-6b} by $G^0$. It is an open subscheme of $G$ in general (\cite{sga3-6b} Th\'eor\`eme 3.10). If $G$ is separated of finite presentation and $G^0$ is reductive in the sense of \cite[D\'efinition 2.7]{sga3-19}, set $\pi_0(G)\coloneqq G/G^0$. This is separated \'etale $S$-scheme of finite presentation (\cite[Proposition 3.1.3]{conrad2014}).
	
	For a reductive group scheme $G$ over $S$, we denote the moduli scheme of Borel subgroups of $G$ by $\Bb_G$ and call it the flag scheme of $G$ (\cite[Corollaire 5.8.3]{sga3-22}). Let $\Pp_G$ denote the moduli scheme of parabolic subgroups of $G$ (\cite[Section 3.2, Th\'eor\`eme 3.3]{sga3-26}). Let $\type G$ be the scheme of parabolic types of $G$ which may be identified with the quotient sheaf $G\backslash\Pp_G$ (\cite[Section 3.2, D\'efinition 3.4]{sga3-26}). Let $t:\Pp_G\to \type G$ denote the canonical quotient map. For a torus $T$ in $G$ (\cite[D\'efinition 1.3]{sga3-9}, cf.~\cite[Corollaire 4.5]{sga3-10}), let $Z_G(T)$ denote the centralizer of $T$ in $G$ (\cite[D\'efinition 6.1]{sga3-6b}, \cite[Corollaire 5.3]{sga3-11}). For a maximal torus $H$ of $G$ (\cite[D\'efinition 1.3]{sga3-12}), let $N_G(H)$ be the normalizer of $H$ in $G$ (\cite[D\'efinition 6.1]{sga3-6b}, \cite[Proposition 2.1.2]{conrad2014}). Define the Weyl group scheme as $W(G,H)\coloneqq N_G(H)/H$ (\cite[Corollaire 5.4]{sga3-12}, \cite[Proposition 3.2.8]{conrad2014}).
	
	Let $S$ be an affine scheme, and $G$ be a reductive group scheme. Then for a cocharacter $\mu$ of $G$, write $P_G(\mu)$ for the parabolic subgroup of $G$ defined by
	\[P_G(\mu)(T)=\{g\in G(T):\ \lim_{a\to 0}\mu(a)g\mu(a)^{-1}\ {\rm exists}\}\]
	for affine $S$-schemes $T$ (see \cite[sections 2.1 and 2.2]{conradpseudo} for the details).

	\subsection*{Lie groups}
	For a real Lie group $G$, its unit component will be denoted by $G_0$. Write $\pi_0(G)=G/G_0$ for the component group of $G$.
	
	\subsection*{Algebraic number theory}
	For a number field $E$, write $\OO_E$ for the ring of its integers. We denote the absolute discriminant of $E/\QQ$ by $d_E$.
	
	\subsection*{Representations}
	
	Let $K$ be an affine group scheme over a commutative ring $k$. Then for representations $V,W$ of $K$, let $\Hom_K(V,W)$ be the $k$-module of homomorphisms from $V$ to $W$ as representations of $K$.
	
	Suppose that we are given a Lie algebra $\lies$ and an affine group scheme $K$ over a commutative ring $k$. Let $V,W$ be $k$-modules, equipped with the (not necessarily compatible) structures of representations of $\lies$ and $K$. Then we denote the $k$-module of $k$-linear maps from $V$ to $W$ which intertwine the actions of $\lies$ and $K$ by $\Hom_{\lies,K}(V,W)$.
	
	\section{Algebras of twisted differential operators and Picard algebroids}\label{sec:tdo_and_pic_alg}
	
	Twisted $\D$-modules are modules over tdos (sheaves of twisted differential operators). In this section, we see that the theory of tdos and Picard algebroids in \cite{kashiwara1989,beilinsonbernstein1993} works over general base schemes possibly after minor changes. In particular, we impose perspectives from Picard algebroids in their pullback towards the construction of the left action on the transfer bimodule (section \ref{sec:bimodules}). As a new aspect, we introduce the base change functor.
	
	We would like to work with the settings of rings and sheaves simultaneously to prove that these are equivalent for affine schemes. This is nontrivial since tdos on a smooth scheme $X$ over a base scheme are not algebras over $\OO_X$.
	
	Throughout this section, $s$ is either a homomorphism $k\to k'$ of commutative rings or a morphism $S'\to S$ of schemes unless specified otherwise. Let $f$ be either a homomorphism $R\to R'$ of smooth $k$-algebras or a morphism $X\to Y$ of smooth $S$-schemes. Write $x$ and $y$ for the structure homomorphisms $k\to R$ and $k\to R'$ (resp.\ the structure morphisms $X\to S$ and $Y\to S$) respectively. Set $X'=X\times_S S'$ and $Y'=Y\times_S S'$. Write $f'$, $x'$, and $y'$ for the base changes of $f$, $x$, and $y$ by $s$ respectively. Let $s_R:R\to k'\otimes_k R$ and $s_{R'}:R'\to R'\otimes_k k'$ be the canonical homomorphisms. Similarly, let $s_X:X'\to X$ and $s_Y:Y'\to Y$ denote the canonical projections.
	
	\subsection{Twisted differential operators}\label{sec:tdo}
	\begin{lemma}\label{lem:commutator}
		Let $A$ be an almost commutative filtered $k$-algebra, and $n\in\NN$.
		\begin{enumerate}
			\renewcommand{\labelenumi}{(\arabic{enumi})}
			\item $F_0A$ is a commutative $k$-subalgebra of $A$.
			\item For $P\in F_n A$,
			$\left[P,-\right]$ determines an element of $\mathrm{Der}_k(F_0A,\Gr^{n-1} A)$.
			\item The map $\sigma:F_n A\to
			\mathrm{Der}_k(F_0A,\Gr^{n-1}A)$ determined by (2) is $F_0A$-linear.
			\item The map $\sigma$ descends to $\sigma:\Gr^n A\to \mathrm{Der}_k(F_0A,\Gr^{n-1}A)$.
			\item The map $\sigma: F_1A\to \Theta_{F_0A/k}$ is a Lie algebra homomorphism. Here $F_1 A$ is a Lie algebra for the commutator (use the assumption that $A$ is almost commutative).
		\end{enumerate}
	\end{lemma}
	
	\begin{proof}
		Part (1) is clear. We prove (2). Since $A$ is almost commutative, $\left[P,a\right]$ belongs to $F_{n-1} A$ for every $a\in F_0A$ and $P\in F_nA$. Apply the quotient map $F_{n-1} A\to \Gr^{n-1} A$ to get a $k$-linear map $F_0A\to \Gr^{n-1} A$. For any $a,a'\in F_0A$ we have
		\[\begin{split}
			\left[P,aa'\right]-a\left[P,a'\right]-a'\left[P,a\right]
			&=Paa'-aa'P-aPa'+aa'P-a'Pa+a'aP\\
			&=\left[\left[P,a\right],a'\right]\\
			&\equiv 0 \pmod{F_{n-2}A}.
		\end{split}\]
		
		It is clear that $\sigma$ is additive. Part (3) follows from
		\[\sigma(aP)(a')=aPa'-a'aP=(a\sigma(P))(a')\]
		for $P\in F_n A$ and $a,a'\in F_0A$. If $P\in F_{n-1}A$ and $a\in F_0A$, we have $\left[P,a\right]=0$ in $\Gr^{n-1}A$. This shows (4).

		To prove (5), pick $e,e'\in F_1A$, and $a\in F_0A$. Then we get
		\[\begin{split}
			\left[\sigma(e),\sigma(e')\right](a)
			&=\sigma(e)\sigma(e')(a)-\sigma(e')\sigma(e)(a)\\
			&=\left[e,\left[e',a\right]\right]
			-\left[e',\left[e,a\right]\right]\\
			&=\sigma(\left[e,e'\right])(a).
		\end{split}\]
		This completes the proof.
	\end{proof}
	
	Likewise, we obtain
	
	\begin{lemma}\label{lem:commutator/scheme}
		Let $\Aa$ be an almost commutative filtered $x^{-1}\OO_S$-algebra, and $n\in\NN$.
		\begin{enumerate}
			\renewcommand{\labelenumi}{(\arabic{enumi})}
			\item $F_0\Aa$ is a commutative $x^{-1}\OO_S$-subalgebra of $\Aa$.
			\item For an open subset $U$ and $P\in F_n \Aa(U)$,
			$\left[P,-\right]$ determines an element of
			\[\Der_S(F_0 \Aa,\Gr^{n-1} \Aa)(U).\]
			\item The map $\sigma:F_n \Aa\to
			\Der_S(F_0\Aa,\Gr^{n-1}\Aa)$ determined by (2) is $F_0\Aa$-linear.
			\item The map $\sigma$ descends to $\sigma:\Gr^n \Aa\to \Der_S(F_0\Aa,\Gr^{n-1}\Aa)$.
			\item The map $\sigma: F_1\Aa\to \Theta_{F_0\Aa/x^{-1}\OO_S}$ is a Lie algebra homomorphism, where $F_1\Aa$ is a Lie algebra for the commutator.
		\end{enumerate}
	\end{lemma}
	
	\begin{definition}[Algebras of twisted differential operators, tdos]\label{defn:tdo}
		\begin{enumerate}
			\renewcommand{\labelenumi}{(\arabic{enumi})}
			\item An almost commutative filtered $k$-algebra $(A,F_\bullet A)$, equipped with a $k$-algebra homomorphism $i:R\to A$ is called a tdo for $R/k$ if the following conditions are satisfied:
			\begin{enumerate}
				\item[(i)] $F_\bullet A$ is exhaustive;
				\item[(ii)] the map $i$ is an isomorphism onto $F_0A$;
				\item[(iii)] the canonical $R$-algebra homomorphism $\Sym_R\Gr^1 A\to\Gr A$ is an isomorphism; 
				\item[(iv)] the map
				$\sigma:\Gr^1 A\to\Theta_{R/k};~e\mapsto \left[e,-\right]$ is an isomorphism.
			\end{enumerate}
			\item A homomorphism $(A,F_\bullet,i)\to (A',F_\bullet,i')$ of tdos is a filter-preserving $k$-algebra homomorphism $\varphi:A\to A'$ satisfying $\varphi\circ i =i'$.
			\item We denote the category of tdos by $\TDO_{R/k}$.
		\end{enumerate}
	\end{definition}
	
	We readily generalize the definition of tdos in \cite{kashiwara1989}:
	
	\begin{definition}[{\cite[Definition 2.3.3]{kashiwara1989}, \cite[Definition 2.1.1]{beilinsonbernstein1993}}]
		\begin{enumerate}
			\renewcommand{\labelenumi}{(\arabic{enumi})}
			\item An almost commutative filtered $x^{-1}\OO_S$-algebra sheaf $\Aa$, equipped with an $x^{-1}\OO_S$-algebra homomorphism $i:\OO_X\to\Aa$ is called a tdo for $X/S$ or simply a tdo on $X$ if the following conditions are satisfied:
			\begin{enumerate}
				\item[(i)] $i$ is an isomorphism onto $F_0\Aa$;
				\item[(ii)] the filtration $F_\bullet \Aa$ on $\Aa$ is exhaustive;
				\item[(iii)] the canonical $\OO_X$-algebra homomorphism $\Sym_{\OO_X}\Gr^1 \Aa\to\Gr \Aa$ is an isomorphism;
				\item[(iv)] the map $\sigma:\Gr^1 \Aa\to\Theta_{X/S};~e\mapsto \left[e,-\right]$ is an isomorphism.
			\end{enumerate}
			\item A homomorphism $(\Aa,F_\bullet,i)\to (\Aa',F_\bullet,i')$ of tdos is a filter-preserving $x^{-1}\OO_S$-algebra homomorphism $\varphi:\Aa\to \Aa'$ such that $\varphi\circ i =i'$.
			\item We denote the category of tdos by $\TDO_{X/S}$.
		\end{enumerate}
	\end{definition}
	
	\begin{remark}\label{rem:tdo_is_groupoid}
		The conditions (iii) and (iv) imply that $\TDO_{R/k}$ and $\TDO_{X/S}$ are groupoids.
	\end{remark}
	
	The following observations are immediate consequences of the definitions.
	
	\begin{lemma}\label{lem:tdozariski}
		Let $x:X\to S$ be a smooth morphism of schemes.
		\begin{enumerate}
			\renewcommand{\labelenumi}{(\arabic{enumi})}
			\item For an open subscheme $U\subset X$, the usual restriction determines a functor
			\[(-)|_U:\TDO_{X/S}\to\TDO_{U/S}.\]
			\item If $f$ factors thorugh an open subscheme $V\subset S$, we have a canonical categorical isomorphism $\TDO_{X/S}\cong\TDO_{X/V}$.
			\item Tdos form a stack over the lisse-Zariski site of $S$.
		\end{enumerate}
	\end{lemma}
	
	\begin{example}
		Put $S=\Spec F$, where $F$ is a field of characteristic zero. Then our definition of tdos agrees with \cite[Definition 2.3.3]{kashiwara1989}.
	\end{example}
	
	\begin{example}[{\cite[Th\'eor\`eme (16.11.2)]{ega44}}]\label{ex:egadiffop}
		Let $X$ be a smooth scheme over a base $S$ on which every nonzero integer $n$ is invertible. Then Grothendieck's sheaf of differential operators (\cite{ega44} Corollaire (16.8.10)) is a tdo.
	\end{example}
	
	\begin{example}[{\cite[Example 4.5, Proposition 4.6, Corollary 4.7]{berthelotogus}}]\label{ex:pddiffop}
		Suppose that there is a positive integer $m$ such that $m\OO_X$ is zero. Then the sheaf of PD differential operators on $X$ is a tdo.
	\end{example}

	\begin{example}\label{ex:optdo}
		Let $\Aa$ be a tdo on $X$. Then its opposite algebra sheaf $\Aa^{\op}$ is a tdo on $X$. Similarly, if we are given a tdo $A$ for $R/k$, $A^{\op}$ is a tdo for $R/k$. 
	\end{example}
	
	\begin{example}
		Let $\Aa$ be a tdo on $X$, and $\Ll$ be a line bundle on $X$. Then $\Ll^\vee\otimes_{\OO_X} \Aa\otimes_{\OO_X}\Ll$ is a tdo since the conditions are Zariski local in $X$.
	\end{example}

	\begin{lemma}\label{lem:Dquasicoherence}
		Let $X\to S$ be a smooth morphism of schemes, and $\Aa$ be a tdo on $X$. Then $F_n\Aa$ and $\Aa$ are quasi-coherent sheaves on $X$ for the multiplication of $\OO_X$ from the left side ($n\geq 0$).
	\end{lemma}
	
	\begin{proof}
		We may prove $F_n\Aa$ is quasi-coherent by induction on $n$. The assertion is clear if $n=0$ since $\OO_X\cong F_0\Aa$. For $n\geq 1$, recall that we have an exact sequence
		\[0\to F_{n-1}\Aa\to F_n\Aa\to \Gr^n \Aa\to 0.\]
		We have an isomorphism
		$\Gr^n\Aa\cong\Sym^n_{\OO_X}\Theta_{X/S}$. Combine it with the induction hypothesis to deduce that $F_n\Aa$ is quasi-coherent.
	\end{proof}
	
	\begin{corollary}
		The global sections determine a functor
		\[\Gamma(X,-):\TDO_{\Spec R/\Spec k}\to\TDO_{R/k}.\]
	\end{corollary}
	
	\begin{proposition}\label{prop:tdofreebasisaffine}
		Let $(A,F_\bullet A,i)$ be a tdo. Regard $A$ as an $R$-module for the multiplication from the left side. Suppose that $\Theta_{R/k}$ is free as $R$-module. Fix elements $e_1,\ldots,e_n\in F_1 A$ which form a free $R$-basis of $\Gr^1 A$.
		Then
		$\{e_1^{i_1}\cdots e_n^{i_n}\in D:~i_1,\ldots,i_n\geq 0\}$ is a free basis of $A$ as an $R$-module. Furthermore, each filtration step $F_dA$ is free of finite rank for any $d\geq 0$.
	\end{proposition}
	
	\begin{proof}
		We consider the free $R$-module
		$\tilde{A}=\oplus_{(i_1,\ldots,i_n)\in\NN} R e_1^{i_1}\cdots e_n^{i_n}$ with the grading by total degree
		$\deg e_1^{i_1}\cdots e_n^{i_n}=\sum_{j=1}^n i_j$. Then $\tilde{A}$ comes with a natural filtration. The canonical map $\tilde{A}\to A$ respects the filtration. Moreover, the induced map $\Gr \tilde{A}\to \Gr A$ is an isomorphism. Hence $\tilde{A}\to A$ is an isomorphism of filtered $R$-modules, which proves the claim.
	\end{proof}
	
	\begin{proposition}\label{prop:tdofreebasis}
		Let $(\Aa,F_\bullet \Aa,i)$ be a tdo on $X$. Then $\Aa$ is a locally free $\OO_X$-module for the multiplication from the left side. Furthermore, each filtration step $F_d\Aa$ is locally free of finite rank for any $d\geq 0$.
	\end{proposition}
	
	In their theory over fields of characteristic zero, Kashiwara, Beilinson, and Bernstein could specify filtrations on tdos.
	
	\begin{definition}[{\cite[section 1.1]{beilinsonbernstein1993}}]\label{defn:D-filtration}
		Let $A$ be a $k$-algebra, equipped with a $k$-algebra homomorphism $i:R\to A$.
		\begin{enumerate}
			\renewcommand{\labelenumi}{(\arabic{enumi})}
			\item For each integer $n\geq 0$, let $F^{\vee}_n A\subset A$ be the $k$-submodule consisting of $P\in A$ with the following property: For all $a_1,a_2,\ldots,a_{n+1}\in R$,
			\[\left[\cdots\left[\left[P,i(a_1)\right],i(a_2)\right],\ldots,i(a_{n+1})\right]=0.\]
			We call $F^{\vee}$ the maximal $D$-filtration on $A$.
			\item We say $A$ is a $D$-algebra (relative to $x$) if the filtration $F^{\vee}_\bullet A$ is exhaustive, i.e., for every element $P\in A$, there exists an integer $n\geq 0$ such that for all $a_1,a_2,\ldots,a_{n+1}\in R$,
			\[\left[\cdots\left[\left[P,a_1\right],a_2\right],\ldots,a_{n+1}\right]=0.\]
		\end{enumerate}
	\end{definition}
	
	They proved that the maximal $D$-filtration is the unique filtration which turns a $k$-algebra $A$ with a $k$-algebra homomorphism $i:R\to A$ into a tdo, provided $k$ is a field of characteristic zero (\cite[Proposition 2.3.2]{kashiwara1989}, \cite[2.1.1.~Remark]{beilinsonbernstein1993}, see also Corollary \ref{cor:tdofilt=maxDfilt}). Moreover, they gave a nice criterion for tdos in terms of the maximal $D$-filtrations (\cite[Proposition 2.3.2]{kashiwara1989}, \cite[2.1.2.~Lemma]{beilinsonbernstein1993}, see also Corollary \ref{cor:tdobymdfilt}). Hence their notion of $D$-algebras is a good approximation of tdos. This idea partly extends to general bases and yields
	
	\begin{lemma}\label{lem:F<Fvee}
		Let $(A,F_\bullet A,i)$ be a tdo. Then for each $n\geq 0$, we have $F_n A\subset F_n^{\vee} A$.
	\end{lemma}
	
	\begin{proof}
		We proceed by induction on $n$. Notice that by definition, $F_0^{\vee} A$ is the centralizer of $i(R)$. Since $R$ is commutative, $F_0D=i(R)\subset F_0^{\vee} A$. Let $n\geq 1$ and $P\in F_n A$. Take an arbitrary element $a\in R$. Since $A$ is almost commutative, $\left[P,i(a)\right]$ belongs to $F_{n-1} A$. The induction hypothesis shows $\left[P,i(a)\right]
		\in F_{n-1}A\subset F_{n-1}^{\vee}A$. Since $a$ is arbitrary, $P\in F_n^{\vee} A$ follows as claimed.
	\end{proof}
	
	However, in general the filtration of a tdo $(A,F_\bullet A,i)$ is not the maximal $D$-filtration when working over general bases.
	
	\begin{example}
		Let $(A,F_\bullet A,i)$ be a tdo, and $e\in F_1A$. Then we have
		\[\begin{split}
			\left[e^2,i(a)\right]
			&=e^2i(a)-i(a)e^2\\
			&=e\left[e,i(a)\right]+ei(a)e-i(a)e^2\\
			&=\left[e,\left[e,i(a)\right]\right]
			+\left[e,i(a)\right]e
			+\left[e,i(a)\right]e+i(a)e^2-i(a)e^2\\
			&=i(\sigma(e)^2(a))+2\left[e,i(a)\right]e
		\end{split}\]
		for all $a\in R$.
		
		Put $k=\FF_2$ and $R=\FF_2\left[t\right]$. Let $A$ be an arbitrary tdo for $\FF_2\left[t\right]/\FF_2$. For example, $A$ can be the ring of PD differential operators on $\FF_2\left[t\right]$ (Example \ref{ex:pddiffop}). In this case, 
		$A$ is equal to $\oplus_{n\geq 0} \FF_2\left[t\right]\frac{d^n}{dt^n}$
		with the usual filtered $\FF_2$-algebra structure. Since $\sigma:\Gr^1 A\to \Theta_{R/k}$ is an isomorphism, one can find $e\in F_1A\setminus F_0A$ such that
		$\sigma(e)=\frac{d}{dt}$. For such $e$, $e^2$ belongs to $F_0^{\vee}A$, because by the above computation,
		\[\left[e^2,i(a)\right]
		=i\left(\frac{d^2a}{dt^2}\right)=0.\]
		Hence $F_0^{\vee} A$ is strictly bigger than $F_0 A$. Moreover, we claim that $e$ is central in $A$. For this, observe that $A$ is generated by $F_1A$ since $\Gr A\cong \Sym_R \Theta_{R/k}$. Hence it suffices to show that $\left[e^2,-\right]=0$ on $F_1A$. Now the isomorphism
		$\Gr^1 A\cong\Theta_{R/k}=
		\FF_2\left[t\right]\frac{d}{dt}$
		shows that every element of $F_1 A$ can be expressed as $i(a)e+i(b)$ for some $a,b\in R$. Therefore, the assertion follows with
		\[\left[e^2,i(a)e+i(b)\right]
		=e^2i(a)e-i(a)e^3
		=\left[e^2,i(a)\right]e+i(a)e^3-i(a)e^3=0.\]
	\end{example}
	
	We will discuss in section \ref{sec:equivtdo} how we can generalize Beilinson and Bernstein's characterization for $A$ flat over $\ZZ$ (cf.\ Corollary \ref{cor:tdobymdfilt} and Corollary \ref{cor:tdofilt=maxDfilt}).

	\subsection{Picard algebroids and $\Omega^{\geq 1}$-torsors}\label{sec:picalg}
	
	In this section, we study Picard algebroids and $\Omega^{\geq 1}$-torsors, key notions for the construction of tdos.
	
	\begin{definition}\label{defn:picalg}
		\begin{enumerate}
			\renewcommand{\labelenumi}{(\arabic{enumi})}
			\item A Picard algebra for $R/k$ or simply a Picard algebra over $R$ is an $R$-module $E$, equipped with an $R$-module homomorphism $\sigma:E\to\Theta_{R/k}$, a $k$-module homomorphism
			$\left[-,-\right]=\left[-,-\right]_E:E\otimes_k E\to E$,
			and an element $1_E\in E$ such that
			\begin{enumerate}
				\item[(i)] $E$ is a Lie algebra over $k$ for $\left[-,-\right]$,
				\item[(ii)] for $e\in E$, $\left[1_E,e\right]=0$,
				\item[(iii)] $\sigma$ is a Lie algebra homomorphism,
				\item[(iv)] we have $\left[e_1,ae_2\right]=a\left[e_1,e_2\right]+\sigma(e_1)(a)e_2$ for $a\in R$ and $e_1,e_2\in E$, and that
				\item[(v)] the sequence
				$
				0\to R\overset{1_E}{\to}E
				\overset{\sigma}{\to}\Theta_{R/k}\to 0
				\label{eq:defnexactseqofpicalg}
				$
				is exact.
			\end{enumerate}
			We will call $1_E$ the unit of $E$. We will put a filtration on a Picard algebra $E$ by $F_0E=R1_E$ and $F_1E=E$.
			\item A homomorphism $(E,\sigma,1_E)\to (E',\sigma',1_{E'})$ is an $R$-linear Lie algebra homomorphism $\varphi:E\to E'$ over $k$ such that $\varphi(1_E)=1_{E'}$.
			\item We denote the category of Picard algebras by $\PA_{A/k}$.
		\end{enumerate}
	\end{definition}
	
	\begin{example}[Elementary twists of Picard algebras: $E^\gamma$]
		Let $(E,[-,-],\sigma,1_E)$ be a Picard algebra. Then for any $\gamma\in k^\times$, $(E,\sigma,[-,-],\gamma^{-1} 1_E)$ is a Picard algebra which we denote by $E^\gamma$. This construction is natural in the following sense: For a homomorphism $\varphi:E_1\to E_2$ of Picard algebras, the map $\varphi$ is also a homomorphism of Picard algebras from $E_1^\gamma$ to $E_2^\gamma$.
	\end{example}
	
	Likewise, we formulate
	\begin{definition}[{\cite[Definition 2.1.3]{beilinsonbernstein1993}}]\label{defn:picalg/sch}
		\begin{enumerate}
			\renewcommand{\labelenumi}{(\arabic{enumi})}
			\item A Picard algebroid on $X$ is an $\OO_X$-module $\Ee$, equipped with an $\OO_X$-module homomorphism $\sigma:\Ee\to\Theta_{X/S}$, an $x^{-1}\OO_S$-module homomorphism
			$\left[-,-\right]=\left[-,-\right]_{\Ee}:\Ee\otimes_{x^{-1}\OO_S} \Ee\to \Ee$,
			and $1_{\Ee}\in \Ee(X)$ such that
			\begin{enumerate}
				\item[(i)] $\Ee$ is a Lie algebra over $x^{-1}\OO_S$ for $\left[-,-\right]$,
				\item[(ii)] for any open subset $U$, any $e\in \Ee(U)$, $\left[1_{\Ee}|_U,e\right]=0$,
				\item[(iii)] $\sigma$ is an $x^{-1}\OO_S$-Lie algebra homomorphism,
				\item[(iv)] we have $\left[e_1,ae_2\right]=a\left[e_1,e_2\right]+\sigma(e_1)(a)e_2$ for any open subset $U\subset X$, $a\in \OO_X(U)$ and $e_1,e_2\in \Ee(U)$,
				\item[(v)] the sequence
				$0\to \OO_X\overset{1_{\Ee}}{\to}\Ee
				\overset{\sigma}{\to}\Theta_{X/S}\to 0$
				is exact.
			\end{enumerate}
			We call $1_{\Ee}$ the unit of $\Ee$ and on any Picard algebroid $\Ee$ we have a finite filtration with $F_0\Ee=\OO_X1_{\Ee}$ and $F_1\Ee=\Ee$.
			\item A homomorphism $(\Ee,\sigma,1_{\Ee})\to (\Ee',\sigma',1_{\Ee'})$ is an $\OO_X$-linear Lie algebra homomorphism $\varphi:\Ee\to \Ee'$ over $x^{-1}\OO_S$ such that
			$\varphi(1_{\Ee})=1_{\Ee'}$.
			\item We denote the category of Picard algebroids by $\PA_{X/S}$.
		\end{enumerate}
	\end{definition}
	
	\begin{remark}\label{rem:PA_is_groupoid}
		The categories $\PA_{R/k}$ and $\PA_{X/S}$ are groupoids by the five lemma.
	\end{remark}
	
	\begin{remark}
		By property (v), any Picard algebroid is a quasi-coherent $\OO_X$-module. This readily implies that in the affine case, passage to global sections defines a functor
		\[\Gamma(X,-):\PA_{\Spec R/\Spec k}\to\PA_{R/k}.\]
	\end{remark}
	
	\begin{example}[Elementary twists of Picard algebroids: $\Ee^\gamma$]
		Let $(\Ee,[-,-],\sigma,1_{\Ee})$ be a Picard algebroid on $X$. Then for any $\gamma\in (x^{-1}\OO_S)(X)^\times$, $(\Ee,[-,-],\sigma,\gamma^{-1} 1_{\Ee})$ is a Picard algebroid which we denote by $\Ee^\gamma$. This construction is functorial in $\Ee$ in the obvious way.
	\end{example}
	
	\begin{example}
		Put $X=\Spec R$ and $S=\Spec k$. Let $\gamma\in k^\times$ be a unit. Regard $\gamma$ as a constant global section of $x^{-1}\OO_S$. Then for a Picard algebroid $\Ee$ on $X$, we have a canonical identification $\Gamma(X,\Ee)^\gamma=\Gamma(X,\Ee^\gamma)$.
	\end{example}
	
	The next result is obvious by definitions.
	
	\begin{proposition}\label{prop:Zariskilocalpropertyofpa}
		Let $X$ be a smooth scheme over $S$.
		\begin{enumerate}
			\renewcommand{\labelenumi}{(\arabic{enumi})}
			\item For an open subscheme $U\subset X$, the usual restriction determines a functor
			\[(-)|_U:\PA_{X/S}\to\PA_{U/S}.\]
			\item If the structure morphism $x:X\to S$ factors through an open subscheme $V\subset S$, we have a canonical categorical isomorphism $\PA_{X/S}\cong\PA_{X/V}$.
			\item Picard algebroids form a stack over the lisse-Zariski site of $S$.
		\end{enumerate}
	\end{proposition}

	We next relate Picard algebroids with certain torsors (see section \ref{sec:torsor} for basic definitions).
	
	\begin{definition}
		Let $\Omega^{\geq 1}_{X/S}$ be the stupid truncation of the de Rham complex $\Omega^\bullet_{X/S}$ (\cite[section 16.6]{ega44}) at the 0th degree, i.e.,
		$0\to\Omega^1_{X/S}\overset{d}{\to} \Omega^2_{X/S}\overset{d}{\to}\cdots$,
		where the left hand $0$ is placed in degree $0$. We define $\Omega^{\geq 1}_{R/k}$ in a similar way.
	\end{definition}

	\begin{construction}\label{cons:torsor}
		Let $\Ee$ be a Picard algebroid for $X/S$. Let $\mathcal{C}(\Ee)$ be the sheaf of $\OO_X$-linear sections of $\sigma$. This is endowed with the group action of $\Omega^{1}_{X/S}$ for $(\omega+\varphi)(\xi)=\varphi(x)+\langle\omega,\xi\rangle 1_\Ee$ ($\omega\in\Omega^1_{X/S}$, $\varphi\in\Cc(\Ee)$, $\xi\in\Theta_{X/S}$). Let $c$ be the curvature map in \cite[1.2]{beilinsonbernstein1993}, i.e., under the identification $(\wedge^2\Theta_{X/S})^\vee\cong \Omega^{2}_{X/S}$, $c:\mathcal{C}(\Ee)\to \Omega^2_{X/S}$ is given by
		\[c(\varphi)(\xi_1,\xi_2)=[\varphi(\xi_1),\varphi(\xi_2)]-\varphi([\xi_1,\xi_2])\in\Ker\sigma\overset{1_{\Ee}}{\cong}\OO_X.\]
	\end{construction}
	
	\begin{proposition}[{\cite[Lemma 2.1.6]{beilinsonbernstein1993}}]\label{prop:classification/sch}
		The correspondence $\Ee\to (\mathcal{C}(\Ee),c)$ determines a categorical equivalence $\PA_{X/S}\simeq \Tor(\Omega^{\geq 1}_{X/S})$ (see section \ref{sec:torsor} for the definition of $\Tor(\Omega^{\geq 1}_{X/S})$). 
	\end{proposition}
	
	One can prove a similar assertion in the ring setting:

	\begin{proposition}\label{prop:classification}
		We have a categorical equivalence $\PA_{R/k}\simeq \Tor(\Omega^{\geq 1}_{R/k})$, which is defined in a similar way to Proposition \ref{prop:classification}.
	\end{proposition}
	
	\begin{corollary}\label{cor:localization_PA}
		Assume $X=\Spec R$ and $S=\Spec k$. Then the global section functor determines a 2-commutative diagram of equivalences of categories:
		\[\begin{tikzcd}
			\PA_{X/S}\ar[r, "{\Gamma(X,-)}"]\ar[d, "\sim"{sloped, below}]&\PA_{R/k}\ar[d, "\sim"{sloped, above}]\\
			\Tor(\Omega^{\geq 1}_{X/S})\ar[r, "{\Gamma(X,-)}"]&\Tor(\Omega^{\geq 1}_{R/k}).
		\end{tikzcd}\]
	\end{corollary}
	
	\begin{proof}
		This is immediate from Propositions \ref{prop:classification}, \ref{prop:classification/sch} and \ref{prop:tor_affine}.
	\end{proof}
	
	We obtain basic examples and constructions through these equivalences:
	
	\begin{example}
		Let $\Bun_1(X)$ denote the groupoid of line bundles on $X$. Then the log differential $d\log:\OO^\times_X\to \Omega^{\geq 1}_{X/S}$ in Example \ref{ex:logdiff} gives rise to a functor
		\[\Bun_1(X)\simeq \Tor(\OO^\times_X)\to \Tor(\Omega^{\geq 1}_{X/S})\simeq \PA_{X/S},\]
		where the first equivalence is the canonical one.
	\end{example}
	
	\begin{example}[The untwisted Picard algebra/algebroid]
		Let $\tilde{\Theta}_{R/k}$ be the Picard algebra attached to the trivial $\Omega^{\geq 1}_{R/k}$-torsor $(\Omega^1_{R/k},0)$. Explicitly, it is described as follows:
		\[\begin{array}{cc}
			\tilde{\Theta}_{R/k}=R\oplus \Theta_{R/k},
			&1_{\tilde{\Theta}_{R/k}}=1\in R\subset \tilde{\Theta}_{R/k},
		\end{array}\]
		\[\left[a_1+\partial_1,a_2+\partial_2\right]
		=\left[\partial_1,\partial_2\right]
		+\partial_1(a_2)-\partial_2(a_1),
		~(a,a'\in R,~\partial,\partial'\in\Theta_{R/k}),\]
		and $\sigma:\tilde{\Theta}_{R/k}\to\Theta_{R/k}$ is the projection. We can define and construct the untwisted Picard algebroid $\tilde{\Theta}_{X/S}$ in a similar way.
	\end{example}
	
	\begin{example}[{\cite[2.5]{kashiwara1989}}]
		For an $\Omega^{\geq 1}_{R/k}$-torsor $(C,c)$ and a closed 2-form $\eta\in \Omega^{2,\cl}_{R/k}$, define a new $\Omega^{\geq 1}_{R/k}$-torsor $(C,c+\eta)$ by $(c+\eta)(\varphi)=c(\varphi)+\eta$. Correspondingly, for a Picard algebra $(E,\sigma,1_E)$, we obtain a new structure of a Picard algebra on the $R$-module $E_\eta=E$ for $1_E$, $\sigma$, and the Lie bracket $[-,-]_\eta$ defined by
		\[\left[e_1,e_2\right]_\eta
		=\left[e_1,e_2\right]
		+\langle\eta,\sigma(e_1)\wedge \sigma(e_2)\rangle 1_E\]
		for $e_1,e_2\in E_\eta$. This is clearly additive in $\eta$, in the sense that
		$(E_{\eta_1})_{\eta_2}=E_{\eta_1+\eta_2}$ for $\eta_1,\eta_2\in \Omega^{2,\cl}_{R/k}$. We can define and construct $\Ee_\eta$ for a Picard algebroid $\Ee$ and a closed 2-form $\eta\in\Omega^{2,\cl}_{X/S}(X)$ in a similar way.
	\end{example} 
	
	\begin{remark}
		By our hypotheses, $\Theta_{R/k}$ is finitely generated and projective over $R$. Therefore, the $k$-module $Z^2(\Theta_{R/k},R)$ of 2-cocycles is canonically identified with $\Omega^{2,\cl}_{R/k}$. From this perspective, our module $\tilde{\Theta}_{R/K,\eta}$ coincides with the Lie-Rinehart algebra $\Theta_{R/k}(\eta)$ constructed by Maakestad in \cite[section 2]{maakestad}.
	\end{remark}
	
	\begin{example}
		Let $\eta\in\Omega^{2,\cl}_{R/k}$ and $\gamma\in k^\times$. Then the map $(a,\partial)\mapsto (\gamma^{-1}a,\partial)$
		determines an isomorphism $\tilde{\Theta}_{R/k,\gamma\eta}\cong (\tilde{\Theta}_{R/k,\eta})^\gamma$ of Picard algebras.
	\end{example}

	\begin{example}[{\cite[2.5]{kashiwara1989}}]
		For an $\Omega^{\geq 1}_{R/k}$-torsor $(C,c)$ and a 1-form $\theta\in \Omega^{1}_{R/k}$, we have a morphism $\varphi_\theta:(C,c)\to (C,c+d\theta);~\varphi\mapsto c(\varphi)+d\theta$. Moreover, it gives rise to a bijection
		\[\Omega^1_{R/k}\cong\End_{\Tor(\Omega^{\geq 1}_{R/k})}((C,c));~
		\theta\mapsto \varphi_\theta.\]
		Correspondingly, for a Picard algebra $(E,\sigma,1_E)$ and $\theta\in \Omega^1_{R/k}$, we have an isomorphism
		\[\varphi_\theta:E\to E_{d\theta};~e\mapsto 
		e-\langle\theta,\sigma(e)\rangle 1_{E_{\eta+d\theta}}.\]
		This is additive in $\theta$ and closed 2-forms $\eta$ in the sense that for 1-forms $\theta_1,\theta_2\Omega^1_{R/k}$ and $\eta\in\Omega^{2,\cl}_{R/k}$, the diagram
		\[\begin{tikzcd}
			E_\eta\ar[rr, "\varphi_{\theta_1}"]
			\ar[rd, "\varphi_{\theta_1+\theta_2}"']
			&&E_{\eta+d\theta_1}
			\ar[ld, "\varphi_{\theta_2}"]\\
			&E_{\eta+d(\theta_1+\theta_2)}
		\end{tikzcd}\]
		commutes. We can define and construct $\varphi_\theta$ for a Picard algebroid $\Ee$ and a 1-form $\eta\in\Omega^{1}_{X/S}(X)$ in a similar way. Moreover, we have a group isomorphism $\Omega^{1,\cl}_{X/S}(X)\cong
		\End_{\PA_{X/S}}(\Ee);~\theta\mapsto \varphi_\theta$.
	\end{example}
	
	\begin{example}
		Let $\Ee$ be a Picard algebroid on $X$. Then we have a group homomorphism $d\log:\OO_X(X)^\times\to\Omega^{1,\cl}_{X/S}(X)\cong\End_{\PA}(\Ee)$. Explicitly, for a unit $a\in\OO_X(X)^\times$, the automorphism is given by
		$e\mapsto e-\frac{\sigma(e)(a)}{a}1_{\Ee}$.
	\end{example}
	
	\begin{example}
		The groupoid $\PA_{\PP^1_S/S}$ is a setoid since $\Omega^1_{\PP^1_S/S}(\PP^1_S)=0$.
	\end{example}

	\subsection{Functorial constructions}\label{sec:functorspicalgtdo}
	
	In this section we establish the equivalence between tdos and Picard algebroids and related functorial constructions.
	
	\subsubsection{Equivalence of $\TDO$ and $\PA$}\label{subsec:lie}
	
	In this section, we aim to construct the equivalence of groupoids of tdos and Picard algeborids. 
	
	\begin{lemma}[{\cite[2.1.3]{beilinsonbernstein1993}}]\label{lem:lieDispicalg}
		Let $(A,F_\bullet A,i)$ be a tdo for $R/k$. Then $F_1A$, the commutator, the Lie algebra homomorphism $\sigma:F_1 A\to \Theta_{R/k}$ in Lemma \ref{lem:commutator}, and $1\in A$ form a Picard algebra. Similarly, for a tdo$(\Aa,F_\bullet \Aa,i)$, $F_1 \Aa$, the commutator, the Lie algebra homomorphism $\sigma:F_1 \Aa\to \Theta_{X/S}$ in Lemma \ref{lem:commutator/scheme}, and $1\in \Aa(X)$ form a Picard algebroid.
	\end{lemma}
	
	The resulting functors $\TDO_{R/k}\to \PA_{R/k}$ and $\TDO_{X/S}\to\PA_{X/S}$ will be denoted by $\Lie$. We next construct its quasi-inverse.
	
	\begin{definition}\label{def:A_E}
		Let $(E,[-,-],\sigma,1_E)$ be a Picard algebra. Then let $A_E$ be the $k$-algebra generated by $E$ with the following relations (in the tensor algebra):
		\begin{enumerate}
			\renewcommand{\labelenumi}{(\roman{enumi})}
			\item for $a\in R$ and $e\in E$, $a1_E\otimes e=ae$;
			\item for $e,e'\in E$, $e\otimes e'-e'\otimes e=\left[e,e'\right]$;
			\item $1_E=1$.
		\end{enumerate}
		We will write $F_\bullet A_E$ for the filtration on $A_E$ generated by $F_\bullet E$.
	\end{definition}
	
	\begin{remark}
		If we fix a section of $\sigma$ to identify $E$ with $\tilde{\Theta}_{R/k,\mu}$, $A_E$ is canonically isomorphic to $U(R,\Theta_{R/k},\mu)$ in \cite[section 3]{maakestad}. Our construction is slightly different from Maakestad's (or originally Rinehart's). In \cite{maakestad}, Maakestad considered the nonunital subalgebra $U^+$ of the enveloping algebra of $E$ generated by $E$ and then defined $U(R,\Theta_{R/k},\mu)$ as the quotient of $U^+$ by the two sided ideal generated by $a\otimes e-ae$ ($a\in R$, $e\in\Theta_{R/k}(\mu)$). The element $1\in R\subset\Theta_{R/k}(\mu)$ eventually becomes the unit of $U(R,\Theta_{R/k},\mu)$ by the added relations.
	\end{remark}
	
	\begin{example}
		Write $D_R=D_{R/k}=A_{\tilde{\Theta}_{A/k}}$. Explicitly, $D_R$ is the $k$-algebra generated by the $k$-module $R\oplus \Theta_{R/k}$ with the following relations:
		\begin{enumerate}
			\renewcommand{\labelenumi}{(\roman{enumi})}
			\item for $a_1,a_2\in R$, $a\otimes a'=aa'$;
			\item for $a\in R$ and $\partial\in \Theta_{R/k}$, $\partial\otimes a-a\otimes \partial=\partial(a)$, and 
			\item for $\partial,\partial'\in\Theta_{R/k}$, $\partial\otimes \partial'-\partial'\otimes \partial=\left[\partial,\partial'\right]$, where the bracket is taken in $\Theta_{R/k}$;
			\item for $a\in R$ and $\partial\in\Theta_{R/k}$, $a\otimes \partial=a\partial$;
			\item the units of $R$ and $D_{R}$ coincide.
		\end{enumerate}
		Regard $R\oplus \Theta_{R/k}$ as a graded $k$-module for $\deg a=0$ and $\deg\partial=1$ ($a\in R$, $\partial\in\Theta_{R/k}$). Define a natural filtration on the tensor $k$-algebra generated by $R\oplus\Theta_{R/k}$ from the canonical grading. This induces the filtration $F_\bullet D_R$ on $D_R$.
	\end{example}
	
	\begin{definition}
		Let $\Ee$ be a Picard algebroid on $X$. Then let $\Aa_{\Ee}$ be the $x^{-1}\OO_S$-algebra generated by $\Ee$ with the following relations (in the tensor algebra):
		\begin{enumerate}
			\renewcommand{\labelenumi}{(\roman{enumi})}
			\item for an open subset $U\subset X$ and sections $a\in \OO_X(U)$, $e\in \Ee(U)$, we have $a1_{\Ee}\otimes e=ae$;
			\item for an open subset $U\subset X$ and $e,e'\in \Ee(U)$, $e\otimes e'-e'\otimes e=\left[e,e'\right]$;
			\item $1_{\Ee}=1$.
		\end{enumerate}
		We will write $F_\bullet \Aa_{\Ee}$ for the filtration on $\Aa_{\Ee}$ generated by $F_\bullet \Ee$.
	\end{definition}
	
	Following \cite[Definition 2.1.3]{beilinsonbernstein1993}, we formulate
	
	\begin{theorem}\label{thm:DEistdo}
		Define $i:A\to A_E$ by $a\mapsto a1_E$. Then $(A_E,i)$ is a tdo. In this way, we obtain a functor $A_{(-)}:\PA_{R/k}\to\TDO_{R/k}$. We can also define a functor $\Aa_{(-)}:\PA_{X/S}\to\TDO_{X/S}$ in a similar way.
	\end{theorem}

	The essential difficulty for these two results is to prove the Poincar\'e-Birkhoff-Witt theorem (locally on $X$ in the sheaf setting). This is achieved by a straightforward modification of \cite[Proof of Theorem 3.1]{rinehart1963}, which we omit in this paper.
	
	\begin{corollary}\label{cor:pavstdo}
		The functors $A_{(-)}:\PA_{R/k}\to \TDO_{R/k}$ (resp.~$\Aa_{(-)}:\PA_{X/S}\to \TDO_{X/S}$) and $\Lie:\TDO_{R/k}\to \PA_{R/k}$ (resp.~$\Lie:\TDO_{X/S}\to \PA_{X/S}$) are quasi-inverse to each other.
	\end{corollary}
	
	\begin{proof}
		We only prove the assertion for the ring setting.
		Let $E$ be a Picard algebra. Then the canonical map $E\to A_E$ factors through $\Lie A_E$ which is a Picard algebra (Lemma \ref{lem:lieDispicalg}). The map $E\to \Lie A_E$ is a homomorphism of Picard algebras. This is an isomorphism (Remark \ref{rem:PA_is_groupoid}). Conversely, let $A$ be a tdo. Then the map $\Lie A\to A$ extends to an isomorphism $A_{\Lie A}\to A$ of tdos. This is an isomorphism (Remark \ref{rem:tdo_is_groupoid}). These constructions naturally give bijections on Hom sets. 
	\end{proof}

	These equivalences lead us to basic examples and constructions of tdos and Picard algebroids.

	\begin{definition}[The sheaf of untwisted differential operators on $X/S$]
		Let $\D_X=\D_{X/S}$ be the tdo attached to the untwisted Picard algebroid on $X$.
	\end{definition}
	
	\begin{example}\label{ex:D_G}
		Let $G$ be an $S$-affine smooth group scheme over a scheme $S$. Let $g:G\to S$ denote the structure morphism. Then the canonical $G$-equivariant Lie algebra homomorphism $\lieg\to g_\ast\Theta_{G/S}$ gives rise to a canonical $G$-equivariant $\OO_G$-module isomorphism $g^\ast\lieg\cong\Theta_{G/S}$. It also attaches a $G$-equivariant $\OO_S$-algebra homomorphism $U(\lieg)\to g_\ast\D_G$, where $U(\lieg)$ is the enveloping algebra of $\lieg$. One can see that the associated $\OO_G$-module homomorphism
		\[g^\ast U(\lieg)\to \D_G\]
		is an isomorphism by the Poincar\'e-Birkhoff-Witt theorem of $U(\lieg)$ (work locally in $S$ if necessary to assume $S$ affine and $\lieg$ free). We remark that the induced structure of a tdo on the left $\OO_{G}$-module $g^\ast U(\lieg)$ is uniquely characterized by the following properties:
		\begin{enumerate}
			\item[(i)] The left $\OO_G$-module structure on $g^\ast U(\lieg)$ coincides with that determined by the left multiplication via the structure homomorphism $\OO_{G}\to g^\ast U(\lieg)$;
			\item[(ii)] As a filtered left $\OO_{G}$-module, $F_\bullet g^\ast U(\lieg)=g^\ast F_\bullet U(\lieg)$, where $F_\bullet U(\lieg)$ is the canonical filtration on $U(\lieg)$;
			\item[(iii)] The canonical map $g^{-1}U(\lieg) \to g^\ast U(\lieg)$ is a homomorphism of sheaves of $g^{-1}\OO_S$-algebras.
		\end{enumerate}
	\end{example}
	
	\begin{example}[Elementary twists of tdos: $\Aa^\gamma$]
		Let $\Aa$ be a tdo on $X$, and $\gamma\in (x^{-1}\OO_S)(X)^\times$. Then define a new tdo $\Aa^\gamma$ by
		$\Aa_{(\Lie\Aa)^\gamma}$.
	\end{example}
	
	\begin{example}[The opposite Picard algebra/algebroid]\label{ex:lievsop}
		Let $(E,\sigma,1_E)$ be a Picard algebra for $R/k$. Then the passage of Example \ref{ex:optdo} through the equivalence $\TDO_{R/k}\simeq \PA_{R/k}$ leads us to a new Picard algebra $(E^{\op},-\sigma,1_E)$, where 
		\begin{enumerate}
			\renewcommand{\labelenumi}{(\roman{enumi})}
			\item $E^{\op}=E$ is a right (= left) $R$-module for $ea\coloneqq ae+\sigma(e)(a)1_E$, and
			\item $E^{\op}$ is the opposite Lie algebra to $E$.
		\end{enumerate}
		Note that $E$ is an $R$-bimodule over $k$ for the given left action and (i). We can define and construct the opposite Picard algebroid $\Ee^{\op}$ in a similar way.
	\end{example}
	
	\begin{remark}
		One could prove that $(E^{\op},-\sigma,1_E)$ is a Picard algebra directly. This is elementary. Likewise, it is rather an elementary fact that $E$ is equipped with the structure of an $R$-bimodule as explained above.
	\end{remark}
	
	\begin{corollary}[{\cite[2.7]{kashiwara1989}}]\label{cor:oppositecocyle}
		For a tdo $\Aa$ on $X$, we have a canonical isomorphism of tdos
		\[\omega_{X/S}\otimes_{\OO_X}\Aa^{-1}\otimes_{\OO_X} \omega^\vee_{X/S}\cong \Aa^{\op}.\]
	\end{corollary}
	
	\begin{proof}
		Write $\Ee=\Lie\Aa$. We define a map $\Phi:\omega_{X/S}\times\Ee^{-1}\times\omega_{X/S}^\vee\to \Ee^{\op}$ by
		\[\Phi(\omega,e,\delta)=-\langle\omega,\delta\rangle e
		-\langle L_{\sigma(e)}\omega,\delta\rangle 1_{\Ee}.\]
		Regard the $\OO_X$-modules $\omega_{X/S}$ and $\omega^\vee_{X/S}$ as $\OO_X$-bimodules for the multiplication
		\[\OO_X\otimes_{x^{-1}\OO_S}\OO_X\to\OO_X\]
		of $\OO_X$. Then the map $\Phi$ descends to an $\OO_X$-bimodule homomorphism
		\[\varphi:\omega_{X/S}\otimes_{\OO_X}\Ee^{-1}\otimes_{\OO_X}\omega_{X/S}^\vee\to \Ee^{\op},\]
		which is easily verified with Example \ref{ex:lievsop}.
		
		Identify the domain of $\varphi$ with $\Lie (\omega_{X/S}\otimes_{\OO_X}\Aa^{-1}\otimes_{\OO_X} \omega_{X/S}^\vee)$ to put the structure of a Picard algebroid on $\omega_{X/S}\otimes_{\OO_X}\Ee^{-1}\otimes_{\OO_X}\omega_{X/S}^\vee$. Then it is easily verified by the standard argument on affine open patches on which $\Theta_{X/S}$ is free of finite rank to assume $\omega_{X/S}$ is free of rank 1 that $\varphi$ is a homomorphism of Picard algebroids.
		
		In view of Corollary \ref{cor:pavstdo} and Example \ref{ex:lievsop}, the assertion follows by applying $\Aa_{(-)}$ to the isomorphism $\varphi$.
	\end{proof}
	
	\subsubsection{The base change functors $k'\otimes_k-$ and $s^\ast_X$}\label{sec:Sbasechange}
	
	In this section, we study base change functors. We define a functor $\Tor(\Omega^{\geq 1}_{X/S})\to \Tor(\Omega^{\geq 1}_{X'/S'})$ by the canonical homomorphism
	\[s^{-1}_X\Omega^{\geq 1}_{X/S}\to \Omega^{\geq 1}_{X'/S'}\]
	(see Construction \ref{cons:deRham}). Correspondingly, we obtain functors
	\[\begin{array}{cc}
		\PA_{X/S}\to \PA_{X'/S'},&\TDO_{X/S}\to\TDO_{X'/S'},
	\end{array}\]
	which we denote by $s^\ast_X$. We call them the base change functors. In fact, for a Picard algebroid $(\Ee,1_{\Ee},\sigma)$, its base change consists of the left module base change $s^\ast_X\Ee$, $1\otimes 1_{\Ee}$, $s^\ast_X\sigma$, and the unique Lie bracket over $(x')^{-1}\OO_{S'}$ with the following properties:
	\begin{enumerate}
		\renewcommand{\labelenumi}{(\roman{enumi})}
		\item $(s_X^\ast\Ee,s^\ast_X\sigma,1\otimes 1_{\Ee})$ is a Picard algebroid under the identification $s^\ast_X\Theta_{X/S}\cong \Theta_{X'/S'}$;
		\item The canonical map $s^{-1}_X\Ee\to s^\ast_X\Ee$ is a Lie algebra homomorphism. Here we put the structure of a Lie algebra on $s^{-1}_X\Ee$ in the canonical way.
	\end{enumerate}
	Similarly, for a tdo $\Aa$ on $X$, $s_X^\ast\Aa$ is the unique tdo with the following properties:
	\begin{enumerate}
		\renewcommand{\labelenumi}{(\roman{enumi})}
		\item The left $\OO_{X'}$-module structure on $s_X^\ast\Aa$ coincides with that determined by the left multiplication via the structure homomorphism $\OO_{X'}\to s_X^\ast\Aa$;
		\item As a filtered left $\OO_{X'}$-module, $F_\bullet(s_X^\ast\Aa)$ is the pullback of the filtered left $\OO_X$-module $F_\bullet\Aa$;
		\item The canonical map $s^{-1}_X\Aa\to s^\ast_X\Aa$ is a homomorphism of $s^{-1}_Xx^{-1}\OO_{S}$-algebras.
	\end{enumerate}
	
	Thanks to our formalism in terms of $\Omega^{\geq 1}_{X/S}$-torsors, we immediately deduce:
	
	\begin{theorem}\label{thm:tdofpqcstack}
		The correspondence $S'\mapsto \TDO_{X\times_S S'/S'}$ determines a stack over the big fpqc site of $S$.
	\end{theorem}

	We next discuss the ring setting. We define the base change functors
	\[\begin{array}{cc}
		k'\otimes_k-:\PA_{R/k}\to \PA_{R\otimes_k k'/k'},&k'\otimes_k-:\TDO_{R/k}\to \TDO_{R\otimes_k k'/k'}
	\end{array}\]
	through the canonical homomorphism $\Omega^\bullet_{R/k}\to k'\otimes_k \Omega^{\bullet}_{R/k}\cong\Omega^{\bullet}_{R\otimes_k k'/k'}$. We would like to describe them explicitly. Since $R$ is smooth over $k$, $\Omega^1_{R/k}$ is finitely generated and projective as an $R$-module. Hence we have a Lie algebra isomorphism
	\[\begin{split}
		k'\otimes_k \Theta_{R/k}
		&\cong(R\otimes_k k')\otimes_R \Theta_{R/k}\\
		&\cong \Hom_{R\otimes_k k'}((R\otimes_k k')\otimes_R \Omega^1_{R/k},R\otimes_k k')\\
		&\cong \Hom_{R\otimes_k k'}(\Omega^1_{R\otimes_k k'/k'},R\otimes_k k')\\
		&\cong\Theta_{R\otimes_k k'/k'}.
	\end{split}\]
	For a Picard algebra $E$ over $R$,
	\[\begin{array}{ccc}
		k'\otimes_k E,
		&k'\otimes_k\sigma:k'\otimes_k E\to k'\otimes_k \Theta_{R/k}
		\cong \Theta_{R\otimes_k k'/k'},
		&1\otimes 1_E
	\end{array}\]
	form the base change of $E$. Similarly, for a tdo $A$ over $R$, $k'\otimes_k A$ is naturally equipped with the structure of a tdo over $R\otimes_k k'$ since the $p$-th symmetric power $\Sym^p\Theta_{R/k}$ is a finitely generated and projective $R$-module for every nonnegative integer $p$. This exhibits the base change of $A$.
	
	The next assertion is evident from constructions.
	
	\begin{lemma}
		Put $X=\Spec R$, $S=\Spec k$, and $S'=\Spec k'$. Then the diagram
		\[\begin{tikzcd}
			&\TDO_{X/S}
			\ar[dl, "{\Gamma(X,-)}"']\ar[rr, "\Lie"]\ar[dd]
			& & \PA_{X/S}
			\ar[dl, "{\Gamma(X,-)}"']\ar[dd, "s_X^\ast"] \\
			\TDO_{R/k}
			\ar[rr,crossing over, "\Lie"]
			\ar[dd, "k'\otimes_k-"'] & & \PA_{R/k}\\
			& \TDO_{X\times_S S'/S'}
			\ar[dl, "{\Gamma(X\times_S S',-)}"]\ar[rr, "\Lie\ \ \ \ \ \ \ \ "']& &
			\PA_{X\times_S S'/S'}
			\ar[dl, "{\Gamma(X\times_S S',-)}"] \\
			\TDO_{k'\otimes_k R/k'}
			\ar[rr, "\Lie"']
			& &\PA_{k'\otimes_k R/k'}
			\ar[from=uu,crossing over, "k'\otimes_k-"]
		\end{tikzcd}\]
		is 2-commutative.
	\end{lemma}

	Finally, we discuss the compatibility of base changes with the globalization $x_\ast$ which generalizes the global section functor $\Gamma(X,-)$ in the classical theory of tdos. Namely, let $\Aa$ be a tdo on a smooth $S$-scheme $X$. Then by general theory of sheaves, we obtain an $\OO_S$-algebra $x_\ast\Aa$ since $\Aa$ is an $x^{-1}\OO_S$-algebra. Hence for a morphism $s:S'\to S$ of schemes, we get the base change map
	\begin{equation}
		s^\ast x_\ast\Aa\to x'_\ast s_X^\ast\Aa
		\label{eq:bcmapfortdo}
	\end{equation}
	as an $\OO_S$-module (cf.\ Example \ref{ex:noflatbcmap}).
	
	\begin{proposition}\label{prop:bcmaptdo}
		The map \eqref{eq:bcmapfortdo} is a homomorphism of $\OO_{S'}$-algebras.
	\end{proposition}
	
	We remark that $s^\ast x_\ast\Aa$ is canonically equipped with the structure of an $\OO_{S'}$-algebra since $x_\ast\Aa$ is an $\OO_S$-algebra.
	
	\begin{proof}
		By the universal property of scalar extension, we may restrict to $s^{-1} x_\ast\Aa$. Then the assertion follows since we have a compatible base change map $s^{-1} x_\ast\Aa\to x'_\ast s_X^{-1}\Aa$ which is a homomorphism of $s^{-1}\OO_S$-algebras.
	\end{proof}

	\subsubsection{Functoriality in $X$: $f_\cdot$ and $f^\cdot$}\label{sec:Xbasechange}
	In this section, we discuss functoriality of Picard algebras and Picard algebroids in $R$ and $X$ respectively, following \cite{beilinsonbernstein1993}.
	
	We define the pullback functor $\Tor(\Omega^{\geq 1}_{Y/S})\to\Tor(\Omega^{\geq 1}_{X/S})$ by the canonical homomorphism $f^{-1}\Omega^{\bullet}_{Y/S}\to \Omega^{\bullet}_{X/S}$. Correspondingly, we obtain functors 
	\[\begin{array}{cc}
		f^\cdot:\PA_{Y/S}\to \PA_{X/S},&f^\cdot:\TDO_{Y/S}\to\TDO_{X/S}.
	\end{array}\]
	We call these functors the pullback.
	
	\begin{example}
		Let $j:U\to X$ be an open immersion of smooth $S$-schemes.
		Then
		\[\begin{array}{cc}
			j^\cdot:\PA_{X/S}\to\PA_{U/S}&(\text{resp.~}j^\cdot:\TDO_{X/S}\to\TDO_{U/S})
		\end{array}\]
		is isomorphic to the restriction functor of
		Proposition \ref{prop:Zariskilocalpropertyofpa} (1) (resp.~Lemma \ref{lem:tdozariski} (1)).
	\end{example}
	
	The following assertion is straightforward:
	
	\begin{theorem}[{\cite[2.2]{beilinsonbernstein1993}}]\label{thm:stackinlisseetaletop}
		The correspondence $X\mapsto \TDO_{X/S}$ determines a stack over the lisse-\'etale site of $S$. In particular, tdos form a stack over the small \'etale site of $X$.
	\end{theorem}
	
	\begin{remark}
		By Theorem \ref{thm:stackinlisseetaletop} and general theory, tdos form a stack over the smooth site of $S$. In particular, tdos form a stack over the smooth site of $X$.
	\end{remark}
	
	To describe $f^\cdot$ for Picard algebroids, define an $\OO_X$-linear map
	\begin{equation}
		f_\ast:\Theta_{X/S}\to f^\ast\Theta_{Y/S}
		\label{eq:pushofvectorfield}
	\end{equation}
	by the diagram
	\[\begin{tikzcd}
		f^\ast\Theta_{Y/S}\ar[r]&\Der_S(f^{-1}\OO_Y,\OO_X)
		&\Theta_{X/S}\ar[l, "-\circ f^\sharp"']\\
		f^\ast\iHom_{\OO_Y}(\Omega^1_{Y/S},\OO_Y)\ar[r, "\sim"]
		\ar[u, "\sim"{sloped, above}]
		&\iHom_{\OO_X}(f^\ast\Omega^1_{Y/S},\OO_X)
		\ar[u, "\sim"{sloped, below}]
		&\iHom_{\OO_X}(\Omega^1_{X/S},\OO_X)\ar[l]
		\ar[u, "\sim"{sloped, below}]
	\end{tikzcd}\]
	of canonical homomorphisms. Unwinding the definitions, we obtain:
	
	\begin{proposition}\label{prop:pullbackofpicalg/sch}
		Let $\Ee$ be a Picard algebroid for $Y/S$. Then the $\OO_X$-module
		\[\Theta_{X/S}\times_{f^\ast\Theta_{Y/S}} f^\ast\Ee\]
		with the global section $1_{f^\cdot\Ee}\coloneqq(0,1\otimes 1_{\Ee})$ and the projection exhibits the pullback of $\Ee$ for the well-defined Lie bracket locally expressed as
		\[\left[(\partial,\tilde{e}),
		(\partial',\tilde{e}')\right]=(\left[\partial,\partial'\right], \sum_j \partial(b_j')\otimes e_j'-\sum_i \partial'(b_i)\otimes e_i+\sum_{i,j} b_ib_j'\otimes \left[e_i,e_j'\right])\]
		for local sections $(\partial,\tilde{e})$ and $(\partial',\tilde{e}')$ of $\Theta_{X/S}\times_{f^\ast\Theta_{Y/S}} f^\ast\Ee$ with presentations
		\[\begin{array}{ccc}
			\tilde{e}=\sum_i b_i\otimes e_i,&
			\tilde{e}'=\sum_j b'_j\otimes e'_j
			&(b_i,b'_j\in\OO_X,~e_i,e_j\in f^{-1}\Ee)
		\end{array}\]
		We remark that $\tilde{e}$ and $\tilde{e}'$ have the last presentations locally in $X$.
	\end{proposition}
	
	\begin{example}
		We have a canonical isomorphism $f^\cdot\D_Y\cong\D_X$. 
	\end{example}
	
	Let us note the structure constraints of the stack of Picard algebroids in terms of this construction:
	\begin{itemize}
		\item For morphisms $f:X\to Y$ and $g:Y\to Z$ of smooth $S$-schemes, $\alpha_{f,g}:(g\circ f)^\cdot\cong f^\cdot\circ g^\cdot$ is given by the canonical maps
		\[\begin{split}
			f^\cdot(g^\cdot(\Ee))
			&= f^\cdot(\Theta_{Y/S}\times_{g^\ast\Theta_{Z/S}}
			g^\ast \Ee)\\
			&=\Theta_{X/S}\times_{f^\ast\Theta_{Y/S}}
			f^\ast(\Theta_{Y/S}\times_{g^\ast \Theta_{Z/S}}g^\ast\Ee)\\
			&\to \Theta_{X/S}\times_{f^\ast\Theta_{Y/S}}
			(f^\ast\Theta_{Y/S}\times_{f^\ast g^\ast\Theta_{Z/S}} f^\ast g^\ast \Ee)\\
			&\cong \Theta_{X/S}
			\times_{f^\ast g^\ast\Theta_{Z/S}}
			f^\ast g^\ast \Ee\\			
			&\cong \Theta_{X/S}
			\times_{(g\circ f)^\ast\Theta_{Z/S}}
			(g\circ f)^\ast \Ee\\
			&= (g\circ f)^\cdot \Ee
		\end{split}\]
		for Picard algebroids $\Ee$ on $Z$, and
		\item $\epsilon_X:\id^\cdot_X\cong \id_{\PA_{X/S}}$ is the projection.
	\end{itemize}
	
	Likewise, we define $f_\cdot:\PA_{R/k}\to \PA_{R'/k}$ and $f_\cdot:\TDO_{R/k}\to \TDO_{R'/k}$. We also define $\alpha_{f,g}$ and $\epsilon_R$ in a similar way to the above. It is evident that the diagram
	\[\begin{tikzcd}
		&\TDO_{Y/S}
		\ar[dl, "{\Gamma(Y,-)}"']\ar[rr, "\Lie"]\ar[dd]
		& & \PA_{Y/S}
		\ar[dl, "{\Gamma(Y,-)}"']\ar[dd, "f^\cdot"] \\
		\TDO_{A/k}
		\ar[rr,crossing over, "\Lie"]
		\ar[dd, "f_\cdot"'] & & \PA_{A/k}\\
		& \TDO_{X/S}
		\ar[dl, "{\Gamma(X,-)}"']\ar[rr]& &
		\PA_{X/S}
		\ar[dl, "{\Gamma(X,-)}"] \\
		\TDO_{B/k}
		\ar[rr, "\Lie"]
		& &\PA_{B/k}.
		\ar[from=uu,crossing over, "\overset{\ }{\overset{\ }{\overset{\ }{f_\cdot}}}"]
	\end{tikzcd}\]
	is 2-commutative for $X=\Spec R'$, $Y=\Spec R$, and $S=\Spec k$. One has a similar description to Proposition \ref{prop:pullbackofpicalg/sch}. Namely, define an $R'$-linear map $f^\ast:\Theta_{R'/k}\to R'\otimes_R \Theta_{R/k}$ in a similar way to \eqref{eq:pushofvectorfield}. Then we have:
	
	\begin{proposition}\label{prop:pullbackofpicalg}
		Let $E$ be a Picard algebroid for $R/k$. Then the $R'$-module
		\[\Theta_{R'/k}\times_{R'\otimes_R\Theta_{R/S}} R'\otimes_RE\]
		with the element $1_{f_\cdot E}\coloneqq(0,1\otimes 1_{E})$ and the projection exhibits $f_\cdot E$ for the well-defined Lie bracket expressed as
		\[\left[(\partial,\tilde{e}),
		(\partial',\tilde{e}')\right]=(\left[\partial,\partial'\right], \sum_j \partial(b_j')\otimes e_j'-\sum_i \partial'(b_i)\otimes e_i+\sum_{i,j} b_ib_j'\otimes \left[e_i,e_j'\right])\]
		for $(\partial,\tilde{e})$ and $(\partial',\tilde{e}')$ in $\Theta_{R'/k}\times_{R'\otimes_R\Theta_{R/S}} R'\otimes_R E$ with presentations
		\[\begin{array}{ccc}
			\tilde{e}=\sum_i b_i\otimes e_i,&
			\tilde{e}'=\sum_j b'_j\otimes e'_j
			&(b_i,b'_j\in R',~e_i,e_j\in E).
		\end{array}\]
	\end{proposition}

	Towards the definition of the transfer bimodule in section \ref{sec:operationswithDmodules} (essentially \eqref{eq:bimod}), we need a better understanding to this complicated Lie bracket.
	
	\begin{construction}\label{cons:ambientfE}
		Let $(A,F_\bullet A,i)$ be a tdo, and $j:E\to \Lie A$ be an isomorphism of Picard algebras. Typical examples are $(\Lie A,A)$ and $(E,A_E)$. Hence such $A$ uniquely exists up to isomorphisms for each $E$. We consider $R'$ as an $R$-module via its $R$-algebra structure. Apply $-\otimes_R A$ to get the structure of a right $A$-module on $R'\otimes_R A$. Let $\End_A(R'\otimes_R A)$ denote the $k$-algebra of right $A$-linear endomorphisms of $R'\otimes_R A$. We also define a $k$-linear map $\iota:f_\cdot E=\Theta_{R'/k}\times_{R'\otimes_{R}\Theta_{R/k}}(R'\otimes_R E))\to \Hom_k(R',R'\otimes_R A)$ by the sum of the two arrows
		\[(\Theta_{R'/k}\times_{R'\otimes_R\Theta_{R/k}}(R'\otimes_R E))\otimes_k R'
		\overset{\pr_1}{\to}\Theta_{R'/k}\otimes_k R'\to R'\overset{R'\otimes 1}{\to} R'\otimes_R A,\]
		\[\begin{split}
			(\Theta_{R'/k}\times_{R'\otimes_R\Theta_{R/k}}(R'\otimes_R E))\otimes_k R'
			&\overset{\pr_2}{\to}(R'\otimes_R E)\otimes_k R'
			\overset{C_{R'\otimes_R E,R'}}{\cong} R'\otimes_k (R'\otimes_R E)\\
			&\overset{m_{R'}}{\to} R'\otimes_R E
			\overset{R'\otimes j}{\to} R'\otimes_R A,
		\end{split}\]
		where the unlabeled middle arrow in the first line is given by the action.
		Namely, we have a $k$-linear map
		\[\iota((\partial,\tilde{e})):R'\to R'\otimes_R A;~b\mapsto \partial(b)\otimes 1+b (R'\otimes j)(\tilde{e})\]
		for $(\partial,\tilde{e})\in f_\cdot E$.
		If we write $\tilde{e}=\sum_i b_i\otimes e_i$, the image is expressed as
		\[\partial(b)\otimes 1+b (R'\otimes j)(\tilde{e})=
		\partial(b)\otimes 1+\sum_i bb_i\otimes j(e_i).\]
		We define $\iota:f^\cdot \Ee\to
		\iHom_{x^{-1}\OO_S}(\OO_X,f^\ast \Aa)$ for a tdo $\Aa$ on $X$ and a homomorphism $j:\Ee\to \Lie \Aa$ of Picard algebroids.
	\end{construction}

	\begin{lemma}\label{lem:R-linear}
		Let $A$, $E$, and $j$ be as above. If we regard $R'\otimes_R A$ as a right $R$-module for multiplication from the right, $\Theta_{R'/k}\times_{R'\otimes_R\Theta_{R/k}}(R'\otimes_R E))\to \Hom_k(R',R'\otimes_R A)$ factors through the set $\Hom_R(R',R'\otimes_R A)$. In particular, we obtain a map
		\[\Theta_{R'/k}\times_{R'\otimes_R\Theta_{R/k}}(R'\otimes_R E))\to \End_A(R'\otimes_R A)\]
		by scalar extension, which we denote by the same symbol $\iota$. For $\Aa$, $\Ee$, and $j$ as above, we obtain $\iota:f^\cdot \Ee\to\iEnd_{f^{-1}\Aa}(f^\ast \Aa)$ in a similar way.
	\end{lemma}
	
	\begin{proof}
		We only prove the assertion for the ring setting.
		Let $(\partial,\tilde{e})\in f_\cdot E$, $b\in R'$, and $a\in R$. We wish to prove
		\[\partial (bf(a))\otimes 1+ bf(a)(R'\otimes j)(\tilde{e})
		=(\partial(b)\otimes 1+b(R'\otimes j)(\tilde{e}))i(a)\]
		(recall that $i$ is the given homomorphism $R\to A$ in Definition \ref{defn:tdo}).
		We have
		\begin{flalign*}
			&\partial (bf(a))\otimes 1+ bf(a)(R'\otimes j)(\tilde{e})
			-(\partial(b)\otimes 1+b(R'\otimes j)(\tilde{e}))i(a)\\
			&=(\partial(b)f(a)+b\partial(f(a)))\otimes 1 +bf(a)(R'\otimes j)(\tilde{e})
			-(\partial(b)\otimes 1+b(R'\otimes j)(\tilde{e}))i(a)\\
			&=b\partial(f(a))\otimes 1 +bf(a)(R'\otimes j)(\tilde{e})
			-b(R'\otimes j)(\tilde{e})i(a).
		\end{flalign*}
		
		The proof will be completed by showing
		\[\partial(f(a))\otimes 1 =
		(R'\otimes j)(\tilde{e})i(a)-f(a)(R'\otimes j)(\tilde{e}).\]
		Observe that we have an equality
		$i(\sigma(-)(a))=[j(-),i(a)]$
		of $R$-linear maps from $E$ to $A$. Apply $R'\otimes_R-$ to this, and use $(R'\otimes_R \sigma)(\tilde{e})=f^\ast(\partial)$ to deduce the the assertion.
	\end{proof}

	\begin{theorem}\label{thm:ambientfcdot}
		Let $E,A,j$ be as above.
		\begin{enumerate}
			\renewcommand{\labelenumi}{(\arabic{enumi})}
			\item The map $\iota:f_\cdot E\to \End_{A}(R'\otimes_R A)$ is an injective Lie algebra homomorphism.
			\item We have $\iota((0,1\otimes 1_E))=\id_{R'\otimes_R A}$.
			\item For $b\in R'$
			and $(\partial,\tilde{e})\in f_\cdot E$, we have
			$\iota((0,b\otimes 1_E))\circ
			\iota((\partial,\tilde{e}))=
			\iota((b\partial,b\tilde{e}))$.
		\end{enumerate}
	\end{theorem}
	
	\begin{proof}
		Choose a retract map from a free $R$-module of finite rank to $\Theta_{R/k}$ as an $R$-module. Take $\Sym_R$ to show that for each integer $n\geq 0$, $\Sym^n_R\Theta_{R/k}$ is a finitely generated and projective $R$-module. Hence the sequence
		\[0\to R'\otimes_R F_n A\to R'\otimes_R F_{n+1} A\to R'\otimes_R\Gr^{n+1} A\to 0\]
		is split exact. Take $\varinjlim\limits_n$ to see that the canonical $R'$-module homomorphism $R'\otimes_R F_n A\to R'\otimes_R A$ is injective, and that $R'\otimes_R F_\bullet A$ is an exhaustive filtration of $R'\otimes_R A$. In particular, 
		\begin{equation}
			R'\cong R'\otimes_R F_0 A\to R'\otimes_R A;~b\mapsto b\otimes 1,\label{eq:bc_1}
		\end{equation}
		\begin{equation}
			R'\otimes j:R'\otimes_R E\cong R'\otimes_R F_1 A\to R'\otimes_R A;
			~b\otimes e\mapsto b\otimes j(e)\label{eq:bc_2}
		\end{equation}
		are injective.
		
		Let $(\partial,\tilde{e})\in\Ker \iota$. Then we have
		$(R'\otimes j)(\tilde{e})
		=\iota((\partial,\tilde{e}))(1\otimes 1)
		=0$ in $R'\otimes_R A$. Since \eqref{eq:bc_1} is injective, $\tilde{e}=0$. For each $b\in R'$, we have $0=\iota((\partial,\tilde{e}))(b\otimes 1)
		=\partial(b)\otimes 1$.
		Since \eqref{eq:bc_2} is injective, $\partial(b)=0$. This shows that $\iota$ is injective.
		
		The assertion that $\iota$ is a Lie algebra homomorphism is straightforward: 
		Pick elements $(\partial,\tilde{e}),(\partial',\tilde{e}')\in f_\cdot E$. Then for $b\in R'$, we have
		\begin{flalign*}
			&\left[\iota((\partial,\tilde{e})),
			\iota((\partial',\tilde{e}'))\right](b)\\
			&=\iota((\partial,\tilde{e}))(\partial'(b)\otimes 1
			+b (R'\otimes j)(\tilde{e}'))-
			\iota((\partial',\tilde{e}'))
			(\partial(b)\otimes 1
			+b(R'\otimes j)(\tilde{e}))\\
			&=\partial(\partial'(b))\otimes 1
			+\partial'(b) (R'\otimes j)(\tilde{e})
			+\iota((\partial,\tilde{e}))(b (R'\otimes j)(\tilde{e}'))\\
			&-\partial'(\partial(b))\otimes 1
			-\partial(b) (R'\otimes j)(\tilde{e}')
			-\iota((\partial',\tilde{e}'))(b (R'\otimes j)(\tilde{e})).
		\end{flalign*}
		Write $\tilde{e}=\sum b_i\otimes e_i$ and $\tilde{e}'=\sum b_j'\otimes e_j'$. Then we have
		\[\partial'(b) (R'\otimes j)(\tilde{e})
		-\partial(b) (R'\otimes j)(\tilde{e}')
		=\partial'(b)b_i\otimes j(e_i)
		-\partial(b)b_j'\otimes j(e_j'),\]
		\begin{flalign*}
			&\iota((\partial,\tilde{e}))(b (R'\otimes j)(\tilde{e}'))
			-\iota((\partial',\tilde{e}'))(b (R'\otimes j)
			(\tilde{e}))\\
			&=\sum\iota((\partial,\tilde{e}))(bb_j'\otimes j(e_j'))
			-\sum\iota((\partial',\tilde{e}'))(bb_i\otimes j(e_i))\\
			&=\sum \iota((\partial,\tilde{e}))(bb_j'\otimes 1)j(e_j')
			-\sum\iota((\partial',\tilde{e}'))(bb_i\otimes 1)j(e_i)\\
			&=\sum (\partial (bb_j')\otimes 1+bb_j'b_i\otimes j(e_i))j(e_j')-\sum (\partial' (bb_i)\otimes 1+bb_ib_j'\otimes j(e_j'))j(e_i)\\
			&=\partial(b)b_j'\otimes j(e_j')
			+\sum \partial(b_j')b\otimes j(e_j')
			+\sum bb_j'b_i\otimes j(e_i)j(e_j')\\
			&-\partial'(b)b_i\otimes j(e_i)
			-\sum \partial'(b_i)b\otimes j(e_i)
			-\sum bb_ib_j'\otimes j(e_j')j(e_i)\\
			&=\partial(b)b_j'\otimes j(e_j')
			+\sum \partial(b_j')b\otimes j(e_j')
			+\sum bb_ib_j'\otimes j(\left[e_i,e_j'\right])\\
			&-\partial'(b)b_i\otimes j(e_i)
			-\sum \partial'(b_i)b\otimes j(e_i).
		\end{flalign*}
		We now deduce
		\begin{flalign*}
			&\left[\iota((\partial,\tilde{e})),
			\iota((\partial',\tilde{e}'))\right](b)\\
			&=\left[\partial,\partial'\right](b)\otimes 1
			+\sum \partial(b_j')b\otimes j(e_j')
			-\sum \partial'(b_i)b\otimes j(e_i)
			+\sum bb_ib_j'\otimes j(\left[e_i,e_j'\right])\\
			&=\iota(\left[(\partial,\tilde{e}),
			(\partial',\tilde{e}')\right])(b).
		\end{flalign*}
		
		Part (2) is evident by definition.
		
		Part (3) follows from
		\[\begin{split}
			(\iota((0,b\otimes 1_E))\circ
			\iota((\partial,\tilde{e})))(b')
			&=\iota((0,b\otimes 1_E))
			(\partial(b')\otimes 1+b'\tilde{e})\\
			&=b\partial(b')\otimes 1+bb'\tilde{e}\\
			&=\iota((b\partial,b\tilde{e}))(b')
		\end{split}\]
		for $b'\in R'$.
	\end{proof}

	Similarly, we have:
	
	\begin{theorem}\label{thm:fcdotispicalg/sch}
		Let $\Ee,\Aa,j$ be as in Construction \ref{cons:ambientfE}.
		\begin{enumerate}
			\renewcommand{\labelenumi}{(\arabic{enumi})}
			\item The map $\iota:f^\cdot \Ee\to \iEnd_{f^{-1} \Aa}(f^\ast \Aa)$ is a monorphism of sheaves of Lie algebras.
			\item We have $\iota(1_{f^\cdot\Ee})=\id_{f^\ast \Aa}$.
			\item We have
			$\iota((0,b\otimes 1_{\Ee}))\circ \iota((\partial,\tilde{e}))=
			\iota((b\partial,b\tilde{e}))$ for local sections $b\in\OO_X$ and $(\partial,\tilde{e})\in f^\cdot \Ee$.
		\end{enumerate}
	\end{theorem}
	
	\begin{corollary}
		The homomorphism $\iota$ extends to
		\[\begin{array}{cc}
			A_{f_\cdot E}\to \End_{A}(R'\otimes_R A)&(\mathrm{resp.~}\Aa_{f^\cdot \Ee}\to \iEnd_{f^{-1} \Aa}(f^\ast \Aa)).
		\end{array}\]
	\end{corollary}
	
	\begin{proof}
		We only prove the assertion for the ring setting. We check the relations of Definition \ref{def:A_E}. Part (i) (resp.~(ii), (iii)) follows from Theorem \ref{thm:ambientfcdot} (3) (resp.~(1), (2)). This completes the proof.
	\end{proof}
	
	We put $E=\Lie A$ (resp.~$\Ee=\Lie\Aa$) and $j=\id$ to obtain 
	\begin{equation}
		\begin{array}{cc}
			f_\cdot A\to \End_A(A_{R'\gets R})&(\mathrm{resp.~}f^\cdot\Aa\to \iEnd_{f^{-1} \Aa}(f^\ast \Aa))	
		\end{array}
		\label{eq:bimod}
	\end{equation}
	through the proof of Corollary \ref{cor:pavstdo}. We will revisit these actions in section \ref{sec:bimodules}.

	Let us see a basic example of $f_\cdot$ and $f^\cdot$.
	
	\begin{proposition}
		Let $\iota_2:R'\to R\otimes_k R'$ be the canonical homomorphism to the second factor. Let $A$ be a tdo for $R'/k$, we have a natural isomorphism $(\iota_2)_\cdot A\cong D_R\otimes_k A$ of $k$-algebras.
	\end{proposition}
	
	\begin{proof}
		We have canonical $k$-algebra homomorphisms to $(\iota_2)_\cdot A$ from $R$ and $A$. Write $E=\Lie A$. Identify $(\iota_2)_\cdot E$ with $\Theta_{R/k}\otimes_k R' \times R\otimes_k E$ by the canonical splitting
		\[\Theta_{R\otimes_k R'/k}\cong \Theta_{R/k}\otimes_k R'\oplus R\otimes_k \Theta_{R'/k}.\]
		Define $\Theta_{R/k}\to (\iota_2)_\cdot E\subset (\iota_2)_\cdot A$ and $E\to (\iota_2)_\cdot E\subset (\iota_2)_\cdot A$ by $\xi\mapsto (\xi\otimes 1,0)$ and $e\mapsto(0, 1\otimes e)$ respectively. They induce filtered $k$-algebra homomorphisms $D_R\to (\iota_2)_\cdot A$ and $A\to (\iota_2)_\cdot A$. Moreover, $D_R$ and $A$ are commutative in $(\iota_2)_\cdot A$. Therefore we obtain a filtered map
		\[D_R\otimes_k A\to (\iota_2)_\cdot A.\]
		One can show that this is an isomorphism by passing to the associated graded algebras.
	\end{proof}
	
	\begin{proposition}\label{prop:algebroidtdoprojection}
		Let $X$ and $Y$ be smooth $S$-schemes. Let $z=x\circ \pr_1=y\circ \pr_2:X\times_S Y\to S$ denote the structure morphism (recall our notations on projections). For a tdo $\Aa$ on $Y$, there is a unique structure of a tdo on
		\[\OO_{X\times_S Y} \otimes_{\pr_1^{-1}\OO_X\otimes_{z^{-1}\OO_S} \pr_2^{-1}\OO_Y} (\pr_1^{-1}\D_X\otimes_{z^{-1}\OO_S} \pr_2^{-1}\Aa)\]
		with the following properties:
		\begin{enumerate}
			\renewcommand{\labelenumi}{(\roman{enumi})}
			\item The left $\OO_{X\times_S Y}$-module structure on
			\[\OO_{X\times_S Y} \otimes_{\pr_1^{-1}\OO_X\otimes_{z^{-1}\OO_S} \pr_2^{-1}\OO_Y} (\pr_1^{-1}\D_X\otimes_{z^{-1}\OO_S} \pr_2^{-1}\Aa)\]
			coincides with that determined by the left multiplication via the structure homomorphism;
			\item We have
			\begin{flalign*}
				&F_\bullet(\OO_{X\times_S Y} \otimes_{\pr_1^{-1}\OO_X\otimes_{z^{-1}\OO_S} \pr_2^{-1}\OO_Y} (\pr_1^{-1}\D_X\otimes_{z^{-1}\OO_S} \pr_2^{-1}\Aa))\\
				&=\OO_{X\times_S Y} \otimes_{\pr_1^{-1}\OO_X\otimes_{z^{-1}\OO_S} \pr_2^{-1}\OO_Y} (\pr_1^{-1}F_\bullet\D_X\otimes_{z^{-1}\OO_S} \pr_2^{-1}F_\bullet\Aa)
			\end{flalign*}
			as a filtered left $\OO_{X\times_S Y}$-module, where $F_\bullet\D_X$ and $F_\bullet\Aa$ are the given filtrations on $\D_X$ and $\Aa$ respectively;
			\item The canonical map 
			\[\pr_1^{-1}\D_X\otimes_{z^{-1}\OO_S} \pr_2^{-1}\Aa
			\to \OO_{X\times_S Y} \otimes_{\pr_1^{-1}\OO_X\otimes_{z^{-1}\OO_S} \pr_2^{-1}\OO_Y} (\pr_1^{-1}\D_X\otimes_{z^{-1}\OO_S} \pr_2^{-1}\Aa)\]
			is a homomorphism of sheaves of rings.
		\end{enumerate}
		Moreover, there is a canonical isomorphism
		\[\OO_{X\times_S Y} \otimes_{\pr_1^{-1}\OO_X\otimes_{z^{-1}\OO_S} \pr_2^{-1}\OO_Y} (\pr_1^{-1}\D_X\otimes_{z^{-1}\OO_S} \pr_2^{-1}\Aa)\cong\pr^\cdot_2\Aa\]
		of tdos. A similar assertion holds for
		\[ (\pr_1^{-1}\D_X\otimes_{z^{-1}\OO_S} \pr_2^{-1}\Aa)\otimes_{\pr_1^{-1}\OO_X\otimes_{z^{-1}\OO_S} \pr_2^{-1}\OO_Y}\OO_{X\times_S Y} .\]
	\end{proposition}

	We also see how the pullback behaves under basic operations of Picard algebroids and tdos:
	\begin{proposition}\label{prop:pullbackvsscalartwist}
		Let $\Ee$ be a Picard algebroid on $Y$, and $\gamma\in(y^{-1}\OO_S)(S)^\times$. We denote the image of $\gamma$ in $(x^{-1}\OO_S)(S)$ by the same symbol $\gamma$. Then we have $f^\cdot \Ee^\gamma =(f^\cdot\Ee)^\gamma$.
	\end{proposition}
	
	\begin{proof}
		This is evident by definitions.
	\end{proof}
	
	\begin{proposition}[{\cite[Lemma 2.11.1]{kashiwara1989}}]\label{prop:pullbacksidechange}
		Let $\Aa$ be a tdo on $Y$. Then there is a natural isomorphism 
		\[\omega_{X/S}^\vee\otimes_{\OO_X} f^\cdot(\omega_{Y/S}^\vee\otimes_{\OO_Y}\Aa^{\op}\otimes_{\OO_Y}
		\omega_{Y/S})^{\op} \otimes_{\OO_X} \omega_{X/S}
		\cong
		f^\cdot\Aa.\]
	\end{proposition}
	
	\begin{proof}
		This follows from Proposition \ref{prop:pullbackvsscalartwist} and the twice use of Corollary \ref{cor:oppositecocyle}:
		\[\begin{split}
			\omega_{X/S}^\vee\otimes_{\OO_X} f^\cdot(\omega_{Y/S}^\vee\otimes_{\OO_Y}\Aa^{\op}\otimes_{\OO_Y}
			\omega_{Y/S})^{\op} \otimes_{\OO_X} \omega_{X/S}
			&\cong \omega_{X/S}^\vee\otimes_{\OO_X} f^\cdot(\Aa^{-1})^{\op} \otimes_{\OO_X} \omega_{X/S}\\
			&\cong f^\cdot(\Aa^{-1})^{-1}\\
			&\cong f^\cdot\Aa.
		\end{split}\]
	\end{proof}
	
	\begin{proposition}\label{prop:pullbacklinebundlecompatible}
		Let $\Aa$ be a tdo on $Y$, and $\Ll$ be a line bundle on $Y$. Then there is a natural isomorphism
		$f^\cdot(\Ll\otimes_{\OO_Y}\Aa\otimes_{\OO_Y}\Ll^\vee)
		\cong
		f^\ast\Ll\otimes_{\OO_X} f^\cdot\Aa
		\otimes_{\OO_X} (f^\ast\Ll)^\vee$.
	\end{proposition}

	\begin{proof}
		Define an action of $(f^\ast\Ll)\otimes_{\OO_X} f^\cdot\Aa
		\otimes_{\OO_X} (f^\ast\Ll)^\vee$ on
		\[f^\ast(\Ll\otimes_{\OO_Y}\Aa\otimes_{\OO_Y}\Ll^\vee)
		\cong f^\ast\Ll\otimes_{\OO_X} f^\ast\Aa 
		\otimes_{f^{-1}\OO_Y} f^{-1}\Ll\]
		by
		\begin{flalign*}
			&(f^\ast\Ll)\otimes_{\OO_X} f^\cdot\Aa
			\otimes_{\OO_X} (f^\ast\Ll)^\vee
			\otimes_{f^{-1}\OO_S} f^\ast\Ll\otimes_{\OO_X} f^\ast\Aa \otimes_{f^{-1}\OO_Y} f^{-1}\Ll\\
			&\to (f^\ast\Ll)\otimes_{\OO_X} f^\cdot\Aa
			\otimes_{f^{-1}\OO_S}
			f^\ast\Aa \otimes_{f^{-1}\OO_Y} f^{-1}\Ll\\
			&\to (f^\ast\Ll)\otimes_{\OO_X} 
			f^\ast\Aa \otimes_{f^{-1}\OO_Y} f^{-1}\Ll,
		\end{flalign*}
		where the first arrow is given by the canonical pairing of $(f^\ast\Ll)^\vee$ and $f^\ast\Ll$, and the second one is obtained by the left action of $f^\cdot\Aa$ on $f^\ast\Aa$ defined by \eqref{eq:bimod}. This induces a Picard algebroid homomorphism
		\[\Lie 
		(f^\ast\Ll\otimes_{\OO_X} f^\cdot\Aa
		\otimes_{\OO_X} (f^\ast\Ll)^\vee)
		\to
		\Lie (f^\cdot(\Ll\otimes_{\OO_Y}\Aa\otimes_{\OO_Y}\Ll^\vee))\]
		(recall the construction of the Lie bracket of $f^\cdot\Lie\Aa$). Pass to the equivalence of Corollary \ref{cor:pavstdo} to deduce the desired isomorphism.
	\end{proof}

	\subsubsection{Compatibility}\label{sec:compatibility}
	
	The aim of this section is to prove the following result:
	
	\begin{theorem}\label{thm:basechangevspullbackfortdo}
		Let $\Aa$ be a tdo on $Y$. Then there is a natural isomorphism
		\[s^\ast_X f^\cdot \Aa\cong (f')^\cdot s^\ast_Y\Aa.\]
	\end{theorem}
	Due to the compatibility of the two pullbacks with the equivalence of tdos and Picard algebroids, we may prove a similar assertion for Picard algebroids. Let $\Ee$ be a Picard algebroid. Then we have a canonical $\OO_{X'}$-module homomorphism
	\begin{equation}
		s^\ast_X f^\cdot \Ee
		=s^\ast_X (\Theta_{X/S}\times_{f^\ast\Theta_{Y/S}}f^\ast\Ee)
		\to s^\ast_X\Theta_{X/S}\times_{s^\ast_X f^\ast\Theta_{Y/S}}
		s^\ast_X f^\ast\Ee.
		\label{eq:canhomofbcpullback}
	\end{equation}
	The target in the above map is naturally identified with
	\begin{equation}
		\Theta_{X'/S'} \times_{(f')^\ast \Theta_{Y'/S'}} 
		(f')^\ast s^\ast_Y \Ee=(f')^\cdot s^\ast_Y\Ee
		\label{eq:identificationwithpullbackbasechange}
	\end{equation}
	(use $s_Y\circ f'=f\circ s_X$). The map \eqref{eq:canhomofbcpullback} is a morphism of Picard algebroids under the identification \eqref{eq:identificationwithpullbackbasechange}. In fact, we may reduce to the ring setting through the affine setting since the constructions are local in $X$ and $S$. Then the assertion is obvious.

	\subsection{Equivariant tdos}\label{sec:equivtdo}
	Suppose that we are given an action $a:G\times_S X\to X$ of a smooth group $S$-scheme on $X$. Let $\pr_1,\pr_2,\pr_{23},m_G,i$ be as in section \ref{sec:equivsheaf}.
	
	\begin{definition}[{\cite[4.6]{kashiwara1989}}]
		An equivariant tdo on $X$ is a tdo $\Aa$ on $X$, equipped with an isomorphism $I:a^\cdot \Aa\cong \pr_2^\cdot \Aa$ such that the equalities
		\[\begin{array}{cc}
			\pr_{23}^\cdot I\circ(\id_G\times a)^\cdot I=(m_G\times \id_X)^\cdot I,
			&i^\cdot I=\id_{\Aa}
		\end{array}\]
		hold.
	\end{definition}
	
	We define equivariant Picard algebroids in a similar way. Since the equivalence of the groupoids of tdos and Picard algebroids respects the pullback functor, this equivalence extends to that of the equivariant objects.
	
	\begin{example}
		The untwisted Picard algebroid $\tilde{\Theta}_{X/S}$ is $G$-equivariant for the trivial action on $\OO_X$ and the action on $\Theta_{X/S}$ in section A.3.
	\end{example}
	
	\begin{example}
		Let $\Ll$ be a line bundle on $X$ with a $G$-equivariant structure $I_{\Ll}:a^\ast\Ll\cong\pr_2^\ast\Ll$, and $\Aa$ be a tdo on $X$ with a $G$-equivariant structure $I_{\Aa}:a^\cdot\Aa\cong\pr_2^\cdot\Aa$. Then $\Ll\otimes_{\OO_X} \Aa\otimes_{\OO_X} \Ll^\vee$ is $G$-equivariant for the composite isomorphism
		\[\begin{split}
			a^\cdot (\Ll\otimes_{\OO_X} \Aa\otimes_{\OO_X} \Ll^\vee)
			&\cong (a^\ast\Ll)\otimes_{\OO_{G\times_SX}} a^\cdot \Aa
			\otimes_{\OO_{G\times_SX}}(a^\ast\Ll^\vee)\\
			&\cong (\pr_2^\ast\Ll)\otimes_{\OO_{G\times_SX}}
			\pr_2^\cdot \Aa
			\otimes_{\OO_{G\times_SX}}(\pr_2^\ast\Ll^\vee)\\
			&\cong 
			\pr_2^\cdot (\Ll\otimes_{\OO_X} \Aa\otimes_{\OO_X} \Ll^\vee).
		\end{split}\]
		Here the first and last isomorphisms are obtained by Corollary \ref{prop:pullbacklinebundlecompatible}. The second isomorphism is given by application of $I_{\Aa}$ and $I_{\Ll}$.
	\end{example}
	
	\begin{example}[Restriction of groups]
		Let $K$ be a smooth group $S$-scheme, equipped with a homomorphism $f:K\to G$, and $\Aa$ be a tdo on $X$ with a $G$-equivariant structure $I$. Define an action of $K$ on $X$ by the restriction of $a$ along $f$. Then $(f\times \id_X)^\cdot I$ determines the structure of a $K$-equivariant tdo on $\Aa$.
	\end{example}

	\subsection{Torsion free case}
	In this section, we explain how we can weaken the definition of tdos on $X$ when $X$ is flat over $\ZZ$ and how we can realize tdos in a reasonable way for such $X$. By the Zariski localization in $X$, we may and do work with the ring setting through Corollaries \ref{cor:localization_PA} and \ref{cor:pavstdo}. Let $k\to R$ be a smooth morphism. In this section, assume that $R$ is flat over $\ZZ$.
	
	\begin{lemma}\label{lem:bracket}
		Let $A$ be an almost commutative filtered $k$-algebra, and $n\geq 1$. Then for $P_1,\ldots, P_n\in F_1 A$ and $a\in F_0A$, we have
		\[\sum_{j=1}^n \left[P_j,a\right]
		P_1\overset{\overset{j}{\vee}}{\cdots}P_n\equiv
		\left[P_1\cdots P_n,a\right]\pmod{F_{n-2} A}.\]
	\end{lemma}
	
	\begin{proof}
		We prove it by induction on $n\geq 1$. The assertion is clear if $n=1$. For $n\geq 2$, we have
		\[\begin{split}
			\left[P_1\cdots P_n,a\right]
			&=P_1\cdots P_n a-aP_1\cdots P_n\\
			&=P_1\cdots P_n a+\left[P_1,a\right]P_2\cdots P_n
			-P_1aP_2\cdots P_n\\
			&=\left[P_1,a\right]P_2\cdots P_n
			+P_1\left[P_2\cdots P_n,a\right]\\
			&\equiv \left[P_1,a\right]P_2\cdots P_n
			+\sum_{j=2}^n P_1\left[P_j,a\right]
			P_2\overset{\overset{j}{\vee}}{\cdots}P_n
			\pmod{F_{n-2} A}
		\end{split}\]
		where the last congruence in the forth row follows from the induction hypothesis. Since $A$ is almost commutative, we also have $P_1\left[P_j,a\right]
		\equiv\left[P_j,a\right]P_1 \pmod{F_0A}$ for $2\leq j\leq n$. Hence we obtain
		\[\left[P_1\cdots P_n,a\right]
		\equiv \left[P_1,a\right]P_2\cdots P_n
		+\sum_{j=2}^n P_1\left[P_j,a\right]
		P_2\overset{\overset{j}{\vee}}{\cdots}P_n\\
		\equiv \sum_{j=1}^n \left[P_j,a\right]
		P_1\overset{\overset{j}{\vee}}{\cdots}P_n
		\pmod{F_{n-2}A}.\]
	\end{proof}
	
	\begin{proposition}[{\cite[Lemma 2.1.2]{beilinsonbernstein1993}}]\label{prop:pssymbol}
		Let $(A,F_\bullet A)$ be an almost commutative filtered $k$-algebra, equipped with a $k$-algebra homomorphism $i:R\to A$. Then $(A,F_\bullet A,i)$ is a tdo if and only if the following conditions are satisfied:
		\begin{enumerate}
			\renewcommand{\labelenumi}{(\roman{enumi})}
			\item $F_\bullet A$ is exhaustive;
			\item the map $i$ is an isomorphism onto $F_0A$;
			\item The map $\sigma:\Gr^1 A\to\Theta_{R/k};~P\mapsto \left[P,-\right]$ is an isomorphism;
			\item $\Gr A$ is generated by $\Gr^1 A$ as an $R$-algebra.
		\end{enumerate}
	\end{proposition}
	
	The ``only if'' direction is clear. For the converse, let us define a graded $R$-algebra homomorphism $\Sigma=(\Sigma^n):\Sym_R \Theta_{R/k}\to \Gr A$
	by
	$\Theta_{R/k}\cong\Gr^1 A\subset \Gr A$. For the proof of the ``if'' direction and applications to latter arguments, let us verify a stronger result:
	
	\begin{lemma}\label{lem:Sigma,Pi}
		Assume that the conditions (i)-(iv) in Proposition \ref{prop:pssymbol} are satisfied. Let $n\geq 0$.
		\begin{enumerate}
			\item The map $\Sigma$ is an isomorphism.
			\item For $P\in F_n A$, define $\tilde{\Pi}^n(P)$ by the image of $\left[P,-\right]$ under the composite map
			\[\begin{split}
				\mathrm{Der}_{k}(R,\Gr^{n-1}A)
				&\cong\Gr^{n-1} A\otimes_R\mathrm{Der}_k(R,R)\\
				&\xcong{(\Sigma^{n-1})^{-1}\otimes_R \Theta_{R/K}}\
				\Sym^{n-1}_R \Theta_{R/k}\otimes_R \Theta_{R/k}\\
				&\to \Sym^n_R \Theta_{R/k},
			\end{split}\]
			where the first isomorphism follows since $R$ is smooth:
			\[\begin{tikzcd}
				\Hom_R(\Omega^1_{R/k},\Gr^{n-1} A)\ar[r, "\sim"]
				&\mathrm{Der}_{k}(R,\Gr^{n-1}A)\\
				\Gr^{n-1} A\otimes_R\Hom_R(\Omega^1_{R/k},R)
				\ar[r,"\sim"]\ar[u, "\sim"{sloped, above}]
				&\Gr^{n-1} A\otimes_R\mathrm{Der}_k(R,R).\ar[u]
			\end{tikzcd}\]
			The right vertical arrow is given by $P\otimes \partial\mapsto \partial(-)P$. The resulting map
			\[\tilde{\Pi}^n:F_nA\to \Sym^n_R \Theta_{R/k}\]
			descends to $\Pi^n:\Gr^n A\to \Sym^n_R \Theta_{R/k}$.
			\item We have $\Pi^n\circ \Sigma^n=n\id_{\Sym^n_{R} \Theta_{R'/k}}$.
		\end{enumerate}		
	\end{lemma}
	
	\begin{proof}
		Since $\Gr A$ is generated by $\Gr^1 A$, $\Sigma^n$ is surjective for each $n\geq 0$ (recall the definition of $\Sigma$). We prove the rest simultaneously by induction on $n$. It follows from the condition (ii) that $\Sigma^0$ is an isomorphism, and that $\tilde{\Pi}^0=0$. In particular, (2) and (3) hold for $n=0$.
		
		Assume $n\geq 1$. Then $\tilde{\Pi}^n$ is well-defined by the induction hypothesis. Part (2) follows since $\left[P,-\right]=0$ as an element of $\mathrm{Der}_k(R,\Gr^{n-1} A)$ for $P\in F_{n-1}A$. Part (3) is a consequence of Lemma \ref{lem:bracket}. Finally, we prove that $\Sigma^n$ is injective. Since $R$ is smooth, $\Theta_{R/k}$ is finitely generated and projective. Therefore $\Sym_R\Theta_{R/k}$ is flat (smooth) over $R$ (check Zariski locally on $\Spec R$). Since $R$ is flat over $\ZZ$, so is $\Sym_R\Theta_{R/k}$. The equality $\Pi^n\circ \Sigma^n=n\id_{\Sym^n_{R} \Theta_{R'/k}}$ now implies that $\Sigma^n$ is injective. This completes the proof.
	\end{proof}

	The ``if'' direction of Proposition \ref{prop:pssymbol} now follows from Lemma \ref{lem:Sigma,Pi} (1) (recall the definition of $\Sigma$).

	\begin{corollary}[{\cite[Proposition 2.3.2]{kashiwara1989}, \cite[2.1.1.~Remark]{beilinsonbernstein1993}}]\label{cor:tdofilt=maxDfilt}
		Let $(A,F_\bullet A,i)$ be a tdo. Then $F_\bullet A$ coincides with the maximal $D$-filtration.
	\end{corollary}
	
	\begin{proof}
		Recall that $F_n A\subset F_n^{\vee} A$ for $n\geq 0$ (Lemma \ref{lem:F<Fvee}). The proof is completed by showing $F_n^{\vee} A\cap F_m A\subset F_nA$ (equivalently, $F_n^{\vee}A\cap F_mA= F_nA$) for $m\geq n\geq 0$. In fact, let $P\in F^{\vee}_n A$. Since $F_\bullet A$ is exhaustive, there exists an element $m\geq n$ such that $P\in F_m A$. We deduce $P\in F_n A$ by $F_n^{\vee} A\cap F_m A= F_n A$.
		
		The assertion is clear if $n=m$. We prove the containment by induction on $s=m+n$. The assertion follows if $s=0$ since $s=0$ implies $n=m=0$. Let $s\geq 1$. We may assume $n<m$. Let $P\in F_n^{\vee} A\cap F_m A$ and $a\in R$. Then we have $\Pi^m(P)=0$. In fact, if $n=0$, we have $\left[P,i(a)\right]=0$; otherwise, i.e., if $n>0$, then $\left[P,i(a)\right]$ belongs to $F_{n-1}^{\vee} A\cap F_{m-1} A$. Since $0\leq n-1\leq m-1$ and $0\leq (n-1)+(m-1)=s-2<s$, the induction hypothesis shows $\left[P,i(a)\right]\in F_{n-1} A\subset F_{m-2} A$.
		
		Observe that
		\[\begin{split}
			0&=\Sigma^m(\Pi^m(P))\\
			&=\Sigma^m(\Pi^m(\Sigma^m((\Sigma^m)^{-1}(P))))\\
			&=m\Sigma^m((\Sigma^m)^{-1}(P))\\
			&=mP\in \Gr^m A
		\end{split}\]
		(recall that we verified in the proof of Proposition \ref{prop:pssymbol} that $\Sigma^m$ is an isomorphism). Since $\Sym_R\Theta_{R/k}$ is flat over $\ZZ$, we get $P=0$ in $\Gr^m A$. Hence $P$ belongs to $F_{m-1} A$. The induction hypothesis implies $P\in F_n^{\vee} A\cap F_{m-1} A=F_n A$ as desired.
	\end{proof}
	
	\begin{corollary}[{\cite[Proposition 2.3.2]{kashiwara1989}, \cite[2.1.2.~Lemma]{beilinsonbernstein1993}}]\label{cor:tdobymdfilt}
		Let $(A,i)$ be a $D$-algebra. Put the maximal $D$-filtration on $A$. Then $A$ is a tdo if and only if the following conditions are satisfied:
		\begin{enumerate}
			\renewcommand{\labelenumi}{(\roman{enumi})}
			\item $i$ is an isomorphism onto $F_0^{\vee} A$;
			\item The map $\sigma:\Gr^1 A\to\Theta_{R/k};~e\mapsto \left[e,-\right]$ is surjective;
			\item $\Gr A$ is generated by $\Gr^1 A$ as an $R$-algebra.
		\end{enumerate}
	\end{corollary}
	
	\begin{proof}
		We have already proved the ``only if'' direction. Suppose that the above conditions (i)-(iii) are satisfied. We prove that $F_\bullet^{\vee}$ is a tdo filtration. The only nontrivial part is to show that $\sigma$ is injective. Let $e\in F_1^{\vee} A$ with $\sigma(e)=0$. This says $\left[e,i(a)\right]=0$ in $A$ for all $a\in R$. Hence $e$ belongs to $F_0^{\vee} A$ as desired.
	\end{proof}
	
	\begin{corollary}
		Let $(A,i), (A',i')$ be tdos. Then a $k$-algebra homomorphism
		$\varphi:A\to A'$
		is a morphism of tdos if and only if $\varphi\circ i=i'$.
	\end{corollary}
	
	\begin{corollary}\label{cor:generatetdo}
		Let $(E,\sigma,1_E)$ be a Picard algebra, $A$ be a (possibly noncommutative) $k$-algebra, equipped with a $k$-algebra homomorphism $j:R\to A$, and $\iota:E\to A$ be a $k$-module homomorphism. Suppose that these satisfy the following conditions:
		\begin{enumerate}
			\renewcommand{\labelenumi}{(\roman{enumi})}
			\item $\iota$ is an $R$-module homomorphism;
			\item $\iota$ is a Lie algebra homomorphism over $k$;
			\item $\iota(1_E)=1$;
			\item $j$ is an isomorphism onto $F^{\vee}_0 A$;
			\item $\iota$ is an isomorphism onto $F^{\vee}_1 A$;
		\end{enumerate}
		Then $A_E$ is isomorphic onto the $k$-subalgebra $A_{E,\iota}$ of $A$ generated by $\iota(E)\subset A$.
	\end{corollary}
	
	\begin{proof}
		The map $\iota$ induces a $k$-algebra homomorphism $\iota:A_E\to A$ from the conditions (i)--(iii). Note that the first relation follows as
		\[\iota(a1_E)\iota(e)=(j(a)\iota(1_E))\iota(e)
		=(j(a)1)\iota(e)=j(a)\iota(e)=\iota(ae)\]
		by definition of the $R$-module structure on $A$.
		
		Put the image filtration $F_\bullet A_{E,\iota}$ on $A_{E,\iota}$ along $\iota$. Then we have commutative diagrams
		\begin{equation}
			\begin{tikzcd}
				R\ar[r, "\sim"]\ar[rrr, bend left, "j"]\ar[d, "1_E"', hook]
				&F_0A_E\ar[r, two heads, "\iota"]\ar[d, hook]
				&F_0A_{E,\iota}\ar[r, hook]\ar[d, hook]&D\\
				E\ar[r, "\sim"]\ar[rrru, bend right=50, "\iota"']&F_1 A_E\ar[r, two heads, "\iota"]
				&F_1A_{E,\iota}\ar[ru, hook]
			\end{tikzcd}\label{diag:compatible_embedding}
		\end{equation}
		(see the proof of Theorem \ref{thm:DEistdo} for the left horizontal arrows). Since $j$ and $\iota$ are injective by (iv) and (v) respectively, the middle horizontal arrows are isomorphisms. Hence $\iota$ induces isomorphisms
		\[\begin{array}{cc}
			F_0 A_E\cong F_0 A_{E,\iota}=\iota(R1_E)=j(R),
			&E\cong F_1 A_E\cong F_1 A_{E,\iota}=\iota(E).
		\end{array}\]
		Since $R$ is commutative, we have $F_0 A_{E,\iota}\subset F_0^{\vee} A_{E,\iota}$. The condition (iii) implies the converse containment
		\[F_0^{\vee} A_{E,\iota}\subset F_0^{\vee} A=j(R)
		=F_0 A_{E,\iota}.\]
		Similarly, (iv) implies
		\[F_1^{\vee} A_{E,\iota}\subset F_1^{\vee} A=\iota(E)
		=F_1 A_{E,\iota}.\]
		The converse containment holds immediately by definition of Picard algebras. In fact, we have
		\[\left[\left[\iota(e),j(a)\right],j(a')\right]
		=\left[\left[\iota(e),\iota(a1_E)\right],\iota(a'1_E)\right]
		=\iota(\left[\left[e,a1_E\right],a'1_E\right])
		=0\]
		for $e\in E$ and $a,a'\in R$. For the first equality, recall that $ A$ is regarded as an $R$-module for the multiplication via $j$. Hence we get $\iota(a1_E)=j(a)\iota(1_E)=j(a)1=j(a)$.
		
		In view of the right square in \eqref{diag:compatible_embedding}, $\Gr^i\iota:\Gr^i A_E\to \Gr^i A_{E,\iota}$ ($i\in\{0,1\}$) are isomorphisms since so are its horizontal arrows.
		Since $\Gr A_E$ is generated by $\Gr^1 A_E$ over $R$,
		$\Gr A_{E,\iota}$ is generated by $\Gr^1 A_{E,\iota}$. It also implies that the map $\Gr^1 A_{E,\iota}\to \Theta_{R/k}$ in Lemma \ref{lem:commutator} is an isomorphism (compare it with $\Gr^1 A_E\cong \Theta_{R/k}$).
		We now deduce that $A_{E,\iota}$ is a tdo (Proposition \ref{prop:pssymbol}). Moreover, $\iota$ is an isomorphism onto $A_{E,\iota}$ by Remark \ref{rem:tdo_is_groupoid}.
	\end{proof}
	
	\begin{corollary}
		Let $\End_k(R)$ be the $k$-algebra of endomorphisms of $R$ as a $k$-module. Then $D_R$ is isomorphic to the $k$-subalgebra of $\End_k(R)$ generated by the canonical actions of $R$ and $\Theta_{R/k}$ on $R$.
	\end{corollary}
	\begin{proposition}[{\cite[section 1.4]{beilinsonbernstein1993}}]
		Let $(A,i)$ be a tdo for $R/k$. Define a $k$-algebra homomorphism $j:R'\to \End_R(R'\otimes_R A)$ by the scalar multiplication from the left side. Then $f_\cdot A$ is isomorphic to the $k$-subalgebra $\End_R(R'\otimes_R A)$ generated by $F_1^{\vee}\End_A(R'\otimes_R A)$.
	\end{proposition}
	\begin{proof}
		We check the conditions of Corollary \ref{cor:generatetdo} for the map $\iota:f_\cdot \Lie A\to\End_A(R'\otimes_R A)$ in Construction \ref{cons:ambientfE}.
		
		The map $\iota$ is a Lie algebra homomorphism by definition of $f_\cdot E$. It follows by construction that $\iota(1_{f_\cdot E})$ is the identity map. To see that $\iota$ is $R'$-linear, let $(\partial,\tilde{e})\in f_\cdot \Lie A$ and $b,b'\in R'$. Then we have
		\[\begin{split}
			j(b)\circ \iota((\partial,\tilde{e}))(b'\otimes 1)
			&=j(b)(\partial(b')\otimes 1+b'\tilde{e})\\
			&=b\partial(b')\otimes 1+bb'\tilde{e}\\
			&=\iota(b(\partial,\tilde{e}))(b'\otimes 1).
		\end{split}\]
		This shows $\iota(b(\partial,\tilde{e}))=j(b)\iota((\partial,\tilde{e}))$.
		
		Unwinding the definitions, $j$ factors as
		\[\begin{tikzcd}
			R'\ar[rr, "j"]\ar[rd, "1_{f_\cdot \Lie A}"']
			&&\End_{A}(R'\otimes_R A)\\
			&f_\cdot\Lie A.\ar[ru, "\iota"']
		\end{tikzcd}\]
		Since the two diagonal arrows are injective (see Theorem \ref{thm:ambientfcdot} and Proposition \ref{prop:pullbackofpicalg}), so is $j$. It is clear by definition that $F^{\vee}_0 \End_{A}(R'\otimes_R A)$ consists of the endomorphisms of $R'\otimes_R A$ which are both $R'$-linear and $A$-linear. Hence $j$ injectively maps into
		$F^{\vee}_0 \End_{A}(R'\otimes_R A)$. Passing to the adjunctions, we can identify
		$F^{\vee}_0 \End_{A}(R'\otimes_R A)$ with the set of elements $\tilde{e}\in R'\otimes_R A$ satisfying $a\tilde{e}=\tilde{e}a$ for all $a\in R$ through $T\mapsto T(1\otimes 1)$. Under this identification, $j$ is given by
		\[b\mapsto b\otimes 1\in R'\otimes_R F_0^{\vee} A\cong R'.\]
		To see that $j$ is onto $F^{\vee}_0 \End_{A}(R'\otimes_R A)$, let $\tilde{e}\in F^{\vee}_0 \End_{A}(R'\otimes_R A)\subset R'\otimes_R A$. Recall that $R'\otimes_R F_\bullet^{\vee} A$ is an exhaustive filtration of $R'\otimes_R A$ by the proof of Theorem \ref{thm:ambientfcdot} (1). There is a nonnegative integer $n\geq 0$ such that $\tilde{e}\in R'\otimes_R F^{\vee}_n A$. We prove by induction on $n\geq 0$ that $\tilde{e}$ belongs to $R'\otimes_R F_0^{\vee} A\cong R'$. The assertion is clear if $n=0$. Let $n\geq 1$. Write $\tilde{e}=\sum b_i\otimes P_i$ with $b_i\in R'$ and $P_i\in F_n^{\vee} A$. Since $\tilde{e}\in F^{\vee}_0 \End_{A}(R'\otimes_R A)$, we have
		\[\sum b_i\otimes \left[P_i,i(a)\right]=0.\]
		The below diagram shows $(R'\otimes_R \Pi^n)(\sum b_i\otimes P_i)=0$:
		\[\begin{tikzcd}
			R'\otimes_R F^{\vee}_n A_E\ar[d]\ar[rd]&\\
			R'\otimes_R\mathrm{Der}(R,\Gr^{n-1} A_E)\ar[r, "\sim"]
			&\mathrm{Der}_k(R,R'\otimes_R \Gr^{n-1} A_E)\\
			R'\otimes_R \Hom_R(\Omega^1_{R/k},\Gr^{n-1} A_E)
			\ar[r, "\sim"]\ar[u, "\sim"{sloped, above}]
			&\Hom_R(\Omega^1_{R/k},R'\otimes_R \Gr^{n-1} A_E).
			\ar[u, "\sim"{sloped, below}]
		\end{tikzcd}\]
		Notice that $R'\otimes_R \Sym_R\Theta_{R/k}\cong\Sym_{R'} (R'\otimes_R \Theta_{R/k})$ is flat over $\ZZ$ by our hypothesis. Since $(R'\otimes_R \Pi^n)\circ (R'\otimes_R \Sigma^n)=n\id$ and $\Sigma^n$ is an isomorphism,
		$R'\otimes_R \Pi^n$ is injective. Therefore $\sum b_i\otimes P_i$ belongs to $R'\otimes_R F^{\vee}_{n-1} A_E$. The induction hypothesis shows $\tilde{e}\in R'\otimes_R F^{\vee}_0 A$. 
		
		Recall again that $\iota$ is injective (Theorem \ref{thm:ambientfcdot} (1)). Let $(\partial,\tilde{e})\in f_\cdot \Lie A$. Then $\left[\iota((\partial,\tilde{e})),j(b')\right]$ belong to $j(R')$ for all $b'\in R'$ by
		\[\begin{split}
			\left[\iota((\partial,\tilde{e})),j(b')\right](b\otimes P)
			&=\iota((\partial,\tilde{e})(b'b\otimes P)-(j(b')\iota((\partial,\tilde{e}))(b\otimes P)\\
			&=\partial(b'b)\otimes P+bb'\tilde{e}P
			-b'\partial(b)\otimes P-b'b\tilde{e}P\\
			&=\partial(b')b\otimes P\\
			&=j(\partial(b'))(b\otimes P).
		\end{split}\]
		for $b\otimes P\in R'\otimes_R A$. This shows $\iota((\partial,\sum b_i\otimes e_i))\in F_1^{\vee}\End_R(R'\otimes_R A)$.
		
		The proof is completed by showing that $\iota$ is onto $F_1^{\vee}\End_A(R'\otimes_R A)$. Suppose that we are given an operator $T\in F_1^{\vee}\End_A(R'\otimes_R A)$. Then we obtain an element $\partial\in\Theta_{R'/k}$ by Lemma \ref{lem:commutator}. Choose a positive integer $n$ such that $T(1\otimes 1)\in R'\otimes_R F^{\vee}_n A$. One can prove $T(1\otimes 1)\in R'\otimes_R \Lie A$ by induction on $n$. In fact, the assertion is clear if $n=1$. Let $n\geq 2$. For $a\in R$, we have
		\[\begin{split}
			T(1\otimes 1)i(a)
			&=T(1\otimes i(a))\\
			&=T(f(a)\otimes 1)\\
			&=(T\circ j(f(a)))(1\otimes 1)\\
			&=\left[T,j(f(a))\right](1\otimes 1)+(j(f(a))T)(1\otimes 1)
		\end{split}\]
		since $T$ is $A$-linear. Since $\left[T,j(f(a))\right]\in F^{\vee}_0 \End_{A}(R'\otimes_R A)=j(R')$, one can write
		\[\left[T,j(f(a))\right](1\otimes 1)=b_a\otimes 1\]
		for some $b_a\in R'$. If we express $T(1\otimes 1)=\sum b_i\otimes P_i$, we have
		\[R'\otimes_R F^{\vee}_0 A\ni b_a\otimes 1=\sum b_i\otimes \left[P_i,i(a)\right].\]
		Since $n\geq 2$, we get $\Pi^n(\sum b_i\otimes P_i)=0$. A similar argument to the above implies
		\[\sum b_i\otimes P_i\in R'\otimes_R F^{\vee}_1 A.\]
		
		Write $T(1\otimes 1)=\sum b_i\otimes e_i$. The pair $(\partial,T(1\otimes 1))$ belongs to $f_\cdot E$ from
		\[(f^\ast\partial)(a)=
		(f^\ast (j^{-1}(\left[T,j(-)\right]))(a)
		=j^{-1}(\left[T,j(f(a))\right]),\]
		\[\begin{split}
			\left[T,j(f(a))\right](1\otimes 1)
			&=T(1\otimes 1)i(a)-aT(1\otimes 1)\\
			&=\sum b_i\otimes \left[e_i,i(a)\right]\\
			&=\sum b_if(i^{-1}(\left[e_i,i(a)\right]))\otimes 1\\
			&=\sum b_i f(\sigma(e_i)(a))\otimes 1
		\end{split}\]
		for $a\in R$. Finally, $\iota((\partial, T(1\otimes 1)))=T$ follows as
		\[\begin{split}
			\iota((\partial, T(1\otimes 1)))(b\otimes 1)
			&=\partial(b)\otimes 1 +bT(1\otimes 1)\\
			&=j^{-1}(\left[T,j(b)\right])\otimes 1+(j(b)T)(1\otimes 1)\\
			&=\left[T,j(b)\right](1\otimes 1)+(j(b)T)(1\otimes 1)\\
			&=(Tj(b))(1\otimes 1)\\
			&=T(b\otimes 1).
		\end{split}\]

	\end{proof}

	\section{Generalities on twisted $\mathcal D$-modules}\label{sec:Dmodules}
	
	Twisted $\D$-modules are literally modules over tdos. In the last section, we introduced the notions of tdos and Picard algebroids. In this section, we see some elementary examples of twisted $\D$-modules. Then we give categorical and homological study of twisted (quasi-coherent and coherent) $\D$-modules.
	
	Let $k\to R$ be a smooth homomorphism of commutative rings, and $x:X\to S$ be a smooth morphism of schemes.
	
	\subsection{Twisted $\D$-modules}\label{sec:twistedd-mod}
	\begin{definition}[Modules over Picard algebras]\label{def:moduleoverpicardalgebra}
		Let $(E,[-,-],\sigma,1_E)$ be a Picard algebra. A left $E$-module is an $R$-module $M$ endowed with a $k$-linear map $E\otimes_k M\to M,\,e\otimes m\mapsto em$ satisfying the four conditions
		\begin{enumerate}
			\renewcommand{\labelenumi}{(\roman{enumi})}
			\item $\forall a\in R,\,e\in E,\,m\in M:$ $a(em)=(ae)m$,
			\item $\forall a\in R,\,e\in E,\,m\in M:$ $e(am)=(ae)m+\sigma(e)(a)m$,
			\item $\forall e,e'\in E,\,m\in M:$ $[e,e']m=e(e'm)-e'(em)$,
			\item $\forall m\in M:$ $1_Em=m$.
		\end{enumerate}
		Likewise, a right $E$-module is an $R$-module $M$ endowed with a $k$-linear map $M\otimes_k E\to M$ satisfying
		\begin{enumerate}
			\renewcommand{\labelenumi}{(\roman{enumi})}
			\item $\forall a\in R,\,e\in E,\,m\in M:$ $(ma)e=m(ae)$,
			\item $\forall a\in R,\,e\in E,\,m\in M:$ $(me)a=m(ae)+m\sigma(e)(a)$,
			\item $\forall e,e'\in E,\,m\in M:$ $m[e,e']=(me)e'-(me')e$,
			\item $\forall m\in M:$ $m=m1_E$.
		\end{enumerate}
		Morphisms of $E$-modules are defined in the usual way.
	\end{definition}
	
	Then the categories of left (resp.\ right) $E$-modules and left (resp.\ right) $A_E$-modules are isomorphic as an easy consequence of our definition.
	
	\begin{definition}[Modules over Picard algebroids]
		Let $\Ee$ be a Picard algebroid on $X$.
		\begin{enumerate}
			\item A left (resp.\ right) $\Ee$-module is a $\OO_X$-module $\Mm$ endowed with an $x^{-1}\OO_S$-linear morphism $\Ee\otimes_{x^{-1}\OO_S} \Mm\to \Mm$ (resp.\ $\Mm\otimes_{x^{-1}\OO_S} \Ee\to \Mm$) whose sections satisfy the respective conditions in Definition \ref{def:moduleoverpicardalgebra}. Let $\Mod(\Ee)$ (resp.~$\Mod_{\mathrm{r}}(\Ee)$) denote the category of left (resp.\ right) $\Ee$-modules.
			\item We say a left or right $\Ee$-module is quasi-coherent if it is so over $\OO_X$. We denote the category of quasi-coherent left (resp.~right) $\Ee$-modules by $\Mod_{\qc}(\Ee)$ (resp.~$\Mod_{{\mathrm{r}},\qc}(\Ee)$).
		\end{enumerate}
	\end{definition}
	
	We again have an isomorphism $\Mod_{\bullet,\ast}(\Ee)\cong\Mod_{\bullet,\ast}(\Aa_{\Ee})$ for $\bullet\in\{\emptyset,\mathrm{r}\}$ and $\ast\in\{\emptyset,\qc\}$.
	
	\begin{example}\label{ex:linebdltwistofmodule}
		Let $\Aa$ be a tdo on $X$, and $\Ll$ be a line bundle on $X$. For a left (resp.\ right) $\Aa$-module $\Mm$, $\Ll\otimes_{\OO_X}\Mm$ (resp.\ $\Mm\otimes_{\OO_X}\Ll^\vee$) is a left (resp.\ right) $\Ll \otimes_{\OO_X} \Aa \otimes_{\OO_X} \Ll^\vee$-module.
	\end{example}
	
	\begin{definition}[Twisted integrable connections]
		Let $\Ee$ (resp.~$\Aa$) be a Picard algebroid (resp.~a tdo) on $X$. Then a locally free $\OO_X$-module of finite rank with the structure of a left $\Ee$-module (resp.~$\Aa$-module) is called an integrable left $\Ee$-connection (resp.~$\Aa$-connection). If $\Ee=\tilde{\Theta}_{X/S}$ (resp.~$\Aa=\D_X$) then we simply say an integrable connection. We define integrable right $\Ee$-connections (resp.~$\Aa$-connections) in a similar way.
	\end{definition}
	
	\begin{example}\label{ex:O_XasDmodule}
		For the untwisted Picard algebroid $\tilde{\Theta}_{X/S}$, $\OO_X$ is a left $\tilde{\Theta}_{X/S}$-module via the canonical actions: $(a,\partial)b=a\partial(b)$ for all $a,b\in\OO_X$, $\partial\in\Theta_{X/S}$.
	\end{example}
	
	\begin{example}
		The canonical bundle $\omega_{X/S}$ is a right $\D_X$-module for the minus of the Lie derivative (see Corollary \ref{cor:omegaisrightDmod}).
	\end{example}

	\begin{example}[Integrability, {\cite[Theorem 1.4.10]{hottaetal2008}}]\label{thm:integrability}
		Put $S=\Spec F$ with $F$ a field of characteristic zero. Then an $\Aa$-module $\mathcal M$ is $\mathcal O_X$-coherent if and only if $\mathcal M$ is a locally free $\mathcal O_X$-module of finite rank. Let us note that for a purely algebraic argument, use Krull's intersection theorem to define $\ord(f_i)$ in \cite[Proof of Theorem 1.4.10]{hottaetal2008}.
	\end{example}
	
	We remark that this integrability property fails for general bases. We have two counterexamples in the ring setting (see Lemma \ref{lem:modulecataffinecase} if necessary).
	
	\begin{example}\label{counterexample1}
		Let $p$ be any prime. Write $D(1)=D_{\FF_p\left[x\right]/\FF_p}$. Consider the canonical left action of $D(1)$ on $\FF_p\left[x\right]$. Then the ideal $(x^p)\subset\FF_p\left[x\right]$ is a $D(1)$-submodule. In particular, we obtain a nonzero finite dimensional $D(1)$-module $\FF_p\left[x\right]/(x^p)$.
	\end{example}
	
	\begin{example}\label{counterexample2}
		Put $k=\CC\left[t\right]\subset\CC\left[t,x\right]=R$ (see the beginning of section \ref{sec:tdo_and_pic_alg}). Then we have $D_{R/k}\cong k\otimes_\CC D(1)$, where $D(1)=D_{\CC\left[x\right]/\CC}$. Define the structure of a $D_{R/k}$-module on $\CC\left[x\right]$ as follows:
		\begin{itemize}
			\item Put the canonical structure of a $D(1)$-module on $\CC\left[x\right]$.
			\item For $f\in \CC\left[x\right]$, $t\cdot f=0$.
		\end{itemize}
		In particular, $\CC\left[x\right]$ is finitely generated as a $D_{R/k}$-module. However, this is not projective as a $\CC\left[t,x\right]$-module.	
	\end{example}
	
	The results below will be used in section \ref{sec:operationswithDmodules}.
	
	\begin{proposition}\label{prop:resolutions}
		Let $\Aa$ be a tdo on $X$. 
		\begin{enumerate}
			\item Every q-injective complex of $\Aa$-modules is q-injective in $K(\OO_X)$.
			\item Every q-flat complex of left $\Aa$-modules is q-flat in $K(\OO_X)$.
		\end{enumerate}
	\end{proposition}
	
	\begin{proof}
		Recall that $\Aa$ is flat as an $\OO_X$-module for both right and left multiplications (apply the dual argument to Proposition \ref{prop:tdofreebasis} for the right-flatness). Notice that $\Aa\otimes_{\OO_X}-$ is left adjoint to the restriction functor $\Mod(\Aa)\to\Mod(\OO_X)$. The dual argument to Proposition \ref{prop:tdofreebasis} implies that $\Aa$ is a flat $\OO_X$-module for the right multiplication. The assertions are now straightforward.
	\end{proof}

	\subsection{Quasi-coherent twisted $\mathcal D$-modules}\label{sec:qcd-mod}
	
	Recall that for a sheaf $\Rr$ of rings on a topological space $Z$, we have a general notion of quasi-coherent $\Rr$-modules:
	
	\begin{definition}[{\cite[(5.1.3)]{ega1}}]\label{defn:qcoh}
		A left or right $\Rr$-module is called quasi-coherent if there is an open covering $Z=\cup_\lambda U_\lambda$ such that for each index $\lambda$, there is a presentation
		\[\Rr|_{U_\lambda}^{I_\lambda}\to
		\Rr|_{U_\lambda}^{J_\lambda}\to \Mm|_{U_\lambda}\to 0,\]
		where $I_\lambda$ and $J_\lambda$ are (small) sets.
	\end{definition}
	
	Let $\Aa$ be a tdo on $X$. In this case, we have no risk of confusion with the notation $\Mod_{\qc}(\Aa)$ by the following fact:
	
	\begin{proposition}\label{prop:defnofqcoh}
		A left $\Aa$-module $\Mm$ is quasi-coherent as an $\Aa$-module if and only if it is so as an $\OO_X$-module.
	\end{proposition}
	
	To prove this, we remark the following observation:
	
	\begin{lemma}\label{lem:modulecataffinecase}
		Suppose that $X$ is affine. Write $A=\Gamma(X,\Aa)$. Then the global section functor $\Gamma(X,-)$ determines a categorical equivalence $\Mod_{\qc}(\Aa)\simeq\Mod(A)$.
	\end{lemma}
	
	\begin{proof}
		Write $R$ for the coordinate ring of $X$. Then we have a canonical isomorphism
		\[\OO_X\otimes_{R_X} A_X\cong \Aa\]
		of $(\OO_X,A_X)$-bimodules since $\Aa$ is a quasi-coherent $\OO_X$-module for the left multiplication (Proposition \ref{prop:tdofreebasis}). Hence the quasi-inverse $\OO_X\otimes_{R_X}(-)_X$ to $\Gamma(X,-):\Mod_{\qc}(\OO_X)\to \Mod(R)$ lifts to that of $\Gamma(X,-):\Mod_{\qc}(\Aa)\to \Mod(A)$. 
	\end{proof}
	
	\begin{proof}[Proof of Proposition \ref{prop:defnofqcoh}]
		To prove the ``only if'' direction, we may assume that $\Mm$ has a presentation
		$\Aa^I\to \Aa^J\to \Mm\to 0$.
		Recall that $\Aa$ is a quasi-coherent $\OO_X$-module for the left multiplication (Lemma \ref{lem:Dquasicoherence}). The assertion now follows since $\Mod_{\qc}(\Aa)$ is closed under formation of colimits in $\Mod(\Aa)$.
		
		Conversely, suppose that $\Mm$ is quasi-coherent as an $\OO_X$-module. We may assume $X$ is affine. Write $A=\Gamma(X,\Aa)$ and $M=\Gamma(X,\Mm)$. Choose a presentation
		$A^I\to A^J\to M\to 0$
		as an $A$-module. It lifts to an exact sequence of $\Aa$-modules
		$\Aa^I\to \Aa^J\to \Mm\to 0$
		as desired by Proposition \ref{lem:modulecataffinecase}.
	\end{proof}
	
	Let us also state a derived analog of Lemma \ref{lem:modulecataffinecase}.
	
	\begin{proposition}\label{prop:derivedglobalsectionforaffineschemes}
		Let $X$ and $A$ be as in Lemma \ref{lem:modulecataffinecase}. Recall that we have defined the functors $R\Gamma(X,-):D(\Aa)\to D(A)$ and $H^i(X,-):D(\Aa)\to\Mod(A)$ ($i\in\ZZ$) in section \ref{sec:homologicalalgebra}. For any complex $\Mm^\bullet\in D_{\qc}(\Aa)$ and any integer $i$, we have a canonical isomorphism
		\[H^i(X,\Mm^\bullet)\cong \Gamma(X,H^i(\Mm^\bullet)).\]
		In particular, $R\Gamma(X,-):D(\Aa)\to D(A)$ restricts to
		\[\begin{array}{ccc}
			D^{\geq 0}_{\qc}(\Aa)\to D^{\geq 0}(A),
			&D^{\leq 0}_{\qc}(\Aa)\to D^{\leq 0}(A),
			&\Gamma(X,-):~\Mod_{\qc}(\Aa)\to \Mod(A).
		\end{array}\]
	\end{proposition}
	
	\begin{proof}
		Without loss of generality by passing to shifts, we may assume $i=0$. For an integer $p$, let $\tau^{\geq p}:D(\Aa)\to D^{\geq p}(\Aa)$ and $\tau^{\leq p}:D(\Aa)\to D^{\leq p}(\Aa)$ be the truncation functors (\cite[section 1.11]{lipman2009}). Then a similar argument to \cite[Proposition (3.9.2)]{lipman2009} implies $H^0(X,\Mm^\bullet)\cong H^0(X,\tau^{\geq 0}\Mm^\bullet)$. In fact, $d$ in the proof of \cite[Proposition (3.9.2)]{lipman2009} is 0 in our case by \cite[Th\'eor\`eme (1.3.1)]{ega31} and the remark below Proposition \ref{prop:independent}. We may therefore assume $\Mm^\bullet\in D^{\geq 0}_{\qc}(\Aa)$. 
		
		The assertion now follows from the distinguished triangle
		$H^0(\Mm_\bullet)\to \Mm^\bullet\to \tau^{\geq 1}\Mm^\bullet\overset{+1}{\to}$
		and the evident equality $H^0(X,\tau^{\geq 1}\Mm^\bullet)=0$.
	\end{proof}
	
	\begin{proposition}\label{prop:(co)limofmodqc}
		The category $\Mod_{\rm qc}(\Aa)$ is Grothendieck abelian, whose small colimits and finite limits are computed in $\Mod(\ZZ_X)$.
	\end{proposition}
	
	\begin{proof}
		Since small colimits and finite limits of $\Mod(\Aa)$ and $\Mod_{\qc}(\OO_X)$ are computed in the category $\Mod(\ZZ_X)$,
		so are those of $\Mod_{\qc}(\Aa)$. It will therefore suffice to show that $\Mod_{\qc}(\Aa)$ admits a generator.
		
		It is clear by Proposition \ref{prop:defnofqcoh} that the functor $\Aa\otimes_{\OO_X}-:\Mod(\OO_X)\to\Mod(\Aa)$ respects quasi-coherent sheaves. The claim now follows from \cite[Lemma 2.1.8]{hayashi2018}. In fact, note that $\Mod_{\qc}(\OO_X)$ admits a generator by a result of Gabber.
	\end{proof}
	
	The adjoint functor theorem now implies
	
	\begin{corollary}[Coherentor]
		The canonical embedding $\Mod_{\qc}(\Aa)\hookrightarrow\Mod(\Aa)$ admits a right adjoint functor.
	\end{corollary}

	\subsection{Coherent twisted $\mathcal D$-modules}\label{sec:coh}
	
	Let $\Aa$ be a tdo on $X$. In this section, assume that $X$ is locally Noetherian\footnote{This is satisfied, for example, if $S$ is locally Noetherian.}.
	
	\begin{definition}[{\cite[Definition 1.4.8]{hottaetal2008}}]
		A left or right $\mathcal A$-module $\mathcal M$ is called $\mathcal A$-{\em coherent} if $\mathcal M$ is a locally finitely generated $\mathcal A$-module and if for {\em every} open $U\subseteq X$, every locally finitely generated submodule of $\mathcal M|_U$ is a locally finitely presented $\mathcal A_U$-module. We denote the full subcategory of $\Mod(\Aa)$ (resp.\ $\Mod_{\mathrm{r}}(\Aa)$) consisting of left (resp.\ right) $\Aa$-coherent modules by $\Mod_c(\Aa)$ (resp.\ $\Mod_{\rm{r,c}}(\Aa)$). 
	\end{definition}

	We discuss the coherence property via filtrations.
	
	\begin{proposition}\label{prop:Dnoetherian}
		Suppose that $X$ is Noetherian\footnote{This is satisfied, for example, if $x$ and $S$ are quasi-compact.}.
		\begin{enumerate}
			\item Every tdo $\Aa$ on $X$ is left (resp.\ right) Noetherian, i.e., $\Aa$ is a Noetherian object in $\Mod_{\qc}(\Aa)$ (resp.\ $\Mod_{{\mathrm{r}},\qc}(\Aa)$).
			\item For any $p\in X$, the stalk $\Aa_{p}$ is left and right Noetherian.
		\end{enumerate}
	\end{proposition}
	
	\begin{proof}
		Part (2) is a routine consequence of (1). For (1), we may assume that $X$ is affine. Moreover, assume that $\Theta_{X/S}$ admits a free basis $\{\partial_1,\ldots,\partial_n\}$. Then $\Gr\Aa$ is isomorphic to the polynomial algebra sheaf $\mathcal O_X[\partial_1,\dots,\partial_n]$ by definition of tdos. In particular, $\Aa$ is left and right Noetherian by Lemma \ref{lem:coherentDgenerators} and \cite[Proposition D.1.4]{hottaetal2008}.
	\end{proof}
	
	For $\Aa$-modules, we begin with an easy observation in
	
	\begin{lemma}\label{lem:coherentDgenerators}
		Assume $X$ is Noetherian. Then a locally finitely generated $\Aa$-module $\mathcal M$ is generated by a coherent $\mathcal O_X$-submodule.
	\end{lemma}
	
	This statement means that there is a coherent $\mathcal O_X$-submodule $\mathcal F$ of $\mathcal M$ with the property that the canonical map
	$
	\Aa\otimes_{\mathcal O_X}\mathcal F\to\mathcal M
	$
	is an epimorphism of sheaves.
	
	\begin{proof}
		Since $\mathcal M$ is locally finitely generated, we find a finite affine covering $U_i$, $i\in I$, of $X$ with the property that on each $U_i$, $\mathcal M|_{U_i}$ is a finitely generated $\Aa|_{U_i}$-module. This is the same to say that there is a coherent $\mathcal O_{U_i}$-submodule $\mathcal F_i$ in $\mathcal M|_{U_i}$ which generates $\mathcal M|_{U_i}$ as $\Aa_{U_i}$-module. Now each $\mathcal F_i$ extends to a coherent $\mathcal O_X$-submodule $\mathcal G_i$ of $\mathcal M$ (cf.~\cite[Th\'eor\`eme (9.4.7)]{ega1}), and the finite sum $\mathcal G$ of the modules $\mathcal G_i$, $i\in I$, is a coherent $\mathcal O_X$-submodule of $\mathcal M$ which generates $\mathcal M$.
	\end{proof}
	
	\begin{proposition}\label{prop:schemegoodfiltration}
		Suppose that $X$ is Noetherian. For a quasi-coherent left $\mathcal A$-module $\Mm$, the following conditions are equivalent:
		\begin{enumerate}
			\renewcommand{\labelenumi}{(\alph{enumi})}
			\item $\mathcal M$ is a locally finitely generated $\mathcal A$-module.
			\item $\mathcal M$ admits the structure of an exhaustive filtered $\Aa$-module such that $\Gr\mathcal M$ is quasi-coherent, and that $\Gr\mathcal M$ is locally finitely generated as a $\Gr\mathcal A$-module.
		\end{enumerate} 
	\end{proposition}
	
	We call a filtration of $\mathcal M$ satisfying (b) {\em good}. To construct a filtration consisting of quasi-coherent subsheaves, we need the following simple observation:
	
	\begin{lemma}\label{lem:bimodtensoqucoherence}
		Let $\mathcal F$ be an $\OO_X$-bimodule such that $\mathcal F$ is quasi-coherent as a left $\OO_X$-module. Let $\mathcal G$ be a quasi-coherent $\OO_X$-module. Then $\Ff \otimes_{\OO_X}\Gg$ is a quasi-coherent $\OO_X$-module\footnote{Here $X$ needs not be locally Noetherian.}.
	\end{lemma}
	\begin{proof}
		Since the assertion is local on $X$, we may assume that $\Gg$ admits a presentation
		\[\OO_X^I\to\OO_X^J\to \Gg\to 0.\]
		We thus get an exact sequence of $\OO_X$-modules
		$\Ff^I\to \Ff^J
		\to \Ff\otimes_{\OO_X}\Gg\to 0$.
		Since $\Ff$ is quasi-coherent for the left action, so are $\Ff^I$ and $\Ff^J$. Therefore $\Ff\otimes_{\OO_X}\Gg$ is quasi-coherent as well.
	\end{proof}
	
	\begin{proof}[Proof of Proposition \ref{prop:schemegoodfiltration}]
		Suppose that (a) is satisfied. Then there is a coherent $\OO_X$-submodule $\mathcal M_0\subseteq\mathcal M$ such that the action map $s:\mathcal A\otimes_{\mathcal O_X}\mathcal M_0\to\mathcal M$ is an epimorphism (Lemma \ref{lem:coherentDgenerators}). Lemma \ref{lem:bimodtensoqucoherence} implies that the canonical images of $F_i\mathcal A\otimes_{\OO_X}\mathcal M_0$ under $s$ define a filtration $F_i\mathcal M$ on $\mathcal M$ satisfying the conditions in (b).
		
		Conversely, suppose that we are given a filtration $F_i\Mm$ in (b). We wish to prove (a). Since the assertion is local in $X$, we may assume that $\Gr\Mm$ is a finitely generated $\Gr\Aa$-module, and that $X$ is affine. Write $A=\Gamma(X,\Aa)$, $F_iA=\Gamma(X,F_i\Aa)$, $M=\Gamma(X,\Mm)$, and $F_iM=\Gamma(X,F_i\Mm)$ for $i\in\NN$. Then $M$ is a filtered $A$-module for $F_\bullet M$ by Lemma \ref{lem:bimodtensoqucoherence}. Moreover, we can identify the induced action of $\Gr A$ on $\Gr M$ with $\Gamma(X,a)$, where $a$ is the action map of $\Gr\Aa$ on $\Gr\Mm$. The assertion now follows from \cite[Corollary D.1.2]{hottaetal2008}. This completes the proof.
	\end{proof}
	
	\begin{proposition}\label{prop:coherentDmodules}
		Let $\mathcal M$ be a quasi-coherent left $\mathcal A$-module.
		Then the following are equivalent:
		\begin{enumerate}
			\renewcommand{\labelenumi}{(\alph{enumi})}
			\item $\Mm$ is $\Aa$-coherent;
			\item $\mathcal M$ is a locally finitely presented $\mathcal A$-module;
			\item $\mathcal M$ is a locally finitely generated $\mathcal A$-module.
		\end{enumerate}
	\end{proposition}
	
	\begin{proof}
		Firstly, we show the equivalence of (b) and (c). The implication (b) $\Rightarrow$ (c) is clear. Suppose that (c) is satisfied. We wish to prove (b). Since the assertion is local in $X$, we may assume that $X$ is affine, and that $\Mm$ is finitely generated as an $\Aa$-module. Then $\Mm$ is finitely presented by Proposition \ref{prop:Dnoetherian} and Lemma \ref{lem:modulecataffinecase}.
		
		The implication (a) $\Rightarrow$ (c) is clear. The converse is deduced by applying (c) $\Rightarrow$ (b) for open subschemes $U\subset X$ and locally finitely generated $\Aa|_U$-submodules of $\Mm|_U$. This is possible since open subschemes of Noetherian schemes are again Noetherian.
	\end{proof}
	
	\begin{corollary}
		The category $\Mod_{\rm c}(\Aa)$ of $\Aa$-coherent $\Aa$-modules is a thick abelian subcategory of $\Mod_{\rm qc}(\Aa)$.
	\end{corollary}
	
	\begin{corollary}
		Any tdo $\Aa$ is an $\Aa$-coherent $\Aa$-module (from the left and from the right actions).
	\end{corollary}

	\section{Operations with twisted $\mathcal D$-modules}\label{sec:operationswithDmodules}
	
	A remarkable aspect of twisted $\D$-modules over fields of characteristic zero is their derived functoriality. This is important in representation theory as a tool of geometric constructions of representations. We aim to establish derived functoriality of twisted $\D$-modules over schemes in this section towards construction of integral models of cohomologically induced modules. To verify that our geometric construction in section \ref{sec:arithmeticmodelsofaq(lambda)} provides us integral forms of cohomologically induced modules, we introduce the base change functor and study the base change property of other functors. In this occasion, we would like to discuss them at the largest reasonable generality.
	
	We adopt the notations of section \ref{sec:tdo_and_pic_alg} for $X,Y,S,S',X',Y',x,x',s,y,y',f,f',s_X,s_Y$ unless specified otherwise.
	
	\subsection{Localization and globalization}\label{sec:localization}
	
	In this section, we introduce the localization and globalization functors for twisted $\D$-modules over general bases. Let $\Aa$ be a tdo on $X$. Then we have a pair of adjoint functors
	
	\begin{equation}
		\Aa\otimes_{x^{-1}x_\ast\Aa} x^{-1}(-):
		\Mod_{\qc}(x_\ast\Aa)\rightleftarrows \Mod_{\qc}(\Aa):x_\ast.
		\label{eq:adjglobloc;abel}
	\end{equation}
	In fact, for the right adjoint functor $x_\ast$, use the counit $x^{-1}x_\ast\Aa\to \Aa$. We also have a pair of adjoint functors
	\begin{equation}
		\Aa\otimes^{L}_{x^{-1}x_\ast\Aa}x^{-1}(-)
		:D(x_\ast\Aa) \rightleftarrows D(\Aa):Rx_\ast
		\label{eq:adjglobloc;derived}
	\end{equation}
	of derived categories by Proposition \ref{prop:rightderivedfunctor} (2). We call $\Aa\otimes^L_{x^{-1}x_\ast\Aa} x^{-1}(-)$ and $R x_\ast$ the (derived) localization and globalization functors respectively. We remark that the derived localization functor $\Aa\otimes^L_{x^{-1}x_\ast\Aa} -: D(x_\ast\Aa)\to D(\Aa)$ can be computed as the left derived functor of the corresponding left adjoint functor $\Aa\otimes_{x^{-1}x_\ast\Aa} x^{-1}(-):\Mod(x_\ast\Aa)\to\Mod(\Aa)$ (see below Proposition \ref{prop:leftderivedfunctor}). By Propositions \ref{prop:resolutions} or \ref{prop:independent}, $Rx_\ast:D(\Aa)\to D(x_\ast\Aa)$ enjoys the 2-commutative diagram
	\begin{equation}
		\begin{tikzcd}
			D(\Aa)\ar[r, "Rx_\ast"]\ar[d]&D(x_\ast\Aa)\ar[d]\\
			D(\OO_X)\ar[r, "Rx_\ast"]&D(\OO_S),
		\end{tikzcd}
		\label{diag:compatibilityofglobalization}
	\end{equation}
	where the vertical arrows are defined by restriction.
	
	To work within cohomologically quasi-coherent complexes, consider
	\begin{condition}\label{cond:preservationofqcoh}
		The functor $x_\ast|_{\Mod_{\qc}(\OO_X)}$ is valued in $\Mod_{\qc}(\OO_S)$.
	\end{condition}

	\begin{remark}[{\cite[Corollaire 9.2.2]{ega1}, \cite[(5.6) Theorem]{altmaneltal}}]\label{rmk:onDaffinei}
		Sufficient conditions on a morphism $x:X\to S$ of schemes in general to preserve quasi-coherent modules are:
		\begin{enumerate}
			\renewcommand{\labelenumi}{(\roman{enumi})}
			\item $x$ is concentrated,
			\item $x$ is quasi-compact and $X$ is locally Noetherian,
			\item $S$ is the spectrum of a field.
		\end{enumerate}
	\end{remark}
	
	In the rest of this small section, we assume Condition \ref{cond:preservationofqcoh}. Then $x_\ast\Aa$ is a quasi-coherent $\OO_S$-algebra, i.e., a monoid object of the symmetric monoidal category $\Mod_{\qc}(\OO_S)$ since $x_\ast\Aa$ is an $\OO_S$-algebra, and is a quasi-coherent $\OO_S$-module for the (left) multiplication. Under the current hypothesis, a similar result to Proposition \ref{prop:defnofqcoh} holds for $x_\ast\Aa$:
	
	\begin{lemma}\label{lem:qcohmodoverqcohalg}
		Let $S$ be an arbitrary scheme, and $\Bb$ be a (possibly noncommutative) quasi-coherent $\OO_S$-algebra. Then a left $\Bb$-module $\Mm$ is quasi-coherent as an $\Bb$-module if and only if it is so as an $\OO_S$-module. In particular, $\Mod_{\qc}(\Bb)$ consists precisely of the quasi-coherent left $\Bb$-modules.
	\end{lemma}
	
	\begin{proof}
		The ``only if'' direction is proved in a similar way to Proposition \ref{prop:defnofqcoh}. For the converse, we may assume $S$ is affine. Then one can find an epimorphism $\OO_S^I\to \Mm$. Take the scalar extension to get an epimorphism $p:\Bb^I\to\Mm$ of left $\Bb$-modules. Apply the same argument to $\Ker p$ to get a presentation $\Bb^J\to \Bb^I\to\Mm\to 0$. This completes the proof.
	\end{proof}
	
	Hence we can use the notations $\Mod_{\qc}(x_\ast\Aa)$ and $D_{\qc}(x_\ast\Aa)$ without any confusion. Indeed, we will never use these notations unless Condition \ref{cond:preservationofqcoh} holds.
	
	The next assertion is now evident by definitions:
	
	\begin{proposition}
		Assume Condition \ref{cond:preservationofqcoh} holds. Then the adjunction \eqref{eq:adjglobloc;abel} restricts to 
		\[\Aa\otimes_{x^{-1}x_\ast\Aa} x^{-1}(-):
		\Mod_{\qc}(x_\ast\Aa)\rightleftarrows \Mod_{\qc}(\Aa):x_\ast.\]
	\end{proposition}
	
	One can prove its derived analog:
	\begin{proposition}\label{prop:x_ast;qc_pres}
		\begin{enumerate}
			\item Assume Condition \ref{cond:preservationofqcoh} holds. Then the derived localization functor $\Aa\otimes^L_{x^{-1}x_\ast\Aa} x^{-1}(-): D(x_\ast\Aa)\to D(\Aa)$ respects cohomologically quasi-coherent complexes.
			\item Suppose that $x$ is concentrated. Then $Rx_\ast$ restricts to $D_{\qc}(\Aa)\to D_{\qc}(x_\ast\Aa)$.
		\end{enumerate}
	\end{proposition}
	
	\begin{proof}
		Part (2) follows from the 2-commutative diagram \eqref{diag:compatibilityofglobalization} and \cite[Proposition (3.9.2)]{lipman2009}.
		
		We prove (1). Since the problem is local in $S$, we may assume $S$ affine. By a similar argument to \cite[Proposition (3.9.1)]{lipman2009}, it will suffice to show $\Aa\otimes^L_{x^{-1}x_\ast\Aa} x^{-1}\Mm\in D_{\qc}(\Aa)$ for $\Mm\in\Mod_{\qc}(x_\ast\Aa)$. Since $S$ is affine, one can take a resolution $\Ff^\bullet\to \Mm\to 0$ by free $x_\ast\Aa$-modules. Then $\Aa\otimes^L_{x^{-1}x_\ast\Aa} x^{-1}\Mm$ is computed by $\Aa\otimes_{x^{-1}x_\ast\Aa} x^{-1}\Ff^\bullet$. Notice that this complex consists of free left $\Aa$-modules. In particular, $\Aa\otimes^L_{x^{-1}x_\ast\Aa} x^{-1}\Ff^\bullet$ is a complex of quasi-coherent $\OO_X$-modules (see Proposition \ref{prop:defnofqcoh}). This shows $\Aa\otimes^L_{x^{-1}x_\ast\Aa} x^{-1}\Mm$ lies in $D_{\qc}(x_\ast\Aa)$ as desired.
	\end{proof}
	
	Finally, let us give a quick review on the twisted $\D$-affinity. Henceforth we put $S=\Spec F$, where $F$ is a field. We identify
	$x_\ast: \Mod_{\qc}(\Aa)\to\Mod(x_\ast\Aa)$
	with the global section functor
	\[\Gamma(X,-): \Mod_{\qc}(\Aa)\to\Mod(\Gamma(X,\Aa))\]
	through the isomorphism $\Mod(x_\ast\Aa)\cong\Mod(\Gamma(X,\Aa))$. 
	
	Recall that $X$ is quasi-separated over $S$ since $X$ is locally Noetherian. In the rest of this small section, we assume that $X$ is quasi-compact. 
	
	\begin{definition}[$\mathcal A$-affine varieties, {\cite[Definition 1.4.2]{hottaetal2008}}]\label{def:Daffine}
		We say that $X$ is {\em $\Aa$-affine} if the following two conditions are satisfied:
		\begin{enumerate}
			\renewcommand{\labelenumi}{(\roman{enumi})}
			\item $H^q(X,\mathcal M)=0$ for every $\mathcal M\in\Mod_{\rm qc}(\Aa)$ and $q>0$.
			\item For every $\mathcal M\in\Mod_{\rm qc}(\Aa)$ with $\Gamma(X,\Mm)=0$, we have $\Mm=0$.
		\end{enumerate}
	\end{definition}
	
	\begin{example}
		If $X$ is affine, $X$ is $\Aa$-affine for any tdo $\Aa$ by \cite[Corollaire (1.3.2)]{ega31}.
	\end{example}
	
	We will see another fundamental example at the end of section \ref{sec:basechangeglobalsection}. The following fact is standard:
	
	\begin{proposition}[{\cite[Propositions 1.4.4 and 1.4.13]{hottaetal2008}}]\label{cor:firstequivalenceoverfields}
		Assume that $X$ is $\Aa$-affine. Then
		\begin{enumerate}
			\item Every $\mathcal M\in\Mod_{\rm qc}(\Aa)$ is generated by its global sections as $\Aa$-module.
			\item The functor
			$
			\Gamma(X,-):
			\Mod_{\rm qc}(\Aa)\to
			\Mod(\Gamma(X,\Aa))
			$
			is an equivalence of categories.
			\item The functor $\Gamma(X,-)$ induces an equivalence between $\Aa$-coherent modules and finitely generated $\Gamma(X,\Aa)$-modules.
		\end{enumerate}
	\end{proposition}
	
	\subsection{Base change}
	
	Let $A$ be a tdo for $R/k$. Then the adjunction
	$k'\otimes_k-: \Mod(k)\rightleftarrows \Mod(k')$
	clearly extends to an adjunction of (left or right) $A$-modules and $k'\otimes_k A$-modules.
	
	We next consider its sheafified analog. Since the functors for left and right modules have similar behavior, we only deal with the left ones for a while. Let $\Aa$ be a tdo on $X$. 
	
	\begin{proposition}\label{prop:bcadj}
		The adjunction 
		$s_X^\ast:\Mod(\OO_X)\rightleftarrows\Mod(\OO_{X'})
		:(s_X)_\ast$
		lifts to
		\[s_X^\ast:\Mod(\Aa)\rightleftarrows
		\Mod(s_X^\ast\Aa):(s_X)_\ast.\]
	\end{proposition}
	
	\begin{proof}
		We have an adjunction
		\[s^\ast_X\Aa\otimes_{s^{-1}_X\Aa} s_X^{-1}(-):
		\Mod(\Aa)\rightleftarrows
		\Mod(s_X^\ast\Aa):(s_X)_\ast\]
		by general theory of sheaves. In fact, use the canonical homomorphisms $s^{-1}_X\Aa\to s^\ast_X\Aa$ and $\Aa\to (s_X)_\ast s_X^\ast\Aa$ of sheaves of rings to define $s^\ast_X\Aa\otimes_{s^{-1}_X\Aa} s_X^{-1}(-)$ and $(s_X)_\ast$ respectively. We can identify $s^\ast\Aa_{s^{-1}_X\Aa}s_X^{-1}(-)$ with $s^\ast_X$. In fact, for a left $\Aa$-module $\Mm$, we have a canonical $\OO_{X'}$-module homomorphism $s_X^\ast \Mm\to s^\ast_X\Aa\otimes_{s^{-1}_X\Aa} s_X^{-1}\Mm$. One can prove that this is an isomorphism. For this, think of $s^{-1}\Aa$ as an $(s^{-1}_X\OO_X,s^{-1}_X\Aa)$-bimodule to regard $s^\ast_X\Aa\otimes_{s^{-1}_X\Aa} s_X^{-1}(-)$ as the triple tensor product $(\OO_{X'}\otimes_{s^{-1}_X\OO_X} s^{-1}_X\Aa)\otimes_{s^{-1}_X\Aa} s^{-1}_X(-)$ which canonically simplifies to $\OO_{X'}\otimes_{s^{-1}_X\OO_X} s^{-1}_X(-)=s^\ast_X$ by the above map. Finally, we remark that the two adjunctions are compatible with this identification. This completes the proof.
	\end{proof}
	
	It follows by our definition (or Proposition \ref{prop:defnofqcoh}) that $s_X^\ast$ respect quasi-coherent sheaves. If $(s_X)_\ast:\Mod(\OO_{X'})\to\Mod(\OO_X)$ restricts to $\Mod_{\qc}(\OO_{X'})\to\Mod_{\qc}(\OO_X)$, so does
	\[(s_X)_\ast:\Mod(s^\ast_X \Aa)\to\Mod(\Aa)\]
	also to $\Mod_{\qc}(s^\ast_X \Aa)\to\Mod_{\qc}(\Aa)$. 
	
	Let us collect basic properties of $s^\ast_X$:
	\begin{proposition}
		The functor $s_X^*$ preserves the following properties of $\Aa$-modules:
		\begin{enumerate}
			\renewcommand{\labelenumi}{(\roman{enumi})}
			\item (local) finite generation,
			\item (local) finite presentation,
			\item $\mathcal O_X$-local freeness.
		\end{enumerate}
	\end{proposition}
	
	\begin{proof}
		The claim (iii) follows by definitions. The claims (i) and (ii) are immediate from the identification $s^\ast_X=s^\ast\Aa\otimes_{s^{-1}_X\Aa}s^{-1}_X(-)$.
	\end{proof}
	
	The next observation will be used in later sections for the study of derived functors:
	
	\begin{lemma}\label{lem:preservationflat}
		Let $\Ff^\bullet\in K(\Aa)$ be a q-flat complex. Then $s_X^\ast\Ff^\bullet$ is q-flat in $K(s_X^\ast\Aa)$.
	\end{lemma}
	
	\begin{proof}
		Think of $s_X^\ast\Ff$ as $s_X^\ast \Aa \otimes_{s_X^{-1}\Aa} s_X^{-1}\Ff$ (see the proof of Proposition \ref{prop:bcadj}). Then the assertion follows from general theory of sheaves. 
	\end{proof}
	
	Finally, we discuss the derived setting. 
	\begin{proposition}\label{prop:derivedadjunction_basechange}
		The adjunction
		$Ls^\ast_X:D(\OO_X)\rightleftarrows D(\OO_{X'}):R(s_X)_\ast$
		lifts to
		\[Ls^\ast_X:D(\Aa)\rightleftarrows D(s^\ast_X\Aa):R(s_X)_\ast.\]
		In other words,
		\begin{enumerate}
			\item There is a left derived functor $Ls^\ast_X:D(\Aa)\to D(s^\ast_X\Aa)$ which enjoys the 2-commutative diagram
			\begin{equation}
				\begin{tikzcd}
					D(\Aa)\ar[r, "Ls_X^\ast"]\ar[d]
					&D(s^\ast_X\Aa)\ar[d]\\
					D(\OO_X)\ar[r, "Ls_X^\ast"]
					&D(\OO_{X'}),
				\end{tikzcd}
				\label{diag:compatibilityofbc}
			\end{equation}
			where the vertical arrows are defined by restriction.
			\item There is a right derived functor $R(s_X)_\ast:D(s^\ast_X\Aa)\to D(\Aa)$ which enjoys the 2-commutative diagram
			\begin{equation}\begin{tikzcd}
					D(s^\ast_X\Aa)\ar[r, "R(s_X)_\ast"]\ar[d]
					&D(\Aa)\ar[d]\\
					D(\OO_{X'})\ar[r, "R(s_X)_\ast"]
					&D(\OO_{X}),
				\end{tikzcd}
				\label{diag:compatibilityofrestriction}
			\end{equation}
			where the vertical arrows are defined by restriction.
			\item The derived functors of (1) and (2) form an adjunction.
			\item The natural transformations of units and counits defined by the adjunctions
			\[\begin{array}{cc}
				Ls^\ast_X:D(\OO_X)\rightleftarrows D(\OO_{X'}):R(s_X)_\ast,
				&Ls^\ast_X:D(\Aa)\rightleftarrows D(s^\ast_X\Aa):R(s_X)_\ast
			\end{array}\]
			commute with the restriction.
		\end{enumerate}
	\end{proposition}
	
	\begin{proof}
		Regard $s^\ast_X=(s^\ast_X\Aa)\otimes_{s^{-1}_X\Aa} s^{-1}_X(-)$. Then the derived adjunction is obtained by Lemma \ref{lem:q-inj}. Recall that the horizontal derived functors in \eqref{diag:compatibilityofbc} (resp.~\eqref{diag:compatibilityofrestriction}) are computed by q-flat (resp.~q-injective) resolutions. Hence the 2-commutativity of the diagrams \eqref{diag:compatibilityofbc} and \eqref{diag:compatibilityofrestriction} follows from Proposition \ref{prop:resolutions}.
		
		Finally, we prove (4). For this, it will suffice to show that the diagram
		\[\begin{tikzcd}
			\Hom_{D(s^\ast_X\Aa)}(Ls^\ast_X\Mm^\bullet,\Nn^\bullet)\ar[r, "\sim"]\ar[d]
			&\Hom_{D(\Aa)}(\Mm^\bullet,R(s_X)_\ast\Nn^\bullet)\ar[d]\\
			\Hom_{D(\OO_{X'})}(Ls^\ast_X\Mm^\bullet,\Nn^\bullet)\ar[r, "\sim"]
			&\Hom_{D(\OO_X)}(\Mm^\bullet,R(s_X)_\ast\Nn^\bullet)
		\end{tikzcd}\]
		commutes for $\Mm^\bullet\in D(\Aa)$ and $\Nn^\bullet \in D(s^\ast_X\Aa)$. The vertical arrows are obtained by the restriction. The horizontal arrows are given by the adjunctions. We may compute the Hom sets and the bijections of the adjunctions by q-flat and q-injective resolutions (see Proposition \ref{prop:rightderivedfunctor} (3) and \cite[Proof of Theorem 14.4.5]{kashiwaraschapira2}). The assertion now follows from Proposition \ref{prop:resolutions}.
	\end{proof}
	
	\begin{corollary}
		\begin{enumerate}
			\item The functor $Ls^\ast_X$ restricts to $D_{\qc}(\Aa)\to D_{\qc}(s^\ast_X\Aa)$.
			\item If $s_X$ is concentrated then $R(s_X)_\ast$ restricts to $D_{\qc}(s^\ast_X\Aa)\to D_{\qc}(\Aa)$
		\end{enumerate}
	\end{corollary}
	\begin{proof}
		Part (1) (resp.~(2)) is a consequence of \cite[Proposition (3.9.1) (resp.~(3.9.2))]{lipman2009}.
	\end{proof}

	\subsection{Flat base change for globalization}\label{sec:basechangeglobalsection}
	
	Let $\Aa$ denote a tdo on $X$. In this section, we study the flat base change formula of the derived globalization functor. For this, notice that the diagram
	\begin{equation}
		\begin{tikzcd}
			&\Mod(\OO_{X'})\ar[rr, "(s_X)_\ast"]&&\Mod(\OO_X)\ar[dd, "x_\ast"]\\
			\Mod(s^\ast_X\Aa)\ar[ru]\ar[dd, "x'_\ast"']\ar[rr, "\ \ \ \ \ \ \ \ \ \ \ (s_X)_\ast"']
			&&\Mod(\Aa)\ar[ru]&\\
			&\Mod(\OO_{S'})\ar[rr, "s_\ast\ \ \ "]\ar[from=uu, crossing over, "\overset{x'_\ast}{\overset{\overset{ \overset{}{}\overset{\overset{}{}\overset{}{}\overset{}{}\overset{}{}}{}}{}}{}}"']&&\Mod(\OO_S)\\
			\Mod(x'_\ast s^\ast_X\Aa)\ar[r]\ar[ru]&\Mod(s^\ast x_\ast \Aa)\ar[r, "s_\ast"]
			&\Mod(x_\ast \Aa)\ar[ru]\ar[from=uu, crossing over, "\underset{x'_\ast}{\underset{\underset{ \underset{}{}\underset{}{}}{}}{}}"]
		\end{tikzcd}\label{diag:globalization_abelian_compatible}
	\end{equation}
	strictly commutes. Here the bottom left horizontal arrow $\Mod(x'_\ast s^\ast_X\Aa)\to\Mod(s^\ast x_\ast \Aa)$ is defined by \eqref{eq:bcmapfortdo}. The remaining unlabeled four arrows from front to back are defined by restriction. In fact, for each $\Mm\in \Mod(s^\ast_X\Aa)$, the two actions of $x_\ast \Aa$ on $x'_\ast(s_X)_\ast \Mm=s_\ast x'_\ast\Mm$ induced from the front rectangle coincide by definition of \eqref{eq:bcmapfortdo}. This implies that the front rectangle commutes. The commutativity of the other rectangles are evident.
	
	Since the left adjoint functor of each functor in \eqref{diag:globalization_abelian_compatible} respects q-flat complexes, \cite[Proposition 14.4.7]{kashiwaraschapira2} implies a 2-commutative diagram
	\begin{equation}
		\begin{tikzcd}
			&D(\OO_{X'})\ar[rr, "R(s_X)_\ast"]&&D(\OO_X)\ar[dd, "Rx_\ast"]\\
			D(s^\ast_X\Aa)\ar[ru]\ar[dd, "Rx'_\ast"']\ar[rr, "\ \ \ \ \ \ \ \ \ \ \ R(s_X)_\ast"']
			&&D(\Aa)\ar[ru]&\\
			&D(\OO_{S'})\ar[rr, "Rs_\ast\ \ \ "]\ar[from=uu, crossing over, "\overset{Rx'_\ast}{\overset{\overset{\overset{\overset{}{}\overset{}{}}{}\overset{}{}}{}}{}}"']&&D(\OO_S)\\
			D(x'_\ast s^\ast_X\Aa)\ar[r]\ar[ru]&D(s^\ast x_\ast \Aa)\ar[r, "Rs_\ast"]
			&D(x_\ast \Aa).\ar[ru]\ar[from=uu, crossing over, "\underset{Rx'_\ast}{\underset{\underset{ \underset{}{}\underset{}{}}{}}{}}"]
		\end{tikzcd}\label{diag:globalization_derived_compatible}
	\end{equation}
	In particular, the front 2-commutative rectangle gives rise to the base change map
	\begin{equation}
		Ls^\ast Rx_\ast(-)\to Rx'_\ast Ls_X^\ast(-).
		\label{eq:bcmapmodule}
	\end{equation}
	
	\begin{proposition}\label{prop:coincidencebcmapofglobalsection}
		The map \eqref{eq:bcmapmodule} is a lift of \eqref{eq:bcmapofderivedfunctor}.
	\end{proposition}
	
	\begin{proof}
		This is straightforward from Proposition \ref{prop:derivedadjunction_basechange}, \eqref{diag:compatibilityofrestriction}, and \eqref{diag:globalization_derived_compatible}.
	\end{proof}

	Now Example \ref{ex:flatbasechange} implies
	
	\begin{theorem}\label{thm:flatbasechangeofglobalsectionI}
		Assume $x$ concentrated and $s$ flat. Then the base change map \eqref{eq:bcmapmodule} is an equivalence on $D_{\qc}(\Aa)$.
	\end{theorem}
	
	\begin{corollary}
		Suppose $x$ concentrated and $s:S'\to S$ faithfully flat. Then for every quasi-coherent left $\Aa$-module $\mathcal M$ and every nonnegative integer $q$,
		$$
		R^qx_*\mathcal M=0\quad\iff\quad(R^qx'_*)s_X^*\mathcal M=0.
		$$
	\end{corollary}
	
	\begin{corollary}\label{cor:A-affine_is_geometric}
		Suppose that we are given an extension $F\subset F'$ of fields. Put $S=\Spec F$ and $S=\Spec F'$. Assume that $X$ is quasi-compact. Then $X$ is $\Aa$-affine if and only if $X\times_S S'$ is $s^\ast_X\Aa$-affine.
	\end{corollary}
	
	\begin{proof}
		The ``if'' direction is immediate from the faithfully flat descent. Assume that $X$ is $\Aa$-affine. The vanishing of higher cohomology groups is verified by
		\[H^q(X\times_S S',\Mm')\cong H^q(X,(s_X)_\ast \Mm')=0\]
		for $\Mm'\in\Mod_{\qc}(s^\ast_X\Aa)$ and $q>0$, where the first isomorphism follows since $s_X$ is affine.
		
		Let $\Mm'$ be a quasi-coherent left $s^\ast_X\Aa$-module with $\Gamma(X\times_S S,\Mm')=0$. The proof will be completed by showing $\Mm'=0$. With $\Gamma(X\times_S S,\Mm')=\Gamma(X,(s_X)_\ast\Mm')$, the $\Aa$-affinity implies $(s_X)_\ast \Mm'=0$. To prove $\Mm'=0$, we may work locally in $X$ to assume that $X$ is affine. The assertion is clear in this case.
	\end{proof}
	
	\cite[Th\'eor\`eme principal]{beilinsonbernstein1981} now implies:
	
	\begin{corollary}[{\cite[Theorem 2.1.3 (2)]{harris2013}}]
		The flag variety of a connected reductive algebraic group over a field $F$ of characteristic zero is $\mathcal D$-affine over $F$.
	\end{corollary}

	\subsection{Descent properties for twisted $\mathcal D$-modules}
	In this section, we note the fpqc descent property of twisted $\D$-modules with respect to the base. Let $\Aa$ be a tdo on $X$.
	\begin{theorem}\label{thm:fpqcDmoduledescent}
		\begin{enumerate}
			\item The correspondence $S'\mapsto \Mod_{\qc}(s_X^\ast\Aa)$ determines a stack of categories over the big fpqc site of $S$.
			\item The following properties are local in the fpqc topology of $S$:
			\begin{itemize}
				\item[(i)] local finite generation,
				\item[(ii)] local finite presentation,
				\item[(iii)] $\mathcal O_X$-flatness,
				\item[(iv)] $\mathcal O_X$-local freeness of finite type,
				\item[(v)] $\mathcal O_X$-local freeness of finite rank $n$.
			\end{itemize}
		\end{enumerate}
	\end{theorem}

	\begin{remark}
		A special case of (2) is the faithfully flat descent for $\Aa$-modules in the case $S=\Spec F$ and $S'=\Spec F'$, $F'/F$ an extension of fields. In particular, properties (i) - (v) over $F$ may be checked over an algebraic or separably closed extension $F'$ of $F$.
	\end{remark}
	
	\begin{proof}
		For (1), we may assume that $S$ is affine. Since the base change is Zariski local in $X$, we may also assume that $X$ is affine to reduce the assertion to a corresponding statement in the ring setting. This is straightforward.
		
		Let us prove (2). First we observe that properties (i) and (ii) are local in the Zariski topology of $X$ and therefore of $S$. We may prove that if $X,S$ are affine and if we are given a faithfully flat morphism $s:S'\to S$ of affine schemes and a quasi-coherent left $\Aa$-module $\Mm$ such that $s^\ast_X \Mm$ satisfies (i) (resp.~(ii)), so is $\Mm$. Notice that $\Gamma(X\times_S S',s^\ast_X\Aa)$ is faithfully flat over $\Gamma(X,\Aa)$. The assertion now follows from \cite[Chap.\ I, \S 3, no.\ 6, Proposition 11]{bourbaki}.
		
		Properties (iii), (iv) and (v), may be verified within the category of quasi-coherent $\mathcal O_X$-modules, where the statements follow from \cite[Proposition (2.5.1) and Proposition (2.5.2), (iii) and (iv)]{ega42}.
	\end{proof}
	
	\begin{corollary}\label{cor:galoisdescent}
		Let $s:S'\to S$ be a Galois covering of Galois group $\Gamma$, and $\Aa$ be a tdo on $X$. Then $\Mod_{\qc}(\Aa)$ is equivalent to the category of quasi-coherent left $s_X^\ast\Aa$-modules with semilinear actions of $\Gamma$.
	\end{corollary}

	\subsection{Transfer bimodules}\label{sec:bimodules}
	
	Let $f$ be either a homomorphism $R\to R'$ of smooth $k$-algebras or a morphism $X\to Y$ of smooth $S$-schemes. Let $A$ be a tdo for $R/k$, and $\Aa$ be a tdo on $Y$.
	
	\begin{definition}
		We set $A_{R'\gets R}\coloneqq R'\otimes_R A$ (resp.~$\Aa_{X\to Y}\coloneqq f^\ast\Aa$). We regard it as an $(f_\cdot A,A)$- (resp.~$(f^\cdot \Aa,f^{-1}\Aa)$-)bimodule over $k$ (resp.~$x^{-1}\OO_S$) for the multiplication of $A$ (resp.~$f^{-1}\Aa$) from the right side and the left action of $f_\cdot A$ (resp.~$f^\cdot\Aa$) defined by \eqref{eq:bimod}. 
	\end{definition}
	
	In fact, recall that for $(\partial,\tilde{e})\in f_\cdot \Lie A$, $\iota((\partial,\tilde{e}))$ commutes with the right action of $A$ and likewise for $f^\cdot\Lie\Aa$ by definition of the target of $\iota$ in \eqref{eq:bimod}. This implies that $A_{R'\gets R}$ and $\Aa_{X\to Y}$ are really bimodules.
	
	\begin{remark}
		The action of $f_\cdot \Lie A$ on $A_{R'\gets R}$ is given by
		\begin{equation}
			(\partial,\sum b_i\otimes e_i)\cdot
			b\otimes P
			=\partial(b)\otimes P+\sum bb_i\otimes e_iP
			\label{eq:actionoffcdotEonbimod}
		\end{equation}
		for $(\partial,\sum b_i\otimes e_i)\in f_\cdot \Lie A$ and $b\otimes P\in A_{R'\gets R}$.
	\end{remark}

	\begin{construction}\label{cons:change_tdo}
		Let $M$ be a left $A$-module. Then $A_{R'\gets R}\otimes_A M$ is endowed with the canonical structure of a left $f_\cdot A$-module. Identify $A_{R'\gets R}\otimes_A M$ with $R'\otimes_R M$ to get the structure of a left $f_\cdot
		A$-module on $R'\otimes_R M$. We define a lift $f^\ast:\Mod(\Aa)\to\Mod(f^\cdot \Aa)$ of the pullback functor $f^\ast:\Mod(\OO_Y)\to\Mod(\OO_X)$ in a similar way.
	\end{construction}
	
	\begin{lemma}
		The action of $f_\cdot \Lie A$ on $R'\otimes_R M$ is given by
		\begin{equation}
			(\partial,\sum b_i\otimes e_i)(b\otimes m)
			=\partial(b)\otimes m+\sum bb_i\otimes e_im.
		\end{equation}
	\end{lemma}
	
	\begin{proof}
		The map
		\[R'\otimes_R M\cong A_{R'\gets R}\otimes_A M
		\xrightarrow{(\partial,\sum b_i\otimes e_i)}
		A_{R'\gets R}\otimes_A M
		\cong R'\otimes_R M\]
		is expressed as
		\[\begin{split}
			b\otimes m
			&\mapsto b\otimes 1\otimes m\\
			&\mapsto (\partial,\sum b_i\otimes e_i)\cdot (b\otimes 1)\otimes m\\
			&=(\partial(b)\otimes 1+\sum bb_i\otimes e_i)
			\otimes m\\
			&\mapsto \partial(b)\otimes m+\sum bb_i\otimes e_im.
		\end{split}\]
	\end{proof}

	We see the preservation of $\Aa$-coherence by $f^\ast$ for a digression:
	
	\begin{proposition}\label{prop:smoothpreimage}
		Let $f:X\to Y$ be a smooth morphism over $S$ of smooth $S$-schemes. Suppose that $X$ and $Y$ are locally Noetherian. Then $f^*$ preserves $\mathcal \Aa$-coherent modules.
	\end{proposition}
	
	\begin{proof}
		Let $\Mm$ be a coherent left $\Aa$-module. Let $U$ be an open subset of $Y$ such that $\Mm|_U$ is a finitely generated $\Aa$-module. In view of Proposition \ref{prop:coherentDmodules}, it will suffice to show that $(f^\ast\Mm)|_{f^{-1}(U)}$ is a coherent left $f^\cdot \Aa$-module. Choose a presentation
		\[(\Aa|_U)^n\to \Mm|_U\to 0,\]
		where $n$ is a nonnegative integer. Apply $f^\ast_U$ to the above exact sequence to get
		\[(f^\ast\Aa)^n|_{f^{-1}(U)}\to (f^\ast\Mm)|_{f^{-1}(U)}
		\to 0 \] 
		since $f_U$ is flat. The left $f^\cdot\Aa$-module $f^\ast\Aa$ is locally finitely generated by the exact sequence of locally free $\OO_X$-modules of finite rank
		\[0\to f^\ast\Omega^1_{Y/S}\to
		\Omega^1_{X/S}\to\Omega^1_{X/Y}\to 0\]
		(\cite[Proposition (17.2.3)]{ega44}). This completes the proof.
	\end{proof}
	
	Let us note the functoriality of the transfer bimodules.
	
	\begin{proposition}\label{prop:naturalityoffcdot}
		\begin{enumerate}
			\item For a Picard algebra $E$ for $R/k$, we have a commutative diagram
			\begin{equation}
				\begin{tikzcd}
					g_\cdot f_\cdot E\ar[d, "\alpha_{f,g}"]\ar[r]&
					g_\cdot f_\cdot A_E\ar[r]
					&\End_{A_E}((A_{f_\cdot E})_{R''\gets R'}\otimes_{A_{f_\cdot E}}(A_{E})_{R'\gets R})
					\ar[d, "\sim"{sloped, above}]\\
					(g\circ f)_\cdot E\ar[rr, hook, "\iota"]
					&&\End_{A_E}(A_{R''\gets R}).
				\end{tikzcd}
				\label{diag:associativityofdot}
			\end{equation}
			The left upper arrow is the canonical embedding. The right upper arrow is defined by application of Construction \ref{cons:change_tdo} to $(A_{E})_{R'\gets R}$.		
			Recall section \ref{sec:Xbasechange} for the definitions of $\alpha_{f,g}$ and $\iota$. The right vertical arrow is obtained by cancellation.
			\item For a Picard algebroid $\Ee$ on $Z$, we have a commutative diagram
			\begin{equation}
				\begin{tikzcd}
					g^\cdot f^\cdot \Ee\ar[d, "\alpha_{f,g}"]\ar[r]&
					g^\cdot f^\cdot \Aa_{\Ee}\ar[r]
					&\iEnd_{f^{-1}g^{-1}\Aa_{\Ee}}((\Aa_{f^\cdot \Ee})_{X\to Y}\otimes_{\Aa_{f^\cdot \Ee}}(\Aa_{\Ee})_{Y\gets Z})
					\ar[d, "\sim"{sloped, above}]\\
					(g\circ f)^\cdot \Ee\ar[rr, hook, "\iota"]
					&&\iEnd_{(g\circ f)^{-1}\Aa_{\Ee}}(\Aa_{X\to Z}).
				\end{tikzcd}
				\label{diag:associativityofdot/sch}
			\end{equation}
		\end{enumerate}
	\end{proposition}
	
	\begin{proof}
		We only prove (1).
		Let $E$ be a Picard algebra over $R$. We compare the left actions of $(g\circ f)_\cdot E$ on
		\[(A_{f_\cdot E})_{R''\gets R'}\otimes_{A_{f_\cdot E}}(A_{E})_{R'\gets R}\cong A_{R''\gets R}=R''\otimes_R A_E.\]
		By passage to the right $A_E$-linear action, we may restrict ourselves to $R''$ (recall Lemma \ref{lem:R-linear}). Then the assertion is evident by definitions. We note that the left action of $f_\cdot E$ at $1\otimes 1\in (A_{E})_{R'\gets R}$ agrees with the composition of the projection $f_\cdot E\to R'\otimes_R E$ and the canonical embedding $R'\otimes_R E\to R'\otimes_R A_E$.
	\end{proof}

	\begin{corollary}
		Let $R\overset{f}{\to} R'\overset{g}{\to} R''$ be a homomorphism of smooth $k$-algebras, and $A$ be a tdo for $R/k$. Then there is a natural isomorphism $R''\otimes_{R'}(R'\otimes_R -)\cong R''\otimes_R-$ of functors from $\Mod(A)$ to $\Mod((g\circ f)_\cdot A)$ via the identification $(g\circ f)_\cdot A\cong g_\cdot f_\cdot A$. Similarly, we have a natural isomorphism $f^\ast \circ g^\ast\cong (g\circ f)^\ast$ of functors from $\Mod(\Aa)$ to $\Mod((g\circ f)^\cdot\Aa)$ 
		for a sequence $X\overset{f}{\to} Y\overset{g}{\to} Z$ of morphisms of smooth $S$-scheme and a tdo $\Aa$ on $Z$.
	\end{corollary}

	We apply the idea of the proof of Proposition \ref{prop:naturalityoffcdot} to the base change functor to obtain:
	
	\begin{proposition}\label{prop:basechangevspullbackformodules}
		Let $k\to k'$ be a homomorphism of commutative rings, and $f:R\to R'$ be a homomorphism of smooth $k$-algebras. Let $A$ be a tdo for $R/k$. 
		\begin{enumerate}
			\item The isomorphism
			\[(k'\otimes_k A)_{k'\otimes_k R'\gets k'\otimes_k R}
			\cong k'\otimes_k A_{R'\gets R}\]
			defined by canceling $k'$ is $(k'\otimes_k f_\cdot A,A)$-bilinear. Here $k'\otimes_k f_\cdot A$ acts on
			\[(k'\otimes_k A)_{k'\otimes_k R'\gets k'\otimes_k R}\]
			via the identification $k'\otimes_k f_\cdot A\cong f'_\cdot(k'\otimes_k A)$ in (the proof of) Theorem \ref{sec:compatibility}.
			\item The isomorphism $(k'\otimes_k R')\otimes_{k'\otimes_k R} (k'\otimes_k-)\cong k'\otimes_k (R'\otimes_R-)$ of functors from $\Mod(R)$ to $\Mod(k'\otimes_k R')$ lifts to that from $\Mod(A)$ to $\Mod(k'\otimes_k f_\cdot A)$. A similar assertion holds for right modules.
		\end{enumerate}
	\end{proposition}
	
	\begin{proof}
		Write $E=\Lie A$. For (1), it will suffice to prove that the diagram
		\[\begin{tikzcd}
			k'\otimes_k f_\cdot E\ar[d]\ar[r]&
			k'\otimes_k f_\cdot A\ar[r]
			&\End_{A}(k'\otimes_k A_{R'\gets R})\\
			f'_\cdot(k'\otimes_k E)\ar[rr, hook, "\iota"]
			&&\End_{k'\otimes_k A}((k'\otimes_k A)_{k'\otimes_k R'\gets k'\otimes_k R}).
			\ar[u, hook]
		\end{tikzcd}\]
		commutes. Here the left upper horizontal arrow is defined by the base change of the left action of $f_\cdot A$ on $A_{R'\gets R}$. The left vertical arrow is as in Theorem \ref{thm:basechangevspullbackfortdo}. The right vertical arrow is defined by the canonical cancellation. To compare the actions, we may restrict ourselves to $k'\otimes_k R'$.
		
		Observe the the map	$k'\otimes_k f_\cdot E\to f'_\cdot(k'\otimes_k E)$
		of Picard algebras in (the proof of) Theorem \ref{thm:basechangevspullbackfortdo} satisfies
		the commutative diagrams
		\[\begin{tikzcd}
			k'\otimes_k f_\cdot E\ar[d]\ar[r, "k'\otimes_k\sigma"]
			&k'\otimes_k \Theta_{R'/k}\ar[d, "\sim"{sloped, above}]\\
			f'_\cdot(k'\otimes_k E)\ar[r, "\sigma"]
			&\Theta_{k'\otimes_k R'/k'},
		\end{tikzcd}\]
		\[\begin{tikzcd}
			k'\otimes_k f_\cdot E\ar[d]\ar[r, "k'\otimes_k\pr_2"]
			&k'\otimes_k (R'\otimes_RE)\ar[d, "\sim"{sloped, above}]\\
			f'_\cdot(k'\otimes_k E)\ar[r, "\pr_2"]
			&(k'\otimes_k R')\otimes_{k'\otimes_k R} (k\otimes_k E),
		\end{tikzcd}\]
		where $\pr_2$ are the canonical projections, and the right vertical arrows in the two diagrams are the canonical isomorphisms. The diagram in the first paragraph commutes by definition of $\iota$.
		
		Part (2) is a consequence of (1).
	\end{proof}
	
	A similar argument implies:
	
	\begin{proposition}\label{prop:basechangevspullbackformodules/sch}
		Let $s:S'\to S$ be a morphism of schemes.
		\begin{enumerate}
			\item The isomorphism of left $\OO_{X'}$-modules
			$(s^\ast_Y\Aa)_{X'\to Y'}\cong s^\ast_X \Aa_{X\to Y}$
			defined by
			$s^\ast_Xf^\ast\cong (f')^\ast s^\ast_Y$
			is $(s^\ast_X f^\cdot \Aa,(f')^{-1}s^{-1}_Y\Aa)$-bilinear. Here $s^\ast_X f^\cdot\Aa$ acts on $(s^\ast_Y\Aa)_{X'\to Y'}$ via the identification $s^\ast_X f^\cdot\Aa\cong (f')^\cdot s^\ast_Y\Aa$ in Theorem \ref{sec:compatibility}. For the right action of $(f')^{-1}s^{-1}_Y\Aa$ on $s^\ast_X \Aa_{X\to Y}$, use $(f')^{-1}s^{-1}_Y\cong s^{-1}_Xf^{-1}$.
			\item The isomorphism $s^\ast_Xf^\ast\cong (f')^\ast s^\ast_Y$ of functors from $\Mod(\OO_Y)$ to $\Mod(\OO_{X'})$ lifts to an isomorphism of functors from $\Mod(\Aa)$ to $\Mod(s^\ast_X f^\cdot\Aa)$. A similar assertion holds for right modules.
		\end{enumerate}
	\end{proposition}

	As Picard algebroids satisfy the \'etale descent property, one can see:
	
	\begin{theorem}\label{thm:etale_descent}
		The correspondence $U\mapsto \Mod_{\qc}(j^\cdot \Ee)$ with $j:U\to X$ \'etale determines a stack of categories over the small \'etale site of $X$.
	\end{theorem}
	
	A possible approach would be to introduce the notion of modules over torsors. Here we prove it directly. To save letters, we denote canonical homomorphisms by `$\can$' if there is no risk of confusion until the end of the proof of Theorem \ref{thm:etale_descent}. Let us also record:
	
	\begin{lemma}\label{lem:descent}
		Let $\{\Spec R_\lambda\to\Spec R\}$ be an fpqc covering, and $(M_\lambda,\phi_{\lambda,\mu})$ be a descent datum of $M$, i.e., a descent datum isomorphic to
		$(R_\lambda\otimes_R M,\can)$.
		Then the canonical map $M\to \prod_\lambda M_\lambda$ is injective.	
	\end{lemma}

	Observe that for an \'etale homomorphism $f:R\to R'$ of smooth $k$-algebras and a Picard algebra $E$ for $R/k$, the projection $f_\cdot E\to R'\otimes_R E$ is an isomorphism of $R'$-modules. We may and do identify them below. 
	
	\begin{proof}[Proof of Theorem \ref{thm:etale_descent}]
		
		We may assume $X$ and $S$ are affine. Write $X=\Spec R$ and $S=\Spec k$. Let $\{f_\lambda:R\to R_\lambda\}$ be a finite set of \'etale homomorphisms of smooth $k$-algebras such that the map $R\to\prod_\lambda R_\lambda$ is faithfully flat. For indices $\lambda,\mu,\nu$, write
		\[\begin{array}{cc}
			R_{\lambda\mu\nu}
			=R_\lambda\otimes_R R_\mu\otimes_R R_\nu,
			&R_{\lambda\mu}=R_\lambda\otimes_R R_\mu.
		\end{array}\]
		
		Let $M,M'$ be left $E$-modules. Let $\varphi,\varphi'$ be $E$-module homomorphisms from $M$ to $M'$. It follows by the usual faithfully flat descent for $R$-modules that if
		$R_\lambda\otimes_R\varphi=R_\lambda\otimes_R\varphi'$ for all $\lambda$ then we have $\varphi=\varphi'$.
		
		Let $\varphi_\lambda:R_\lambda\otimes_RM\to R_\lambda\otimes_R M'$ be isomorphisms of $(f_\lambda)_\cdot E$-modules such that
		\[R_{\lambda\mu}\otimes_{R_\lambda} \varphi_\lambda
		=R_{\lambda\mu}\otimes_{R_\mu}\varphi_\mu.\]
		Then we obtain an $R$-linear map $\varphi:M \to M'$ by the faithfully flat descent. We wish to prove that it is $E$-linear. This is verified by localizing the actions. That is, compare the commutative diagrams
		\[\begin{tikzcd}
			E\times M\ar[d, "\mathrm{action}"']\ar[r, "\can"]
			&\prod_\lambda((R_\lambda\otimes_R E)\times (R_\lambda\otimes_RM))
			\ar[d, "\mathrm{action}"]\\
			M\ar[d, "\varphi"']\ar[r, "\can"]&\prod_\lambda R_\lambda\otimes_R M
			\ar[d, "{(\varphi_\lambda)}"]\\
			M'\ar[r, hook, "\can"]&\prod_\lambda R_\lambda\otimes_R M',
		\end{tikzcd}\]
		\[\begin{tikzcd}
			E\times M\ar[r, "\can"]\ar[d, "E\times\varphi"']
			&\prod_\lambda((R_\lambda\otimes_R E)\times (R_\lambda\otimes_RM))
			\ar[d, "{((R_\lambda\otimes_R E)\times\varphi_\lambda)}"]\\
			E\times M'\ar[r, "\can"]\ar[d, "{\mathrm{action}}"']
			&\prod_\lambda((R_\lambda\otimes_R E)\times (R_\lambda\otimes_RM))
			\ar[d, "{\mathrm{action}}"]\\
			M'\ar[r, hook, "\can"]&\prod_\lambda R_\lambda\otimes_R M'.
		\end{tikzcd}\]
		The clockwise arrows coincide by our hypothesis on $\varphi_\lambda$. We conclude by Lemma \ref{lem:descent} that $\varphi$ respects the action of $E$.
		
		Finally, suppose that we are given a descent datum $(M_\lambda,\phi_{\lambda,\mu})$. Then we obtain an $R$-module $M$ by the faithfully flat descent of $R$-modules. Let $\phi_\lambda:M\to M_\lambda$ be the canonical homomorphism for each $\lambda$. For $m\in M$ and an index $\lambda$, we define an $R_\lambda$-linear map $\alpha_m:R_\lambda\otimes_R E\to M_\lambda$ by
		$\alpha_{m,\lambda}(a\otimes e)=(a\otimes e)\phi_\lambda(m)$. One can check in a similar way to the previous paragraph that $(\phi_\lambda)$ satisfy the descent condition. Let $\alpha_m:E\to M$ be the corresponding homomorphism. We now define $em=\alpha_m(e)$. It is easy to show that $M$ is an $E$-module for this action, and that the canonical descent datum attached to $M$ is isomorphic to $(M_\lambda,\phi_{\lambda,\mu})$. This completes the proof.
	\end{proof}

	The idea of the simplification $j^\cdot\Ee\cong j^\ast\Ee$ by the projection is also useful for constructing homomorphisms from Lie algebras to Picard algebroids:
	
	\begin{construction}[{\cite[4.6]{kashiwara1989}}]\label{cons:differentialaction}
		Let $G$ be a smooth $S$-affine group scheme over $S$, $X$ be a smooth $G$-scheme over $S$, and $\Ee$ be a $G$-equivariant Picard algebroid on $X$. Let $g:G\to S$ denote the structure morphism. Let $\pr_1$, $\pr_2$, and $i$ be as in section \ref{sec:equivsheaf}. Observe that we have a map
		\[\pr^\ast_1g^\ast\lieg\to\pr^\ast_1\Theta_{X/S}\to \pr^\cdot_2\Ee
		\cong a^\cdot\Ee\overset{\pr_2}{\to}a^\ast\Ee.\]
		Here for the first (resp.~second) arrow, see Example \ref{ex:D_G} (resp.~Proposition \ref{prop:algebroidtdoprojection}). The isomorphism $\pr^\cdot_2\Ee
		\cong a^\cdot\Ee$ is the structure constraint of the equivariant structure. The last map is the projection. Apply $i^\ast$ to the above sequence, and pass to the adjunction of $(x^\ast,x_\ast)$ to get a Lie algebra homomorphism $\lieg\to x_\ast\Ee$.
	\end{construction}

	\subsection{Derived inverse image}\label{sec:derivedinverse}
	
	Let $\Aa$ be a tdo on $Y$. 
	
	\begin{theorem}\label{thm:inverseimage}
		\begin{enumerate}
			\item The functor $f^\ast:\Mod(\Aa)\to \Mod(f^\cdot \Aa)$ admits a left derived functor 
			\[Lf^\ast:D(\Aa)\to D(f^\cdot\Aa),\]
			which satisfies the 2-commutative diagram
			\begin{equation}
				\begin{tikzcd}
					D(\Aa)\ar[r, "Lf^\ast"]\ar[d]
					&D(f^\cdot\Aa)\ar[d]\\
					D(\OO_Y)\ar[r, "Lf^\ast"]
					&D(\OO_X),\label{diag:compatibilitypullback}
				\end{tikzcd}
			\end{equation}
			where the vertical arrows are defined by restriction. We call it the derived inverse image of twisted left $\D$-modules.
			\item The functor $Lf^\ast$ restricts to
			$Lf^\ast:D_{\qc}(\Aa)\to D_{\qc}(f^\cdot\Aa)$.
			\item For another morphism $g:Y\to Z$ of smooth $S$-schemes and a tdo $\mathcal B$ on $Z$, 
			there is an equivalence
			$L(g\circ f)^\ast\simeq Lf^\ast\circ Lg^\ast$
			of functors from $D(\mathcal B)$ to $D(f^\cdot g^\cdot \Bb)$ which lifts the same relation of functors from $D(\OO_Z)$ to $D(\OO_X)$ (\cite[(3.6.4)$^\ast$]{lipman2009}).
			\item There exists an equivalence $Ls^\ast_X\circ Lf^\ast\simeq L(f')^\ast\circ Ls^\ast_Y$.
		\end{enumerate}
	\end{theorem}
	
	\begin{proof}
		The left derived functor $Lf^\ast$ exists by a similar argument to Corollary \ref{cor:derivedtensorproduct} (2). Moreover, it is computed by q-flat resolutions by its proof. The 2-commutativity of \eqref{diag:compatibilitypullback} then follows from Proposition \ref{prop:resolutions}.
		
		Part (2) follows from (1) and \cite[Proposition (3.9.1)]{lipman2009}.
		
		For (3), let $\Ff^\bullet$ be a q-flat complex of left $\Bb$-modules. Then $L(g\circ f)^\ast\Ff^\bullet \simeq (g\circ f)^\ast\Ff^\bullet$ and $Lg^\ast\Ff^\bullet \simeq g^\ast\Ff^\bullet$. Let $\Gg^\bullet\to g^\ast\Ff^\bullet$ be a q-flat resolution as a complex of left $g^\cdot\Bb$-modules. Since $g^\ast\Ff^\bullet$ is q-flat as a complex of left $\OO_Y$-modules, we get
		\[L(g\circ f)^\ast \Ff^\bullet\simeq 
		(g\circ f)^\ast \Ff^\bullet \cong f^\ast g^\ast \Ff^\bullet \simeq f^\ast\Gg^\bullet \simeq Lf^\ast (g^\ast\Ff^\bullet)
		\simeq Lf^\ast (Lg^\ast\Ff^\bullet).\]
		
		Part (4) follows by a similar argument to (3).
	\end{proof}

	We also define the derived inverse image of twisted right $\D$-modules by
	\[\begin{split}
		D_{\mathrm{r}}(\Aa)
		&\xrightarrow{-\otimes_{\OO_Y}\omega^\vee_{Y/S}}
		D_{\mathrm{r}}(\omega_{Y/S}\otimes_{\OO_Y}\Aa\otimes_{\OO_Y} \omega^\vee_{Y/S})\\
		&=D((\omega_{Y/S}\otimes_{\OO_Y}\Aa\otimes_{\OO_Y} \omega^\vee_{Y/S})^{\op})\\
		&\overset{Lf^\ast}{\to}
		D(f^\cdot((\omega_{Y/S}\otimes_{\OO_Y}\Aa\otimes_{\OO_Y} \omega^\vee_{Y/S})^{\op}))\\
		&\cong D((\omega_{X/S}\otimes_{\OO_X}f^\cdot(\Aa)\otimes_{\OO_X} \omega^\vee_{X/S})^{\op})\\
		&\cong D_{\mathrm{r}}(\omega_{X/S}\otimes_{\OO_X}f^\cdot(\Aa)\otimes_{\OO_X} \omega^\vee_{X/S})\\
		&\xrightarrow{-\otimes_{\OO_X}\omega_{X/S}} D_{\mathrm{r}}(f^\cdot\Aa)
	\end{split}\]
	which we denote by the same symbol $Lf^\ast$. One can easily show that similar properties to the above hold for the derived inverse image of twisted right $\D$-modules.

	\subsection{Derived direct image}\label{sec:deriveddirectimage}
	The aim of this section is to study the derived direct image functor over general bases.
	
	\subsubsection{Definition}\label{sec:deff+}
	
	We start with the ring setting. Let $f:R\to R'$ be a homomorphism of smooth $k$-algebras, and $A$ be a tdo for $R/k$. Recall that we introduced the $(f_\cdot A,A)$-bimodule $A_{R'\gets R}$. 
	
	\begin{definition}
		Set $f^+\coloneqq -\otimes^L_{f_\cdot A} A_{R'\gets R}:D_{{\mathrm{r}}}(f_\cdot A)\to D_{{\mathrm{r}}}(A)$.
	\end{definition}
	
	We next work with the sheaf setting. Let $f:X\to Y$ be a morphism of smooth $S$-schemes, and $\Bb$ be a tdo on $Y$. We define the derived direct image functor $f_+:D(f^\cdot\Bb)\to D(\Bb)$ as $f_+=Rf_\ast(-\otimes^L_{f^\cdot\Bb}\Bb_{X\to Y})$, where
	
	\begin{equation}
		-\otimes^L_{f^\cdot\Bb}\Bb_{X\to Y}:
		D(f^\cdot\Bb)\to
		D(f^{-1}\Bb)
		\label{eq:derivedtensorforderivedpush}
	\end{equation}
	\begin{equation}
		Rf_*: D(f^{-1}\Bb)\to
		D(\Bb).
		\label{eq:ordinaryderivedpushforderivedpush}
	\end{equation}
	We also define the direct image functor for left modules by
	\[\begin{split}
		D(f^\cdot\Bb)
		&=D_{\mathrm{r}} ((f^\cdot \Bb)^{\op})\\
		&\xrightarrow{-\otimes_{\OO_X}\omega_{X/S}}
		D_{\mathrm{r}} (\omega^\vee_{X/S}\otimes_{\OO_X}(f^\cdot \Bb)^{\op}\otimes_{\OO_X}\omega_{X/S})\\
		&\cong D_{\mathrm{r}}(f^\cdot (\omega_{Y/S}^\vee\otimes_{\OO_Y}\Bb^{\op} 
		\otimes_{\OO_Y}\omega_{Y/S}))\\
		&\overset{f_+}{\to}D_{\mathrm{r}}(\omega_{Y/S}^\vee\otimes_{\OO_Y}\Bb^{\op} 
		\otimes_{\OO_Y}\omega_{Y/S})\\
		&\xrightarrow{-\otimes_{\OO_Y}\omega^\vee_{Y/S}}
		D_{\mathrm{r}}(\Bb^{\op})\\
		&=D(\Bb),
	\end{split}\]
	which we denote by the same symbol $f_+$. For the isomorphism in the third line, see Proposition \ref{prop:pullbacksidechange}. Define an $(f^{-1}\Bb,f^\cdot\Bb)$-bimodule $\Bb_{Y\gets X}$ by
	\begin{equation}
		\omega_{X/S}\otimes_{\OO_X}(\omega^\vee_{Y/S}\otimes_{\OO_Y}\Bb^{\op}\otimes_{\OO_Y}\omega_{Y/S})_{X\to Y}\otimes_{f^{-1}\OO_Y} f^{-1}\omega^\vee_{Y/S},
	\end{equation}
	which may be also expressed as
	\[\omega_{X/S}\otimes_{f^{-1}\OO_Y} f^{-1}\omega_{Y/S}^\vee\otimes_{f^{-1}\OO_Y}f^{-1}\Bb^{\op}\cong f^{-1}\Bb \otimes_{f^{-1}\OO_Y} f^{-1}\omega_{Y/S}^\vee \otimes_{f^{-1}\OO_Y} \omega_{X/S}.\]
	We caution that the above isomorphism is just as bimodules. One can then identify
	\[f_+:D(f^\cdot\Bb)\to D(\Bb)\]
	with $Rf_\ast(\Bb_{Y\gets X}\otimes_{f^\cdot\Bb} \Mm)$ (cf.\ \cite[section 1.3, (C.2.15)]{hottaetal2008}). The direct image $f_+$ is local in $Y$ in the following sense:
	
	\begin{proposition}
		For $U$ an arbitrary open subscheme of $Y$, there is a canonical equivalence
		\begin{equation}
			(-)|_U \circ f_+\simeq (f_U)_+(-|_{f^{-1}(U)}).
			\label{eq:localityoff+}
		\end{equation}
	\end{proposition}
	
	\begin{proof}
		Use Example \ref{ex:bcthmforopimm}.
	\end{proof}

	To see the naturality of composition, we need an observation. Let $g:Y\to Z$ be a morphism of smooth $S$-schemes. Let $\Aa$ be a tdo on $Z$. Put $\Bb=g^\cdot\Aa$.
	
	\begin{lemma}\label{lem:acyclicityofbimodchain}
		We have $(g^\cdot\Aa)_{X\to Y}
		\otimes^L_{f^{-1}g^\cdot\Aa} f^{-1}\Aa_{Y\to Z}\simeq \Aa_{X\to Z}$ as $(f^\cdot g^\cdot\Aa,f^{-1}g^{-1}\Aa)$-bimodules.
	\end{lemma}
	
	\begin{proof}
		Recall that $\Aa_{Y\to Z}=g^\ast\Aa$ is flat as a left $\OO_Y$-module by Lemma \ref{prop:tdofreebasis}. Identify
		\[(g^\cdot\Aa)_{X\to Y}
		\otimes_{f^{-1}g^\cdot\Aa}f^{-1}(-)\]
		with $f^\ast$ to see that $\Aa_{Y\to Z}$ is $(g^\cdot\Aa)_{X\to Y}
		\otimes_{f^{-1}g^\cdot\Aa}f^{-1}(-)$-acyclic. Therefore we get
		\[(g^\cdot\Aa)_{X\to Y}
		\otimes^L_{f^{-1}g^\cdot\Aa} f^{-1}\Aa_{Y\to Z}
		\simeq (g^\cdot\Aa)_{X\to Y}
		\otimes_{f^{-1}g^\cdot\Aa} f^{-1}\Aa_{Y\to Z}.\]
		We have an isomorphism
		$(g^\cdot\Aa)_{X\to Y}
		\otimes_{f^{-1}g^\cdot\Aa} f^{-1}\Aa_{Y\to Z}\cong \Aa_{X\to Z}$
		of $(f^\cdot g^\cdot\Aa, f^{-1}g^{-1}\Aa)$-bimodules by \eqref{diag:associativityofdot/sch}. This completes the proof.
	\end{proof}
	
	\begin{construction}\label{cons:compositionlaw}
		Let $f:X\to Y$ and $g:Y\to Z$ morphisms of smooth $S$-schemes. Let $\Aa$ be a tdo on $Z$. Then we define a natural transformation $g_+\circ f_+\to (g\circ f)_+$ by
		\[\begin{split}
			g_+(f_+(-))
			&=Rg_\ast(f_+(-)
			\otimes^L_{g^\cdot\Aa} \Aa_{Y\to Z})\\
			&=Rg_\ast(Rf_\ast(-
			\otimes^L_{f^\cdot g^\cdot \Aa} (g^\cdot\Aa)_{X\to Y})	
			\otimes^L_{g^\cdot\Aa}\Aa_{Y\to Z})\\
			&\to Rg_\ast Rf_\ast((-
			\otimes^L_{f^\cdot g^\cdot \Aa} (g^\cdot\Aa)_{X\to Y})	
			\otimes^L_{f^{-1}g^\cdot\Aa} f^{-1}\Aa_{Y \to Z})\\
			&\simeq R(g\circ f)_\ast(-
			\otimes^L_{f^\cdot g^\cdot \Aa} ((g^\cdot\Aa)_{X\to Y}
			\otimes^L_{f^{-1}g^\cdot\Aa} f^{-1}\Aa_{Y \to Z}))\\
			&\simeq R(g\circ f)_\ast(-
			\otimes^L_{f^\cdot g^\cdot \Aa} \Aa_{X\to Z})\\
			&=(g\circ f)_+(-).
		\end{split}\]
	\end{construction}
	
	\begin{theorem}\label{thm:fpluscompositions1}
		The natural transformation $g_+\circ f_+\to (g\circ f)_+$ in Construction \ref{cons:compositionlaw} is an equivalence if $f$ is a closed immersion.
	\end{theorem}
	
	\begin{proof}
		This is immediate from Theorem \ref{thm:projformula}.
	\end{proof}
	
	Recall that every morphism $f:X\to Y$ of schemes over a base scheme $S$ is decomposed into the graph morphism $i:X\to X\times_S Y$ and the projection $\pr_2:X\times_S Y\to Y$. Notice that $i$ is an immersion. Moreover, if the structure morphism $y:Y\to S$ is separated then $i$ is a closed immersion. This happens, for instance, if we localize $Y$ and $S$ so that they are affine. Hence Theorem \ref{thm:fpluscompositions1} reduces the study of $f_+$ to the cases of closed immersions and projections. For this reason, we will revisit $f_+$ after the case-by-case study. We will also prove the equivalence $g_+\circ f_+\simeq (g\circ f)_+$ under easy hypotheses.
	
	To end this small section, we compare $f^+$ and $f_+$. Put $S=\Spec k$, $X=\Spec R'$, and $Y=\Spec R$. We denote the structure (homo)morphisms $R\to R'$ and $X\to Y$ by the same symbol $f$. Let $\Aa$ be a tdo on $Y$. Set $A=\Gamma(Y,\Aa)$.
	
	\begin{construction}\label{cons:comparef+}
		We define a natural transformation $f^+\circ R\Gamma(X,-)\to R\Gamma(Y,f_+(-))$ of the functors from $D_{\mathrm{r}}(f^\cdot\Aa)$ to $D_{\mathrm{r}}(A)$ as follows:
		\[\begin{split}
			f^+\circ R\Gamma(X,-)&\simeq R\Gamma(X,-)\otimes^L_{f_\cdot A} R\Gamma(X,\Aa_{X\to Y})\\
			&\to R\Gamma(X,-\otimes^L_{f^\cdot \Aa} \Aa_{X\to Y})\\
			&\simeq R\Gamma(Y,Rf_\ast(-\otimes^L_{f^\cdot \Aa} \Aa_{X\to Y}))\\
			&=R\Gamma(Y,f_+(-)).
		\end{split}\]
		For the first equivalence, we may regard $\Aa_{X\to Y}$ as a left $\OO_X$-module. Then we get a canonical equivalence $R\Gamma(X,\Aa_{X\to Y})\simeq A_{R'\gets R}$ since $\Aa_{X\to Y}$ is a quasi-coherent $\OO_X$-module. The arrow in the second line arises from the canonical lax monoidal structure of $R\Gamma(X,-)$.
	\end{construction}
	
	\begin{theorem}\label{thm:comparisonf_+f^+}
		Suppose $\dim^+(f,\Aa)<\infty$ (see Definition \ref{defn:dim^+(f,R)}). Then the functor
		\[f_+:D_{{\mathrm{r}}}(f^\cdot \Aa)\to D_{{\mathrm{r}}}(\Aa)\]
		restricts to $D^-_{{\mathrm{r}},\qc}(f^\cdot \Aa)\to D^-_{{\mathrm{r}},\qc}(\Aa)$. Moreover, the natural transformation of Construction \ref{cons:comparef+} is an equivalence on $D_{{\mathrm{r}},\qc}(f^\cdot\Aa)$. In particular, the diagram
		\begin{equation}
			\begin{tikzcd}
				D^-_{{\mathrm{r}},\qc}(f^\cdot \Aa)\ar[r, "f_+"]
				\ar[d, "{R\Gamma(X,-)}"']
				&D^-_{{\mathrm{r}},\qc}(\Aa)\ar[d, "{R\Gamma(Y,-)}"]\\
				D_{{\mathrm{r}}}(f_\cdot A)\ar[r, "f^+"]
				&D_{{\mathrm{r}}}(A)
			\end{tikzcd}
			\label{diag:comparisonf_+f^+}
		\end{equation}
		is 2-commutative.
	\end{theorem}
	
	The preservation of quasi-coherence is nontrivial since we do not have a natural structure of $\OO_Y$-module on $\Mm\otimes^L_{f^\cdot\Aa}\Aa_{X\to Y}$.

	\begin{proof}
		First, let $\Ff$ be a free right $f^\cdot\Aa$-module. Then we have a canonical equivalence
		\[f_+\Ff\simeq f_\ast(\Ff\otimes_{f^\cdot\Aa} \Aa_{X\to Y})\]
		since the sheaf $\Ff\otimes_{f^\cdot\Aa} \Aa_{X\to Y}$ of abelian groups admits the structure of a quasi-coherent $\OO_X$-module (see below Proposition \ref{prop:independent}). As a right $f^{-1}\Aa$-module, $\Ff\otimes_{f^\cdot\Aa} \Aa_{X\to Y}$ can be identified with the coproduct of copies of $\Aa_{X\to Y}$. Since $f$ is quasi-compact, $f_\ast(\Ff\otimes_{f^\cdot\Aa} \Aa_{X\to Y})$ is the coproduct of copies of the right $\Aa$-module $f_\ast\Aa_{X\to Y}$.
		
		We wish to prove that $f_\ast\Aa_{X\to Y}$ is a quasi-coherent right $\Aa$-module. Set $F_\bullet \Aa_{X\to Y}=f^\ast F_\bullet\Aa$. Then $F_\bullet \Aa_{X\to Y}$ is an exhaustive filtration on $\Aa_{X\to Y}$ with $\Gr^\bullet\Aa_{X\to Y}\cong f^\ast\Gr^\bullet\Aa$ as an $(\OO_X,f^{-1}\OO_Y)$-bimodule since $\Gr^p \Aa$ is a locally free $\OO_Y$-module of finite rank for every integer $p$. In particular, the induced structures of left and right $\OO_X$-modules on $\Gr^\bullet\Aa_{X\to Y}$ coincide. Since $F_p \Aa_{X\to Y}$ are quasi-coherent left $\OO_X$-modules for all $p\geq 0$, $f_\ast F_\bullet \Aa_{X\to Y}$ is an exhaustive filtration of $\OO_X$-bimodules with $\Gr^\bullet f_\ast \Aa_{X\to Y}\cong f_\ast\Gr^\bullet\Aa_{X\to Y}\cong f_\ast f^\ast \Gr^\bullet\Aa$. The assertion now follows since $f_\ast f^\ast \Gr^p\Aa$ is a quasi-coherent (left or right) $\OO_X$-module for every integer $p$.
		
		We conclude by the former two paragraphs that if $\Ff$ is a free right $f^\cdot\Aa$-module, $f_+\Aa$ lies in $\Mod_{{\mathrm{r}},\qc}(\Aa)$ (up to equivalences in the derived category). We next prove the canonical homomorphism $f^+(R\Gamma(X,\Ff))\to R\Gamma(Y,f_+\Ff)$ is an equivalence. In view of Proposition \ref{prop:derivedglobalsectionforaffineschemes}, we may drop the derived operations, i.e., we are reduced to showing that the canonical homomorphism $\Gamma(X,\Ff)\otimes_{f_\cdot A} A_{R'\gets R}\to \Gamma(Y,f_\ast(\Ff\otimes_{f^\cdot\Aa} \Aa_{X\to Y}))$ obtained as an abelian analog of Construction \ref{cons:comparef+} is an isomorphism. We may again pass to coproducts to replace $\Ff$ by $f^\cdot\Aa$. Then the assertion follows from Lemma \ref{lem:Dquasicoherence}.
		
		The general statement follows from Proposition \ref{prop:derivedglobalsectionforaffineschemes} and the proofs of \cite[Lemma (1.11.3) (iii) and Complement (1.11.3.1)]{lipman2009}.
	\end{proof}
	
	\begin{variant}\label{var:comparisonf_+f^+;imm}
		Put $S=\Spec k$ again. Let $i:Y\to X$ be an immersion of smooth affine $S$-schemes. Write $R$ and $R'$ for the coordinate rings of $X$ and $Y$ respectively. We denote the corresponding $k$-algebra homomorphism $R\to R'$ by $p$. Let $\Aa$ be a tdo on $X$. Write $A=\Gamma(X,\Aa)$. Suppose that $\dim(i,\Aa)<\infty$. 
		\begin{enumerate}
			\item The functor $i_+$ restricts to $i_\ast(-\otimes_{i^\cdot\Aa} \Aa_{Y\to X}):\Mod_{{\mathrm{r}},\qc}(i^\cdot \Aa)\to \Mod_{{\mathrm{r}},\qc}(\Aa)\subset D^-_{{\mathrm{r}},\qc}(\Aa)$.
			\item As a left $p_\cdot A$-module, $A_{R'\gets R}$ is flat.
			\item The diagram \eqref{diag:comparisonf_+f^+} restricts to
			\begin{equation}
				\begin{tikzcd}
					\Mod_{{\mathrm{r}},\qc}(i^\cdot \Aa)\ar[rr, "{i_\ast(-\otimes_{i^\cdot\Aa} \Aa_{Y\to X})}"]
					\ar[d, "{\Gamma(X,-)}"']
					&&\Mod_{{\mathrm{r}},\qc}(\Aa)\ar[d, "{\Gamma(Y,-)}"]\\
					\Mod_{{\mathrm{r}}}(f_\cdot A)\ar[rr, "{-\otimes_{p_\cdot A} A_{R'\gets R}}"]
					&&\Mod_{{\mathrm{r}}}(A)
				\end{tikzcd}
				\label{diag:comparisoni_+i^+}
			\end{equation}
			\item The functor $i_\ast(-\otimes_{i^\cdot\Aa} \Aa_{Y\to X}):\Mod_{{\mathrm{r}},\qc}(i^\cdot \Aa)\to \Mod_{{\mathrm{r}},\qc}(\Aa)$ respects small colimits and finite limits.
		\end{enumerate}
	\end{variant}
	
	We will prove it in the next section. The assertion below is immediate from Variant \ref{var:comparisonf_+f^+;imm}:
	
	\begin{corollary}\label{cor:directimageforqcptimmersions}
		Let $S$ be a scheme, $i:Y\to X$ be an affine immersion of smooth $S$-schemes, and $\Aa$ be a tdo on $X$. Assume that $Ri_\ast:D_{\mathrm{r}}(i^{-1}\Aa)\to D_{\mathrm{r}}(\Aa)$ is locally bounded. Then $i_+$ restricts to an exact and cocontinuous functor
		\[i_\ast(-\otimes_{i^\cdot\Aa} \Aa_{Y\to X}):
		\Mod_{{\mathrm{r}},\qc}(i^\cdot\Aa)\to\Mod_{{\mathrm{r}},\qc}(\Aa)\subset D_{\mathrm{r}}(\Aa).\]
	\end{corollary}
	
	Note that the hypotheses on $i$ hold if $i$ is a closed immersion.
	
	\subsubsection{Direct image for immersion}
	
	In this section, we plan to prove Variant \ref{var:comparisonf_+f^+;imm}. We start with a more general setting: Let $i:Y\hookrightarrow X$ be an immersion of smooth schemes over an arbitrary base scheme $S$, and $\Aa$ be a tdo on $X$. We study the direct image functor for $i$. Define a map $\act_1:i^\cdot\Aa\to \Aa_{Y\to X}$ by the left action at $1\otimes 1$.

	\begin{lemma}\label{lem:bimodlocfreeimmcase}
		\begin{enumerate}
			\item As a left $i^\cdot\Aa$-module, $\Aa_{Y\to X}$ is locally free.
			\item The map $\act_1:i^\cdot\Aa\to \Aa_{Y\to X}$ locally admits an $i^\cdot\Aa$-linear retract.
		\end{enumerate}
	\end{lemma}
	
	Part (2) will be used in the next section.
	
	\begin{proof}
		We may replace $X$ and $Y$ by $U$ and $U\cap Y$ in Theorem \ref{thm:relativecoordinate} to assume that $X$ admits a relative coordinate $((x_1,x_2,\ldots,x_r),(x_{r+1},\ldots,x_n))$. Let $\partial_p$ be as in Theorem \ref{thm:relativecoordinate}.
		
		To compute the left action of $i^\cdot\Aa$ on $\Aa_{Y\to X}$, let us set $F_p\Aa_{Y\to X}= i^\ast F_p\Aa$ for each nonnegative integer $p$. Then $F_\bullet \Aa_{Y\to X}$ is an exhaustive filtration on $\Aa_{Y\to X}$ with $\Gr^\bullet\Aa_{Y\to X}\cong i^\ast\Gr^\bullet\Aa$ since $\Gr^p\Aa$ is locally free of finite rank as an $\OO_X$-module for every $p$. It is easy to show that the left action of the filtered sheaf $i^\cdot\Aa$ of rings on $\Aa_{Y\to X}$ respects the filtration. Then $\partial_i$ acts on $\Gr^\bullet\Aa_{Y\to X}$ by the multiplication for $1\leq i\leq r$ under the identifications $\Gr^1i^\cdot\Aa\cong\Theta_{Y/S}$ and $\Gr^\bullet\Aa_{Y\to X}\cong i^\ast\Gr^\bullet\Aa\cong \Sym_{\OO_Y}i^\ast\Theta_{X/S}$. Therefore the $i^\cdot\Aa$-linear map
		$\oplus_{i_{r+1},\ldots,i_n\in\NN}
		i^\cdot\Aa \to \Aa_{Y\to X}$ defined componentwisely by $1\mapsto 1\otimes e_{r+1}^{i_{r+1}}\cdots e_n^{i_n}$ is an isomorphism by passing to the associated graded sheaves and by comparing the $\OO_Y$-bases (see Proposition \ref{prop:tdofreebasisaffine}), where $e_i$ are lifts of $\partial_i$. In particular, $\{1\otimes e_{r+1}^{i_{r+1}}\cdots e_n^{i_n}\in\Aa_{Y\to X}(Y):~i_{r+1},\ldots,i_n\in\NN \}$
		is a free basis of $\Aa_{Y\to X}$ as a left $i^\cdot\Aa$-module. The projection to the component of $i_{r+1}=\cdots=i_n=0$ is an $i^\cdot\Aa$-linear retract of $\act_1$.
	\end{proof}
	
	\begin{proof}[Proof of Variant \ref{var:comparisonf_+f^+;imm}]
		Part (4) follows from (2), (3), and Lemma \ref{lem:modulecataffinecase}. Part (3) follows from (1), (2), and Proposition \ref{prop:derivedglobalsectionforaffineschemes}.
		
		To prove (1) and (2), let $\Mm$ be a quasi-coherent right $i^\cdot\Aa$-module. Write $M=\Gamma(Y,\Mm)$. Then by Theorem \ref{thm:comparisonf_+f^+}, we have an equivalence
		\begin{equation}
			R\Gamma(X,i_+\Mm)\simeq M\otimes^L_{p_\cdot A} A_{R'\gets R}.
			\label{eq:comparison1}	
		\end{equation}
		Since $i_+\Mm$ lies in $D^{\geq 0}(\Aa)$ by Lemma \ref{lem:bimodlocfreeimmcase} (1), $H^q(R\Gamma(X,i_+\Mm))$ vanishes if $q<0$. It is evident that $H^q(M\otimes^L_{p_\cdot A} A_{R'\gets R})$ vanishes if $q>0$. Therefore all the cohomology groups but in the zeroth degree of $R\Gamma(X,i_+\Mm)$ and $M\otimes^L_{p_\cdot A} A_{R'\gets R}$ vanish. Part (2) now follows since
		\[\Gamma(X,-):\Mod_{{\mathrm{r}},\qc}(i^\cdot\Aa)\to \Mod_{\mathrm{r}}(p_\cdot A)\]
		is essentially surjective (Lemma \ref{lem:modulecataffinecase}). We may also identify the equivalence \eqref{eq:comparison1} with
		\begin{equation}
			H^0(X,i_+ \Mm)\cong M\otimes_{p_\cdot A} A_{R'\gets R}
			\label{eq:comparison2}	
		\end{equation}
		by the vanishing of cohomology groups.
		
		Recall that for each integer $q$, we have an isomorphism 
		\[H^q(X,i_+\Mm)\cong \Gamma(X,H^q(i_+\Mm))\]
		by Theorem \ref{thm:comparisonf_+f^+} and Proposition \ref{prop:derivedglobalsectionforaffineschemes}. Part (1) now follows from Theorem \ref{thm:comparisonf_+f^+} and Lemma \ref{lem:modulecataffinecase}.
	\end{proof}

	\subsubsection{Kashiwara's theorem}\label{sec:kashiwara}
	
	Let $i:Y\hookrightarrow X$ be a closed immersion of smooth $S$-schemes, and $\Aa$ be a tdo on $X$. In this section, we aim to analyze Kashiwara's theorem on $i_+$. We caution that Kashiwara's theorem fails in our general setting:
	
	\begin{example}\label{ex:counterexofkashiwaraeq}
		Put $S=\Spec\FF_p$ with $p$ a positive prime. Let $i:Y=\Spec\FF_p\hookrightarrow\Spec \FF_p\left[x\right]=X$ be the closed immersion defined by the ideal $(x)\subset\FF_p\left[x\right]$. Let us work with $\Aa=\D_{X/S}$. Write
		$D(1)=D_{\FF_p\left[x\right]/\FF_p}$
		and $\partial= \frac{\partial}{\partial x}
		\in\Theta_{\FF_p\left[x\right]/\FF_p}\subset D(1)$.
		Then the structure of the right $D(1)$-module $M\coloneqq\Gamma(X,i_+\OO_S)$ is described as
		\[\begin{array}{ccc}
			M=\oplus_{n\geq 0} \FF_p\partial^n,
			&\partial^n\cdot x=n\partial^{n-1},
			&\partial^n\cdot \partial=\partial^{n+1}.
		\end{array}\]
		We have a nontrivial $D(1)$-linear endomorphism $\varphi$ determined by $\varphi(1)=\partial^p$. In particular, $i_+$ is not fully faithful on $\Mod_{{\mathrm{r}},\qc}(\D_{Y/S})$.
		
		Notice also that this computation implies that every nonzero right $D(1)$-module obtained by $i_+$ is of infinite dimension over $\FF_p$. On the other hand, observe that $\FF_p\left[x\right]$ is a right $D(1)$-module for
		$x^n\cdot x=x^{n+1}$ and 
		$x^n\cdot\partial=-nx^{n-1}$.
		It is clear that the ideal $(x^p)\subset\FF_p\left[x\right]$ is a $D(1)$-submodule. In particular, we obtain a nonzero finite dimensional right $D(1)$-module $N=\FF_p\left[x\right]/(x^p)$. Write $\mathcal N$ be the corresponding $\D_{X/S}$-module.
		Although $\mathcal N$ is supported on $Y$, it does not lie in the essential image of $i_+$ since $\dim_{\FF_p} N<\infty$.
	\end{example}
	
	\begin{remark}
		By Lemma \ref{lem:bimodlocfreeimmcase} (1), $i_+$ is still conservative (i.e., reflects isomorphisms). In fact, it is easy to show that $\Aa_{Y\to X}$ is nonzero unless $Y=\emptyset$.
	\end{remark}
	
	An evident obstruction for the failure of Kashiwara's theorem is due to $\ZZ$-torsion.
	
	\begin{definition}
		Let $Z$ be a topological space, and $\Rr$ be a sheaf of rings on $Z$. Then we write $\Mod_{{\mathrm{r}}}^{\ZZ\textrm{-}\mathrm{flat}}(\Rr)$ for the full subcategory of $\Mod_{{\mathrm{r}}}(\Rr)$ consisting of right $\Rr$-modules which are flat over $\ZZ_Z$.
	\end{definition}

	\begin{theorem}\label{thm:kashiwara}
		\begin{enumerate}
			\item The functor $-\otimes_{i^\cdot\Aa} \Aa_{Y\to X}:\Mod_{{\mathrm{r}}}(i^\cdot\Aa)\to \Mod_{{\mathrm{r}}}(i^{-1}\Aa)$ is faithful.
			\item The functor in (1) restricts to a fully faithful functor $\Mod_{{\mathrm{r}}}^{\ZZ\textrm{-}\mathrm{flat}}(i^\cdot\Aa)\to \Mod_{{\mathrm{r}}}(i^{-1}\Aa)$.
		\end{enumerate}
		
	\end{theorem}

	\begin{corollary}\label{cor:kashiwara}
		\begin{enumerate}
			\item The functor $i_+:\Mod_{{\mathrm{r}}}(i^\cdot\Aa)\to\Mod_{{\mathrm{r}}}(\Aa)$ is faithful.
			\item The functor in (1) restricts to a fully faithful functor $i_+:\Mod_{{\mathrm{r}}}^{\ZZ\textrm{-}\mathrm{flat}}(i^\cdot\Aa)\to\Mod_{{\mathrm{r}}}^{\ZZ\textrm{-}\mathrm{flat}}(\Aa)$.
		\end{enumerate}
	\end{corollary}
	
	\begin{proof}
		For a right $i^\cdot\Aa$-module $\Mm$, $\Mm\otimes_{i^\cdot\Aa} \Aa_{Y\to X}$ is locally a direct sum of copies of $\Mm$. In particular, if $\Mm$ is $\ZZ$-flat, so is $\Mm\otimes_{i^\cdot\Aa} \Aa_{Y\to X}$. Since $i_\ast$ respects $\ZZ$-flat sheaves, so does $i_+$.
		
		Since the counit of the adjunction
		$i^{-1}:\Mod_{{\mathrm{r}}}(\Aa)\rightleftarrows
		\Mod_{{\mathrm{r}}}(i^{-1}\Aa):i_\ast$
		is an isomorphism, $i_\ast$ is fully faithful. Combine it with Theorem \ref{thm:kashiwara} to deduce the assertions.
	\end{proof}

	\begin{proof}[Proof of Theorem \ref{thm:kashiwara}]
		Observe that for a right $i^\cdot\Aa$-module $\mathcal M$, \[\Mm\otimes_{i^\cdot\Aa}\act_1:\Mm\cong \Mm\otimes_{i^\cdot\Aa} i^\cdot\Aa\to \Mm\otimes_{i^\cdot\Aa}\Aa_{Y\to X}\]
		is monic by Lemma \ref{lem:bimodlocfreeimmcase} (2). In particular, (1) follows.
		
		Let $\Mm$ and $\mathcal N$ be right $i^\cdot\Aa$-modules, and $\varphi:\Mm\otimes_{i^\cdot\Aa} \Aa_{Y\to X}\to\mathcal N \otimes_{i^\cdot\Aa} \Aa_{Y\to X}$ be a right $i^{-1}\Aa$-linear map. We prove that $\varphi\circ (\Mm\otimes_{i^\cdot\Aa}\act_1)$ factors through $\mathcal{N}\otimes_{i^\cdot\Aa}\act_1:\mathcal N\hookrightarrow \mathcal N \otimes_{i^\cdot\Aa} \Aa_{Y\to X}$. We may replace $X$ and $Y$ with $U$ and $U\cap Y$ respectively in Theorem \ref{thm:relativecoordinate} to assume that $X$ admits a relative coordinate $((x_1,x_2,\ldots,x_r),(x_{r+1},\ldots,x_n))$. 
		Since
		\[\varphi(u\otimes 1\otimes 1)x_p
		=\varphi(u\otimes 1\otimes x_p)
		=\varphi(u\otimes x_p\otimes 1)
		=0\]
		for a local section $u\in\Mm$ and $r+1\leq p\leq n$, it will suffice to show that every local section $v\in \mathcal{N}\otimes_{i^\cdot\Aa} \Aa_{Y\to X}$ with $vx_p=0$ for $r+1\leq p\leq n$ belongs to $\mathcal N$.
		
		Since the assertion is local, we may express
		\[v=\sum_{I=(i_{r+1},i_{r+2},\ldots,i_n)} v_I\otimes \partial^{i_{r+1}}_{r+1}\cdots\partial^{i_{n}}_{n}.\]
		We wish to prove $v_I=0$ unless $i_{r+1}=i_{r+2}=\cdots =i_n=0$; otherwise, set
		\[q\coloneqq \max\{p\in\{r+1,\ldots,n\}:~
		\exists I=(i_{r+1},i_{r+2},\ldots,i_n)\
		\mathrm{s.t.}\ v_I\neq 0,\ i_p\neq 0\}\geq r+1.\]
		Then we have
		\begin{equation}
			\begin{split}
				vx_q
				&=\sum_{\overset{I=(i_{r+1},i_{r+2},\ldots,i_q,0,0,\ldots,0)}{i_q= 0}}
				v_I\otimes \partial^{i_{r+1}}_{r+1}\cdots\partial^{i_{q-1}}_{q-1}
				x_q\\
				&+\sum_{\overset{I=(i_{r+1},i_{r+2},\ldots,i_q,0,0,\ldots,0)}{i_q\neq 0}}
				v_I\otimes \partial^{i_{r+1}}_{r+1}\cdots\partial^{i_q}_qx_q\\
				&=\sum_{\overset{I=(i_{r+1},i_{r+2},\ldots,i_q,0,0,\ldots,0)}{i_q= 0}}
				v_I\otimes x_q
				\partial^{i_{r+1}}_{r+1}\cdots\partial^{i_{q-1}}_{q-1}\\
				&+\sum_{\overset{I=(i_{r+1},i_{r+2},\ldots,i_q,0,0,\ldots,0)}{i_q\neq 0}}
				v_I\otimes \partial^{i_{r+1}}_{r+1}\cdots
				\partial^{i_{q-1}}_{q-1}
				([\partial_q^{i_q},x_q]+x_q\partial_q^{i_q})\\
				&=\sum_{\overset{I=(i_{r+1},i_{r+2},\ldots,i_q,0,0,\ldots,0)}{i_q\neq 0}}
				i_q v_I\otimes \partial^{i_{r+1}}_{r+1}\cdots\partial^{i_{q-1}}_{q-1}
				\partial^{i_{q}-1}_{q}.
			\end{split}
			\label{eq:keycomputationtokashiwarathm}
		\end{equation}
		Since $\Mm$ is $\ZZ$-flat, $v_I=0$ for $I=(i_{r+1},i_{r+2},\ldots,i_q,0,0,\ldots,0)$ with $i_q\neq 0$. This contradicts to the definition of $q$.
	\end{proof}
	
	We next wish to know the essential image of $(i_+)|_{\Mod_{{\mathrm{r}}}^{\ZZ\textrm{-}\mathrm{flat}}(i^\cdot\Aa)}$. The problem is still delicate (see Example \ref{ex:counterexofessimg}). For this, let us see the adjointness property of $i_+$.
	
	\begin{lemma}\label{lem:i+cloimm}
		\begin{enumerate}
			\item For $\Mm\in \Mod_{\mathrm{r}}(i^\cdot\Aa)$, we have
			\[i_+\Mm\simeq i_\ast(\Mm\otimes_{i^\cdot\Aa}\Aa_{X\to Y})\in \Mod_{\mathrm{r}}(\Aa)\subset D_{\mathrm{r}}(\Aa).\]
			\item The functor $i_+$ restricts to an exact functor $\Mod_{\mathrm{r}}(i^\cdot\Aa)\to\Mod_{\mathrm{r}}(\Aa)$ which will be denoted by the same symbol $i_+$.
			\item The functor $i_+$ restricts to $D_{{\mathrm{r}},\qc}(i^\cdot\Aa)\to D_{{\mathrm{r}},\qc}(\Aa)$.
		\end{enumerate}
	\end{lemma}
	
	\begin{proof}
		This is evident since $i_\ast:\Mod(i^{-1}\OO_X)\to \Mod(\OO_X)$ is exact.
	\end{proof}

	Define a functor $i^\natural:\Mod_{\mathrm{r}}(\Aa)\to\Mod_{\mathrm{r}}(i^\cdot\Aa)$ by $\iHom_{i^{-1}\Aa}(\Aa_{Y\to X},i^{-1}(-))$. 
	
	Let $\Ii\subset \OO_X$ denote the ideal corresponding to $Y$. 
	
	\begin{lemma}[{\cite[Proposition 1.5.16]{hottaetal2008}}]\label{lem:Gamma_Y_trivial}
		We have
		\[i^\natural=\iHom_{i^{-1}\Aa}(\Aa_{Y\to X},i^{-1}\Gamma_Y(-)).\]
	\end{lemma}
	
	\begin{proof}
		Pick any local section $\varphi\in\iHom_{i^{-1}\Aa} (\Aa_{Y\to X},i^{-1}\Nn)$. We wish to prove that $\varphi$
		factors through $i^{-1}\Gamma_Y \Nn$. Since the structure homomorphism $i^\sharp:i^{-1}\OO_X\to \OO_Y$ is epic, $\Aa_{Y\to X}$ is generated by $1\otimes 1$ as a right $i^{-1}\Aa$-module. Therefore it suffices to check $\varphi(1\otimes 1)\in i^{-1}\Gamma_Y \Nn$.
		
		Define $a:\Nn\to \iHom_{\OO_X}(\Ii,\Nn)$ by the action of $\Ii\subset\OO_X$ on $\Nn$. Write $V\subset X$ for the complementary open subscheme to $Y$. Since $\Ii|_V=\OO_V$, we have $\Ker a\subset\Gamma_Y\Nn$. Apply $i^{-1}$ to get $\Ker i^{-1}a\subset i^{-1}\Gamma_Y\Nn$. Since $\Ii$ is locally finitely presented as an $\OO_X$-module around $Y$ by the proof of Theorem \ref{thm:relativecoordinate}, the canonical map
		$b:i^{-1}\iHom_{\OO_X}(\Ii,\Nn)\to \iHom_{i^{-1}\OO_X}(i^{-1}\Ii,i^{-1}\Nn)$
		is an isomorphism. We thus get
		$\Ker b\circ i^{-1}a=\Ker i^{-1}a\subset i^{-1}\Gamma_Y\Nn$.
		We now deduce
		\[\varphi(1\otimes 1)\in \Ker b\circ i^{-1}a\subset i^{-1}\Gamma_Y\Nn\]
		since $i^{-1}\Ii$ annihilates $\varphi(1\otimes 1)$. In fact, we have
		$\varphi(1\otimes 1)r=\varphi(1\otimes r)
		=\varphi(r\otimes 1)=0$ for $r\in i^{-1}\Ii$.
		
	\end{proof}
	
	\begin{proposition}[{\cite[Propositions 1.5.25 (i), 1.5.14 (ii)]{hottaetal2008}}]\label{prop:basisofinatural}
		\begin{enumerate}
			\item The functor $i^\natural$ is right adjoint to $i_+$.
			\item Let $U$ be an open subscheme of $X$. Write $i_U:Y\cap U\hookrightarrow U$ for the attached closed immersion. Then there is a natural isomorphism
			$(-)|_{U\cap Y}\circ i^{\natural}\cong i_U^\natural((-)|_U)$.
			\item The functor $i^\natural$ respects quasi-coherent sheaves.
			\item The functor $i^\natural$ respects $\ZZ$-flat sheaves.
		\end{enumerate}
		
	\end{proposition}
	
	\begin{proof}
		First, we prove (1). Let $\Mm\in \Mod_{\mathrm{r}}(i^\cdot\Aa)$ and $\mathcal N\in \Mod_{\mathrm{r}}(\Aa)$. Then we have an adjunction
		\[\Hom_{\Aa}(i_+\Mm,\Nn)
		\cong\Hom_{i^{-1}\Aa}(\Mm\otimes_{i^\cdot\Aa} \Aa_{Y\to X},
		i^{-1}\Gamma_Y\Nn)\]
		by \cite[Chapter II, Proof of Proposition 6.6]{iversen1986} (cf.~\cite[Lemma C.1.13]{hottaetal2008}). Pass to the adjunction of $(-\otimes_{i^\cdot\Aa} \Aa_{Y\to X},\iHom_{i^{-1}\Aa} (\Aa_{Y\to X},-))$ and use Lemma \ref{lem:Gamma_Y_trivial} to get
		\[\begin{split}
			\Hom_{i^{-1}\Aa}(\Mm\otimes_{i^\cdot\Aa} \Aa_{Y\to X},
			i^{-1}\Gamma_Y\Nn)
			&\cong \Hom_{i^\cdot\Aa}(\Mm,
			\iHom_{i^{-1}\Aa} (\Aa_{Y\to X},i^{-1}\Gamma_Y\Nn))\\
			&=\Hom_{i^\cdot\Aa}(\Mm,
			i^\natural\Nn).
		\end{split}\]
		This shows (1).
		
		Part (2) is evident by definition. For (3) and (4), we may simplify $i^\natural$ to $\iHom_{i^{-1}\OO_Y}(\OO_X,i^{-1}(-))$ since we are only interested in the structure of the underlying $\OO_X$-module. Then (4) is clear. It remains to prove (3). In view of (2), we may replace $X$ and $Y$ by open subschemes $U$ and $Y\cap U$ in Theorem \ref{thm:relativecoordinate}. In particular, $Y$ is defined by the ideal generated by $x_{r+1},\ldots,x_n$. Then $i^\natural$ can be identified with $L_{n-r}i^\ast$ by using the Koszul complex (see above \cite[Proposition 1.5.14]{hottaetal2008}). The assertion now follows from \cite[II. Proposition 4.4]{hartshorneresidue}.
	\end{proof}
	
	We now see that Kashiwara's theorem on the essential image of $i_+$ fails even if we restrict ourselves to $\ZZ$-flat sheaves:
	
	\begin{example}\label{ex:counterexofessimg}
		Put $S=\Spec\ZZ$ and $X=\bfA^1_S$. Define a closed reduced subsceheme $i:Y\hookrightarrow X$ by $0\in\bfA^1_S(S)$. Write $D(1)=\Gamma(X,\D_{X/S})$. Define a $\ZZ$-flat right $D(1)$-module $M'$ by
		\[\begin{array}{ccc}
			M'=\oplus_{n\geq 0} \ZZ\frac{\partial^n}{n!},
			&\frac{\partial^n}{n!}\cdot x=\frac{\partial^{n-1}}{(n-1)!},
			&\frac{\partial^n}{n!}\cdot \partial=
			(n+1)\frac{\partial^{n+1}}{(n+1)!}.
		\end{array}\]
		It is easy to show that $M'$ is not finitely generated as a $D(1)$-module.
		
		Let $\Mm'$ be the corresponding $\D_{X/S}$-module. Then $\Mm'$ is clearly supported in $Y$. If $\Mm'$ lies in the essential image of $i_+$, the counit $i_+i^\natural\Mm'\to\Mm'$ is an isomorphism. However, it is easy to show that $i^\natural\Mm'\cong\OO_S$. Hence $M\coloneqq \Gamma(X,i_+i^\natural\Mm')$ is given by
		\[\begin{array}{ccc}
			M=\oplus_{n\geq 0} \ZZ\partial^n,
			&\partial^n\cdot x=n\partial^{n-1},
			&\partial^n\cdot \partial=\partial^{n+1}.
		\end{array}\]
		Since $M$ is finitely generated, $M$ and $M'$ cannot be isomorphic to each other.
	\end{example}

	\begin{definition}
		Let $\Mm$ be a right $\OO_X$-module. Then define their subsheaves $\Gamma_{\left[Y,p\right]}\Mm$ ($p\in\NN$) and $\Gamma_{\left[Y\right]}\Mm$ by
		$\Gamma_{\left[Y,p\right]}\Mm=\{v\in\Mm:~v\Ii^p=0\}$ and
		$\Gamma_{\left[Y\right]}\Mm=
		\varinjlim_p \Gamma_{\left[Y,p\right]}\Mm$.
	\end{definition}
	
	It is easy to show that if $\Mm$ is a right $\Aa$-module, then $\Gamma_{\left[Y,p\right]}\Mm$ ($p\in\NN$) and $\Gamma_{\left[Y\right]}\Mm$ are $\Aa$-submodules of $\Mm$. 
	
	\begin{theorem}\label{thm:kashiwaraequiv}
		\begin{enumerate}
			\item For $\mathcal M\in \Mod_{{\mathrm{r}}}(i^\cdot\Aa)$, $\Gamma_{\left[Y\right]}i_+\Mm=i_+\Mm$.
			\item For $\Nn\in \Mod_{{\mathrm{r}}}^{\ZZ\textrm{-}\mathrm{flat}}(\Aa)$, the counit $u:i_+ i^\natural\Nn\to \Nn$ is a monomorphism into $\Gamma_{\left[Y\right]}\Nn$.
			\item Suppose that $X$ is a $\QQ$-scheme. Then we have a canonical isomorphism $i_+\circ i^\natural\cong\Gamma_{\left[Y\right]}$ on $\Mod_{{\mathrm{r}}}(i^\cdot\Aa)$.
		\end{enumerate}
	\end{theorem}
	
	\begin{proof}
		Iterated use of a computation similar to \eqref{eq:keycomputationtokashiwarathm} shows (1).
		
		For (2), let $\Nn\in \Mod_{{\mathrm{r}}}^{\ZZ\textrm{-}\mathrm{flat}}(\Aa)$. The counit $u$ is valued in $\Gamma_{\left[Y\right]}\Nn$ by (1). We next prove that $u$ is monic. Let
		$\ev:i^\natural\Nn\otimes_{i^\cdot\Aa}\Aa_{Y\to X}
		\to i^{-1}\Gamma_Y\Nn$
		be the evaluation map (recall Lemma \ref{lem:Gamma_Y_trivial}). Then $u$ is equal to the composite map
		$i_+i^\natural \Nn\overset{i_\ast u}{\to} i_\ast i^{-1}\Gamma_Y\Nn\cong
		\Gamma_Y\Nn\subset \Nn$.
		Therefore it will suffice to show that $\ev$ is monic.
		
		We may replace $X$ and $Y$ by $U$ and $U\cap Y$ in Theorem \ref{thm:relativecoordinate} to assume that $X$ admits a relative coordinate $((x_1,x_2,\ldots,x_r),(x_{r+1},\ldots,x_n))$. We may prove that $\ev_p$ is injective at each $p\in Y$. Recall that we have a decomposition
		\[(\iHom_{i^{-1}\Aa}(\Aa_{Y\to X},i^{-1}\Nn)\otimes_{i^\cdot\Aa}\Aa_{Y\to X})_p
		=\oplus_{I=(i_{r+1},\ldots,i_n)\in\NN^{n-r}}
		\iHom_{i^{-1}\Aa}(\Aa_{Y\to X},i^{-1}\Nn)\otimes 1\otimes
		\partial^{i_{r+1}}_{r+1}\cdots\partial^{i_{q-1}}_{q-1}.\]
		For
		$\Phi=\sum\varphi_I\otimes \partial^{i_{r+1}}_{r+1}\cdots\partial^{i_n}_n
		\in (\iHom_{i^{-1}\Aa}(\Aa_{Y\to X},i^{-1}\Nn)\otimes_{i^\cdot\Aa}\Aa_{Y\to X})_p$,
		set
		\[q(\Phi)\coloneqq \max\{p\in\{r+1,\ldots,n\}:~
		\exists I=(i_{r+1},i_{r+2},\ldots,i_n)\
		\mathrm{s.t.}\ \varphi_I\neq 0,\ i_p\neq 0\}\cup\{r\}
		\geq r.\]
		We prove by induction on $q(\Phi)$ that if $\ev_p(\Phi)=0$ then $\Phi=0$. If $q(\Phi)=r$ then we have an equality
		$\Phi=\varphi_{(0,0,\ldots,0)}\otimes 1 \otimes 1$.
		Suppose that $\ev_p(\Phi)=0$ then $\varphi_{(0,0,\ldots,0)}(1\otimes 1)=0$. Since $\Aa_{Y\to X}$ is generated by $1\otimes 1$ as a right $i^{-1}\Aa$-module, $\varphi_{(0,0,\ldots,0)}=0$. This shows $\Phi=0$. Suppose $q(\Phi)\leq q$ with $q\geq r+1$. Set
		\[j=\max\{i_{q}\in\NN:~
		\exists I=(i_{r+1},i_{r+2},\ldots,i_n)\
		\mathrm{s.t.}\ \varphi_I\neq 0\}\cup\{0\}.\]
		Suppose that $j=0$ then $q(\Phi)\leq q-1$ by definition of $j$. Hence the induction hypothesis shows $\Phi=0$. Suppose $j\geq 1$. Iterated use of \eqref{eq:keycomputationtokashiwarathm} implies
		\[\Phi x_{q}^j
		=\sum_{\overset{I=(i_{r+1},i_{r+2},\ldots,i_{q},0,0,\ldots,0)}{i_{q}=j}}
		j! \varphi_I\otimes \partial^{i_{r+1}}_{r+1}\cdots\partial^{i_{q-1}}_{q-1}.\]
		Since $0=\ev(\Phi)x^j_{q}=\ev(\Phi x^j_{q})$,
		the induction hypothesis implies $j!\varphi_I=0$ for $I=(i_{r+1},i_{r+2},\ldots,i_n)$ with $i_{q}=j$. Since $i^\natural \Nn$ is $\ZZ$-flat, $\varphi_I=0$. This contradicts to the definition of $j$. This completes the proof of (2).
		
		For (3), assume that $X$ is a $\QQ$-scheme. We may regard $u$ as a homomorphism into $\Gamma_{\left[Y\right]}\Nn$ by (2). We prove that $u$ is an epimorphism. Write $V\subset X$ for the complementary open subscheme to $Y$. Since $\Ii|_V=\OO_V$, $(\Gamma_{\left[Y\right]}\Nn)|_V=0$. In particular, $u|_V$ is an epimorphism. We next see that $u$ is epic around $Y$. We may again assume that $X$ admits a relative coordinate. It will suffice to prove that every local section $v\in\Gamma_{\left[Y,p\right]}\Nn$ lies in the image of $u$ by induction on $p$ (see the stalks). If $p=1$, $\Ii$ annihilates $v$. Hence $v$ determines a local section of $\iHom_{\OO_X}(i_\ast\OO_Y,\Nn)$. Apply $i^{-1}$ to get $\varphi\in \iHom_{i^{-1}\OO_X}(\OO_Y,i^{-1}\Nn)\cong \iHom_{i^{-1}\Aa}(\Aa_{Y\to X},i^{-1}\Nn)$ through the canonical homomorphism $i^{-1}\iHom_{\OO_X}(i_\ast\OO_Y,\Nn)\to \iHom_{i^{-1}\OO_X}(\OO_Y,i^{-1}\Nn)$. We have $u(\varphi\otimes 1\otimes 1)=v$ by definitions.
		
		Let $p\geq 2$. Write $\Image u$ for the image of $u$. Set $E=\sum_{l=r+1}^n x_l\partial_l$. Then it is straightforward that
		$v(E-p)x_q=vx_q(E-(p-1))$
		for $r+1\leq q\leq n$. Iterated use of this equality shows $v(E-p)\Ii^{p-1}=0$. The induction hypothesis implies that $v(E-p)\in\Image u$. We also have $vx_q\Ii^{n-1}=0$ for $r+1\leq q\leq n$. The induction hypothesis implies $vx_q\in\Image u$. Since $\Image u$ is an $\Aa$-submodule of $\Nn$, $vE\in\Image u$. Therefore
		$v=\frac{1}{p}E-\frac{1}{p}v(E-p)\in\Image u$.
		This completes the proof.
	\end{proof}

	\begin{corollary}\label{cor:kashiwara_ess_img}
		Suppose that $X$ is a $\QQ$-scheme.
		\begin{enumerate}
			\item A right $\Aa$-module $\Mm$ lies in the essential image of $i_+$ if and only if $\Gamma_{\left[Y\right]}\Mm=\Mm$.
			\item Suppose that $X$ is locally Noetherian. Then $\Mm\in\Mod_{{\mathrm{r}},\qc}(\Aa)$ lies in the essential image of $i_+$ if and only if $\Mm$ is supported in $Y$.
		\end{enumerate}
	\end{corollary}

	\subsubsection{Derived direct image for projections}\label{sec:derivedpushforp2}
	
	In this section, we aim to study the derived direct image functor for projections. First, we state
	
	\begin{theorem}[Spencer and de Rham resolutions]\label{thm:spencerresolution}
		Let $X$ be smooth and of constant relative dimension $n$ over $S$. Then we have two locally free resolutions of left and right $\D_X$-modules
		\begin{align}
			&0\to\D_X\otimes_{\OO_X}\bigwedge^{n+1}\Theta_{X/S}\to\dots\to\D_X\otimes_{\OO_X}\bigwedge^0\Theta_{X/S}\to\OO_X\to 0,\label{eq:Spencercomplex}\\
			&0\to\Omega_{X/S}^0\otimes_{\OO_X}\D_X\to\dots\to\Omega_{X/S}^n\otimes_{\OO_X}\D_X\to\omega_{X/S}\to 0,\label{eq:deRhamcomplex}
		\end{align}
		where
		$\D_X\otimes_{\OO_X}\bigwedge^0\Theta_{X/S}\to\OO_X$
		is the canonical action map (cf.\ Example \ref{ex:O_XasDmodule}) and the map
		$\Omega_{X/S}^n\otimes_{\OO_X}\D_X\to\omega_{X/S}$
		is given by the minus of the Lie derivative (see Corollary \ref{cor:omegaisrightDmod}).
		In the other degrees, the differentials are given by the formulas
		\begin{align*}
			d:\;&\D_X\otimes_{\OO_X}\bigwedge^q\Theta_{X/S}
			\to\D_X\otimes_{\OO_X}\bigwedge^{q-1}\Theta_{X/S},\\
			&P\otimes e_1\wedge\dots\wedge e_q\;\mapsto\\
			&\sum_{i=1}^q (-1)^{i+1} Pe_i\otimes e_1\wedge\dots\wedge\widehat{e_i}\wedge\dots\wedge e_q\\
			&+\sum_{i<j} (-1)^{i+j} P\otimes[e_i,e_j]\otimes e_1\wedge\dots\wedge\widehat{e_i}\wedge\dots\wedge\widehat{e_j}\wedge\dots\wedge e_q,\\[1em]
			d:\;
			&\Omega_{X/S}^q\otimes_{\OO_X}\D_X\to\Omega_{X/S}^{q+1}\otimes_{\OO_X}\D_X,\\
			&\omega\otimes P\;\mapsto d\omega\otimes P\;+\;\sum_{i=1}^q \partial_i^\vee\wedge\omega\otimes \partial_i P,
		\end{align*}
		where $\{\partial_1,\ldots,\partial_n\}$ is a local free basis of $\Theta_{X/S}$ and $\{\partial^\vee_1,\ldots,\partial^\vee_n\}$ is its dual local basis in $\Omega^1_{X/S}$.
		The actions of $\D_X$ are given by left and right multiplications respectively.
	\end{theorem}
	
	The sequence \eqref{eq:Spencercomplex} (resp.~\eqref{eq:deRhamcomplex}) is called the Spencer (resp.~de Rham) resolution.

	\begin{proof}
		One can readily apply the arguments of \cite[Section 4]{rinehart1963} or \cite[Appendix A, Proposition (A.19)]{macarro2015} to our sheafified setting to get the Spencer resolution \eqref{eq:Spencercomplex}.
		
		We obtain the second resolution from the first one by switching the side of actions via Example \ref{ex:linebdltwistofmodule}. In fact, note that $\D_{X/S}^{\op}\cong \omega_{X/S}\otimes_{\OO_X}\D_{X/S}\otimes_{\OO_X} \omega_{X/S}^\vee$ by Corollary \ref{cor:oppositecocyle}. See also \cite[Appendix A, Proposition (A.24)]{macarro2015}.
	\end{proof}
	
	\begin{remark}
		We can derive \eqref{eq:deRhamcomplex} from the classical de Rham complex of \cite[section 16.6]{ega44} as follows. Let $q\in\NN$. The identity section of $\iHom_{\OO_X}(\Omega^1_{X/S},\Omega^1_{X/S})\cong \Omega^1_{X/S}\otimes_{\OO_X}\Theta_{X/S}$ gives rise to a map
		\[\begin{split}
			\Omega^q_{X/S}\to\Omega^q_{X/S}\otimes_{\OO_X}
			\Omega^1_{X/S}\otimes_{\OO_X}\Theta_{X/S}
			&\overset{C_{\Omega^q_{X/S},	\Omega^1_{X/S}}}{\cong}
			\Omega^1_{X/S}\otimes_{\OO_X}
			\Omega^q_{X/S}\otimes_{\OO_X}\Theta_{X/S}\\
			&\to \Omega^{q+1}_{X/S}\otimes_{\OO_X}\Theta_{X/S}\\
			&\subset \Omega^{q+1}_{X/S}\otimes_{\OO_X}\D_{X/S},
		\end{split}
		\]	
		which is locally expressed as $\omega\mapsto \sum_{i=1}^q \partial_i^\vee\wedge\omega\otimes \partial_i$. We also have
		\[d(-)\otimes 1:\Omega^q_{X/S}\overset{d}{\to} \Omega^{q+1}_{X/S}\hookrightarrow
		\Omega^{q+1}_{X/S}\otimes_{\OO_X}\D_{X/S}.\]
		It is easy to show that the sum of these two maps is right $\OO_X$-linear. Hence we can extend it $\D_{X/S}$-linearly to get the differential map of \eqref{eq:deRhamcomplex}. 
	\end{remark}
	
	Let $X,Y$ be smooth $S$-schemes. Let $\pr_1:X\times_S Y\to X$ and $\pr_2:X\times_SY \to Y$ be the canonical projections. Let $\rank \pr^\ast_1\Theta_{X/S}$ denote the function on $X\times_S Y$ assigning the free rank of $\pr^\ast_1\Theta_{X/S}$ at each $q\in X\times_S Y$.
	
	\begin{theorem}\label{thm:directimageforprojections}
		\begin{enumerate}
			\item For a complex $\Mm^\bullet$ of right $\pr_2^\cdot\Aa$-module, the tensor product
			\[\Sp\Mm^\bullet\coloneqq 
			\Mm^\bullet\otimes_{\pr^{-1}_1\OO_X}\bigwedge^\bullet\Theta_{X/S}\]
			is a graded right $\pr_2^{-1}\Aa$-module for
			$(m\otimes \partial_1\wedge \cdots\wedge \partial_p)a
			=ma\otimes \partial_1\wedge \cdots\wedge \partial_p$.
			Here $a\in \pr^{-1}_2\Aa$ acts on $\Mm^\bullet$ via the canonical map $\pr^{-1}_2\Aa\to \pr^\cdot_2\Aa$ in Proposition \ref{prop:algebroidtdoprojection}. Moreover, $\Sp\Mm^\bullet$ is a cochain complex of $\pr_2^{-1}\Aa$-modules for the well-define map
			\begin{flalign*}
				&m\otimes \partial_1\wedge \cdots\wedge \partial_p\\
				&\mapsto dm\otimes \partial_1\wedge \cdots\wedge \partial_p\\
				&+\sum_{i=1}^p (-1)^{i+q+1} m \partial_i\otimes \partial_1\wedge\dots\wedge\widehat{\partial_i}\wedge\dots\wedge \partial_p\;\\
				&+\sum_{i<j} (-1)^{i+j+q} m[\partial_i,\partial_j]\otimes \partial_1\wedge\dots\wedge\widehat{\partial_i}\wedge\dots\wedge\widehat{\partial_j}\wedge\dots\wedge \partial_p,
			\end{flalign*}
			where $m\otimes \partial_1\wedge \cdots\wedge \partial_p\in \Mm^q\otimes_{\pr^{-1}_1\OO_X}\bigwedge^p\Theta_{X/S}$.
			\item We have $(\pr_2)_+\simeq R(\pr_2)_\ast\circ\Sp$.
		\end{enumerate}
		Suppose that $\pr_2$ is concentrated\footnote{If $x$ is concentrated, so is $\pr_2$.}.
		\begin{enumerate}
			\setcounter{enumi}{2}
			\item The functor $(\pr_2)_+$ restricts to $D_{{\mathrm{r}},\qc}(\pr^\cdot_2 \Aa)\to D_{{\mathrm{r}},\qc}(\Aa)$.
			\item If $Y$ is quasi-compact, $(\pr_2)_+$ is bounded on $D_{{\mathrm{r}},\qc}(\pr^\cdot_2\Aa)$ in the sense of \cite[Definition (1.11.1)]{lipman2009}.
		\end{enumerate}
	\end{theorem}
	
	\begin{proof}
		Write $z=x\circ \pr_1$. Apply $\pr^\cdot_2\Aa\otimes_{\pr^{-1}_1\D_X}$ to \eqref{eq:Spencercomplex} to get a locally free resolution
		\begin{equation}
			\pr^\cdot_2\Aa\otimes_{\pr^{-1}_1\OO_X}\pr^{-1}_1\bigwedge^{\bullet}\Theta_{X/S}
			\to \Aa_{X\times_S Y\to Y}\to 0
			\label{eq:nonabeliankoszulresolution}
		\end{equation}
		of $\pr^\cdot_2\Aa\otimes_{\pr^{-1}_1\D_X}\OO_X\cong \Aa_{X\times_S Y\to Y}$ as a left $\pr^\cdot_2\Aa$-module by Proposition \ref{prop:algebroidtdoprojection}. Moreover, it is an exact sequence of $(\pr^\cdot_2\Aa,\pr_2^{-1}\Aa)$-bimodules for multiplying $\pr_2^{-1}\Aa$ to $\pr^\cdot_2\Aa$ from the right side. Let $\Mm^\bullet$ be as in (1). Apply $\Mm^\bullet\otimes^L_{\pr^\cdot_2\Aa}(-)$ to this resolution to get
		\[\Mm^\bullet\otimes^L_{\pr^\cdot_2\Aa}\Aa_{X\times_S Y\to Y} \simeq \Mm^\bullet\otimes_{\pr^\cdot_2\Aa}
		\left(\pr^\cdot_2\Aa\otimes_{\pr^{-1}_1\OO_X}\pr^{-1}_1\bigwedge^{\bullet}\Theta_{X/S}\right),\]
		which in turn simplifies to $\Sp\Mm^\bullet$. We now deduce (1) and (2) by unwinding the definitions. 
		
		We next prove (3) and (4). Notice that we may assume $Y$ quasi-compact. In fact, (3) is local in $Y$. Recall also that we assumed $Y$ quasi-compact in (4). Since $\pr_2$ is concentrated, $X\times_S Y$ is quasi-compact. Since $\pr^\ast_1\Theta_{X/S}$ is locally free of finite rank, $\rank \pr^\ast_1\Theta_{X/S}$ is bounded.
		
		Identify the resolution \eqref{eq:nonabeliankoszulresolution} with 
		\begin{equation}
			\pr^\cdot_2\Aa\otimes_{\OO_{X\times_S Y}}\pr^\ast_1\bigwedge^{\bullet}\Theta_{X/S}
			\to \Aa_{X\times_S Y\to Y}\to 0
			\label{eq:Spencerresolution}
		\end{equation}
		by extending the right action of $\pr^{-1}_1\OO_{X}$ on $\pr^\cdot_2\Aa$ to that of $\OO_{X\times_S Y}$. Since $\rank \pr^\ast_1\Theta_{X/S}$ is bounded, \eqref{eq:Spencerresolution} is of finite length. Hence we can pass to stupid truncations of the complex $\pr^\cdot_2\Aa\otimes_{\OO_{X\times_S Y}}\pr^\ast_1\bigwedge^{\bullet}\Theta_{X/S}$ to reduce (3) and (4) to the following assertions for $n\in\ZZ$:
		\begin{enumerate}
			\item[(3)'] $R(\pr_2)_\ast\left(-\otimes_{\OO_{X\times_S Y}}\pr^\ast_1\bigwedge^n\Theta_{X/S}\right)$ determines a functor $D_{{\mathrm{r}},\qc}(\pr^\cdot_2 \Aa)\to D_{{\mathrm{r}},\qc}(\Aa)$.
			\item[(4)'] Suppose that $Y$ is quasi-compact. Then $R(\pr_2)_\ast\left(-\otimes_{\OO_{X\times_S Y}}\pr^\ast_1\bigwedge^n\Theta_{X/S}\right)$ is bounded on $D_{{\mathrm{r}},\qc}(\pr^\cdot_2\Aa)$.
		\end{enumerate}
		These are consequences of \cite[Proposition (3.9.2)]{lipman2009}.
	\end{proof}
	
	\subsubsection{$S$-base change theorem}\label{sec:flatbcthm}
	In this section, we establish the (flat) base change theorem of the (derived) direct image in $S$ for closed immersions and projections.
	
	\begin{construction}\label{cons:sbcmap}
		Let $s:S'\to S$ be a morphism of schemes, $f:X\to Y$ be a morphism of smooth $S$-schemes, and $\Aa$ be a tdo on $Y$. Let $\Mm$ be a right $f^\cdot\Aa$-module. Then we have an isomorphism
		$s^{-1}_X(\Mm\otimes_{f^\cdot\Aa}\Aa_{X\to Y})
		\cong
		s^{-1}_X\Mm\otimes_{s^{-1}_Xf^\cdot\Aa}s^{-1}_X\Aa_{X\to Y}$
		of right $s^{-1}_Xf^{-1}\Aa$-modules since $s^{-1}_X$ is symmetric monoidal. By Proposition \ref{prop:basechangevspullbackformodules/sch} (1), we obtain
		a homomorphism
		\[\begin{split}
			s^{-1}_X(\Mm\otimes_{f^\cdot\Aa}\Aa_{X\to Y})
			&\to 
			s^{-1}_X\Mm\otimes_{s^{-1}_Xf^\cdot\Aa}s^\ast_X\Aa_{X\to Y}\\
			&\cong s^\ast_X\Mm\otimes_{s^\ast_Xf^\cdot\Aa}
			s^\ast_X\Aa_{X\to Y}\\
			&\cong s^\ast_X\Mm\otimes_{s^\ast_Xf^\cdot\Aa}
			(s^\ast_Y\Aa)_{X'\to Y'}\\
			&\cong s^\ast_X\Mm\otimes_{(f')^\cdot s^\ast_Y\Aa}
			(s^\ast_Y\Aa)_{X'\to Y'}
		\end{split}\]
		of right $s^{-1}_Xf^{-1}\Aa$-modules. Moreover, the target
		$s^\ast_X\Mm\otimes_{(f')^\cdot s^\ast_Y\Aa}
		(s^\ast_Y\Aa)_{X'\to Y'}$ is equipped with the structure of a right $(f')^{-1}s^\ast_Y\Aa$-module. Apply $f_\ast$ to get a right $s^{-1}_Y\Aa$-linear map
		\[f'_\ast s^{-1}_X(\Mm\otimes_{f^\cdot\Aa}\Aa_{X\to Y})
		\to f'_\ast (s^\ast_X\Mm\otimes_{(f')^\cdot s^\ast_Y\Aa}
		(s^\ast_Y\Aa)_{X'\to Y'}).\]
		By the base change construction, we obtain
		\[s^{-1}_Y f_\ast(\Mm\otimes_{f^\cdot\Aa}\Aa_{X\to Y})
		\to f'_\ast (s^\ast_X\Mm\otimes_{(f')^\cdot s^\ast_Y\Aa}
		(s^\ast_Y\Aa)_{X'\to Y'}).\]
		Take the scalar extension to get
		\[s^\ast_Y f_\ast(\Mm\otimes_{f^\cdot\Aa}\Aa_{X\to Y})
		\to f'_\ast (s^\ast_X\Mm\otimes_{(f')^\cdot s^\ast_Y\Aa}
		(s^\ast_Y\Aa)_{X'\to Y'}).\]
	\end{construction}
	
	One can obtain the flat base change map in the derived setting in a similar way (omitted).
	
	\begin{theorem}\label{thm:Sbasechange1}
		Let $s:S'\to S$ be a morphism of schemes. Let $i:Y\to X$ be a closed immersion of smooth $S$-schemes, $i':Y'\to X'$ be the base change of $i$ to $S'$, and $\Aa$ be a tdo on $X$. Then the natural transformation $s^\ast_X i_+\to i'_+s^\ast_Y$ obtained by Construction \ref{cons:sbcmap} is an isomorphism on $\Mod_{{\mathrm{r}},\qc}(i^\cdot\Aa)$.
	\end{theorem}
	
	\begin{proof}
		Since the assertion is local in $X$, we may assume that $X$, and therefore $Y$ are affine. One can also pass to small colimits in $\Mod_{{\mathrm{r}},\qc}(i^\cdot\Aa)$ by Variant \ref{var:comparisonf_+f^+;imm}. Hence we may prove the isomorphism $s^\ast_Xi_+i^\cdot\Aa\cong i'_+s^\ast_Y i^\cdot\Aa$. This follows from Proposition \ref{prop:basechangevspullbackformodules/sch} (use the left module structure on $\Aa_{Y\to X}$).
	\end{proof}

	\begin{theorem}\label{thm:Sbasechange2}
		Let $X,Y$ be smooth $S$-schemes, and $s:S'\to S$ be a flat morphism, and $\Aa$ be a tdo on $Y$. Let $\pr_2:X\times_S Y\to Y$ denote the projection. Suppose that $\pr_2$ is concentrated. Then the canonical natural transformation
		$s^\ast_Y R(\pr_2)_+\to R(\pr'_2)_+s^\ast_{X\times_S Y}$ is an equivalence on $D_{{\mathrm{r}},\qc}(\pr_2^\cdot\Aa)$.
	\end{theorem}
	
	\begin{proof}
		Since the assertion is local in $Y$, we may assume that $Y$ is affine. In particular, $Y$ is quasi-compact. Then $\rank \pr^\ast_1\Theta_{X/S}$ and $\rank (\pr'_1)^\ast\Theta_{X'/S'}$ are bounded, where
		\[\begin{array}{cc}
			\pr_1:X\times_S Y\to X,&\pr'_1:X'\times_{S'} Y'\to X'
		\end{array}\]
		are the projections. In fact, use the isomorphism $(\pr'_1)^\ast\Theta_{X'/S'}\cong s^\ast_{X\times_S Y}\pr^\ast_1\Theta_{X/S}$
		for the boundedness of $\rank (\pr'_1)^\ast\Theta_{X'/S'}$.
		
		Let $\Mm^\bullet\in D_{\mathrm{r}}(f^\cdot\Aa)$. Recall that $s^\ast_Y (\pr_2)_+\Mm^\bullet$ and $(\pr'_2)_+s^\ast_{X\times_S Y}$ are computed by
		\[s^\ast_Y R(\pr_2)_+\Mm_\bullet
		\simeq s^\ast_Y R(\pr_2)_\ast\Sp\Mm^\bullet,\]
		\[(\pr'_2)_+s^\ast_{X\times_S Y}\Mm^\bullet
		\simeq R(\pr'_2)_\ast \Sp(s^\ast_{X\times_S Y}\Mm^\bullet)\]
		respectively. Since $\rank \pr^\ast_1\Theta_{X/S}$ and $\rank (\pr'_1)^\ast\Theta_{X'/S'}$ are bounded, we can pass to the stupid truncations to reduce the assertion to proving the canonical homomorphism
		\[s^\ast_Y R(\pr_2)_\ast(\Mm^\bullet\otimes_{\OO_{X\times_S Y}}
		\pr^\ast_1\bigwedge^q\Theta_{X/S})
		\to R(\pr'_2)_\ast (s^\ast_{X\times_S Y}
		\Mm^\bullet\otimes_{\OO_{X'\times_{S'}Y'}}
		(\pr'_1)^\ast\bigwedge^q\Theta_{X'/S'})\]
		is an equivalence (cf.\ the proof of Theorem \ref{thm:directimageforprojections} (3) and (4)). This follows from \cite[Proposition (3.9.5)]{lipman2009}.
	\end{proof}
	
	\begin{remark}
		It is clear by definitions that similar assertions in section \ref{sec:derivedpushforp2} and section \ref{sec:flatbcthm} hold for left modules. Use the de Rham resolution instead of the Spencer resolution in the projection setting. In particular, $(\pr_2)_+$ may be expressed as $R(\pr_2)_\ast\circ \DR(\Mm^\bullet)$, where $\DR(\Mm^\bullet)=\pr_1^{-1}\Omega^\bullet_{X/S}\otimes_{\pr^{-1}_1\OO_X} \Mm^\bullet$. The differential and the action of $\pr^{-1}_2\Aa$ are given by
		\[\pr_1^{-1}\Omega^q_{X/S}\otimes_{\pr^{-1}_1\OO_X} \Mm^\bullet\ni\omega\otimes m\mapsto d\omega\otimes m 
		+\sum_{i=1}^q \partial^\vee_i\wedge \omega\otimes m
		+(-1)^q \omega\otimes dm,\]
		\[a(\omega\otimes m)=\omega\otimes am,\]
		where $\{\partial_1,\ldots,\partial_n\}$ are free local generators of $\Theta_{X/S}$ and $\{\partial^\vee_1,\ldots,\partial^\vee_n\}$ are the dual sections in $\Omega^1_{X/S}$.
	\end{remark}
	
	\subsubsection{Properties of $f_+$}
	
	The aim of this section is to see consequences of the case-by-case study of the direct image functor in former sections.

	\begin{theorem}\label{thm:qcohpreservation}
		Let $f:X\to Y$ be a concentrated morphism of smooth $S$-schemes, $\Bb$ be a tdo on $Y$. Then
		\begin{enumerate}
			\item The functor $f_+:D_{{\mathrm{r}}}(f^\cdot\Bb)\to D_{{\mathrm{r}}}(\Bb)$ restricts to $D_{{\mathrm{r}},\qc}(f^\cdot\Bb)\to D_{{\mathrm{r}},\qc}(\Bb)$.
			\item For a flat morphism $s:S'\to S$ of schemes, there is an equivalence $s^\ast\circ f_+\simeq f'_+\circ s^\ast_X$ on $D_{{\mathrm{r}},\qc}(f^\cdot\Bb)$.
		\end{enumerate}
	\end{theorem}

	\begin{proof}
		Since the assertion is local in $S$, we may assume $S$ affine. In this case, we may also assume $Y$ affine. In particular, $y$ is quasi-compact and separated. Since $f$ is concentrated, so is $x$ in this case. Notice also that the graph morphism of $f$ is a closed immersion.
		
		Combine Theorems \ref{thm:fpluscompositions1}, \ref{thm:directimageforprojections} and Lemma \ref{lem:i+cloimm} to get (1). Combine Theorems \ref{thm:fpluscompositions1}, \ref{thm:Sbasechange1}, \ref{thm:Sbasechange2} and Lemma \ref{lem:i+cloimm} for (2).
	\end{proof}
	
	We have to work carefully with the boundedness on $D_{{\mathrm{r}},\qc}(f^\cdot\Bb)$. In fact, the notion of local boundedness does not work for cohomologically quasi-coherent complexes since the functor of extension by zero does not respect quasi-coherent sheaves.
	
	\begin{theorem}
		Let $f:X\to Y$ be a morphism of smooth $S$-schemes, $\Bb$ be a tdo on $Y$. Suppose:
		\begin{enumerate}
			\item[(i)] $y$ is separated.
			\item[(ii)] The projection morphism $X\times_S Y\to Y$ is concentrated.
			\item[(iii)] $Y$ is quasi-compact. 
		\end{enumerate}
		Then $f_+$ is bounded on $D_{{\mathrm{r}},\qc}(f^\cdot\Bb)$.
	\end{theorem}
	
	\begin{proof}
		This follows from Theorems \ref{thm:fpluscompositions1}, \ref{thm:directimageforprojections} and Lemmas \ref{lem:i+cloimm}, \ref{lem:bimodlocfreeimmcase}.
	\end{proof}

	We next prove the composition law. Let $f:X\to Y$ and $g:Y\to Z$ be morphisms of smooth $S$-schemes, $\Aa$ be a tdo on $Z$. The key will be versions of the projection formula. Henceforth assume the following conditions:
	\begin{enumerate}
		\renewcommand{\labelenumi}{(\roman{enumi})}
		\item $f$ is concentrated.
		\item $Rf_\ast:D(x^{-1}\OO_S)\to D(y^{-1}\OO_S)$ is locally bounded.
	\end{enumerate}

	\begin{lemma}\label{lem:projformula}
		\begin{enumerate}
			\item For $\Nn^\bullet\in D^{\rm -}_{\mathrm{r}}(f^{-1}g^\cdot \Aa)$, the canonical natural transformation
			\[R f_\ast(\Nn^\bullet)\otimes^L_{g^\cdot \Aa} (-)\to R f_\ast(\Nn^\bullet\otimes^L_{f^{-1}g^\cdot \Aa} f^{-1}(-))\]
			is an equivalence on $D^{\rm -}_{\qc}(g^\cdot \Aa)$.
			\item The canonical natural transformation
			\[R f_\ast(-)\otimes^L_{g^\cdot \Aa} \Aa_{Y\to Z}\to R f_\ast(-\otimes^L_{f^{-1}g^\cdot \Aa} f^{-1}\Aa_{Y\to Z})\]
			is an equivalence on $D_{\mathrm{r}}(f^{-1}g^\cdot \Aa)$.
		\end{enumerate}
		
	\end{lemma}
	
	\begin{proof}
		We may assume $S$, $Y$, and $Z$ are affine.
		
		For (1), notice that the functors $R f_\ast(\Nn^\bullet)\otimes^L_{g^\cdot \Aa} (-)$ and $R f_\ast(\Nn^\bullet\otimes^L_{f^{-1}g^\cdot \Aa} f^{-1}(-))$ are bounded above by (ii). Hence in view of \cite[Lemma (1.11.3) and Complement (1.11.13.1)]{lipman2009}, we may prove that the canonical morphism $R f_\ast(\Nn^\bullet)\otimes^L_{\Bb} \Gg\to R f_\ast(\Nn^\bullet\otimes^L_{f^{-1}\Bb} f^{-1}\Gg)$ is an equivalence if $\Gg$ is a free $\Bb$-module. The assertion now follows from Lemma \ref{lem:partialomegaaccessibility}.
		
		For (2), we may ignore the right action on $\Aa_{Y\to Z}$. Then $\Aa_{Y\to Z}$ admits a locally free resolution of finite length as a left $g^\cdot \Aa$-module by Lemmas \ref{lem:acyclicityofbimodchain}, \ref{lem:bimodlocfreeimmcase} (1), and the resolution \eqref{eq:nonabeliankoszulresolution}. Hence the functors $R f_\ast(-)\otimes^L_{g^\cdot \Aa} \Aa_{Y\to Z}$ and $R f_\ast(-\otimes^L_{f^{-1}g^\cdot \Aa} f^{-1}\Aa_{Y\to Z})$ are bounded. We may therefore work on $\Mod_{\mathrm{r}}(f^{-1}g^\cdot \Aa)$ by \cite[Lemma (1.11.3)]{lipman2009}. Then the assertion follows from (1).
	\end{proof}
	
	We now deduce:
	\begin{theorem}\label{thm:compositionlaw}
		The natural transformation $g_+\circ f_+\to (g\circ f)_+$ in Construction \ref{cons:compositionlaw} is an equivalence.
	\end{theorem}
	
	The proof of Lemma \ref{lem:projformula} (2) also implies:
	
	\begin{theorem}
		Let $\Bb$ be a tdo on $Y$. Then the functor $f_+:D_{\mathrm{r}}(f^\cdot \Bb)\to D_{\mathrm{r}}(\Bb)$ is locally bounded.
	\end{theorem}

	\subsection{$X$-Base Change Theorem}\label{sec:xbcthm}
	In this section, we aim to establish the classical base change theorem of twisted $\D$-modules over schemes. Due to the failure of Kashiwara's equivalence, we restrict ourselves to limited cases which are enough for applications to equivariant twisted $\D$-modules. In fact, we have the following counterexample:
	
	\begin{example}\label{ex:counterexofbcthm}
		Consider the pullback diagram
		\[\begin{tikzcd}
			S\ar[r, equal]\ar[d, equal]
			&S\ar[d, hook, "i"]\\
			S\ar[r, hook, "i"]&\bfA^1_S,
		\end{tikzcd}\]
		where $i$ is the closed immersion attached to $0\in\bfA^1_S(S)$. If $S=\Spec\CC$, the classical base change theorem asserts $L_1i^\ast\circ i_+\cong\id_{\Mod(\OO_S)}$. Note also that $i^\ast\circ i_+\not\cong\id_{\Mod(\OO_S)}$. On the other hand, if $S=\Spec\FF_p$ with $p$ a positive prime, then $L_1i^\ast\circ i_+\not\cong\id_{\Mod(\OO_S)}$ in general by a similar argument to section \ref{sec:kashiwara}. In fact, recall that we have a Koszul resolution
		\[0\to i^{-1}\OO_{\bfA^1_S}\to i^{-1}\OO_{\bfA^1_S}\overset{i^\sharp}{\to} \OO_S\to 0,\]
		where the arrow $i^{-1}\OO_{\bfA^1_S}\to i^{-1}\OO_{\bfA^1_S}$ is defined by multiplying the coordinate function $x$ of $\bfA^1_S$. Use $L i^\ast=\OO_S\otimes^L_{i^{-1}\OO_{\bfA^1_S}} i^{-1}(-)$ to identify 
		$L_1 i^\ast(-)$ with the kernel of the endomorphism of $i^{-1}(-)$ defined by the multiplication with $x$. Put $\Aa=\D_{\bfA^1_S}$. Then the left $D(1)$-module $M\coloneqq\Gamma(X,i_+\OO_S)$ is given by
		\[\begin{array}{ccc}
			M=\oplus_{n\geq 0} \FF_p\partial^n,&x\cdot \partial^n=-n\partial^{n-1},
			&\partial\cdot\partial^n=\partial^{n+1}
		\end{array}\]
		(see Example \ref{ex:counterexofkashiwaraeq}). We thus get
		$L_1 i^\ast(i_+\OO_S)=\oplus_{n\geq 0} \OO_S\partial^{pn}\neq \OO_S$.
		We remark if necessary that $L_1i^\ast$ is left adjoint to $i_+$ by a similar argument to the proof of Proposition \ref{prop:basisofinatural} (3).
	\end{example}

	\begin{theorem}\label{thm:xbc}
		Suppose that we are given a Cartesian diagram of $S$-schemes
		\[\begin{tikzcd}
			\tilde{X}\ar[r, "\tilde{p}"]\ar[d, "\tilde{f}"']
			&\tilde{Y}\ar[d, "f"]\\
			X\ar[r, "p"]&Y.
		\end{tikzcd}\]
		Assume that the following conditions are satisfied:
		\begin{enumerate}
			\renewcommand{\labelenumi}{(\roman{enumi})}
			\item All of $X$, $Y$, $\tilde{X}$, and $\tilde{Y}$ are smooth over $S$.
			\item $p$ is isomorphic to a projection morphism in the category of morphisms of $S$-schemes.
			\item $f$ is concentrated.
		\end{enumerate}
		Let $\Aa$ be a tdo on $Y$. Then there is an equivalence
		$p^\ast\circ f_+\simeq \tilde{f}_+\circ \tilde{p}^\ast$
		between functors from $D_{\qc}(f^\cdot\Aa)$ to $D_{\qc}(p^\cdot \Aa)$.
	\end{theorem}
	
	\begin{remark}
		Although the shift of the pullback plays a crucial role in the classical base change theorem due to the relation of the two pullbacks $Li^\ast$ and $i^\natural$ for a closed immersion $i$ as in Example \ref{ex:counterexofbcthm}, we do not need to take care of shifts since the horizontal arrows are smooth in our setting.
	\end{remark}
	
	\begin{example}
		For a smooth $S$-scheme $X$ with an action $a:G\times_S X\to X$ of a smooth group $S$-scheme $G$, $a$ is isomorphic to the projection $\pr_2:G\times_S X\to X$. 
	\end{example}
	
	For the proof, observe that for an isomorphism $g:W\to Z$ of smooth $S$-schemes, $g^\ast$ is quasi-inverse to $g_+$. Therefore in view of the naturality of the pullback and direct image functors, we may assume that $f$ is a projection. In this way, Theorem \ref{thm:xbc} is reduced to:
	
	\begin{claim}\label{claim:xbc}
		Let $X,Y,Z$ be smooth $S$-schemes. Consider the pullback diagram
		\[\begin{tikzcd}
			X\times_S Z\ar[r, "\tilde{\pr}_1"]\ar[d, "\tilde{f}"']
			&X\ar[d, "f"]\\
			Y\times_S Z\ar[r, "\pr_1"]&Y,
		\end{tikzcd}\]
		where $\pr_1$ is the canonical projection. Suppose that $f$ is concentrated. Let $\Aa$ be a tdo on $Y$. Then there is an equivalence
		$\pr^\ast_1\circ f_+\simeq \tilde{f}_+\circ \tilde{\pr}^\ast_1$
		on $D_{\qc}(f^\cdot\Aa)$.
	\end{claim}
	
	To find the canonical natural transformation, observe
	
	\begin{lemma}
		There is a canonical $(\tilde{f}^{-1}\pr^\cdot_1\Aa,\tilde{\pr}^{-1}_1f^\cdot\Aa)$-bimodule isomorphism
		\begin{equation}
			\tilde{f}^{-1}\Aa_{Y\times_S Z\to Y}
			\otimes_{\tilde{f}^{-1}\pr^{-1}_1\Aa}
			\tilde{\pr}^{-1}_1\Aa_{Y\gets X}
			\cong (\pr^\cdot_1\Aa)_{X\times_S Z\gets Y\times_S Z}
			\otimes_{\tilde{\pr}^\cdot_1f^\cdot\Aa}(f^\cdot\Aa)_{X\times_S Z\to X}.
			\label{eq:bimodhom}
		\end{equation} 
	\end{lemma}
	
	\begin{proof}
		The left hand side of \eqref{eq:bimodhom} is computed as
		\begin{flalign*}
			&\tilde{f}^{-1}\Aa_{Y\times_S Z\to Y}
			\otimes_{\tilde{f}^{-1}\pr^{-1}_1\Aa}
			\tilde{\pr}^{-1}_1\Aa_{Y\gets X}\\
			&\cong \tilde{f}^{-1}\Aa_{Y\times_S Z\to Y}
			\otimes_{\tilde{f}^{-1}\pr^{-1}_1\Aa}
			\tilde{\pr}^{-1}_1(f^{-1}\Aa\otimes_{f^{-1}\OO_Y}
			f^{-1}\omega_{Y/S}^\vee\otimes_{f^{-1}\OO_Y}
			\omega_{X/S})\\
			&\cong\tilde{f}^{-1}\Aa_{Y\times_S Z\to Y}
			\otimes_{\tilde{\pr}^{-1}_1f^{-1}\OO_Y}
			\tilde{\pr}^{-1}_1f^{-1}\omega_{Y/S}^\vee\otimes_{\tilde{\pr}^{-1}_1f^{-1}\OO_Y}
			\tilde{\pr}^{-1}_1\omega_{X/S}\\
		\end{flalign*}
		by the canonical cancellations.
		
		Notice that the left action of $\tilde{\pr}^\cdot_1f^\cdot\Aa$ on $(f^\cdot\Aa)_{X\times_S Z\to X}$ restricts to that of \[\tilde{\pr}^{-1}_1f^\cdot\Aa\otimes_{\tilde{\pr}^{-1}_1x^{-1}\OO_S} \tilde{\pr}^{-1}_2 \D_Z\]
		on $\tilde{\pr}^{-1}_1f^\cdot\Aa\otimes_{\tilde{\pr}^{-1}_1x^{-1}\OO_S} \tilde{\pr}^{-1}_2 \OO_Z$,
		where $\tilde{\pr}_2:X\times_SZ\to Z$ denotes the canonical projection (cf.\ Proposition \ref{prop:algebroidtdoprojection}). Hence we can cancel $\OO_{X\times_S Z}$ to get an identification
		\begin{flalign*}
			&(\pr^\cdot_1\Aa)_{X\times_S Z\gets Y\times_S Z}
			\otimes_{\tilde{\pr}^\cdot_1f^\cdot\Aa}(f^\cdot\Aa)_{X\times_S Z\to X}\\
			&\cong (\pr^\cdot_1\Aa)_{X\times_S Z\gets Y\times_S Z}
			\otimes_{\tilde{\pr}^{-1}_1f^\cdot\Aa\otimes_{\tilde{\pr}^{-1}_1x^{-1}\OO_S} \tilde{\pr}^{-1}_2 \D_Z}(\tilde{\pr}^{-1}_1f^\cdot\Aa\otimes_{\tilde{\pr}^{-1}_1x^{-1}\OO_S} \tilde{\pr}^{-1}_2 \OO_Z)\\
			&\cong (\pr^\cdot_1\Aa)_{X\times_S Z\gets Y\times_S Z}
			\otimes_{\tilde{\pr}^{-1}_2 \D_Z}\tilde{\pr}^{-1}_2 \OO_Z.
		\end{flalign*}
		We remark that the induced right action on $\tilde{\pr}^{-1}_1f^\cdot\Aa$ on $(\pr^\cdot_1\Aa)_{X\times_S Z\gets Y\times_S Z}
		\otimes_{\tilde{\pr}^{-1}_2 \D_Z}\tilde{\pr}^{-1}_2 \OO_Z$ arises from that on $(\pr^\cdot_1\Aa)_{X\times_S Z\gets Y\times_S Z}$. Similarly, one can cancel $\D_Z$:
		\[(\pr^\cdot_1\Aa)_{X\times_S Z\gets Y\times_S Z}
		\otimes_{\tilde{\pr}^{-1}_2 \D_Z}\tilde{\pr}^{-1}_2 \OO_Z
		\cong (\pr^\cdot_1\Aa)_{X\times_S Z\gets Y\times_S Z}
		\otimes_{\tilde{\pr}^\cdot_1f^\cdot\Aa}(f^\cdot\Aa)_{X\times_S Z\to X}.\]
		This completes the proof.
	\end{proof}
	
	\begin{construction}
		Let $\Mm^\bullet\in D(f^\cdot\Aa)$. Apply $R\tilde{f}_\ast(-\otimes_{\tilde{\pr}^{-1}_1f^\cdot\Aa} \tilde{\pr}^{-1}_1\Mm^\bullet)$ to \eqref{eq:bimodhom} to get
		\[R\tilde{f}_\ast
		(\tilde{f}^{-1}\Aa_{Y\times_S Z\simeq Y}\otimes_{\tilde{f}^{-1}\pr^{-1}_1\Aa}
		\tilde{\pr}^{-1}_1\Aa_{Y\gets X}
		\otimes_{\tilde{\pr}^{-1}_1f^\cdot\Aa} \tilde{\pr}^{-1}_1\Mm^\bullet)
		\to \tilde{f}_+(\tilde{\pr}^\ast_1\Mm^\bullet).\]
		Compose it with the base change and projection maps to obtain
		\[\begin{split}
			\pr^\ast f_+\Mm^\bullet
			&=\Aa_{Y\times_S Z\to Y}\otimes_{\pr^{-1}_1\Aa}
			\pr^{-1}_1Rf_\ast(\Aa_{Y\gets X}
			\otimes_{f^\cdot\Aa} \Mm^\bullet)\\
			&\to \Aa_{Y\times_S Z\to Y}\otimes_{\pr^{-1}_1\Aa}
			R\tilde{f}_\ast\tilde{\pr}^{-1}_1(\Aa_{Y\gets X}
			\otimes_{f^\cdot\Aa} \Mm^\bullet)\\
			&\to R\tilde{f}_\ast
			(\tilde{f}^{-1}\Aa_{Y\times_S Z\to Y}\otimes_{\tilde{f}^{-1}\pr^{-1}_1\Aa}
			\tilde{\pr}^{-1}_1\Aa_{Y\gets X}
			\otimes_{\tilde{\pr}^{-1}_1f^\cdot\Aa} \tilde{\pr}^{-1}_1\Mm^\bullet)\\
			&\simeq \tilde{f}_+(\tilde{\pr}^\ast_1\Mm^\bullet).
		\end{split}\]	
		In this way, we obtain a natural transformation
		\begin{equation}
			\pr^\ast_1\circ f_+\to \tilde{f}_+\circ \tilde{\pr}^\ast_1.
			\label{eq:xbcmap}
		\end{equation}
	\end{construction}

	We wish to prove that \eqref{eq:xbcmap} is an equivalence on $D_{\qc}(f^\cdot\Aa)$. We may assume $S$, $Y$, and $Z$ are affine. In particular, we may assume the struture morphisms $Y\to S$ and $Y\times_S Z\to S$ to be separated. In this case, we may also assume $f$ is either a projection or a closed immersion. The proof is therefore reduced to:
	
	\begin{claim}\label{claim:xbcthm}
		\begin{enumerate}
			\item Let $X,Y,Z$ be smooth $S$-schemes. Consider the pullback diagram
			\[\begin{tikzcd}
				X\times_S Z\ar[r, "\tilde{\pr}_1"]\ar[d, "\tilde{i}"']
				&X\ar[d, "i"]\\
				Y\times_S Z\ar[r, "\pr_1"]&Y,
			\end{tikzcd}\]
			where $\pr_1$ is the canonical projection. Suppose that 
			\begin{enumerate}
				\item[(i)] $i$ is affine,
				\item[(ii)] $Ri_\ast:D(i^{-1}\Aa)\to D(\Aa)$ is locally bounded,
				\item[(iii)] $R\tilde{i}_\ast:D(\tilde{i}^{-1}\pr_1^\cdot\Aa)\to D(\pr_1^\cdot\Aa)$ is locally bounded.
			\end{enumerate}
			Then $\pr^\ast_1\circ i_+\cong \tilde{i}_+\circ \tilde{\pr}^\ast_1$ on $\Mod_{\qc}(i^\cdot\Aa)$.
			\item Let $X,Y,Z$ be smooth $S$-schemes. Consider the pullback diagram
			\[\begin{tikzcd}
				X\times_S Y\times_S Z\ar[r, "\pr_{23}"]\ar[d, "\pr_{12}"']
				&Y\times_S Z\ar[d, "\pr_1"]\\
				X\times_S Y\ar[r, "\pr_2"]&Y,
			\end{tikzcd}\]
			whose arrows are defined by the canonical projections.
			If $\pr_1$ is concentrated, then $\pr^\ast_2\circ (\pr_1)_+\simeq (\pr_{12})_+\circ \pr_{23}^\ast$ on $D_{\qc}(\pr_1^\cdot\Aa)$.
		\end{enumerate}
	\end{claim}
	
	\begin{proof}[Proof of Claim \ref{claim:xbcthm} (1)]
		The assertion is local in $Y$ and $Z$ by \eqref{eq:localityoff+}. Hence we may pass to the ring setting to show that the map is just the combination of the isomorphism canceling $\D_Z$ and taking the canonical isomorphism $\omega_{X/S}\otimes_{\OO_X}\omega^\vee_{X/S}\cong{\mathcal O}_X$.
	\end{proof}

	\begin{proof}[Proof of Claim \ref{claim:xbcthm} (2)]
		We may assume $Y$ affine by \eqref{eq:localityoff+}. In particular, $Y\times_S Z$ is quasi-compact since $\pr_1$ is quasi-compact. Let $\pr_2:Y\times_S Z\to Z$ denote the canonical projection. Then $\rank \pr^\ast_2\Theta_{Z/S}$ is bounded. One can therefore reduce the assertion to \cite[Proposition (3.9.5)]{lipman2009} by passing to the stupid truncations.
	\end{proof}
	
	\subsection{Equivariant twisted $\D$-modules}
	
	Suppose that $X$ is equipped with an action
	\[a:G\times_S X\to X\]
	of a smooth $S$-affine group scheme $G$. Let $\Aa$ be a $G$-equivariant tdo on $X$. Let $\pr_2$ be the projection $G\times_S X\to X$.
	
	\begin{definition}[{\cite[4.7]{kashiwara1989}}]
		A $G$-equivariant quasi-coherent left $\Aa$-module on $X$ is a quasi-coherent left $\Aa$-module, equipped with the structure $\beta:a^\ast\Mm\cong \pr_2^\ast\Mm$ of a $G$-equivariant $\OO_X$-module such that $\beta$ is $\pr_2^\cdot\Aa$-linear via the isomorphism $\pr^\cdot_2\Aa \cong a^\cdot\Aa$. We denote the category of $G$-equivariant quasi-coherent left $\Aa$-modules by $\Mod_{\qc}^G(\Aa)$. One can also define $G$-equivariant quasi-coherent right $\Aa$-modules in a similar way.
	\end{definition}
	
	A formal consequence of Theorem \ref{thm:xbc} and its formation by isomorphisms implies

	\begin{theorem}\label{thm:inductionofequivtwistedDmod}
		Let $i:Y\hookrightarrow X$ be an affine immersion of smooth $S$-schemes, equipped with actions of respective smooth $S$-affine group schemes $K,G$ over $S$. Suppose that we are given a homomorphism $K\to G$ for which $i$ intertwines the actions. Let $\Aa$ be a $G$-equivariant tdo on $X$. Assume that the functors 
		\[Ri_\ast:D(i^{-1}\Aa)\to D(\Aa),\]
		\[R(G\times_S i)_\ast:D((G\times_S i)^{-1}\pr_2^\cdot\Aa)\to D(\pr_2^\cdot\Aa)\]
		are locally bounded. Then the direct image functor lifts to $i_+:\Mod^G_{\qc}(i^\cdot \Aa)\to \Mod^G_{\qc}(\Aa)$.
	\end{theorem}
	
	\begin{remark}
		Let $G$ be a connected complex reductive algebraic group, and $K\subset G$ be the fixed point subgroup of an involution of $G$. Then any $K$-orbit on $\mathcal{B}_G$ is affinely embedded to $\mathcal{B}_G$ (\cite[Proposition 4.1]{hechtetal}, \cite[II Lemma 1.2]{bien1990}). See \cite{hayashikgb} for its $\ZZ\left[1/2\right]$-analog in the case $G=\SL_3$. Its generalization to general $G$ is in progress (\cite{hayashiuniform}).
	\end{remark}
	
	\subsection{Application to $(\lieg,K)$-modules}\label{sec:(g,K)-module}
	
	Let $i:K\to G$ be a homomorphism of smooth $S$-affine group schemes over a scheme $S$. Then one can define the notion of $(\lieg,K)$-modules over $S$ in the same line as \cite[section 2.1]{hayashi2019} just by replacing $k\cmod$ (= $\Mod(k)$ in the present paper) with $\Mod_{\qc}(\OO_S)$ (cf.\ \cite[section 2.2]{bernsteinetala}, \cite[section 2.4]{bernsteinetalb}). Let us denote the category of $(\lieg,K)$-modules by $\Mod(\lieg,K)$. For a morphism $s:S'\to S$ of schemes, one can define the base change functor $s^\ast:\Mod(\lieg,K)\to \Mod(s^\ast\lieg,K\times_S S')$ (cf.\ \cite[section 3.1]{hayashi2018}). 
	
	\begin{theorem}\label{thm:globalization}
		Let $X$ be a smooth and concentrated scheme over $S$, equipped with an action of $G$. Let $\Aa$ be a $G$-equivariant tdo on $X$. Then the functor $x_\ast:\Mod_{\qc}(\OO_X)\to \Mod_{\qc}(\OO_S)$ lifts to
		$\Mod^K_{\qc}(\Aa)\to \Mod(\lieg,K)$.
		Moreover, it commutes with flat base change functors.
	\end{theorem}
	
	The base change for $(\lieg,K)$-modules is defined in a similar way to \cite[Proposition 3.1.1]{hayashi2018} (omitted).
	
	\begin{proof}
		The construction of the structure of $(\lieg,K)$-modules is obtained by a minor modification of the proof of \cite[Theorem 11.5.3]{hottaetal2008}. In fact, recall Construction \ref{cons:differentialaction} for the Lie algebra action. One can define the structure of $K$-modules by using the flat base change theorem (\cite[Theorem (6.7)]{altmaneltal}) and the projection formula for affine morphisms. The latter assertion follows from Example \ref{ex:flatbasechange}.
	\end{proof}
	
	As we saw in the former section, we can construct $K$-equivariant $\Aa$-modules with the direct image functor. We will revisit applications to arithmetic models of cohomologically induced modules later.

	\section{$K$-conjugacy classes of stable parabolic subgroups}\label{sec:descentoforbit}
	
	Towards geometric construction of arithmetic forms of cohomologically induced modules, forms of their corresponding orbits are key ingredients. Exactly, the orbits consist of $\theta\otimes_\RR\CC$-stable parabolic subgroups for a connected real reductive algebraic group $G$ with a Cartan involution $\theta$. For a reductive group scheme $G$ and a smooth closed subgroup scheme $K\subset G$, we introduce a certain class of parabolic subgroups of $G$ in this section, which generalizes $\theta\otimes_\RR\CC$-stable parabolic subgroups. We study their moduli spaces. In particular, we give its \'etale local $K$-orbit decomposition under mild conditions on $K$ to obtain forms of $K$-orbits over \'etale local bases. Finally, we compute examples for reductive group schemes over $\ZZ\left[1/2\right]$ given in \cite{hayashilinebdl}.
	
	In this section, we follow the notations and conventions of \cite[sections 1.5 and 3.3]{hayashilinebdl} on the standard $\ZZ$- and $\ZZ\left[1/2\right]$-forms of classical Lie groups and related matrices . In particular, for these forms $G$, we will regard $G$ as group schemes over ($\ZZ\left[1/2\right]$-)schemes by base change without changing the symbols $G$ if there is no risk of confusion.
	
	If we are given a double Galois covering $S'\to S$ of schemes, we denote associated actions of the nontrivial element of the Galois group $\ZZ/2\ZZ$ on objects by $\bar{}$.

	\subsection{Stable parabolic subgroups}\label{sec:stablepsg}
	
	Let $G$ be a reductive group scheme over a scheme $S$, and $K$ be a smooth closed subgroup scheme of $G$. We remark that a reductive subgroup of $G$ is automatically closed by \cite[Theorem 5.3.5]{conrad2014}.
	
	\begin{definition}\label{defn:stable}
		\begin{enumerate}
			\item We call a parabolic subgroup $P\subset G$ is stable relative to $K$ if $P$ is \'etale locally of the form $P_G(\mu)$\footnote{In \cite{conradpseudo}, the group scheme $P_G(\mu)$ is defined only over a commutative ring. Since the definition is local, our definition is valid.} for a cocharacter $\mu$ of $G$ factoring through $K$.
			\item A Borel subgroup is called stable if it is so as a parabolic subgroup.
			\item For a stable parabolic subgroup $P$, we will denote $P_K=P\cap K$.
		\end{enumerate}
	\end{definition}
	
	\begin{remark}\label{rem:unitreduction}
		It is evident by definition that every cocharacter of $K$ factors through the unit component $K^0$. Hence the stability depends only on $K^0$.
	\end{remark}
	
	\begin{remark}[{\cite[Lemma 3.6]{hayashilinebdl}}]
		A Borel subgroup is stable if and only if it is \'etale locally of the form $P_G(\mu)$ for a regular cocharacter $\mu$ of $G$ in the sense of \cite[Definition 3.4]{hayashilinebdl} such that $\mu$ factors through $K$.
	\end{remark}

	\begin{definition}\label{def:stableparabolicsubgroupsrelativetoK}
		Let $\Pp_{G}^{K-\mathrm{st}}$ be the moduli space of stable parabolic subgroups relative to $K$. That is, $\Pp_{G}^{K-\mathrm{st}}$ is the $S$-space defined by
		\[\Pp_{G}^{K-\mathrm{st}}(T)=\{\mathrm{stable\ parabolic\ subgroups\ of\ }G\times_S T\ \mathrm{relative\ to\ } K\times_S T\}.\]
		Likewise, we define the moduli space $\mathcal{B}_{G}^{K-\mathrm{st}}$ of stable Borel subgroups relative to $K$ as
		\[\mathcal{B}_{G}^{K-\mathrm{st}}(T)=\{\mathrm{stable\ Borel\ subgroups\ of\ }G\times_S T\ \mathrm{relative\ to\ } K\times_S T\}.\]
		It is clear by definition that $\Pp_{G}^{K-\mathrm{st}}$ (resp.~$\mathcal{B}_{G}^{K-\mathrm{st}}$) a subsheaf of $\Pp_{G}$ (resp.~$\mathcal{B}_G$) over the big \'etale site of $S$.
	\end{definition}
	
	The following result is our basic motivation to introduce this class:
	
	\begin{proposition}\label{prop:stability}
		Let $P\subset G$ be a stable parabolic subgroup relative to $K$.
		\begin{enumerate}
			\renewcommand{\labelenumi}{(\arabic{enumi})}
			\item The subgroup scheme $P_K$ is smooth.
			\item Suppose that $K$ is reductive. Then $P_K$ is a parabolic subgroup of $K$. Moreover, if $P$ is Borel, so is $P_K$.
		\end{enumerate}
	\end{proposition}
	
	\begin{proof}
		Work locally in the \'etale topology of $S$ to assume $S$ affine and $P=P_G(\mu)$ for a cocharacter $\mu$ of $K$. Then the assertions follow from \cite[Proposition 4.1.10 1]{conrad2014} and \cite[Corollary 3.1.7]{hayashilinebdl}.
	\end{proof}
	
	We wish to classify geometric $K$-conjugacy classes of stable parabolic subgroups. For a milestone, let us introduce the following notions:
	
	\begin{definition}\label{defn:rtype}
		\begin{enumerate}
			\renewcommand{\labelenumi}{(\arabic{enumi})}
			\item Define $\rtype(G,K)$ (resp.~$\rtype_\emptyset(G,K)$) as the \'etale quotient sheaf $K\backslash\Pp^{K\mathrm{-st}}_G$ (resp.~$K\backslash\mathcal{B}^{K\mathrm{-st}}_G$). A section of $\rtype(G,K)$ is called a relative (parabolic) type.
			\item Let $rt:\Pp^{K\mathrm{-st}}_G\to \rtype(G,K)$ denote the canonical projection.
			\item For an $S$-scheme $S'$ and $x'\in\rtype(G,K)(S')$, set \[\Pp^{K\mathrm{-st}}_{G,x'}\coloneqq rt^{-1}(x')=\Pp^{K\mathrm{-st}}_G\times_{\rtype(G,K)} S'.\]
			\item Let $gt:\rtype(G,K)\to\type G$ be the morphism of \'etale $S$-sheaves induced from the $K$-equivariant map $\Pp^{K\mathrm{-st}}_G\subset\Pp_G\overset{t}{\to}\type G$.
			\item For each $x\in \rtype(G,K)(S)$, we denote the induced monomorphism $\Pp^{K\mathrm{-st}}_{G,x}\hookrightarrow \Pp_{G,gt(x)}$
			by $i_x$.
			\item If $K$ is reductive, let $kt:\rtype(G,K)\to\type K$ be the morphism of \'etale $S$-sheaves induced from the $K$-invariant map $\Pp^{K\mathrm{-st}}_G\overset{-\cap K}{\to}
			\Pp_K\overset{t}{\to}\type K$.
		\end{enumerate}
	\end{definition}
	
	In other words, we aim to study $\rtype(G,K)$. This object should be thought of for our purpose since the descent problem of bases of closed $K$-orbits attached to stable parabolic subgroups is reduced to that of their relative types. Recall that $\rtype(G,G)=\type G$, and that it is represented by a finite \'etale $S$-scheme. The goal of this section is to prove its relative analog, i.e., to prove that $\rtype(G,K)$ is represented by a finite \'etale $S$-scheme under a certain assumption. To be precise, the statement is as follows:
	
	\begin{theorem}\label{thm:orbit}
		Let $G$ be a reductive group scheme, and $K$ be a smooth closed subgroup scheme with reductive unit component such that $\pi_0(K)$ is finite \'etale. Let $x$ be a relative parabolic type over $S$.
		\begin{enumerate}
			\renewcommand{\labelenumi}{(\arabic{enumi})}
			\item The $S$-space $\Pp_{G}^{K\mathrm{-st}}$ is represented by a smooth closed subscheme of $\Pp_G$.
			\item The $S$-space $\rtype(G,K)$ is represented by a finite \'etale $S$-scheme.
			\item Stable parabolic subgroups $P$ and $P'$ are \'etale locally $K$-conjugate to each other if and only if $rt(P)=rt(P')$.
			\item The morphism $rt:\Pp^{K\mathrm{-st}}_G\to \rtype(G,K)$ is surjective.
			\item The morphism $rt:\Pp^{K\mathrm{-st}}_G\to \rtype(G,K)$ is smooth projective.
			\item The morphism $i_x$ is a closed immersion.
			\item If $gt(x)=\emptyset$ and $K=K^0$ then we have $kt(x)=\emptyset$.
			\item If $K=K^0$ then the $K$-equivariant map $\Pp^{K\mathrm{-st}}_{G,x}\to\Pp_{K,kt(x)}$ defined by $P\mapsto P\cap K$ is an isomorphism.
		\end{enumerate}
	\end{theorem}
	
	There are elementary examples where $\pi_0(K)$ are finite \'etale:
	
	\begin{example}\label{ex:c1}
		Put $S=\Spec\ZZ\left[1/2\right]$. Let $n$ be a positive integer. Then we have
		\[\SO(2n)\rtimes
		\{I_{2n},\diag(I_{2n-1},-1)\}_S
		\cong\Oo(2n)\]
		\[\SO(2n-1)\rtimes
		\{\pm I_{2n-1}\}_S
		\cong\Oo(2n-1).\]
	\end{example}
	
	\begin{example}\label{ex:c2}
		Put $S=\Spec\ZZ\left[1/2\right]$. Let $p,q$ be positive integers. Then we have
		\[(\SO(2p)\times\SO(2q))\rtimes\{I_{2p+2q},\diag(I_{2p-1},-1,I_{2q-1},-1)\}_S
		\cong\mathrm{S}(\Oo(2p)\times\Oo(2q)),\]
		\[(\SO(2p-1)\times\SO(2q-1))\rtimes\{\pm I_{2p+2q-2}\}_S
		\cong\mathrm{S}(\Oo(2p-1)\times\Oo(2q-1)),\]
		\[(\SO(2p)\times\SO(2q-1))\rtimes\{I_{2p+2q-1},\diag(I_{2p-1},-I_{2q})\}_S
		\cong\mathrm{S}(\Oo(2p)\times\Oo(2q-1)).\]
	\end{example}
	
	\begin{proof}[Proof of: (1) and (2) $\Rightarrow$ (5) - (7) in Theorem \ref{thm:orbit}]
		Suppose that (1) and (2) are verified. Then (5) follows from \cite[Corollaire 3.5]{sga3-26}, \cite[Proposisition 2.4.1 (iv)]{fu2011}, and \cite[Remarques (5.5.4) (v)]{ega2}.

		For (6), observe that the section $x$ is a closed immersion into $\rtype(G,K)$ since $\rtype(G,K)$ is separated over $S$. Similarly, $gt(x)$ is a closed immersion into $\type G$ (recall \cite[sections 3.1 and 3.2]{sga3-26}). Hence their base change maps $i'_x:\Pp^{K\mathrm{-st}}_{G,x}\to \Pp^{K\mathrm{-st}}_{G}$ and $\Pp_{G,gt(x)}\to \Pp_G$ are closed immersions. The assertion now follows since $i'_x$ factors through $i_x$ by definition.
		
		Finally, we prove (7). Notice that $-\cap K:gt^{-1}(\emptyset)\to \Pp_{K}$ factors through $\Pp_{K,\emptyset}=\Bb_K$ if $K=K^0$ (Proposition \ref{prop:stability} (2)). This completes the proof.
	\end{proof}
	
	Part (8) is a consequence of the following observation:
	
	\begin{lemma}\label{lem:conjugacy}
		Assume $K=K^0$. Let $P,P'$ be stable parabolic subgroups which are \'etale locally $K$-conjugate to each other. Then we have $P=P'$ if and only if $P_K=P_{K'}$.
	\end{lemma}
	
	\begin{proof}
		The ``only if'' direction is trivial. To prove the converse, we may assume that $P'=gPg^{-1}$ for a certain element $g\in K(S)$. We take $(-)_K$ to obtain
		$P_K=P'_K=gP_Kg^{-1}$. This implies $k\in P_K(S)\subset P(S)$ and therefore $P'=gPg^{-1}=P$ (Proposition \ref{prop:stability} (2), \cite[Proposition 1.2]{sga3-26}).
	\end{proof}
	
	\begin{proof}[Proof of: (1) - (4) $\Rightarrow$ (8) in Theorem \ref{thm:orbit}]
		Suppose that (1) - (4) are verified. Assume $K=K^0$. We tentatively write $f$ for the map $\Pp^{K\mathrm{-st}}_{G,x}\to\Pp_{K,kt(x)}$ in (8). Then $f$ is a monomorphism by (3) and Lemma \ref{lem:conjugacy}.
		
		We next prove that this map is section-wisely surjective. We may assume that $S$ is nonempty. In this case, $\Pp^{K\mathrm{-st}}_{G,x}$ is nonempty by (4). Working locally in the \'etale topology of $S$, we may assume that $\Pp^{K\mathrm{-st}}_{G,x}(S)$ is nonempty (\cite[Proposition 1.10]{sga3-11}). The assertion now follows from \cite[Th\'eor\`eme 3.3 (i)]{sga3-26}. This completes the proof.
	\end{proof}
	
	It remains to prove (1) - (4). Our strategy is to locally decompose $\Pp_{G}^{K\mathrm{-st}}$ into the disjoint union of $K$-orbits consisting of stable parabolic subgroups in the \'etale topology of $S$. We will achieve this by finding a combinatorial description of $\rtype(G,K)$. We also find a combinatorial description for the set of the parabolic subgroups $P$ indexing the disjoint union. For a while, we forget the assumption that $\pi_0(K)$ is finite \'etale.
	
	Since the statements are local in the \'etale topology of $S$, we may assume that $S$ is nonempty affine, and that $K^0$ has a splitting $(T,N,\Delta(K^0,T))$ (\cite[D\'efinition 1.13]{sga3-22}). Let $W(\Delta(K^0,T))$ denote the Weyl group of the root datum of this splitting. We also choose a simple system $\Pi$ of $\Delta(K^0,T)$. Let $\Delta^+(K^0,T)$ denote the corresponding positive system. Write $N^\vee_+$ for the set of dominant constant cocharacters of $T$. If $S$ is connected, we will denote $N^\vee_+=X_\ast(T)_+$. 
	
	\begin{proposition}\label{prop:reduction}
		\begin{enumerate}
			\item A stable parabolic subgroup $P\subset G$ relative to $K^0$ is locally $K^0$-conjugate to $P_G(\mu)$ in the \'etale topology of $S$ for $\mu\in N^\vee_+$.
			\item For $\mu_1,\mu_2\in N^\vee_+$, $P_G(\mu_1)$ and $P_G(\mu_2)$ are \'etale locally $K^0$-conjugate to each other if and only if $P_G(\mu_1)=P_G(\mu_2)$.
		\end{enumerate}
	\end{proposition}

	\begin{proof}
		For (1), we may assume $P=P_G(\mu)$ for a cocharacter $\mu$ of $K^0$. Pass to an \'etale covering to let $\mu$ factor through a maximal torus of $K$ by \cite[Lemma 5.3.6]{conrad2014}. Since maximal tori of $K$ are \'etale locally conjugate to each other, we may assume that $\mu$ factors through $T$ (\cite[Theorem 3.2.6]{conrad2014}). Then work locally in the Zariski topology of $S$ to get $\mu\in N^\vee$. Finally, we may assume $\mu$ to be dominant by the passage to the action of the Weyl group and the surjectivity of the canonical map
		$N_{K^0(S)}(T)\to W(K^0,T)(S)\supset W(\Delta(K^0,T))$
		(\cite[Corollary 5.1.11]{conrad2014}).
		
		The ``if'' direction of (2) is obvious. Suppose that $P_G(\mu_1)$ and $P_G(\mu_2)$ are \'etale locally $K^0$-conjugate to each other. We may only prove that
		$P_{K^0}(\mu_1)=P_{K^0}(\mu_2)$ (Lemma \ref{lem:conjugacy}). This follows from 
		\cite[Proposition 1.17]{sga3-26} since $\mu_1$ and $\mu_2$ are dominant. In fact, they contain the Borel subgroup attached to $\Delta^+(K^0,T)$.
	\end{proof}

	Without constancy of the cocharacters in $G$, the set of stable parabolic subgroups of the form $P_G(\mu)$ ($\mu\in N^\vee_+$) may not be invariant under base change. To resolve this issue, we give a splitting of $G$ which is compatible with that of $K^0$:
	
	\begin{lemma}\label{lem:relativesplitting}
		There exists a splitting $(H,M,\Delta(G,H))$ of $G$ locally in the \'etale topology of $S$ such that $T\subset H$ and $N^\vee\subset M^\vee$.
	\end{lemma}
	
	\begin{proof}
		Notice that $Z_G(T)$ is a reductive group scheme over $S$ (see \cite[Proof of Theorem 3.1.3]{hayashilinebdl}). Choose a maximal torus $H$ of $Z_G(T)$. Then $H$ contains $T$ by $T\subset Z_{Z_G(T)}(H)=H$. It is easy to show that $H$ is a maximal torus of $G$ (see it at each geometric fiber). We may pass to an \'etale covering to assume that $(G,H)$ is split for a free abelian group $M$ of finite rank. Pass to a Zariski covering to assume $N^\vee \subset M^\vee$. In fact, choose a free basis $\{\mu_i\}$ of $N^\vee$. Then $\mu_i$ is locally constant as a cocharacter of $H$ for each $i$. This completes the proof.
	\end{proof}
	
	Henceforth assume that a splitting $(H,M)$ in Lemma \ref{lem:relativesplitting} is given for the proof of Theorem \ref{thm:orbit}. We denote the canonical pairing of $M^\vee$ and $M$ by $\langle-,-\rangle$. For $\mu\in N^\vee\subset M^\vee$, set
	\[\Phi_\mu=\{\alpha\in \Delta(G,H):
	~\langle\mu,\alpha\rangle\geq 0\}.\]
	Let $\Ss(G,K^0)$ (resp.~$\Ss^+(G,K^0)$) be the set of parabolic subsets of the form $\Phi_\mu$ with $\mu\in N^\vee$ (resp.~$\mu\in N^\vee_+$). The set $\Ss^+(G,K^0)$ is nonempty since it contains $\Delta(G,H)$. For $\Phi\in\Ss(G,K^0)$, we denote the corresponding parabolic subgroup by $P_{\Phi}$. We define $\Ss(G,K^0)$ here because we will use it after a short while.
	
	\begin{proof}[First proof of Theorem \ref{thm:orbit} (1) - (4)]
		Assume first that $K=K^0$. Proposition \ref{prop:reduction} implies $\Pp_{G}^{K-\mathrm{st}}\cong\coprod_{\Phi\in\Ss^+(G,K^0)} K^0/P_{\Phi,K^0}$ (use \cite[Th\'eor\`eme 10.1.2]{sga3-5} for the representability of $K^0/P_{\Phi,K^0}$). The assertions are its immediate consequences.
		
		We next discuss general $K$ with $\pi_0(K)$ finite \'etale over $S$. Part (1) follows from Remark \ref{rem:unitreduction}. The remaining assertions are verified by $\rtype(G,K)\cong\pi_0(K)\backslash \rtype(G,K^0)$. In fact, this is \'etale locally the constant sheaf attached to the corresponding coset. This shows (2). In virtue of this construction, the canonical map $\rtype(G,K^0)\to \rtype(G,K)$ is surjective, which implies (4). For (3), we may write $\pi_0(K)=C_S$ for a finite set $C$. Work locally in the \'etale topology of $S$ to pick an element of $K(S)$ representing $\{c\}_S$ for each $c\in C$ (see \cite[Proposition 1.10]{sga3-11} for existence). Then the assertion is straightforward.
	\end{proof}
	
	\begin{example}
		Let $S'\to S$ be a Galois covering of connected schemes with Galois group $\Gamma=\ZZ/2\ZZ$. Suppose that we are given maximal tori $T\subset K$ and $H\subset G$ such that $T\subset H$, and that $(K\times_S S',T\times_S S')$, $(G\times_S S',H\times_S S')$ are split reductive. Fix a simple system $\Pi$ of $(K\times_S S',T\times_S S')$. Let $w\in W(K,T)(S')$ be the unique element such that $\bar{\Pi}=w\Pi$. Define an action of $\Gamma$ on $\Ss^+(G\times_S S',H\times_S S')$ by $\Phi\mapsto \overline{w\Phi}$. Then we have
		\[\rtype(G,K)=
		\coprod_{\overset{\OO\in\Gamma\backslash\Ss^+(G\times_S S',H\times_S S')}{|\OO|=1}} S
		\coprod
		\coprod_{\overset{\OO\in\Gamma\backslash\Ss^+(G\times_S S',H\times_S S')}{|\OO|=2}} S'\]
	\end{example}
	
	Let us note a way to compute $kt$ in terms of roots of $(G,H)$. For this, it suffices to describe the set of roots in $P_{K^0}(\mu)$ for $\mu\in N^\vee$. Here we do not restrict ourselves to dominant cocharacters because of latter applications (Lemma \ref{lem:theta-stable_parabolic_subset}, Theorem \ref{thm:solutiontotheclassificationproblem}). Pick a splitting $(H,M)$ in Lemma \ref{lem:relativesplitting}. 
	Assume that $\alpha|_T\in N$ for all $\alpha\in\Delta(G,H)$.
	This holds after Zarikski localization by a similar argument to Lemma \ref{lem:relativesplitting}, or if $S$ is connected. Then for a subset $\Phi\subset\Delta(G,H)$, define
	\begin{equation}
		\Phi|_T\coloneqq\{\beta\in N:~
		\mathrm{there\ exists\ }\alpha\in\Phi
		\mathrm{\ such\ that\ }\beta=\alpha|_T\}.
		\label{eq:Phi|_T}
	\end{equation}
	
	\begin{proposition}\label{prop:kt}
		Let $\Phi\in \Ss(G,K^0)$.
		\begin{enumerate}
			\item The subset $\Phi|_T\cap\Delta(K^0,T)$ is a parabolic subset of $\Delta(K^0,T)$. Moreover, the attached parabolic subgroup of $K^0$ is $P_{\Phi,K^0}$.
			\item The following conditions are equivalent:
			\begin{enumerate}
				\item[(a)] $\Phi$ belongs to $\Ss^+(G,K^0)$;
				\item[(b)] Any $\mu\in N^\vee$ satisfying $\Phi=\Phi_\mu$ is dominant;
				\item[(c)] $\Phi|_T$ contains $\Delta^+(K^0,T)$.
			\end{enumerate}
			\item Assume that $\Phi$ is a positive system of $\Delta(G,H)$. Then $\Phi|_T\cap \Delta(K^0,T)$ is a positive system of $\Delta(K^0,T)$. In particular, the following conditions are equivalent:
			\begin{enumerate}
				\item[(a)] $\Phi|_T\cap \Delta(K^0,T)=\Delta^+(K^0,T)$;
				\item[(b)] $\Phi|_T\supset\Delta^+(K^0,T)$;
				\item[(c)] $\Phi|_T\cap \Delta(K^0,T)\subset\Delta^+(K^0,T)$;
				\item[(d)] $\Phi|_T\cap (-\Delta^+(K^0,T))=\emptyset$.
			\end{enumerate}
		\end{enumerate}
		.
	\end{proposition} 
	
	\begin{proof}
		Write $\Phi=\Phi_\mu$ for $\mu\in N^\vee$. Let $\beta\in\Delta(K^0,T)$. Then one can find $\alpha\in \Delta(G,H)$ such that $\beta=\alpha|_T$ since $\alpha|_T\in N$ for $\alpha\in\Delta(G,H)$. For such $\alpha$, the following conditions are equivalent:
		\begin{enumerate}
			\renewcommand{\labelenumi}{(\roman{enumi})}
			\item $\langle \mu,\beta\rangle\geq 0$;
			\item $\langle \mu,\alpha\rangle\geq 0$;
			\item $\alpha\in\Phi$;
			\item $\beta\in\Phi|_T$
		\end{enumerate}
		by the equality $\langle \mu,\beta\rangle=\langle \mu,\alpha\rangle$ and the definition of $\Phi_\mu$. In fact, the equivalences (i) $\iff$ (ii) $\iff$ (iii) and the implication (iii) $\Rightarrow$ (iv) are straightforward. To prove (iv) $\Rightarrow$ (i), we choose $\alpha'\in\Phi$ satisfying $\beta=\alpha'|_T$. Then apply (i) $\iff$ (iii) to $(\alpha',\beta)$ to show (i). 
		
		We now deduce $\Phi|_T\cap \Delta(K^0,T)=\{\beta\in\Delta(K^0,T):~\langle\mu,\beta\rangle\geq 0\}$. Part (1) is now clear. The implications (a) $\Rightarrow$ (c) $\Rightarrow$ (b) in (2) are also consequences of the equivalence in the last paragraph. Condition (b) implies (a) by definition of $\Ss(G,K^0)$.
		Part (3) follows from (1) and Proposition \ref{prop:stability} (2). 
	\end{proof}
	
	Let us revisit Theorem \ref{thm:orbit}. One can still prove (2) as a consequence of a version of Proposition \ref{prop:reduction} for $K$. In fact, $\rtype(G,K)$ will be described in terms of roots. To do this, let us introduce:
	
	\begin{condition}\label{cond:splitting}
		Let $K$ be a smooth closed subgroup scheme of $G$ with $K^0$ reductive. Suppose that we are given a splitting $(T,N,\Delta(K^0,T))$ of $K^0$. Then we say that $K$ satisfies Condition (S) if there exists a constant subscheme $C_S\subset N_K(T)$ attached to a finite subset $C\subset N_{K(S)}(T)$ containing the unit such that:
		\begin{enumerate}
			\renewcommand{\labelenumi}{(\roman{enumi})}
			\item the multiplication map $K^0\times_S C_S\to K$ is an isomorphism of $S$-schemes, and
			\item the action of each element $c\in C$ on the cocharacter group of $T$ respects $N^\vee$.
		\end{enumerate}
	\end{condition}
	
	\begin{example}
		Condition (S) holds if $K=K^0$ (put $C=\{e\}$).
	\end{example}
	
	\begin{example}
		Condition (ii) of (S) holds if $(K,T,N)$ and $C$ are obtained by base change of those over a connected base.
	\end{example}

	\begin{proposition}
		The $S$-scheme $\pi_0(K)$ is finite \'etale if and only if the triple $(K,T,N)$ \'etale locally satisfies Condition (S).
	\end{proposition}
	
	For $c\in N_{K(S)}(T)$ and a cocharacter $\mu$ of $T$, we define a new cocharacter $c\mu$ of $T$ as
	\[(c\mu)(-)=c\mu(-)c^{-1}.\]
	
	\begin{proof}
		The ``only if'' direction is clear. For the converse, we may pick a finite subset $C\subset K(S)$ with $C_S\cong \pi_0(K)$ for $C\ni c\mapsto cK^0\in \pi_0(K)(S)$ (see First proof of Theorem \ref{thm:orbit} (1) - (4)). In particular, one may regard $C_S$ as a closed subscheme of $K$, and the multiplication map $K^0\times_S C_S\to K$ is an isomorphism.
		
		Pick $c\in C$. Since $cTc^{-1}$ is a maximal torus of $K^0$, $cTc^{-1}$ is \'etale locally $K^0$-conjugate to $T$ (\cite[Theorem 3.2.6]{conrad2014}). Replacing $S$, we may assume that there exists $g_c\in K^0(S)$ such that $cTc^{-1}=g_cTg^{-1}_c$. Therefore we may replace $c$ with $g^{-1}_c c$ so that $c$ normalizes $T$.
		
		Fix a basis $\{\mu_i\}$ of $N^\vee$. For each $c\in C$ and $\mu_i\in N^\vee$, $c\mu_i$ is locally constant. Work locally in the Zariski topology of $S$ so that $c\mu_i\in N^\vee$ for all $c,\mu_i$. Since $\mu_i$ form a basis of $N^\vee$, each $c\in C$ normalizes $N^\vee$. This completes the proof.
	\end{proof}
	
	Henceforth we fix $C$ in Condition (S) for the second proof of Theorem \ref{thm:orbit} (2). That is, $K$ is a smooth closed subgroup scheme of $G$ with $K^0$ reductive. The base $S$ is nonempty and affine. A splitting $(T,N,\Delta(K^0,T))$ of $K^0$ satisfying Condition (S) is given. Let $C\subset K(S)$ be a corresponding subset. We also assume that a splitting $(H,M,\Delta(G,H))$ of $G$ such that $T\subset H$ and $N^\vee\subset M^\vee$ is given.
	
	\begin{proposition}\label{prop:classification1}
		For $\mu_1,\mu_2\in N^\vee_+$, the following conditions are equivalent:
		\begin{enumerate}
			\renewcommand{\labelenumi}{(\alph{enumi})}
			\item The parabolic subgroups $P_G(\mu_1)$ and $P_G(\mu_2)$ are \'etale locally $K$-conjugate to each other;
			\item there exist $c\in C$ and $w\in W(\Delta(K^0,T))$ such that $wc\mu_1\in N^\vee_+$ and $P_G(wc\mu_1)=P_G(\mu_2)$;
			\item there exist $c\in C$ and $w\in W(\Delta(K^0,T))$ such that $wc\mu_1\in N^\vee_+$, and that $P_G(wc\mu_1),P_G(\mu_2)$ have the same Lie subalgebras of $\lieg$;
			\item there exist $c\in C$ and $w\in W(\Delta(K^0,T))$ such that $wc\mu_1\in N^\vee_+$ and $\Phi_{wc\mu_1}=\Phi_{\mu_2}$.
		\end{enumerate}
		Moreover, the equalities in (b), (c), and (d) are independent of the choice of $w$ satisfying $wc\mu_1\in N^\vee_+$.
	\end{proposition}
	
	\begin{proof}
		In general, let $\mu\in N^\vee$. Then $P_G(\mu)$ contains the maximal torus $H$ by definition. Moreover, the set of roots of $P_G(\mu)$ with respect to $H$ is exactly $\Phi_\mu$ by \cite[Proposition 2.1.8]{conradpseudo}. The equivalence of (b), (c), and (d) now follows from \cite[Corollaire 5.3.5]{sga3-22} (see also \cite[Proposition 5.1.3, Exemples 5.2.3 b)]{sga3-22}).
		
		Before proving the remaining equivalence, let us prove the latter assertion. From the last paragraph, we may only see the assertion for (b). Suppose that we are given $c\in C$ and
		\[w,w'\in W(\Delta(K^0,T))\]
		such that $P_G(wc\mu_1)=P_G(\mu_2)$ and $wc\mu_1,w'c\mu_1\in N^\vee_+$. Since $P_G(w'c\mu_1)$ is $K^0(S)$-conjugate to $P_G(wc\mu_1)$, we get
		\[P_G(w'c\mu_1)=P_G(wc\mu_1)=P_G(\mu_2)\]
		by Proposition \ref{prop:reduction} (2).
		
		It is evident that (b) implies (a). The proof is completed by showing that (a) implies (d). For this, we may pass to a geometric fiber to assume that $S=\Spec F$ for some algebraically closed field $F$. In this case, one can find $(g,c)\in (K^0\times C)(S)\cong K(S)$ such that
		\[g cP_G(\mu_1)c^{-1} g^{-1}= P_G(\mu_2).\]
		Choose $w\in W(\Delta(K^0,T))$ so that $wc\mu_1\in N^\vee_+$. Then Proposition \ref{prop:reduction} (2) implies
		\[P_{G}(w c\mu_1)=P_G(\mu_2)\]
		since $w$ is represented by an element of $K^0(S)$.
	\end{proof}
	
	Let $\Ss^+(G,K)$ denote the set of equivalence classes of parabolic subsets of the form $\Phi_\mu$ with $\mu\in N^\vee_+$, where the equivalence relation is defined by the equivalent conditions of Proposition \ref{prop:classification1}.
	
	\begin{example}\label{ex:K=K^0}
		If $K=K^0$, the equivalence relation is trivial. In particular, the notations of $\Ss^+(G,K^0)$ do not conflict.
	\end{example}

	\begin{proof}[Second proof of Theorem \ref{thm:orbit} (2)]
		Consider the setting of Proposition \ref{prop:classification1}. Fix a complete system $\{\Phi\}$ of representatives of $\Ss^+(G,K)$. Then we have a $K$-equivariant isomorphism
		\[\Pp_{G}^{K-\mathrm{st}}\cong\coprod_{\Phi\in\Ss^+(G,K)} K/P_{\Phi,K}.\]
		This implies $\rtype(G,K)\cong \Ss^+(G,K)_S$.
	\end{proof}
	
	Let us note that the preceding argument naturally raises the following question towards local computation of $\rtype(G,K)$ (cf.~Theorem \ref{thm:solutiontotheclassificationproblem}):
	
	\begin{problem}\label{prob:classification}
		Consider the setting of Proposition \ref{prop:classification1}. Then compute $\Ss^+(G,K)$.
	\end{problem}
	
	This is independent of the base $S$ in the following sense:
	
	\begin{proposition}\label{prop:basefree}
		Consider the setting of Proposition \ref{prop:classification1}. Let $S'\to S$ be a morphism of affine schemes. Then the base change determines a bijection $\Ss^+(G,K)\to \Ss^+(G\times_SS',K\times_SS')$.
	\end{proposition}
	
	For instance, this enables us to apply the theory of algebraic groups over algebraically closed fields to the classification problem if necessary by taking a geometric point of $S$.
	
	\begin{proof}
		This is immediate by definitions.
	\end{proof}

	We return to the general setting: Let $G$ be a reductive group scheme over $S$, and $K$ be a smooth closed subgroup scheme of $G$ with $K^0$ reductive and $\pi_0(K)$ finite \'etale. A similar argument implies:
	
	\begin{variant}
		\begin{enumerate}
			\item The moduli space $\mathcal{B}^{K\mathrm{-st}}_G$ is represented by a smooth projective $S$-scheme.
			\item We have canonical isomorphisms
			\[\begin{array}{cc}
				\rtype_\emptyset(G,K)\cong gt^{-1}(\emptyset),
				&\mathcal{B}^{K\mathrm{-st}}_G
				\cong \mathcal{P}^{K\mathrm{-st}}_G\times_{\rtype(G,K)}
				\rtype_\emptyset(G,K).
			\end{array}\]
			In particular, $\rtype_\emptyset(G,K)$ is a finite \'etale $S$-scheme.
		\end{enumerate}
	\end{variant}
	
	We see two ``examples''.

	\begin{corollary}\label{cor:restrictionsetting}
		Let $S'\to S$ be a finite \'etale and surjective morphism of schemes, and $G$ be a reductive group scheme over $S$. Regard $G$ as a subgroup scheme of $\res_{S'/S} (G\times_S S')$ for the unit. Then intersection with $G$ determines an isomorphism $\Pp^{G\mathrm{-st}}_{\res_{S'/S} (G\times_S S')}\cong \Pp_G$.
	\end{corollary}
	\begin{proof}
		In view of \cite[Corollaire (17.9.5)]{ega44}, we may pass to geometric fibers to assume that $S$ is the spectrum of an algebraically closed field. Then the map $G\to \res_{S'/S} (G\times_S S')$ is identified with the diagonal map $G\to G^n$ for some positive integer $n$. The assertion now follows since the stable parabolic subgroups relative to $G$ are the diagonal parabolic subgroups.
	\end{proof}
	
	Our constructions in Definition \ref{defn:rtype} commute with the restriction functor:
	
	\begin{theorem}\label{thm:res_vs_Kst}
		Let $S'\to S$ be a finite \'etale morphism. Let $G'$ be a reductive group scheme over $S'$, and $K'$ be a smooth closed subgroup scheme of $G'$ with $\pi_0(K')$ finite \'etale.
		\begin{enumerate}
			\item The group scheme $\res_{S'/S} K'$ is a smooth closed subgroup scheme of $\res_{S'/S} G'$ with reductive unit component and $\pi_0(\res_{S'/S} K')$ finite \'etale.
			\item There are canonical isomorphisms
			\[\res_{S'/S}\Pp^{K'\mathrm{-st}}_{G'}\cong \Pp^{\res_{S'/S} K'\mathrm{-st}}_{\res_{S'/S} G'}\]
			\[\res_{S'/S} \rtype(G',K')\cong \rtype(\res_{S'/S} G',\res_{S'/S} K').\]
		\end{enumerate}
	\end{theorem}
	
	We will use this in section \ref{sec:gl_n}. In fact, the restriction of $\GL_n$ over number fields to $\QQ$ is a fundamental object in automorphic representation theory.
	
	For a proof of Theorem \ref{thm:res_vs_Kst}, let us note general results on the commutation of the restriction functor with the quotient:
	
	\begin{lemma}\label{lem:quotient_vs_res}
		Let $S'\to S$ be a finite \'etale morphism of schemes. Let $X'$ be an $S'$-scheme, equipped with an action of a group $S'$-scheme $G'$. Assume the following conditions:
		\begin{enumerate}
			\renewcommand{\labelenumi}{(\roman{enumi})}
			\item $\res_{S'/S} (X'/G')$ is represented by a locally finitely presented $S$-scheme.
			\item $(\res_{S'/S} X')/(\res_{S'/S} G')$ is represented by a flat locally finitely presented $S$-scheme.
		\end{enumerate}
		Then the canonical morphism
		$(\res_{S'/S} X')/(\res_{S'/S} G')\to \res_{S'/S} (X'/G')$
		is an isomorphism.
	\end{lemma}
	
	\begin{proof}
		In view of \cite[Corollaire (17.9.5)]{ega44}, we may work geometrically fiberwisely to assume \[\begin{array}{cc}
			S=\Spec F,&S'=\Spec F^n
		\end{array}\]
		for some $n\geq 0$ ($F^0=\{0\}$), where $F$ is an algebraically closed field. In this case, $X'$ and $G'$ are identified with disjoint unions of $n$ $F$-schemes $X_i$ and group $F$-schemes $G_i$ respectively. We now get
		\[\res_{S'/S} (X'/G')
		\cong \prod_{1\leq i\leq n} X_i/G_i
		\cong \left(\prod_{1\leq i\leq n} X_i\right)/\left(\prod_{1\leq i\leq n} G_i\right)
		\cong (\res_{S'/S} X')/(\res_{S'/S} G').
		\]
	\end{proof}

	\begin{proof}[Proof of Theorem \ref{thm:res_vs_Kst}]
		Recall that for every quasi-projective $S'$-scheme $X'$, $\res_{S'/S} X'$ is represented by an $S$-scheme by \cite[section 7.6, Theorem 4]{boschetal}. In particular, $\res_{S'/S} X'$ is representable if $X'$ is affine of finite type or projective over $S'$ (\cite[Proposition (5.3.4)]{ega2}).
		
		For (1), notice that $\res_{S'/S} K'$ and $\res_{S'/S} (K')^0$ are representable. Moreover, they are smooth by \cite[section 7.6, Proposition 5]{boschetal}. Therefore a similar argument to Lemma \ref{lem:quotient_vs_res} implies $\res_{S'/S} (K')^0\cong (\res_{S'/S} K')^0$. Since $(K')^0$ is reductive, so is $(\res_{S'/S} K')^0$.
		
		We next prove
		\[\pi_0(\res_{S'/S} K')\cong \res_{S'/S} \pi_0(K').\]
		Since $\pi_0(K')$ is finite \'etale, so is $\res_{S'/S} \pi_0(K')$ by the first paragraph and \cite[section 7.6, Proposition 5]{boschetal} (recall that a morphism of schemes is finite if and only if it is affine and proper, \cite[\href{https://stacks.math.columbia.edu/tag/01WN}{Tag 01WN}]{stacks-project}). In particular, it is locally of finite presentation over $S$. Recall that $\pi_0(\res_{S'/S} K')$ is \'etale over $S$. The isomorphism now follows from Lemma \ref{lem:quotient_vs_res}.
		
		For the first isomorphism of (2), observe that for an $S$-scheme $U$ and a parabolic subgroup $P'$ of $G'\times_{S'} (S'\times_S U)$, we have a parabolic subgroup
		\[\res_{S'\times_S U/U} P'\subset \res_{S'\times_S U/U} (G'\times_{S'} (S'\times_S U))
		\cong (\res_{S'/S} G' )\times_S U.\]
		Suppose that $P'$ is a stable parabolic subgroup relative to $K'$. We wish to prove that
		\[\res_{S'\times_S U/U} P'\]
		is a stable parabolic subgroup relative to $\res_{S'\times_S U/U} (K'\times_S'(S'\times_S U))$. Since the assertion is \'etale local in $U$, we may replace $U$ and $S$ to assume $S=U$ and $P'=P_{G'}(\mu')$ for a cocharacter $\mu'$ of $K'$. Let $\mu$ be the cocharacter of $\res_{S'/S} K'$ defined by the composition of $\res_{S'/S}\mu'$ and the unit $\GG_{m,S}\to \res_{S'/S} \GG_{m,S'}$. Then we have
		\[\res_{S'/S} P_{G'}(\mu')=P_{\res_{S'/S} G'}(\mu)\]
		by the right adjointness of the Weil restriction and the canonical isomorphism
		\[(\res_{S'/S} G')\times_S T\cong \res_{S'\times_ST/T} G'\times_{S'} (S'\times_S T)\]
		for every $S$-scheme $T$.
		
		In this way, we obtain a $\res_{S'/S} G'$-equivariant morphism
		\[\res_{S'/S}\Pp^{K'\mathrm{-st}}_{G'}\to \Pp^{\res_{S'/S} K'\mathrm{-st}}_{\res_{S'/S} G'}.\]
		One can prove that it is an isomorphism in a similar way to Lemma \ref{lem:quotient_vs_res} (recall the first paragraph for the representability of the left hand side).
		
		Finally, we prove the second isomorphism of (2). Since $\rtype(G',K')$ is finite \'etale over $S'$, so is $\res_{S'/S}\rtype(G',K')$ over $S$ (recall the argument in the third paragraph). In view of (1) and Theorem \ref{thm:orbit}, $\rtype(\res_{S'/S} G',\res_{S'/S} K')$ is finite \'etale over $S$. The assertion now follows from Lemma \ref{lem:quotient_vs_res}. This completes the proof.
	\end{proof}

	\subsection{Remarks on the symmetric setting}\label{sec:rem}
	
	Let $S$ be a $\ZZ\left[1/2\right]$-scheme, $G$ be a reductive group scheme over $S$, and $\theta$ be an involution of $G$. Let $K$ be an open and closed subgroup of $G^\theta$ (see \cite[Lemma 3.1.1]{hayashilinebdl}). 
	
	\begin{proposition}\label{prop:stable_vs_theta_stable}
		A parabolic subgroup is stable relative to $K$ if and only if it is $\theta$-stable\footnote{Although $\theta$-stable parabolic subalgebras are required to contain $\theta$-stable \textit{real} Cartan subalgebras in \cite{knappvogan} (see the beginning of \cite[Chapter V, section 1]{knappvogan} for example), we call a parabolic subgroup $P$ $\theta$-stable just if $\theta(P)\subset P$ (equivalently, $\theta(P)=P$) in this paper since currently we are only interested in $K$-conjugacy classes.}.
	\end{proposition}
	
	\begin{proof}
		The ``only if'' direction follows by seeing it locally. Let $P\subset G$ be a parabolic subgroup such that $\theta(P)=P$. We wish to show that $P$ is stable relative to $K$. We may assume that $S$ is affine since the assertion is local in the Zariski topology of $S$. Let $\Tor_{P/S}$ be the moduli scheme of maximal tori of $P$ (\cite[Corollaire 1.10]{sga3-12}).
		Then $\theta$ induces an involution on $\Tor_{P/S}$ which we denote by the same symbol $\theta$. The $S$-scheme $\Tor_{P/S}^\theta$ is smooth by \cite[Theorem 3.2.6]{conrad2014} and \cite[Lemma 3.1.1]{hayashilinebdl}. In particular, we may assume that $P$ admits a $\theta$-stable maximal torus $H$ by passing to an \'etale covering (use the argument of \cite[Lemma 3.1.8]{hayashilinebdl} and \cite[Proposition 4.2]{sga3-14}). Moreover, we may assume that $(G,H)$ is split and that $P=P_G(\mu)$ for a constant cocharacter $\mu$ of $H$ by \cite[Proposition 5.2.3]{conrad2014}. Work locally in the Zariski topology of $S$ to assume that $\theta(\mu)$ is also constant. We may replace $\mu$ by $\mu+\theta(\mu)$ to assume $\mu$ is $\theta$-invariant. In this case, $\mu$ factors through $K$.
	\end{proof}
	
	Recall that we introduced $\Ss^+(G,K^0)$ and $\Ss(G,K^0)$ below Lemma \ref{lem:relativesplitting} as a local parameter set of $K^0$-orbits in $\Pp^{K\mathrm{-st}}_{G}$ and its variant respectively. We plan to compute $\Ss^+(G,K^0)$ locally for the present $K$. Henceforth assume $K=K^0$. In view of Example \ref{ex:K=K^0} and Proposition \ref{prop:basefree}, we may assume $S$ connected. To sum up the setting, suppose that we are given a connected affine $\ZZ\left[1/2\right]$-scheme $S$. Let $G$ be a reductive group scheme over $S$, equipped with an involution $\theta$. Set $K=(G^\theta)^0$. Let $T$ be a maximal torus of $K$, and $H\coloneqq Z_G(T)$. Assume that $(G,H)$ and $(K,T)$ are split. We also fix a positive system $\Delta^+(K,T)$ of the root system $\Delta(K,T)$ of $(K,T)$.

	\begin{lemma}\label{lem:theta-stable_parabolic_subset}
		The set $\Ss(G,K)$ exactly consists of the $\theta$-stable parabolic subsets. Moreover, $\Phi\in\Ss(G,K)$ belongs to $\Ss^+(G,K)$ if and only if $\Delta^+(K,T)\subset\Phi|_T$ (see \eqref{eq:Phi|_T} for the notation of $\Phi|_T$).
	\end{lemma}
	
	\begin{proof}
		We prove that a parabolic subset $\Phi$ of $\Delta(G,H)$ is $\theta$-stable if and only if there exists a cocharacter $\mu\in X_\ast(T)$ such that $\Phi=\Phi_\mu$.
		The ``if'' direction is clear. Suppose that $\Phi$ is $\theta$-stable. Then one can find $\mu\in X_\ast(T)$ such that $\Phi=\Phi_\mu$ by \cite[Proof of Proposition 2.2.9]{conradpseudo} (if necessary, pass to a geometric fiber to assume that $S$ is the spectrum of an algebraically closed field of characteristic not two). It is easy to show $\Phi=\Phi_{\mu+\theta(\mu)}$. Hence we may assume $\theta(\mu)=\mu$. Then $\mu$ factors through $G^\theta$. Since $\GG_{m,S}$ is connected, $\mu$ factors through $K\cap H=Z_K(T)=T$ as desired.
		
		The latter statement follows from Proposition \ref{prop:kt} (2).
	\end{proof} 
	
	Therefore the problem is reduced to classifying $\theta$-stable parabolic subsets $\Phi$ of $\Delta(G,H)$ satisfying $\Delta^+(K,T)\subset\Phi|_T$. Write $W=W(G,H)(S)$. Since $S$ is connected, it is identified with the Weyl group $W$ of the root system $\Delta(G,H)$ of $(G,H)$. We denote the involution on $W$ induced from $\theta$ by the same symbol $\theta$.
	
	\begin{theorem}\label{thm:solutiontotheclassificationproblem}
		
		\begin{enumerate}
			\renewcommand{\labelenumi}{(\arabic{enumi})}
			\item Let $w\in W$ with $\theta(w)=w$. If $\Phi$ is a $\theta$-stable parabolic subset of $\Delta(G,H)$, so is $w\Phi$. In particular, if $\Phi$ is a $\theta$-stable positive system of $\Delta(G,H)$ then so is $w\Phi$.
			\item There exists a $\theta$-stable positive system $\Delta^+(G,H)$ of $\Delta(G,H)$ such that
			\[\Delta^+(G,H)|_T\cap\Delta(K,T)=\Delta^+(K,T)\]
			(recall \eqref{eq:Phi|_T} for the definition of $\Delta^+(G,H)|_T$).
			\item Fix a $\theta$-stable positive system $\Delta^+(G,H)$ in (2). Let $\SPS^\theta$ be the set of $\theta$-stable parabolic subsets of $\Delta(G,H)$ containing $\Delta^+(G,H)$.
			Set
			\[W^\theta_+=\{w\in W:~\theta(w)=w,~(w\Delta^+(G,H))|_T\cap \Delta(K,T)\subset \Delta^+(K,T)\}.\]
			Define an equivalence relation on $W^\theta_+\times\SPS^\theta$ by
			\[(w_1,\Phi_1)\sim (w_2,\Phi_2)\iff \Phi_1=\Phi_2~\mathrm{and}~w_1\Phi_1=w_2\Phi_2.\] 
			Then $(w,\Phi)\mapsto w\Phi$ determines a well-defined bijection
			\begin{equation}
				W^\theta_+\times\SPS^\theta/\sim\cong\Ss^+(G,K).\label{eq:bijectivegeneral}
			\end{equation}
			Moreover, it restricts to
			\begin{equation}
				W^\theta_+\cong\Ss^+_\emptyset(G,K);~w\mapsto w\Delta^+(G,H),\label{eq:bijectionregular}
			\end{equation}
			where $\Ss^+_\emptyset(G,K)$ be the subset of $\Ss^+(G,K)$ consisting of positive systems.
		\end{enumerate}
	\end{theorem}
	
	This is a reformulation of \cite[Proposition 9]{brionhelminck} when $S=\Spec\CC$. 
	
	\begin{proof}
		We may pass to a geometric fiber to assume $S$ to be the spectrum of an algebraically closed field of characteristic not two.
		
		Part (1) is straightforward:
		$\theta(w\Phi)=\theta(w)\theta(\Phi)=w\Phi$.
		
		Assume $\Delta^+(G,H)$ in (2) exists in this and next paragraphs. For $(w,\Phi)\in W^\theta_+\times\SPS^\theta$, (1) and the equality
		\[(w\Phi)|_T \cap \Delta(K,T)
		\supset (w\Delta^+(G,H))|_T \cap \Delta(K,T)
		=\Delta^+(K,T)
		\]
		(use Lemma \ref{lem:theta-stable_parabolic_subset} and Proposition \ref{prop:kt} (3) for the last equality) imply $w\Phi\in\Ss^+(G,K)$. It then follows by definition of the equivalence relation that the map \eqref{eq:bijectivegeneral} from left to right is well-defined.
		
		The induced equivalence relation on the subset
		$\{(w,\Phi)\in W^\theta_+\times\SPS^\theta:~\Phi=\Delta^+(G,H)\}$
		is trivial, and the map \eqref{eq:bijectivegeneral} restricts to an injective map
		$W^\theta_+\hookrightarrow \Ss^+_\emptyset(G,K)$ of \eqref{eq:bijectionregular} since $W$ acts simply on the set of positive systems of $\Delta(G,H)$. We wish to prove that \eqref{eq:bijectionregular} is surjective. Let $\Psi$ be a $\theta$-stable positive system of $\Delta(G,H)$ such that $\Delta^+(K,T)\subset \Psi|_T$. Then we have $\Delta^+(K,T)= \Psi|_T\cap \Delta(K,T)$ by Proposition \ref{prop:kt} (3). Choose $w\in W$ such that $w\Delta^+(G,H)=\Psi$. We wish to prove $w\in W^\theta_+$. The equality
		$\theta(w)\Delta^+(G,H)=\theta(w\Delta^+(G,H))=\theta\Psi=\Psi=w\Delta^+(G,H)$
		implies $w=\theta(w)$. The equality
		$(w\Delta^+(G,H))|_T \cap \Delta(K,T)
		=\Psi|_T \cap \Delta(K,T)
		=\Delta^+(K,T)$
		then shows $w\in W^\theta_+$.
		
		Finally, we verify (2) and the bijection \eqref{eq:bijectivegeneral}. Let $\Phi\in\Ss^+(G,K)$. We wish to find a $\theta$-stable positive system $\Psi$ contained in $\Phi$ such that $\Psi|_T\cap \Delta(K,T)=\Delta^+(K,T)$. In particular, put $\Phi=\Delta(G,H)$ to deduce (2). Let $P$ be the parabolic subgroup of $G$ corresponding to $\Phi$. Choose a $\theta$-stable Borel subgroup $B$ of $P$ and a maximal torus $T'$ of $B\cap K$ (see \cite[Theorem 7.5]{steinberg1968}). We may pass to a $K$-conjugation to assume that $T'=T$, and that the positive system of $\Delta(K,T)$ attached to $B\cap K$ coincides with $\Delta^+(K,T)$. In this case, $H$ is contained in $B$. In fact, a maximal torus $H'$ of $B$ containing $T$ which always exits by definition is contained in $H$ since $H'$ centralizes $T$. Since any maximal torus of $B$ is a maximal torus of $G$ (\cite[Proposition 4.2]{sga3-14}), $H'$ is equal to $H$. In other words, $H$ is a maximal torus of $G$ contained in $B$. The positive system attached to $(B,H)$ now satisfies the condition. 
		
		Fix a positive system $\Delta^+(G,H)$ in (2). Let $\Phi$ and $\Psi$ be as above. Let $w\in W^\theta_+$ be the element satisfying $w\Delta^+(G,H)=\Psi$. Then $w^{-1}\Phi$ is $\theta$-stable by $\theta(w^{-1})=w^{-1}$ and (1). We also have
		\[w^{-1}\Phi\supset w^{-1}\Psi=\Delta^+(G,H).\]
		Hence $(w,w^{-1}\Phi)\in W^\theta_+\times\SPS^\theta$ maps to $\Phi$. This proves the surjecitivity of \eqref{eq:bijectivegeneral}.
		
		In general, let $\Phi_1$ and $\Phi_2$ be parabolic subsets of $\Delta(G,H)$ containing $\Delta^+(G,H)$. We have $\Phi_1=\Phi_2$ if there exists $w\in W$ such that $w\Phi_1=\Phi_2$ (\cite[Corollary of Chap.\ VI, \S 1, no.\ 7, Proposition 21]{bourbakilie}). This shows the injectivity of the map \eqref{eq:bijectivegeneral} (rewrite the equality
		$w_1\Phi_1=w_2\Phi_2$
		as $w^{-1}_2w_1 \Phi_1=\Phi_2$).
	\end{proof}

	\begin{remark}
		Let $\Phi$ be a parabolic subset of $\Delta(G,H)$. Let $P$ be the parabolic subgroup of $G$ corresponding to $\Phi$, and $L$ be the Levi subgroup of $P$ containing $H$. Then for $w\in W$, $w\Phi=\Phi$ if and only if $w\in W(L,H)(S)$.
	\end{remark}
	
	\begin{remark}\label{rem:SPS^theta}
		Choose a $\theta$-stable positive system $\Delta^+(G,H)$. Then $\theta$ acts on the subset $\Pi(G,H)$ of simple roots. Recall also that
		$\Phi\mapsto \{\alpha\in \Pi(G,H):~-\alpha\in\Phi\}$
		determines a bijection between the set of parabolic subsets of $\Delta(G,H)$ containing $\Delta^+(G,H)$ and the power set of $\Pi(G,H)$. It is easy to show that a parabolic subsets of $\Delta(G,H)$ containing $\Delta^+(G,H)$ is $\theta$-stable if and only if the corresponding subset of $\Pi(G,H)$ is $\theta$-stable. In particular, there is a bijection between $\SPS^\theta$ and the set of $\theta$-invariant subsets of $\Pi(G,H)$.
	\end{remark}
	
	\begin{remark}\label{rem:gt}
		Let $x\in \rtype(G,K)(S)\cong\Ss^+(G,K)$ and $(w,\Phi)\in W^\theta_+\times\SPS^\theta$ be a corresponding pair. Identify $\type(G,K)(S)$ with the set of parabolic subsets of $\Delta(G,H)$ containing $\Delta^+(G,H)$. Then in view of the construction of the bijection \eqref{eq:bijectivegeneral}, we have $gt(x)=\Phi$.	
	\end{remark}
	
	\begin{remark}[{\cite[Corollary 5.100]{knapp2002}}]\label{rem:kt}
		For $\Phi\in\Ss(G,K)$,
		$2\rho^\vee(\Phi)\coloneqq \sum_{\alpha\in\Phi} \alpha^\vee=\sum_{\alpha\in\Phi\setminus(-\Phi)}\alpha^\vee$
		belongs to $X_\ast(T)$ and $\Phi=\Phi_{2\rho^\vee(\Phi)}$ (here think of $2\rho^\vee(\Phi)$ just as a single symbol).
		Hence the set of roots of $P_{\Phi,K}$ is given by
		\begin{equation}
			\{\beta\in\Delta(K,T):~
			\langle 2\rho^\vee(\Phi),\beta\rangle\geq 0\}.\label{eq:kt}
		\end{equation}
		In particular, for $x\in\rtype(G,K)(S)$, $kt(x)$ is given by \eqref{eq:kt} for $\Phi\in\Ss^+(G,K)$ corresponding to $x$ under the identification of $(\type K)(S)$ with the set of parabolic subsets of $\Delta(K,T)$ containing $\Delta^+(K,T)$.
		
		In the picture of \eqref{eq:bijectivegeneral}, $kt$ is computed as follows: pick a pair $(w,\Phi)\in W^\theta_+\times\SPS^\theta$ corresponding to $x$. Then we have
		$w\Phi=\Phi_{w (2\rho^\vee(\Phi))}$ and
		\[kt(x)=\{\beta\in\Delta(K,T):~
		\langle w(2\rho^\vee(\Phi)),\beta\rangle\geq 0\}.\]
		
	\end{remark}
	
	\subsection{Absolute descent}\label{sec:absolutedescent}
	
	Let $G$ be a reductive group scheme over a scheme $S$, and $K$ be a smooth closed subgroup scheme of $G$.
	
	\begin{proposition}\label{prop:cocharver}
		Let $S'\to S$ be a morphism of schemes, and $\sigma$ be an automorphism of $S'$ over $S$. Let $\mu$ be a cocharacter of $G\times_S S'$, and $w_\sigma\in G(S')$. If $\mu^\sigma=w_\sigma^{-1}\mu$, then we have $P_{G\times_S S'}(\mu)^\sigma=w_\sigma^{-1} P_{G\times_S S'}(\mu) w_\sigma$.
	\end{proposition}
	\begin{proof}
		Recall that for an $S'$-scheme $f:T'\to S'$, $(T')^\sigma$ can be identified with the $S'$-scheme $\sigma^{-1}\circ f:T'\to S'$. Under this identification, $\mu^\sigma$ is given by
		\[\begin{split}
			\GG_{m,S}(T')
			&=\GG_{m,S}((T')^{\sigma^{-1}})\\
			&=\GG_{m,S'}((T')^{\sigma^{-1}})\\
			&\xrightarrow{\mu_{(T')^{\sigma^{-1}}}}
			(G\times_S S')((T')^{\sigma^{-1}})\\
			&=G((T')^{\sigma^{-1}})\\
			&=G(T').
		\end{split}\]
		The assertion now follows from
		\[\begin{split}
			P_{G\times_S S'}(\mu)^\sigma(T')
			&=P_{G\times_S S'}(\mu)((T')^{\sigma^{-1}})\\
			&=\{g\in (G\times_{S} S')((T')^{\sigma^{-1}}):\ 
			\lim_{a\to 0}
			\mu_{(T')^{\sigma^{-1}}}(a)g\mu_{(T')^{\sigma^{-1}}}
			(a)^{-1}\ {\rm exists}\}\\
			&=\{g\in G(T'):\ \lim_{a\to 0}
			\mu^\sigma_{T'}(a)g\mu^\sigma_{T'}(a)^{-1}\ {\rm exists}\}\\
			&=\{g\in (G\times_S S')(T'):\ \lim_{a\to 0}
			\mu^\sigma_{T'}(a)g\mu^\sigma_{T'}(a)^{-1}\ {\rm exists}\}\\
			&=P_{G\times_S S'}(\mu^\sigma)(T'),
		\end{split}\]
		\[\begin{split}
			P_{G\times_S S'}(\mu^\sigma)(T')
			&=\{g\in G(T'):\ \lim_{a\to 0} w_\sigma^{-1} \mu_{T'}(a) w_\sigma g w_\sigma^{-1} \mu_{T'}(a)^{-1} w_\sigma
			\ {\rm exists}\}\\
			&=\{g\in G(T'):\ \lim_{a\to 0} \mu_{T'}(a) w_\sigma g w_\sigma^{-1} \mu_{T'}(a)^{-1}\ {\rm exists}\}\\
			&=w_\sigma^{-1} \{g\in G(T'):\ \lim_{a\to 0}\mu_{T'}(a)g\mu_{T'}(a)^{-1}\ {\rm exists}\}w_\sigma\\
			&=w_\sigma^{-1} P_{G\times_S S'}(\mu)(T') w_\sigma
		\end{split}\]
		for $S'$-schemes $T'$. 
	\end{proof}
	
	\begin{corollary}
		Let $S'\to S$ be a morphism of schemes, $\sigma$ be an automorphism of $S'$ over $S$, $\mu$ be a cocharacter to $K\times_S S'$, and $w_\sigma\in K(S')$ such that $\mu^\sigma=w_\sigma^{-1}\mu$. Then we have $rt(P_{G\times_S S'}(\mu)^\sigma)=rt(P_{G\times_S S'}(\mu))$.
	\end{corollary}
	
	The remarkable point of the above result is that the condition for $w_\sigma$ is independent of $G$.
	
	For applications, particularly to the standard models in \cite{hayashilinebdl}, we may assume the following conditions:
	
	\begin{enumerate}
		\renewcommand{\labelenumi}{(\roman{enumi})}
		\item $S'\to S$ is a Galois covering of connected schemes of Galois group $\Gamma=\ZZ/2\ZZ$,
		\item $K^0$ is reductive,
		\item $\pi_0(K)$ is finite \'etale over $S$,
		\item there exists a maximal torus $T\subset K$ such that $\mu$ factors through $T'\coloneqq T\times_S S'$, and
		\item $(K^0\times_S S',T')$ is split.
	\end{enumerate} 
	In fact, put $S'=\Spec\ZZ\left[1/2,\sqrt{-1}\right]$ and $S=\Spec\ZZ\left[1/2\right]$ for (i). Let $(G,K)$ be any of the examples in \cite[section 3.3]{hayashilinebdl}. In particular, (ii) holds. Let $T$ be as in \cite[section 3.4]{hayashilinebdl}. We saw in \cite{hayashilinebdl} that (iii)-(v) hold. Recall that $W(K^0\times_S S',T')(S')$ is isomorphic to the Weyl group of the root datum of $(K\times_S S',T\times_S S')$ and that the quotient map
	\[N_{K^0(S')}(T\times_S S')\to W(K^0\times_S S',T')(S')\]
	is surjective. Therefore if we think of finding an element $w\in N_{K^0(S')}(T\times_S S')\subset K^0(S')$ such that $\bar{\mu}=w\mu$, the problem is reduced to the following combinatorics on root data up to the computation of $\bar{\mu}$:
	\begin{problem}
		Let $(M,\Delta,M^\vee,\Delta^\vee)$ be a reduced root datum, $W(\Delta)$ be its Weyl group, and $\bar{\ }$ be an involution of $M^\vee$. Let $\mu\in M^\vee$. When does $w\in W(\Delta)$ such that $\bar{\mu}=w\mu$ exist?
	\end{problem}
	
	Henceforth we put $S'=\Spec\ZZ\left[1/2,\sqrt{-1}\right]$ and $S=\Spec\ZZ\left[1/2\right]$.
	In the case of the standard models, we have $\bar{\mu}=-\mu$ by the dual argument to \cite[section 3.2]{hayashilinebdl}. If the root datum is of type $B_n$, $C_n$, $D_{2n}$, or their direct sums, there always exists an element $w\in W(\Delta)$ such that $w\mu=-\mu$ since their Weyl groups contain $-1$. To see other cases, consider a positive integer $n$.
	\begin{example}[type $A_{n-1}$]
		The diagonal matrices of $\U(n)$ form a maximal torus $T$ of $\U(n)$ which is isomorphic to $\U(1)^n$. We can identify the action of the Weyl group on $X_\ast(T\times_S S')$ with the standard action of the $n$-th symmetric group $\mathfrak{S}_n$ on the free $\ZZ$-module $\ZZ^n$. Since the assertion for $a=(a_i)\in\ZZ^n$ that there exists an element $w\in\mathfrak{S}_n$ such that $wa=-a$ is closed under the action of $\mathfrak{S}_n$ on $\ZZ^n$, we may assume $a_1\geq a_2\geq \cdots \geq a_n$. For such $a$, there exists an element $w\in\mathfrak{S}_n$ such that $wa=-a$ if and only if $a_i+a_{n-i+1}=0$ for all $1\leq i\leq n$. If these conditions are satisfied, the element of the Weyl group represented by $K_n$ in \cite[section 1.5]{hayashilinebdl} satisfies $wa=-a$.
	\end{example}
	\begin{example}[type $D_{n}$]\label{ex:typeD_n}
		A maximal torus of $\SO(2n)$ is given by 
		\[T=\diag(\SO(2),\SO(2),\cdots,\SO(2)).\]
		We can identify the action of the Weyl group on $X_\ast(T\times_S S')$ with the standard action
		\[((\epsilon_i),\sigma)(a_i)=(\epsilon_i a_{\sigma^{-1}(i)})\]
		or equivalently
		\[\begin{array}{cc}
			((\epsilon_i),\sigma)e_j=\epsilon_{\sigma(j)} e_{\sigma(j)}&(1\leq j\leq n)
		\end{array}\]
		of $W=\{((\epsilon_i),\sigma)\in\{\pm 1\}^{n}\rtimes\mathfrak{S}_{n}:~\prod_{i=1}^n\epsilon_i=1\}$ on $\ZZ^n$ through
		\[\ZZ^n\ni(a_1,\cdots,a_n)\mapsto \diag(\mu_2^{a_1},\mu_2^{a_2},\cdots,\mu_2^{a_n}),\]
		where $e_i$ ($1\leq i\leq n$) denotes the standard basis of $\ZZ^n$. For $a=(a_i)\in\ZZ^n$, there exists an element $((\epsilon_i),\sigma)\in W$ such that $((\epsilon_i),\sigma)a=-a$ if and only if either $n$ is even or $\prod_{i=1}^n a_i=0$. In fact, if $n$ is even, we have
		\[(-1,-1,-1,\cdots,-1)\in
		\{(\epsilon_i)\in\{\pm 1\}^n:\ \prod_{i=1}^n\epsilon_i=1\}\subset W\]
		and $(-1,-1,-1,\cdots,-1)a=-a$. Note that $(-1,-1,-1,\cdots,-1)$ is represented by 
		\[w=\diag(w_2,w_2,\cdots,w_2)\in\SO(2n,\ZZ),\]
		where $w_2=\diag(1,-1)$.
		If $n$ is odd, and there exists $i\in\{1,2,\cdots,n\}$ such that $a_i=0$ then
		\[\begin{array}{cc}
			(-1,-1,\cdots,-1,\overset{\overset{i}{\vee}}{1},-1,\cdots,-1)\in W,
			&(-1,-1,\cdots,-1,\overset{\overset{i}{\vee}}{1},-1,\cdots,-1)a=-a.
		\end{array}\]
		Note that $(-1,-1,\cdots,-1,\overset{\overset{i}{\vee}}{1},-1,\cdots,-1)$ is represented by
		\[\diag(\overbrace{w_2,\cdots,w_2}^{i-1},1,1,\overbrace{w_2,\cdots,w_2}^{n-i})\in\SO(2n,\ZZ).\]
		Conversely, suppose that there exists $((\epsilon_i),\sigma)\in W$ such that $((\epsilon_i),\sigma)a=-a$. Then we have $\epsilon_i a_{\sigma^{-1}(i)}=-a_i$ for $1\leq i\leq n$. If $\prod_{i=1}^n a_i\neq 0$, the equality $1=\prod_{i=1}^n\epsilon_i=\prod_{i=1}^n\left(-\frac{a_i}{a_{\sigma^{-1}(i)}}\right)=(-1)^n$ implies that $n$ is even.
	\end{example}
	
	\subsection{Classification of relative parabolic types: Standard setting}\label{sec:halfintegralformofrelativetype}
	
	Put
	\[\begin{array}{cc}
		S=\Spec\mathbb{Z}\left[1/2\right]
		&S'=\Spec\mathbb{Z}\left[1/2,\sqrt{-1}\right].
	\end{array}\]
	Recall that we introduced the standard models $(G,\theta,K)$ of classical (connected) Lie groups, their Cartan involutions, and the corresponding maximal connected compact subgroups in \cite{hayashilinebdl}. The aim of this section is to study $\rtype(G,K)$ for these models. Recall that we also constructed maximal tori $T\subset K$ and $H=Z_G(T)\subset G$ in \cite[sections 3.3 and 3.4]{hayashilinebdl}. We know $(G\times_S S',H\times_S S')$ and $(K\times_S S',T\times_S S')$ are split. We gave $\theta$-stable positive systems $\Delta^+(G\times_S S',H\times_S S')$ in \cite[section 4.1]{hayashilinebdl}. Hence we can determine $\rtype(G,K)$ by computing $\mathscr{S}(G\times_SS',K\times_S S')$ and the Galois actions on them. To achieve this, we are based on Theorem \ref{thm:solutiontotheclassificationproblem}. That is, we give combinatorial descriptions of $W^\theta_+$ and $\SPS^\theta$ in Theorem \ref{thm:solutiontotheclassificationproblem} for $(G\times_S S',K\times_S S')$ to determine $\mathscr{S}(G\times_SS',K\times_S S')$. Then we compute the conjugate actions under the bijection \eqref{eq:bijectivegeneral}.
	
	We also study $\rtype(G,K)$ for $(G,K)=(\GL_{n},\Oo(n)),(\SO(p,q),\mathrm{S}(\Oo(p)\times\Oo(q)))$. Moreover, for each element $\Phi\in\Ss^+(G\times_S S',K^0\times_S S')$, we judge the $K\times_S S'$-orbit attached to $\Phi$ consists of either a single $K^0\times_S S'$-orbit or two $K^0\times_SS'$-orbits.
	
	To achieve these, we will identify $X^\ast(H\times_S S')$ and $X^\ast(T\times_SS')$ with free abelian groups in a natural way via \cite[sections 3.2 and 3.4]{hayashilinebdl}. We will denote the standard basis by $\{e_i\}$.
	
	As for computation of $\theta$, take the base change to the complex numbers to obtain $\theta\lambda=-\bar{\lambda}$ for $\lambda\in X^\ast(H\times_S S')$. The conjugate action on $X^\ast(H\times_S S')$ was computed explicitly in \cite[section 4.1]{hayashilinebdl}.
	
	For each standard model $(G,K)$, let $\Delta^+(G\times_S S',H\times_S S')$ be the positive system given in \cite[section 4.1]{hayashilinebdl}. Set
	\[\Delta^+(K\times_S S',T\times_S S')
	=\Delta^+(G\times_S S',H\times_S S')|_{T\times_S S'}\cap \Delta(K\times_S S',T\times_S S').\]
	We remark that $\Delta^+(K\times_S S',T\times_S S')$ is the direct sum of the positive systems in \cite[Appendix C]{knapp2002}. In fact, recall that $\Delta^+(K\times_S S',T\times_S S')$ was defined by the ordered bases of
	\[X_\ast(T\times_S S')\subset X_\ast(H\times_S S')\]
	given in \cite[section 4.1]{hayashilinebdl}. One can find that they are the standard ordered basis of the free abelian group $X_\ast(T\times_S S')$ (use \cite[sections 3.2 and 3.3]{hayashilinebdl} to identify $X_\ast(T\times_S S')\subset X_\ast(H\times_S S')$ with embeddings of free abelian groups).
	
	In this section, let $p,q,n$ be nonnegative integers with $p+q=n$.
	
	\subsubsection{Description of $W^\theta_+$}
	
	In this section, we try to give a simple combinatorial description of $W^\theta_+$. This is equivalent to classifying the closed $K(\CC)$-orbits in the complex flag variety of $G$. This was done by Matsuki and Oshima in \cite[\S 4]{matsukioshima1990}. We follow up their work to determine the set $W^\theta_+$ explicitly.
	
	Let $e$ denote the unit element of $W$. Before we demonstrate explicit computations for each of the standard models, let us explain our general scheme to describe $W^\theta_+$:
	
	\begin{strategy}\label{strategy}
		\begin{itemize}
			\item According to \cite{matsukioshima1990}, $W^\theta_+$ is trivial if $G=\GL_{2n-1},\U^\ast(2n)$ ($n\geq 1$). In the rest, assume
			$G\neq \GL_{2n-1},\U^\ast(2n)$.
			\item For $G=\GL_{2n}$ (resp.~$\SO(2p+1,2q+1)$) with $n\geq 1$ (resp,~$p,q\geq 0$), we find $w_\theta\in \lieS_{2n}$ (resp.~$w_\theta\in\{\pm 1\}^{n+1}\rtimes\lieS_{n+1}$) such that 
			\begin{equation}
				\theta\lambda=-w_\theta\lambda\label{eq:w_theta}
			\end{equation}
			for all $\lambda \in X^\ast(H\times_S S')$.
			For this, use the equality $\theta\lambda=-\bar{\lambda}$ and the formulas of $\bar{\lambda}$ in \cite[section 4.1]{hayashilinebdl}. It follows by definitions that for $w\in W$ and $\lambda\in X^\ast(H\times_S S')$, we have $\theta(w)\lambda=\theta(w (\theta\lambda))$. Hence for $w\in W$, $w$ belongs to $W^\theta$ if and only if 
			\begin{equation}
				ww_\theta=w_\theta w.\label{eq:w,w_theta}
			\end{equation}
			\item According to \cite[Table 1']{matsukioshima1990}, we have $|W^\theta_+|=2$ if $G=\GL_{2n}$. We find candidates of the elements of $W^\theta_+$ directly in our choice of the positive system, and prove that they belong to $W^\theta_+$. Henceforth we assume $G\neq \GL_{2n}$.
			\item For $G=\SO(2p+1,2q+1)$, we characterize the equality \eqref{eq:w,w_theta} on $w$ by a simple combinatorial condition. In the remaining case, $\theta$ acts trivially on $W$ since we have $H=T$. Therefore we finish computation of $W^\theta$ (except $G=\GL_n,\U^\ast(2n)$).
			\item We compute
			$(w\Delta^+(G\times_S S',H\times_S S'))|_{T\times_S S'}$ for $w\in W^\theta$.
			\item In view of Proposition \ref{prop:kt} (3), the following conditions for $w\in W^\theta$ are equivalent:
			\begin{enumerate}
				\item[(a)] $w$ belongs to $W^\theta_+$;
				\item[(b)] we have $\Delta^+(K\times_S S',T\times_SS')=(w\Delta^+(G\times_S S',H\times_S S'))|_{T\times_S S'}\cap \Delta(K\times_S S',T\times_SS')$;
				\item[(c)] $\alpha\in(w\Delta^+(G\times_S S',H\times_S S'))|_{T\times_S S'}$ for all simple roots $\alpha$ of $\Delta^+(K\times_S S',T\times_SS')$.
			\end{enumerate}
			We give a combinatorial criterion for (c) to determine $W^\theta_+$ by (a). To find a a necessary condition for (c) in an efficient way, we will eventually test
			\[\alpha\in (w\Delta^+(G\times_S S',H\times_S S'))|_{T\times_S S'}\cap \Delta(K\times_S S',T\times_SS')\]
			for other positive roots $\alpha\in \Delta(K\times_S S',T\times_SS')$ (see (b)).
		\end{itemize}
	\end{strategy}

	\begin{proposition}\label{prop:orbit_GL_n}
		Put $(G,K)=(\GL_n,\SO(n))$ with $n\geq 1$. Identify $W$ with $\mathfrak{S}_n$. Then
		we have
		\[W^\theta_+=\begin{cases}
			\{e,~(n-1~n)\}&(n~\mathrm{is~even})\\
			\{e\}&(n~\mathrm{is~odd}).
		\end{cases}\]
	\end{proposition}
	
	\begin{proof}
		See \cite[Table 1']{matsukioshima1990} for the odd case. Henceforth we may assume $n$ to be even. In this case, the assertion is essentially due to \cite{matsukioshima1990}. Here we only show that the $(n-1~n)$ is the nontrivial element of $W^\theta_+$ in our choice of positive system.
		
		Set $w_\theta=(1~2)(3~4)\cdots(n-1~n)$.
		This satisfies the condition \eqref{eq:w_theta} by \cite[Example 4.1.2]{hayashilinebdl}. It is evident by this expression that $(n-1~n)$ belongs to $W^\theta$. For $1\leq i<j\leq \frac{n}{2}$, we have
		\[((n-1~n)(e_{2i-1}-e_{2j-1}))|_{T\times_S S'}=\left\{\begin{array}{cc}
			e_i-e_j&(j\neq \frac{n}{2})\\
			e_i+e_j&(j= \frac{n}{2}),
		\end{array}\right.
		\]
		\[((n-1~n)(e_{2i-1}-e_{2j}))|_{T\times_S S'}=\left\{\begin{array}{cc}
			e_i+e_j&(j\neq \frac{n}{2})\\
			e_i-e_j&(j= \frac{n}{2})
		\end{array}\right.
		\]
		since $e_{2i-1}|_{T\times_SS'}=e_i$ and $e_{2i}|_{T\times_SS'}=-e_i$
		for $1\leq i\leq \frac{n}{2}$ by the description of the embedding $T\hookrightarrow H$ (cf.~\cite[Example 4.1.2]{hayashilinebdl}).
	\end{proof}

	\begin{proposition}[{\cite[Table 1']{matsukioshima1990}}]
		Put $(G,K)=(\U^\ast(2n),\Sp(n))$ with $n\geq 1$. Then we have
		$W^\theta_+=\{e\}$.
	\end{proposition}

	\begin{proposition}\label{prop:w^theta_+_SO(2p+1,2q+1)}
		Put $(G,K)=(\SO(2p+1,2q+1),\SO(2p+1)\times\SO(2q+1))$.
		Identify $W$ with
		$\{((\epsilon_i),\sigma)\in\{\pm 1\}^{n+1}\rtimes\mathfrak{S}_{n+1}:~\prod_{i=1}^{n+1} \epsilon_i=1\}$.
		Then $w=((\epsilon_i),\sigma)\in W$ belongs to $W^\theta_+$ if and only if
		\begin{enumerate}
			\renewcommand{\labelenumi}{(\roman{enumi})}
			\item $\sigma(p+1)=p+1$.
			\item We have $\epsilon_i=1$ for $1\leq i\leq n+1$.
			\item We have $\sigma^{-1}(i)<\sigma^{-1}(i+1)$ for $1\leq i\leq n$ with $i\neq p, p+1$.
		\end{enumerate}
	\end{proposition}

	\begin{proof}
		Set
		\[w_\theta=(\overbrace{1,1,\ldots,1}^{p},-1,\overbrace{1,1,\cdots,1}^{q})
		\]
		This satisfies the condition \eqref{eq:w_theta} by \cite[Examples 4.1.10 and 4.1.11]{hayashilinebdl}.
		
		Let $w=((\epsilon_i),\sigma)\in W$. Then we have
		\[\begin{array}{cc}
			ww_\theta=
			((\epsilon_1,\epsilon_2,\ldots,-\epsilon_{\sigma(p+1)},\ldots,\epsilon_n),\sigma),
			&w_\theta w=
			((\epsilon_1,\epsilon_2,\ldots,-\epsilon_{p+1},\ldots,\epsilon_n),\sigma).
		\end{array}\]
		In particular, $w$ belongs to $W^\theta$ if and only if $\sigma(p+1)=p+1$. Henceforth we assume these equivalent conditions.
		
		Observe that
		\begin{flalign*}
			&w\Delta^+(G\times_S S',H\times_S S')\\
			&=\{\epsilon_{\sigma(i)}e_{\sigma(i)}\pm \epsilon_{\sigma(j)}e_{\sigma(j)}:
			~1\leq i<j\leq n+1,i\neq p+1\}\\
			&\cup\{\pm \sigma(p+1)e_{\sigma(p+1)}
			+\epsilon_{\sigma(j)}e_{\sigma(j)}:~p+1<j\leq n+1\}\\
			&=\{\epsilon_{\sigma(i)}e_{\sigma(i)}\pm e_{\sigma(j)}:
			~1\leq i<j\leq n+1,i\neq p+1\}
			\cup\{\pm e_{\sigma(p+1)} +\epsilon_{\sigma(j)}e_{\sigma(j)}:~p+1<j\leq n+1\}\\
			&=\{\epsilon_{\sigma(i)}e_{\sigma(i)}\pm e_{\sigma(j)}:
			~1\leq i<j\leq n+1,i\neq p+1\}
			\cup\{\pm e_{p+1} +\epsilon_{\sigma(j)}e_{\sigma(j)}:~p+1<j\leq n+1\},
		\end{flalign*}
		where the final equality follows from $\sigma(p+1)=p+1$.
		We note that
		\[e_i|_{T\times_S S'}=\left\{\begin{array}{cc}
			e_i&(1\leq i\leq p)\\
			0&(i=p+1)\\
			e_{i-1}&(p+2\leq i\leq n+1)
		\end{array}\right.\]
		for $1\leq i\leq n+1$ (see \cite[Example 3.4.9]{hayashilinebdl}). It is then easy to show that the fibers of the restriction map
		\[\Delta(G\times_S S',H\times_S S')\to X^\ast(T\times_S S');~
		\alpha\mapsto \alpha|_{T\times_S S'}\]
		at the simple roots and short roots $\beta\in \Delta^+(K\times_S S',T\times_S S')$ are given by
		\[\begin{array}{cc}
			\{e_i-e_{i+1}\}&(\beta=e_i-e_{i+1},~1\leq i\leq p-1)\\
			\{e_i\pm e_{p+1}\}&(\beta=e_i,~1\leq i\leq p)\\
			\{e_{i+1}-e_{i+2}\}&(\beta=e_i-e_{i+1},~p+1\leq i\leq n-1)\\
			\{\pm e_{p+1}+e_{i+1}\}&(\beta=e_i,~p+1\leq i\leq n).
		\end{array}\]
		
		Define a function $u:\{1,2,\ldots,n\}\to \{0,1\}$ by
		\[u(i)=\left\{\begin{array}{cc}
			i&(1\leq i\leq p)\\
			1+1&(p+1\leq i\leq n).
		\end{array}\right.\]
		Let $1\leq i\leq n$.	 
		Since
		\[w\Delta^+(G\times_S S',H\times_S S')
		\supset\{\epsilon_{u(i)} e_{u(i)}\pm e_{p+1}\},\]
		the following conditions are equivalent:
		\begin{enumerate}
			\renewcommand{\labelenumi}{(\roman{enumi})}
			\item $(w\Delta^+(G\times_S S',H\times_S S'))|_{T\times_SS'}$ contains $e_i$;
			\item $\{\epsilon_{u(i)} e_{u(i)}\pm e_{p+1}\}\cap \{e_{u(i)}\pm e_{p+1}\}\neq \emptyset$;
			\item $\epsilon_{u(i)}=1$.
		\end{enumerate}
		In particular, $(w\Delta^+(G\times_S S',H\times_S S'))|_{T\times_SS'}$ contains $e_i$ for $1\leq i\leq n$ if and only if $\epsilon_i=1$ for $1\leq i\leq n+1$ with $i\neq p+1$. 
		
		Henceforth we may assume these equivalent conditions. Since $\prod_{i=1}^{n+1}\epsilon_i=1$, they are also equivalent to the condition that $\epsilon_i=1$ for $1\leq i\leq n+1$. In this case, for $1\leq i\leq n$ with $i\neq p$, $e_i-e_{i+1}$ belongs to $(w\Delta^+(G\times_S S',H\times_S S'))|_{T\times_SS'}$ if and only if $\sigma^{-1}(u(i))<\sigma^{-1}(u(i+1))$. Run through all $i$ to deduce the assertion.
	\end{proof}

	In the rest, we work with the case $H=T$. Then we have $W^\theta=W$ since $\theta$ acts trivially on $X^\ast(H\times_S S')$. It remains to work with the last two items in Strategy \ref{strategy}.
	
	\begin{proposition}
		Put $(G,K)=(\U(p,q),\U(p)\times\U(q))$ with $n=p+q\geq 1$. Identify $W$ with $\mathfrak{S}_n$. Then $w\in W$ belongs to $W^\theta_+$ if and only if $w^{-1}(i)<w^{-1}(i+1)$ for $1\leq i\leq n-1$ with $i\neq p$.
	\end{proposition}
	
	\begin{proof}
		Let $w\in W$. The assertion follows since for $1\leq i\leq n-1$ with $i\neq p$, $e_i-e_{i+1}\in w\Delta^+(G,H)$ if and only if $w^{-1}(i)<w^{-1}(i+1)$ by
		\[w\Delta^+(G\times_S S',H\times_S S')
		=\{e_{w(i)}-e_{w(j)}\in X^\ast(H\times_S S'):~1\leq i<j\leq n\}.\]
	\end{proof}
	
	\begin{proposition}\label{prop:w^theta_+_SO(2p,2q+1)}
		Put $(G,K)=(\SO(2p,2q+1),\SO(2p)\times\SO(2q+1))$ with $n=p+q\geq 1$. Identify $W$ with $\{\pm 1\}^{n}\rtimes\mathfrak{S}_{n}$. Then $w=((\epsilon_i),\sigma)\in W$ belongs to $W^\theta_+$ if and only if
		\begin{enumerate}
			\renewcommand{\labelenumi}{(\roman{enumi})}
			\item $\sigma^{-1}(i)<\sigma^{-1}(i+1)$ for $1\leq i\leq n-1$ with $i\neq p$, and
			\item $\epsilon_i=1$ for $i\neq p$.
		\end{enumerate}
	\end{proposition}
	
	\begin{proof}
		Let $((\epsilon_i),\sigma)\in W$. Then we have
		\begin{flalign*}
			&w\Delta^+(G\times_S S',H\times_S S')\\
			&=\{\epsilon_{\sigma(i)} e_{\sigma(i)}\pm \epsilon_{\sigma(j)}e_{\sigma(j)}\in X^\ast(H\times_S S'):~1\leq i<j\leq n\}
			\cup\{\epsilon_{\sigma(i)} e_{\sigma(i)}\in X^\ast(H\times_S S'):~1\leq i\leq n\}\\
			&=\{\epsilon_i e_{\sigma(i)}\pm e_{\sigma(j)}\in X^\ast(H\times_S S'):~1\leq i<j\leq n\}
			\cup\{\epsilon_i e_i\in X^\ast(H\times_S S'):~1\leq i\leq n\}.
		\end{flalign*}
		Hence for $1\leq i\leq n-1$ with $i\neq p$, $w\Delta^+(G\times_S S',H\times_S S')$ contains $e_i\pm e_{i+1}$ if and only if $\sigma^{-1}(i)<\sigma^{-1}(i+1)$ and $\epsilon_i=1$. In particular,
		\[w\Delta^+(G\times_S S',H\times_S S')\supset\Delta^+(K\times_S S',H\times_S S')\]
		if and only if (i) holds and $\epsilon_i=1$ for $i\neq p,n$.
		
		The remaining simple root of $\Delta^+(K\times_S S',T\times_S S')$ is $e_n$. We have $e_n\in w\Delta^+(G\times_S S',H\times_S S')$ if and only if $\epsilon_n=1$. This completes the proof.
	\end{proof}
	
	\begin{proposition}\label{prop:w^theta_+_Sp_n}
		Put $(G,K)=(\Sp_n,\U(n))$ with $n\geq 1$. Identify $W$ with $\{\pm 1\}^{n}\rtimes\mathfrak{S}_{n}$. Then $w=((\epsilon_i),\sigma)\in W$ belongs to $W^\theta_+$ if and only if the following conditions are satisfied:
		\begin{enumerate}
			\renewcommand{\labelenumi}{(\roman{enumi})}
			\item $\sigma^{-1}(i)<\sigma^{-1}(i+1)$ for $1\leq i\leq \sigma(n)-1$, and that $\sigma^{-1}(i)>\sigma^{-1}(i+1)$ for $\sigma(n)\leq i\leq n-1$.
			\item We have
			\[\epsilon_{i}=\left\{\begin{array}{cc}
				1&(1\leq i\leq \sigma(n)-1)\\
				-1&(\sigma(n)+1\leq i\leq n).
			\end{array}\right.\]
		\end{enumerate}
	\end{proposition}

	\begin{proof}
		Let $((\epsilon_i),\sigma)\in W$. Then we have
		\begin{flalign*}
			&w\Delta^+(G\times_S S',H\times_S S')\\
			&=\{\epsilon_{\sigma(i)} e_{\sigma(i)}\pm \epsilon_{\sigma(j)} e_{\sigma(j)}\in X^\ast(H\times_S S'):~1\leq i<j\leq n\}
			\cup\{2\epsilon_{\sigma(i)} e_{\sigma(i)}\in X^\ast(H\times_S S'):~1\leq i\leq n\}\\
			&=\{\epsilon_{\sigma(i)} e_{\sigma(i)}\pm e_{\sigma(j)}\in X^\ast(H\times_S S'):~1\leq i<j\leq n\}
			\cup\{2\epsilon_i e_i\in X^\ast(H\times_S S'):~1\leq i\leq n\}.
		\end{flalign*}
		
		Let $1\leq i\leq n-1$. Then:
		\begin{description}
			\item[Case 1] Suppose $\sigma^{-1}(i)<\sigma^{-1}(i+1)$. Then $e_i-e_{i+1}\in w\Delta^+(G\times_S S',H\times_S S')$ if and only if $\epsilon_i=1$.
			\item[Case 2] Suppose $\sigma^{-1}(i)>\sigma^{-1}(i+1)$. Then $e_i-e_{i+1}\in w\Delta^+(G\times_S S',H\times_S S')$ if and only if $\epsilon_{i+1}=-1$.
		\end{description}
		
		In particular, $w$ belongs to $W^\theta_+$ if and only if exactly one of the following conditions holds for each $1\leq i\leq n-1$:
		\begin{enumerate}
			\renewcommand{\labelenumi}{(\Roman{enumi})}
			\item $\sigma^{-1}(i)<\sigma^{-1}(i+1)$ and $\epsilon_i=1$;
			\item $\sigma^{-1}(i)>\sigma^{-1}(i+1)$ and $\epsilon_{i+1}=-1$.
		\end{enumerate}
		In particular, the ``if'' direction follows. We next prove the converse. Assume $w\in W^\theta_+$. Condition (ii) is immediate from (i) and the conditions (I), (II). We verify $\sigma^{-1}(i)<\sigma^{-1}(i+1)$ for $1\leq i\leq \sigma(n)-1$ by descending induction on $i$. We clearly have $\sigma^{-1}(i)<\sigma^{-1}(i+1)$ for $i=\sigma^{-1}(n)-1$. Let $k\in\{1,2,\ldots,n\}$. Suppose that $\sigma^{-1}(k)<\sigma^{-1}(k+1)$ for an integer $2\leq k\leq \sigma(n)-1$. Then we have $\epsilon_k=1$ by (I). Then in view of (II), the inequality $\sigma^{-1}(k-1)<\sigma^{-1}(k)$ must hold, and the induction proceeds. Similarly, the ascending induction on $i$ implies that $\sigma^{-1}(i)>\sigma^{-1}(i+1)$ for $\sigma(n)\leq i\leq n-1$. This shows (i), and the proof is completed.
	\end{proof}
	
	\begin{proposition}
		Put $(G,K)=(\Sp(p,q),\Sp(p)\times\Sp(q))$ with $n=p+q\geq 1$. Identify $W$ with $\{\pm 1\}^{n}\rtimes\mathfrak{S}_{n}$. Then $w=((\epsilon_i),\sigma)\in W$ belongs to $W^\theta_+$ if and only if
		\begin{enumerate}
			\renewcommand{\labelenumi}{(\roman{enumi})}
			\item $\sigma^{-1}(i)<\sigma^{-1}(i+1)$ for $1\leq i\leq n-1$ with $i\neq p$, and
			\item $\epsilon_i=1$ for $1\leq i\leq n$.
		\end{enumerate}
	\end{proposition}
	
	\begin{proof}
		Let $((\epsilon_i),\sigma)\in W$. Then we have
		\begin{flalign*}
			&w\Delta^+(G\times_S S',H\times_S S')\\
			&=\{\epsilon_{\sigma(i)} e_{\sigma(i)}\pm \epsilon_{\sigma(j)} e_{\sigma(j)}\in X^\ast(H\times_S S'):~1\leq i<j\leq n\}
			\cup\{2\epsilon_{\sigma(i)} e_{\sigma(i)}\in X^\ast(H\times_S S'):~1\leq i\leq n\}\\
			&=\{\epsilon_i e_{\sigma(i)}\pm e_{\sigma(j)}\in X^\ast(H\times_S S'):~1\leq i<j\leq n\}
			\cup\{2\epsilon_i e_i\in X^\ast(H\times_S S'):~1\leq i\leq n\}.
		\end{flalign*}
		Hence $\{2e_i\in X^\ast(H\times_S S'):~1\leq i\leq n\}\subset w\Delta^+(G\times_S S',H\times_S S')$ if and only if $\epsilon_i=1$ for all $i$. In the rest, assume these equivalent conditions.
		
		For $1\leq i\leq n-1$ with $i\neq p$, $e_i-e_{i+1}\subset w\Delta^+(G\times_S S',H\times_S S')$ if and only if $\sigma^{-1}(i)<\sigma^{-1}(i+1)$. This shows the equivalence.
	\end{proof}
	
	A similar argument to the case of $\SO(2p,2q+1)$ implies:
	
	\begin{proposition}\label{prop:w^theta_+_SO(2p,2q)}
		Put $(G,K)=(\SO(2p,2q),\SO(2p)\times\SO(2q))$ with $n=p+q\geq 1$. Identify $W$ with
		$\{((\epsilon_i),\sigma)\in\{\pm 1\}^{n}\rtimes\mathfrak{S}_{n}:~\prod_{i=1}^n \epsilon_i=1\}$.
		Then $((\epsilon_i),\sigma)\in W$ belongs to $W^\theta_+$ if and only if 
		\begin{enumerate}
			\renewcommand{\labelenumi}{(\roman{enumi})}
			\item $\sigma^{-1}(i)<\sigma^{-1}(i+1)$ for $1\leq i\leq n-1$ with $i\neq p$, and
			\item $\epsilon_i=1$ for $i\neq p,n$.
		\end{enumerate}
		We remark that if these equivalent conditions are satisfied, we obtain $\epsilon_p=\epsilon_n$ from the condition
		$\prod_{i=1}^n \epsilon_i=1$.
	\end{proposition}
	
	A similar argument to the case of $\Sp_n$ implies:
	
	\begin{proposition}\label{prop:w^theta_+_SO^ast(2n)}
		Put $(G,K)=(\SO^\ast(2n),\U(n))$ with $n\geq 1$. Identify $W$ with
		\[\{((\epsilon_i),\sigma)\in\{\pm 1\}^{n}\rtimes\mathfrak{S}_{n}:~\prod_{i=1}^n \epsilon_i=1\}.\]
		Then $w=((\epsilon_i),\sigma)\in W$ belongs to $W^\theta_+$ if and only if the following conditions are satisfied:
		\begin{enumerate}
			\renewcommand{\labelenumi}{(\roman{enumi})}
			\item $\sigma^{-1}(i)<\sigma^{-1}(i+1)$ for $1\leq i\leq \sigma(n)-1$, and that $\sigma^{-1}(i)>\sigma^{-1}(i+1)$ for $\sigma(n)\leq i\leq n-1$.
			\item We have
			\[\epsilon_{i}=\left\{\begin{array}{cc}
				1&(1\leq i\leq \sigma(n)-1)\\
				-1&(\sigma(n)+1\leq i\leq n).
			\end{array}\right.\]
		\end{enumerate}
		We remark that $\epsilon_{\sigma(n)}$ is determined by (ii) since $\prod_{i=1}^{n}\epsilon_i=1$.
	\end{proposition}
	
	\begin{remark}
		Following \cite{matsukioshima1990}, one can parameterize $W^\theta_+$ in a symbolic way: In fact, $W^\theta_+$ is bijective to the set of sequences of $p$ copies of $+$, $q$ copies of $-$, and $\circ$ such that $\circ$ stands at $(p+1)$-th term. If we are given a such sequence, the corresponding element $((\epsilon_i),\sigma)$ is given as follows: put $\epsilon_i=1$ for all $i$. We write the numbers $1,2,\ldots,p$ (resp.~$p+2,p+3,\ldots,n+1$) below $+$ (resp.~$-$) in order. Write $p+1$ below $\circ$.
		If we denote the resulting sequence of numbers by $(a_i)$, $\sigma$ is defined by $\sigma(i)=a_i$. For example, if the sequence is $+-++-\circ+-$ ($p=4$, $q=3$), we obtain
		\[\begin{array}{ccccccccc}
			+&-&+&+&\circ&-&+&-\\
			1&6&2&3&5&7&4&8
		\end{array}\]
		\[\sigma=\left(\begin{array}{cccccccc}
			1 & 2 & 3 & 4 & 5 & 6 & 7 & 8 \\
			1 & 6 & 2 & 3 & 5 & 7 & 4 & 8
		\end{array}\right).\]
		Similarly, $W^\theta_+$ is bijective to the set of sequences of $p$ copies of $+$ and $q$ copies of $-$ for $G=\U(p,q),\Sp(p,q)$. For $G=\SO(2p,2q+1),\SO(2p,2q)$, one can parameterize $(\sigma,(\epsilon_i))$ by the pairs of the same sequences and elements of $\{\pm1\}$. 
		
		For $G=\Sp_n$, $W^\theta_+$ is bijective to the set of pairs of $n$-term sequences $(\delta_i)$ of $\pm$ with $\delta_n=+$ and elements $\epsilon\in\{\pm 1\}$. In fact, define a function $\delta:\{1,2,\ldots,n\}\to \{\pm\}$ by $i\mapsto \delta_i$. Write
		$\delta^{-1}(+)=\{a_1,a_2,\ldots,a_{i_0}\}$ and $\delta^{-1}(-)=\{a_{i_0+1},\ldots,a_n\}$ with
		\[a_1<a_2<\cdots<a_{i_0}=n>a_{i_0+1}>\cdots>a_n.\]
		Then we set $\sigma^{-1}(i)=a_i$ for $1\leq i\leq n$ and 
		\[\epsilon_i=\begin{cases}
			1&(1\leq i<|\delta^{-1}(+)|)\\
			\epsilon&(i=|\delta^{-1}(+)|)\\
			-1&(i>|\delta^{-1}(+)|).
		\end{cases}\]
		For $G=\SO^\ast(2n)$, $W^\theta_+$ is bijective to the set of $n$-term sequences $(\delta_i)$ of $\pm$ with $\delta_n=+$. We can define $\sigma$ in a similar way to the case of $G=\Sp_n$. We set
		\[\epsilon_i=\begin{cases}
			1&(1\leq i\leq |\delta^{-1}(+)|-1)\\
			(-1)^{|\delta^{-1}(-)|}&(i=|\delta^{-1}(+)|)\\
			-1&(|\delta^{-1}(+)|+1\leq i\leq n).
		\end{cases}\]
	\end{remark}
	
	\subsubsection{Classification of standard $\theta$-stable parabolic subsets}\label{sec:stdthetastablepss}
	
	To determine $\SPS^\theta$, it will suffice to compute the action of $\theta$ on the set of simple roots of $\Delta(G\times_SS',H\times_S S')$ by Remark \ref{rem:SPS^theta}. As we explained at the beginning of section \ref{sec:halfintegralformofrelativetype}, we can compute the action of $\theta$ on $\Delta(G\times_SS',H\times_S S')$ from \cite[section 4]{hayashilinebdl}. This is a part of drawing the Vogan diagrams of \cite[Chapter VI, section 10]{knapp2002}. Therefore we only note the action and depict the corresponding Dynkin automorphism below without any more explanations. We may skip noting it if $H=T$ since $\theta$ acts trivially on $X^\ast(H\times_S S')$.
	
	\begin{proposition} Put $G=\GL_{2n+1}$.
		Then we have
		\[\begin{array}{cc}
			\theta(e_{2i-1}-e_{2i+1})=-e_{2i}+e_{2i+2}&(1\leq i\leq n-1),\\
			\theta(e_{2n-1}-e_{2n+1})=-e_{2n}+e_{2n+1},\\
			\theta(-e_{2i}+e_{2i+2})=e_{2i-1}-e_{2i+1}&(1\leq i\leq n-1),\\
			\theta(-e_{2n}+e_{2n+1})=e_{2n-1}-e_{2n+1}.
		\end{array}\]
		If we set
		\[\alpha_i=\left\{\begin{array}{cc}
			e_{2i-1}-e_{2i+1}&(1\leq i\leq n-1)\\
			e_{2n-1}-e_{2n+1}&(i=n)\\
			-e_{2n}+e_{2n+1}&(i=n+1)\\
			-e_{4n-2i+2}+e_{4n-2i+4}&(n+2\leq i\leq 2n),
		\end{array}\right.\]
		the action of $\theta$ is depicted as
		\[\begin{tikzcd}
			\overset{\alpha_1}{\bullet} \arrow[r,dash]\ar[rrrr, bend right]&\overset{\alpha_2}{\bullet} \ar[rr, bend right]\arrow[r,dash]&\cdots\arrow[r,dash]
			&\overset{\alpha_{2n-1}}{\bullet} \arrow[r,dash]\ar[ll, bend left]
			&\overset{\alpha_{2n}.}{\bullet}\ar[llll, bend left]
		\end{tikzcd}\]
	\end{proposition}
	
	\begin{proposition} 
		Put $G=\GL_{2n}$. Then we have
		\[\begin{array}{cc}
			\theta(e_{2i-1}-e_{2i+1})=-e_{2i}+e_{2i+2}&(1\leq i\leq n-1),\\
			\theta(-e_{2i}+e_{2i+2})=e_{2i-1}-e_{2i+1}&(1\leq i\leq n-1),\\
			\theta(e_{2n-1}-e_{2n})=-e_{2n}+e_{2n-1};
		\end{array}\]
		\[\begin{tikzcd}
			\overset{e_1-e_3}{\bullet} \arrow[r,dash]\ar[rrrrrr, bend right]&\overset{e_3-e_5}{\bullet} \ar[rrrr, bend right]\arrow[r,dash]&\cdots\arrow[r,dash]
			&\overset{e_{2n-1}-e_{2n}}{\bullet}\arrow[r,dash]
			&\cdots\arrow[r,dash]
			&\overset{-e_{2n-4}+e_{2n-2}}{\bullet} \arrow[r,dash]\ar[llll, bend left]
			&\overset{-e_{2n-2}+e_{2n}.}{\bullet}\ar[llllll, bend left]
		\end{tikzcd}\]
	\end{proposition}
	
	\begin{proposition}
		Put $G=\U^\ast(2n)$. Then we have
		\[\begin{array}{cc}
			\theta(e_i-e_{i+1})=-e_{n+i}+e_{n+i+1}&(1\leq i\leq n-1),\\
			\theta(e_n-e_{2n})=e_{n}-e_{2n},\\
			\theta(-e_i+e_{i+1})=e_{i-n}-e_{i-n+1}&(n+1\leq i\leq 2n-1).
		\end{array}\]
		If we set
		\[\alpha_i=\left\{\begin{array}{cc}
			e_{i}-e_{i+1}&(1\leq i\leq n-1)\\
			e_{n}-e_{2n}&(i=n)\\
			-e_{3n-i}+e_{3n-i+1}&(n+1\leq i\leq 2n-1),
		\end{array}\right.\]
		the action of $\theta$ is depicted as
		\[\begin{tikzcd}
			\overset{\alpha_1}{\bullet} \arrow[r,dash]\ar[rrrrrr, bend right]&\overset{\alpha_2}{\bullet} \ar[rrrr, bend right]\arrow[r,dash]&\cdots\arrow[r,dash]
			&\overset{\alpha_n}{\bullet}\arrow[r,dash]
			&\cdots\arrow[r,dash]
			&\overset{\alpha_{2n-2}}{\bullet} \arrow[r,dash]\ar[llll, bend left]
			&\overset{\alpha_{2n-1}.}{\bullet}\ar[llllll, bend left]
		\end{tikzcd}\]
	\end{proposition}
	
	\begin{proposition}
		Put $G=\SO(2p+1,2q+1)$ with $p\neq n$. Then we have
		\[\begin{array}{cc}
			\theta(e_i-e_{i+1})=e_i-e_{i+1}&(1\leq i\leq n,~i\neq p,~p+1),\\
			\theta(e_p-e_{p+2})=e_p-e_{p+2},\\
			\theta(e_{n+1}\pm e_{p+1})=e_{n+1}\mp e_{p+1}.
		\end{array}\]
		If we set
		\[\alpha_i=\left\{\begin{array}{cc}
			e_{i}-e_{i+1}&(1\leq i\leq p-1)\\
			e_p-e_{p+2}&(i=p)\\
			e_{i+1}-e_{i+2}&(p+1\leq i\leq n-1)\\
			e_{n+1}+e_{p+1}&(i=n)\\
			e_{n+1}-e_{p+1}&(i=n+1),
		\end{array}\right.\]
		the action of $\theta$ is depicted as
		\[\begin{tikzcd}
			&&&&\overset{\alpha_{n}}{\bullet}\ar[dd, bend left]\\
			\overset{\alpha_1}{\bullet}\ar[r, dash]&\overset{\alpha_2}{\bullet}\ar[r, dash]
			&\cdots\ar[r, dash]&\overset{\alpha_{n-1}}{\bullet}\ar[ru, dash]\ar[rd, dash]\\
			&&&&\overset{\alpha_{n+1}.}{\bullet}\ar[uu, bend right]
		\end{tikzcd}\]
	\end{proposition}

	\begin{proposition}
		Put $G=\SO(1,2n+1)$. Then we have
		\[\begin{array}{cc}
			\theta(e_i-e_{i+1})=e_i-e_{i+1}&(2\leq i\leq n),\\
			\theta(e_{n+1}\pm e_1)=e_{n+1}\mp e_1;
		\end{array}\]
		\[\begin{tikzcd}
			&&&&\overset{e_{n+1}+e_1}{\bullet}\ar[dd, bend left]\\
			\overset{e_2-e_3}{\bullet}\ar[r, dash]&\overset{e_3-e_4}{\bullet}\ar[r, dash]
			&\cdots\ar[r, dash]&\overset{e_{n}-e_{n+1}}{\bullet}\ar[ru, dash]\ar[rd, dash]\\
			&&&&\overset{e_{n+1}-e_1.}{\bullet}\ar[uu, bend right]
		\end{tikzcd}\]
	\end{proposition}

	\subsubsection{Conjugate action}\label{sec:conjaction}
	
	The canonical conjugate action on $\rtype(G,K)(S')$ induces an involution $\bar{\ }$ on $W^\theta_+\times\SPS^\theta/\sim$. Let $w_{0,G}$ (resp.~$w_{0,K}$) be the longest element of the Weyl group $W(G,H)(S')$ (resp.~$W(K,T)(S')$) with respect to our positive system $\Delta^+(G\times_S S',H\times_S S')$ (resp.~$\Delta^+(K\times_S S',T\times_S S')$). We will regard $W(K,T)(S')$ as a subgroup of $W(G,H)(S')$ in the canonical way.
	
	\begin{lemma}
		For $(w,\Phi)\in W^\theta_+\times\SPS^\theta$, we have
		\begin{equation}
			\overline{(w,\Phi)}=(w_{0,K}w w_{0,G},-w_{0,G}\Phi).\label{eq:conjaction}
		\end{equation}
		
	\end{lemma}
	
	\begin{proof}
		The pair $(w,\Phi)$ gives rise to a parabolic subset $w\Phi\in \Ss^+(G\times_S S',K\times_S S')$. Then we have $\overline{w\Phi}=-w\Phi$ since the conjugate action on $X_\ast(T\times_S S')$ is equal to $-1$. Since we are interested in $K\times_S S'$-orbits, we may replace $-w\Phi$ with $-w_{0,K}w\Phi$. Notice that $-w_{0,K}w\Phi$ belongs to $\Ss^+(G\times_S S',K\times_S S')$.
		
		We follow the construction of the inverse map of \eqref{eq:bijectivegeneral} in the proof of Theorem \ref{thm:solutiontotheclassificationproblem} to obtain a representative corresponding to $-w_{0,K}w\Phi$. Notice that $-w_{0,K}w\Delta^+(G,H)$ is a $\theta$-stable positive system contained in $-w_{0,K}w\Phi$ such that
		\[(-w_{0,K}w\Delta^+(G\times_S S',H\times_S S'))|_{T\times_SS'}
		\cap \Delta(K\times_S S',T\times_S S')=\Delta^+(K\times_S S',T\times_S S')\]
		since $(w\Delta^+(G\times_S S',H\times_S S'))|_{T\times_SS'}\cap \Delta(K\times_S S',T\times_S S')=\Delta^+(K\times_S S',T\times_S S')$. We have
		\[w_{0,K}w w_{0,G}\Delta^+(G\times_S S',H\times_S S')
		=-w_{0,K}w\Delta^+(G\times_S S',H\times_S S').\]
		Therefore a lift of the preimage of $-w_{0,K}w\Phi$ is given by
		\[(w_{0,K}w w_{0,G},(w_{0,K}w w_{0,G})^{-1}(-w_{0,K}w\Phi))
		=(w_{0,K}w w_{0,G},-w_{0,G}\Phi).\]
		This completes the proof.
	\end{proof}
	
	Before we apply this lemma to our models, let us restrict our attention. In fact, we may skip the case $w_{0,K}=-1$ in $W(K,T)(S')$ by section \ref{sec:absolutedescent}. We can compute the right hand side in \eqref{eq:conjaction} explicitly if we know $w_{0,G}$ and $w_{0,K}$. We can find that if $H=T$ then the positive systems $\Delta^+(G\times_S S',H\times_S S')$ and $\Delta^+(K\times_S S',T\times_S S')$ are (direct sums of) the positive systems in \cite[Appendix C]{knapp2002}. The corresponding longest elements are given by
	\[\begin{array}{cc}
		\left(\begin{array}{cccc}
			1 & 2 & \cdots & n+1 \\
			n+1 & n & \cdots & 1
		\end{array}\right)&(\mathrm{type}\ A_n)\\
		(-1,-1,\ldots,-1,1)&(\mathrm{type}\ D_{2n+1})\\
		-1&(\mathrm{otherwise}).
	\end{array}\]
	
	The only remaining case is $\GL_{4n+2}$:
	
	\begin{proposition}
		Let $n\geq 0$. Put $(G,K)=(\GL_{4n+2},\SO(4n+2))$. Then
		\[w_{0,G},w_{0,K}\in W(G,H)(S')\cong\mathfrak{S}_{4n+2}\]
		are given by
		\[\begin{array}{c}
			w_{0,G}=(1~2)(3~4)\cdots(4n+1\ 4n+2),\\
			w_{0,K}=(1~2)(3~4)\cdots(4n-1\ 4n).
		\end{array}\]
	\end{proposition}
	
	\begin{proof}
		One can easily check
		\[(1~2)(3~4)\cdots(4n+1\ 4n+2)\Delta^+(\GL_{4n+2},H\times_S S')
		=-\Delta^+(\GL_{4n+2},H\times_S S')\] 
		(cf.~\cite[Example 4.1.2]{hayashilinebdl}).
		
		The longest element $w_{0,K}$ is represented by 
		$\diag(w_2,w_2,\ldots,w_2,1,1)$
		(see Example \ref{ex:typeD_n}). It is straightforward that this element acts on $X^\ast(H\times_S S')\cong \ZZ^{4n+2}$ by $(1~2)(3~4)\cdots(4n-1\ 4n)$.
	\end{proof}

	\subsubsection{Disconnected case}\label{sec:disconn}
	
	In this section, we work with the standard models $(G,K)$ of disconnected classical Lie groups. Following Examples \ref{ex:c1} and \ref{ex:c2}, set
	\[c=\left\{\begin{array}{cc}
		\diag(I_{2n-1},-1)&(G=\GL_{2n},~n\geq 1)\\
		-I_{2n+1}&(G=\GL_{2n+1},~n\geq 0)\\
		\diag(I_{2p-1},-I_{2q+2})&(G=\SO(2p,2q+1))\\
		I_{2p+2q+2}&(G=\SO(2p+1,2q+1))\\
		\diag(I_{2p-1},-1,I_{2q-1},-1)&(G=\SO(2p,2q),~n=p+q\geq 1).
	\end{array}\right.\]
	Then $c$ belongs to $K(S)\setminus K^0(S)$. It follows by definition that $c$ acts on $\rtype(G,K^0)$ as an involution. We compute it under the identification
	\[S'\times_S \rtype(G,K^0)\cong\Ss^+(G,K^0)_{S'}\]
	and the combinatorial description \eqref{eq:bijectivegeneral}. For this, notice that 
	the matrix $c$ normalizes $T$ and therefore $H$. Under the identification $X^\ast(T\times_S S')$ with a free abelian group in \cite[section 3.2]{hayashilinebdl}, the induced action on $X^\ast(T\times_S S')$ is expressed as
	\[\left\{\begin{array}{cc}
		\diag(I_{n-1},-1)&(G=\GL_{2n},~n\geq 1)\\
		I_{n}&(G=\GL_{2n+1},~n\geq 0)\\
		\diag(I_{p-1},-1,I_{q})&(G=\SO(2p,2q+1))\\
		I_{n+1}&(G=\SO(2p+1,2q+1))\\
		\diag(I_{p-1},-1,I_{q-1},-1)&(G=\SO(2p,2q),~n=p+q\geq 1).
	\end{array}\right.\]
	In particular, one can check that this action respects $\Delta^+(K^0\times_S S',T\times_S S')$.
	
	The conjugate action of $c$ on $H$ determines an element $w_d\in W(G,H)(S')$, where
	\[w_d=\left\{\begin{array}{cc}
		(2n-1~2n)&(G=\GL_{2n},~n\geq 1)\\
		e&(G=\GL_{2n+1},~n\geq 0)\\
		(\overbrace{1,1,\ldots,1}^{p-1},-1,
		\overbrace{1,1,\ldots,1}^{q})&(G=\SO(2p,2q+1))\\
		e&(G=\SO(2p+1,2q+1))\\
		(\overbrace{1,1,\ldots,1}^{p-1},-1,
		\overbrace{1,1,\ldots,1}^{q-1},-1)&(G=\SO(2p,2q),~n=p+q\geq 1).
	\end{array}\right.\]

	\begin{proposition}\label{prop:actionofc}
		The involution on $\Ss^+(G\times_S S',K^0\times_S S')$ attached to $c$ is given by
		\[(w,\Phi)\mapsto (w_dw,\Phi)\]
		under the bijection \eqref{eq:bijectivegeneral}.
	\end{proposition}
	
	\begin{proof}
		Let $(w,\Phi)\in W^\theta_+\times \SPS^\theta$. Then we have $w_dw\in W^\theta_+$. In fact, 
		\[\begin{split}
			(w_d w\Delta^+(G\times_S S',H\times_S S'))|_{T\times_SS'}
			&=w_d (w\Delta^+(G\times_S S',H\times_S S'))|_{T\times_SS'}\\
			&= w_d \Delta^+(K^0\times_S S',T\times_S S')\\
			&=\Delta^+(K^0\times_S S',T\times_S S').
		\end{split}
		\]
		Since $\theta(c)=c$, we have $\theta(w_d)=w_d$. This implies $\theta(w_d w)=w_d w$.
		
		The assertion is now evident by construction of the inverse map of \eqref{eq:bijectivegeneral}.
	\end{proof}
	
	\begin{corollary}\label{cor:numberoforbits}
		For $(w,\Phi)\in W^\theta_+\times \SPS^\theta$, the corresponding $K\times_S S'$-orbit consists of two $K^0\times_S S'$-orbits if and only if $w^{-1}w_dw\Phi=\Phi$.
	\end{corollary}
	
	\begin{proposition}\label{prop:numberoforbits1}
		Suppose that $(G,K)$ is not of the form
		\[(\SO(4p,4q+2),\mathrm{S}(\Oo(4p)\times\Oo(4q+2))).\]
		Then we have $\rtype(G,K)\cong \Ss^+(G\times_S S',K\times_S S')_{S}$. Moreover, take \[(w,\Phi)\in \Ss^+(G\times_S S',K^0\times_S S').\]
		Let $x'$ and $y$ be the corresponding elements of $\rtype(G,K^0)(S')$ and $\rtype(G,K)(S)$ respectively. Then the following assertions hold:
		\begin{enumerate}
			\item If $w^{-1}w_dw\Phi=\Phi$, then $x'$ descends to a relative parabolic type $x\in \rtype(G,K^0)(S)$ and $\Pp^{K\mathrm{-st}}_{G,y}=\Pp^{K^0\mathrm{-st}}_{G,x}$.
			\item Suppose $w^{-1}w_dw\Phi\neq \Phi$. Then we have $\Pp^{K\mathrm{-st}}_{G,y}=\Pp^{(S'\times_S K^0)\mathrm{-st}}_{S'\times_S G,x'}$.
		\end{enumerate}
	\end{proposition}
	
	\begin{proof}
		The first assertion holds since the matrices
		\[\left\{\begin{array}{cc}
			\diag(1,-1,1,-1,\ldots,1,-1)&(G=\GL_{2n},~n\geq 1)\\
			\diag(1,-1,1,-1,\ldots,1,-1,1)&(G=\GL_{2n+1},~n\geq 0)\\
			\diag(\overbrace{1,-1,\ldots,1,-1}^{2p},
			\overbrace{1,-1,\ldots,1,-1}^{2q},(-1)^{p+q})
			&(G=\SO(2p,2q+1))\\
			\diag(\overbrace{1,-1,\ldots,1,-1}^{2p},(-1)^p,
			\overbrace{1,-1,\ldots,1,-1}^{2q},(-1)^q)
			&(G=\SO(2p+1,2q+1))\\
			\diag(1,-1,1,-1,\ldots,1,-1)&(G=\SO(2p,2q),~n=p+q\in 2+2\NN).
		\end{array}\right.\]
		in $K(S)$ act on $X^\ast(T\times_S S')$ by $-1$. This also implies that $\overline{(w,\Phi)}$ is equivalent to $(w,\Phi)$ or $(w_dw,\Phi)$.
		
		Suppose $w^{-1}w_dw\Phi=\Phi$. Then $x'$ descends to $x\in \rtype(G,K^0)(S)$ by $\overline{(w,\Phi)}\sim (w,\Phi)$. One can see the equality
		$S'\times_S\Pp^{K\mathrm{-st}}_{G,y}=S'\times_S\Pp^{K^0\mathrm{-st}}_{G,x}$
		since both consists of the same single $K^0\times_S S'$-orbit. The faithfully flat descent implies
		$\Pp^{K\mathrm{-st}}_{G,y}=\Pp^{K^0\mathrm{-st}}_{G,x}$.
		
		Suppose $w^{-1}w_dw\Phi\neq \Phi$. Then $\overline{(w,\Phi)}\sim (w_dw,\Phi)$. In view of Corollary \ref{cor:numberoforbits}, we obtain
		\[S'\times_S\Pp^{K\mathrm{-st}}_{G,y}=\Pp^{(S'\times_S K^0)\mathrm{-st}}_{S'\times_S G,x'}
		\coprod \Pp^{(S'\times_S K^0)\mathrm{-st}}_{S'\times_S G,\bar{x}'}
		\cong S'\times_S\Pp^{(S'\times_S K^0)\mathrm{-st}}_{S'\times_S G,x'}.\]
		The desired equality in (2) now follows by the faithfully flat descent. 
	\end{proof}
	
	It remains to study the case 
	$(G,K)=(\SO(4p,4q+2),\mathrm{S}(\Oo(4p)\times\Oo(4q+2)))$.
	The situation is more complicated this time. To explain it in examples later, let us give easy criteria on coincidence among a given $K^0\times_S S'$-orbit and its translations by the Galois involution and $c$ in terms of cocharacters $\mu$:
	
	\begin{lemma}\label{lem:numerical_criterion}
		Identify $X_\ast(T\times_S S')_+$ with
		\[\{\mu=(a_i)\in \ZZ^{2n+1}:~
		a_1\geq\cdots\geq a_{2p-1}\geq |a_{2p}|,~a_{2p+1}\geq\cdots\geq a_{2p+2q}\geq |a_{2p+2q+1}|\}\subset\ZZ^{2n+1}\]
		by \cite[Example 3.4.8]{hayashilinebdl} and \cite[Appendix C]{knapp2002}\footnote{If $q=0$ then we do not have the inequality condition for $a_{2p+1}$ so that $a_{2p+1}$ can be any integer. }. Let $\mu=(a_i)\in X_\ast(T\times_S S')_+$.
		\begin{enumerate}
			\item We have $\Phi_\mu=-w_{0,K^0}\Phi_\mu$ if and only if the following conditions are satisfied:
			\begin{enumerate}
				\item[(i)] $|a_{2p}|\geq |a_{2n+1}|$.
				\item[(ii)] If $|a_i|=|a_{2n+1}|$ then $a_i=0$ ($1\leq i\leq 2n$).
			\end{enumerate}
			\item We have $w_d\Phi_\mu=-w_{0,K^0}\Phi_\mu$ if and only if the following conditions are satisfied:
			\begin{enumerate}
				\item[(i)] $|a_{2p}|\leq |a_{2n+1}|$.
				\item[(ii)] If $|a_i|=|a_{2p}|$ then $a_i=0$ ($1\leq i\leq 2n+1$ with $i\neq 2p$).
			\end{enumerate} 
			\item We have $w_d\Phi_\mu=\Phi_\mu$ if and only if $a_{2p}=a_{2n+1}=0$.
		\end{enumerate}
		In particular, if the equality $w_d\Phi_\mu=-w_{0,K^0}\Phi_\mu$ holds then we have $\Phi_\mu=-w_{0,K^0}\Phi_\mu=w_d\Phi_\mu$.
	\end{lemma}
	
	\begin{proof}
		We prove (1). Write $-w_{0,K^0}\mu=(b_i)$. Then
		\[b_i=\left\{\begin{array}{cc}
			a_i&(1\leq i\leq 2n)\\
			-a_i&(i=2n+1).
		\end{array}\right.\]
		Hence for $1\leq i<j\leq 2n$, $a_i-a_j=b_i-b_j$ and $a_i+a_j=b_i+b_j$. We may therefore restrict ourselves to $\pm e_i\pm e_{2n+1}$. Without loss of generality, we may assume $a_{2n+1}\geq 0$. 
		
		Suppose that the conditions (i) and (ii) are satisfied. For $1\leq i\leq 2p-1$, we have $e_i\pm e_{2n+1}\in \Phi_\mu$ by $a_i\geq |a_{2p}|\geq a_{2n+1}$. Similarly, for $2p+1\leq i\leq 2n$, we have $e_i\pm e_{2n+1}\in \Phi_\mu$ by $a_i\geq a_{2n+1}$. Let $1\leq i\leq 2n$ with $i\neq 2p$. If $a_i>a_{2n+1}$ then $-e_{i}\pm e_{2n+1}\not\in\Phi_\mu$; Otherwise, $-e_{i}\pm e_{2n+1}\in\Phi_\mu$ since $a_i=a_{2n+1}=0$.
		\begin{description}
			\item[Case 1] Assume $a_{2p}\geq 0$. Then we have $e_{2p}\pm e_{2n+1}\in\Phi_\mu$. If $a_{2p}>a_{2n+1}$ then
			\[-e_{2p}\pm e_{2n+1}\not\in\Phi_\mu;\]
			Otherwise, $-e_{2p}\pm e_{2n+1}\in\Phi_\mu$ since $a_{2p}=a_{2n+1}=0$.
			\item[Case 2] Assume $a_{2p}<0$. Then we have $-e_{2p}\pm e_{2n+1}\in\Phi_\mu$. If $-a_{2p}>a_{2n+1}$ then
			\[e_{2p}\pm e_{2n+1}\not\in\Phi_\mu;\]
			Otherwise, $e_{2p}\pm e_{2n+1}\in\Phi_\mu$ since $a_{2p}=a_{2n+1}=0$.
		\end{description}
		
		Conversely, suppose $\Phi_\mu=-w_{0,K^0}\Phi_\mu$. Condition (i) holds since $\{\pm e_{2p}+e_{2n+1}\}\cap \Phi_\mu\neq \emptyset$ by $a_{2n+1}\geq 0$. Suppose $|a_i|=a_{2n+1}$ for some $1\leq i\leq 2n$. Observe that
		\[\begin{array}{cc}
			|\{-e_i\pm e_{2n+1}\}\cap \Phi_\mu|=1&(a_i>0)\\
			|\{e_i\pm e_{2n+1}\}\cap \Phi_\mu|=1&(a_i<0).
		\end{array}\]
		Hence $a_i$ must be zero.
		
		Observe that 
		$-w_dw_{0,K^0}=\diag(I_{2p-1},-1,I_{2q+1})$
		as an action of $X^\ast(T\times_S S')_+$. Hence (2) follows from a similar argument to (1).
		
		Finally, we prove (3). If $a_p=a_n=0$ then $w_d\mu=\mu$. This shows the ``if'' direction. Conversely, suppose $\Phi_\mu=w_d\Phi_\mu$. Observe that $w^2_d=e$ and $w_d(e_{2p}-e_{2n+1})=-(e_{2p}-e_{2n+1})$. Since $\Phi_\mu$ is a parabolic subset, at least one of $\pm(e_{2p}-e_{2n+1})$ belongs to $\Phi_\mu$. We thus deduce \[\{\pm(e_{2p}-e_{2n+1})\}\subset\Phi_\mu.\]
		It now follows by definition of $\Phi_\mu$ that $a_{2p}=a_{2n+1}$. Apply a similar argument to $\pm(e_{2p}+e_{2n+1})$ to deduce $a_{2p}+a_{2n+1}=0$. We now conclude by the two equalities that $a_{2p}=a_{2n+1}=0$.
	\end{proof}
	
	\begin{proposition}
		Take $(w,\Phi)\in \Ss^+(G\times_S S',K^0\times_S S')$. Let $x'\in \rtype(G,K^0)(S')$ and $y'\in \rtype(G,K^0)(S')$ denote the relative parabolic types corresponding to $(w,\Phi)$ respectively.
		\begin{enumerate}
			\item If $(w,\Phi)\sim \overline{(w,\Phi)}$ then $x'$ and $y'$ descend to
			$x,y\in \rtype(G,K^0)(S)$ respectively. Moreover, we have
			\[\Pp^{K\mathrm{-st}}_{G,y}=\begin{cases}
				\Pp^{K^0\mathrm{-st}}_{G,x}&(\overline{(w,\Phi)}\sim (w_dw,\Phi))\\
				\Pp^{K^0\mathrm{-st}}_{G,x}\coprod \Pp^{K^0\mathrm{-st}}_{G,cx}
				&(\overline{(w,\Phi)}\not\sim (w_dw,\Phi)).
			\end{cases}\]
			\item Suppose $(w,\Phi)\not\sim \overline{(w,\Phi)}\sim (w_dw,\Phi)$. Then we have $cx'=x'$, and $y'$ descends to
			\[y\in \rtype(G,K^0)(S).\]
			Moreover, the equality
			$\Pp^{K\mathrm{-st}}_{G,y}=\Pp^{(S'\times_S K^0)\mathrm{-st}}_{S'\times_S G,x'}$ holds.
			\item If $(w,\Phi)\not\sim \overline{(w,\Phi)}\not\sim (w_dw,\Phi)$, we have $\bar{y}'\neq y'$. Moreover, the equality
			\[\Pp^{(S'\times_S K)\mathrm{-st}}_{S'\times_S G,y'}=
			\Pp^{(S'\times_S K^0)\mathrm{-st}}_{S'\times_S G,x'}\coprod \Pp^{(S'\times_S K^0)\mathrm{-st}}_{S'\times_S G,cx'}\]
			holds in $\Pp^{K\mathrm{-st}}_G$.
			\item The pair $(w,\Phi)$ satisfies exactly one of the assumptions of (1), (2), and (3).
		\end{enumerate}
	\end{proposition}
	
	\begin{proof}
		The final assertion of Lemma \ref{lem:numerical_criterion} shows (4). This also implies $(w,\Phi)\not\sim(w_dw,\Phi)$ under the assumption of (3). One can then prove (1)-(3) by a similar argument to Proposition \ref{prop:numberoforbits1}.
	\end{proof}

	We remark that all of the four cases above (two for (1)) actually occur:
	
	\begin{example}
		Put $p=1$ and $q=0$. Put
		\[\begin{array}{cccc}
			\mu_1=(0,0,0),
			&\mu_2=(0,0,1),
			&\mu_3=(2,1,0),
			&\mu_4=(1,1,1).
		\end{array}\]
		Then we have
		\[\begin{array}{cc}
			\Phi_{\mu_1}=-w_{0,K^0}\Phi_{\mu_1}=w_d\Phi_{\mu_1},
			&\Phi_{\mu_2}\neq -w_{0,K^0}\Phi_{\mu_2}=w_d\Phi_{\mu_2},\\
			\Phi_{\mu_3}=-w_{0,K^0}\Phi_{\mu_3}\neq w_d\Phi_{\mu_3},
			&\Phi_{\mu_4}\neq-w_{0,K^0}\Phi_{\mu_4}\neq w_d\Phi_{\mu_4}.	
		\end{array}\]
	\end{example}

	\section{Arithmetic models of cohomologically induced modules}\label{sec:arithmeticmodelsofaq(lambda)}
	
	In this section, we apply our theory developed in this paper to the theory of $(\lieg,K)$-modules. That is, we explain how one can obtain arithmetic models of the cohomologically induced modules.
	
	If we are given a double Galois covering $S'\to S$ of schemes, we denote associated actions of the nontrivial element of the Galois group $\ZZ/2\ZZ$ on objects by $\bar{}$.
	For a homomorphism $k\to k'$ of commutative rings, we denote the base change functors $k'\otimes_k -$ (of schemes, Lie algebras, modules, and so on) by $(-)_{k'}$.
	
	\subsection{General construction}\label{sec:general_construction}
	
	Let $S$ be a scheme. Let $G$ be a reductive group scheme over $S$, $K$ be a smooth closed subgroup scheme of $G$ with $K^0$ reductive and $\pi_0(K)$ finite \'etale. Let $x\in \rtype (G,K)(S)$, $\Aa$ be a $G$-equivariant tdo on $\Pp_{G,gt(x)}$, and $\Mm$ be a $K$-equivariant left $i^\cdot_x\Aa$-module on $\Pp^{K\mathrm{-st}}_{G,x}$. Write $p_{G,x}:\Pp_{G,x}\to S$ for the structure morphism. Then $R^\bullet (p_{G,x})_\ast (i_x)_+\Mm$ is a $(\lieg,K)$-module for $\bullet\in\NN$. For a flat morphism $s:S'\to S$, there is a canonical isomorphism
	\[s^\ast R^\bullet (p_{G,x})_\ast (i_x)_+\Mm\cong R^\bullet (p_{G\times_S S',x|_{S'}})_\ast (i_x\times_S S')_+s^\ast_{K,x}\Mm,\]
	where $s_{K,x}$ be the projection map $S'\times_S \Pp^{K\mathrm{-st}}_{G,x}\to \Pp^{K\mathrm{-st}}_{G,x}$. Let $p^{K\mathrm{-st}}_{G,x}$ denote the structure morphism of $\Pp^{K\mathrm{-st}}_{G,x}$.

	On the course of writing this section, the first author established a general result on finiteness of the module of global sections of the direct image of certain twisted $\mathcal D$-modules, which is based on the contents of section \ref{sec:tdo_and_pic_alg}-\ref{sec:operationswithDmodules}. Without any circular arguments, we obtain:
	
	\begin{theorem}[{\cite[Theorem 1.1, the end of section 1]{hayashifil}}]\label{thm:filtration}
		Assume that $S$ is Dedekind, and that $\Mm$ is a twisted integrable left connection. Then there exists a natural exhaustive $K$-invariant filtration on
		$(p_{G,x})_\ast(i_x)_+\Mm$
		such that its $p$-th associated graded sheaf is a locally free $\OO_S$-module of finite rank for every integer $p$. Moreover, the 0th step is given by
		\begin{flalign*}
			&(p_{G,x})_\ast (i_\ast(\Mm\otimes_{\OO_{\Pp^{K\mathrm{-st}}_{G,x}}} \omega_{\Pp^{K\mathrm{-st}}_{G,x}/S})
			\otimes_{\OO_{\Pp_{G,gt(x)}}} \omega^{\vee}_{\Pp_{G,gt(x)}/S}
			)\\
			&\cong 
			(p^{K\mathrm{-st}}_{G,x})_\ast (\Mm\otimes_{\OO_{\Pp^{K\mathrm{-st}}_{G,x}}} \omega_{\Pp^{K\mathrm{-st}}_{G,x}/S}\otimes_{\OO_{\Pp^{K\mathrm{-st}}_{G,x}}}
			i^\ast_x \omega^{\vee}_{\Pp_{G,gt(x)}/S}).
		\end{flalign*}

	\end{theorem}
	
	\subsection{Relation to cohomological induction}\label{sec:kitchen}

	Let us recall an algebraic description of the cohomology group $H^\bullet(\Pp_{G,x},i_+\Mm)$
	over fields of characteristic zero. We follow the notations in \cite{knappvogan} for the functors $\mathcal{F}^{\lieq,L}_{\lieg,K}$, $I^{\lieg,K}_{\lieq,L}$, $P^{\lieg,K}_{\lieq,L}$, $H^0(\lieg,-)$, and $\ind^{\lieg}_{\lieq}$ under the algebraic setting (cf.~\cite{januszewskirationality, hayashi2019}). To save space, for a finite dimensional vector space $V$, we denote $\wedge^{\dim_F V} V=\wedge^{\mathrm{top}} V$.
	
	\begin{theorem}\label{thm:geometric_induction}
		Let $G$ be a connected reductive algebraic group over a field $F$ of characteristic zero, $K$ be a (possibly disconnected) reductive algebraic subgroup of $G$, and $Q$ be a stable parabolic subgroup relative to $K$. Set
		$x=rt(Q)\in\rtype(G,K)(F)$.
		Let $u$ be the dimension of the unipotent radical of $Q\cap K$. Choose a Levi subgroup $L$ of $Q\cap K$ (\cite[Chapter VIII, Theorem 4.3]{hochschild1981}).
		
		Let $\Aa$ be a $G$-equivariant tdo on $\Pp_{G,gt(x)}$, and $\Mm$ be a $K$-equivariant integrable left $i^\cdot_{x} \Aa$-connection on $\Pp^{K\mathrm{-st}}_{G,x}$. Let $M$ be the fiber of $\Mm$ at $Q\in \Pp^{K\mathrm{-st}}_{G,z}(F)$. 
		
		\begin{enumerate}
			\item The vector space $M$ is naturally equipped with the structure of a $(\lieq,Q\cap K)$-module. Moreover, its restriction to the unipotent radical of $\lieq$ is trivial.
			\item There exist isomorphisms
			\begin{equation}
				H^\bullet(\Pp_{G_,gt(x)},(i_x)_+\Mm)
				\cong R^{\bullet+u} I^{\lieg,K}_{\lieg,L}\ind^{\lieg}_{\lieq} 
				(M\otimes_F \wedge^{\mathrm{top}} \lieg/\lieq)\label{eq:geometric_induction}
			\end{equation}
			of $(\lieg,K)$-modules.	
		\end{enumerate}
	\end{theorem}
	
	\begin{proof}
		For (1), let $\lambda$ be the linear functional on $\lieq$ corresponding to $\Aa$ (see \cite[Theorem 4.9.2]{kashiwara1989}). The fiber $M$ is naturally equipped with the structure of a representation of $L\cap K$, whose differential action is given by the restriction of $\lambda$ (\cite[Proposition 4.11.1]{kashiwara1989}). It extends to a $(\lieq,L\cap K)$-module by letting $\liel$ act on $M$ by $\lambda$ since $\lambda$ is $Q$-invariant. The character $\lambda$ is trivial on the unipotent radical of $\lieq$ by the same reason. This shows (1).
		
		Part (2) is verified along the same line as \cite[Theorem 5.4 and Corollary 5.5]{kitchen2012} (see also \cite[1.13. Theorem]{milicicpandzic1998} and \cite[3.3.2. Theorem]{pandzic2007}).
	\end{proof}
	
	To see the consequences in representation theory of real reductive Lie groups, recall that a Cartan involution of a (possibly disconnected) real reductive algebraic group $G$ is a holomorphic involution $\theta$ on $G(\CC)$ commuting with the complex conjugation such that the fixed point subgroup of $G(\CC)$ by $\sigma_c$ (= the composition of the complex conjugation and $\theta$) is compact and meets every component of $G(\CC)$ (\cite[Definitions 3.12 and 3.3]{adamstaibi2018}). 
	
	\begin{remark}[{\cite[Theorem 3.13]{adamstaibi2018}}]\label{rem:exists_Cartan}
		Every real reductive algebraic group $G$ admits a Cartan involution. It is unique up to a conjugation by an inner automorphism from $G(\RR)_0$. 
	\end{remark}
	
	As mentioned below \cite[Definition 3.12]{adamstaibi2018}, $\theta$ is complex algebraic. Moreover, it descends to an involution of $G$ since $\theta$ commutes with the complex conjugation (embed the coordinate ring of $G_\CC$ to the ring of holomorphic functions on $G(\CC)$). Henceforth we think of Cartan involutions as involutions of real reductive algebraic groups. Let us also record the preservation property of Cartan involutions by restriction:
	
	\begin{lemma}\label{lem:stabilityofCartan}
		Let $G$ be a possibly disconnected real reductive algebraic group with a Cartan involution $\theta$, and $H$ be a $\theta$-stable real algebraic subgroup. Then $H$ is reductive, and $\theta|_H$ is a Cartan involution on $H$.
	\end{lemma}
	
	\begin{proof}
		See the proof of \cite[Lemma 3.11]{adamstaibi2018} or apply it to $\sigma=\sigma_c$.
	\end{proof}
	
	\begin{example}
		Let $n\geq 1$. Then $((-)^{T})^{-1}$ is a Cartan involution of $\SL_n$, where $(-)^T$ is the transpose.
	\end{example}
	
	\begin{example}\label{ex:mu_n}
		For a nonnegative integer $n$, put $G=\Spec\RR[x]/(x^n-1)$. It is isomorphic to the center of $\SL_n$. Hence $x\mapsto x^{-1}$ is the unique Cartan involution of $G$.
	\end{example}
	
	\begin{example}
		The trivial involution is a Cartan involution for every real algebraic subgroup $K$ of a reductive algebraic group $G$ fixed by a Cartan involution of $G$. In particular, $K(\RR)$ meets every connected component $K(\CC)$ for such $K$. More strongly, we have a Cartan decomposition $K(\RR)\times \sqrt{-1}\liek\cong K(\CC)$ (\cite[Lemma 3.6]{adamstaibi2018}). This property is crucial because of the equivalence of the categories of complex representations of $K(\RR)$ and $K_\CC$ (the unitary trick).
	\end{example}
	
	In relation with the above, let us note:
	
	\begin{lemma}
		For a real reductive algebraic group $G$, the following conditions are equivalent:
		\begin{enumerate}
			\item[(a)] $\id_G$ is a Cartan involution;
			\item[(b)] $G(\RR)$ is compact, and $\pi_0(G(\RR))\cong \pi_0(G(\CC))$.
		\end{enumerate}
	\end{lemma}
	
	\begin{proof}
		It is clear by definition that (b) implies (a). The converse follows from \cite[Lemma 3.6]{adamstaibi2018}.
	\end{proof}
	
	\begin{corollary}\label{cor:theta=id}
		Let $G$ be a real reductive algebraic group with trivial Cartan involution.
		\begin{enumerate}
			\item We have $\pi_0(G)\cong\pi_0(G(\RR))_{\Spec\RR}$.
			\item The open subgroups of $G$ and $G(\RR)$ are in one-to-one correspondence for taking the real points.
		\end{enumerate}
	\end{corollary}

	\begin{example}[Weyl algebraic group]
		Put $G=N_{\SL_2}(\SO(2))$. Then we have $G^0=\SO(2)$, $G(\RR)=\SO(2,\RR)=G^0(\RR)$. In particular, $\pi_0(G(\RR))$ is trivial. On the other hand, we have $\pi_0(G)\cong(\ZZ/2\ZZ)_{\Spec\RR}$. Notice that $\pi_0(G)$ is still constant (cf.~Example \ref{ex:mu_n}).
	\end{example}
	
	\begin{proposition}
		Let $G$ be a real reductive algebraic group with a Cartan involution $\theta$. Then every central torus of $G$ is $\theta$-stable.
	\end{proposition}
	
	\begin{proof}
		Since any involution of a given linear algebraic group respects the unit component and the center, we may replace $G$ with the unit component of the center of $G$ to assume that $G$ is a torus. 
		
		In virtue of Remark \ref{rem:exists_Cartan}, any torus has a unique Cartan involution since its inner automorphisms are trivial. It will therefore suffice to show that every homomorphism between tori respects the Cartan involution. Recall that every torus is the product of copies of $\GG_m$, $\U(1)$, and $\res_{\CC/\RR}\GG_m$ (\cite[Lemma 1.5]{moser-jauslin}). We may only see the homomorphisms among them. Then the assertion is straightforward.
	\end{proof}
	
	\begin{corollary}\label{cor:centraltorus}
		Let $G$ be a real reductive algebraic group with a Cartan involution $\theta$. Let $\liep\subset\lieg$ be the $-1$-eigenspace of the differential of $\theta$.
		\begin{enumerate}
			\item Every compact central torus is contained in $(G^\theta)^0$.
			\item For any split central torus $S$, its Lie algebra $\lies$ is contained in $\liep$.
		\end{enumerate}
	\end{corollary}
	
	Part (2) will be used in the description of the cohomology groups of cohomologically induced modules after Vogan--Zuckerman (see Theorem \ref{thm:a_q->cohomological}).

	For applications of \cite{knappvogan}, we recall that $G(\RR)$ is a real reductive Lie group of Harish-Chandra class in the sense of \cite[Definition 4.29]{knappvogan} for the restriction of $\theta$ and an existing choice of a $G(\RR)$-invariant non-degenerate bilinear form on $\lieg$ if $G$ is connected (\cite[Examples 4) of Definition 4.29]{knappvogan}, see also \cite[Remark 3.16]{adamstaibi2018}). If $K$ is an open subgroup of $G^\theta$ (cf.~Corollary \ref{cor:theta=id}), $(\lieg,K(\RR))$ is a real reductive pair in the sense of \cite[Definition 4.30]{knappvogan}. The corresponding real reductive Lie group (see \cite[Proposition 4.31]{knappvogan}) is realized as an open subgroup of $G(\RR)$ in the obvious way.

	To rewrite the right hand side of \eqref{eq:geometric_induction}, we discuss a real form of the Levi subgroup in a specific setting. Let us also introduce a subgroup of $K$ concerning the connected components of the $K$-orbit in Theorem \ref{thm:geometric_induction}:
	
	\begin{lemma}\label{lem:realformofK^1_C}
		Let $G$ be a connected real reductive algebraic group. Fix a Cartan involution $\theta$. Let $K$ be an open subgroup of $G^\theta$. Choose a maximal torus $T$ of $K^0$. Write $H=Z_G(T)$. Let $Q'$ be a $\theta_{\CC}$-stable parabolic subgroup containing $H_{\CC}$, and $L'$ be the Levi subgroup of $\bar{Q}'$ containing $H_\CC$.
		\begin{enumerate}
			\item The intersection $L'\cap K_\CC$ is a Levi subgroup of $\bar{Q}'\cap K_\CC$.
			\item The subgroup $L'$ is defined over the real numbers. We denote the corresponding real form by $L$.
			\item Let $G(\RR)'$ be the open subgroup of $G(\RR)$ corresponding to $K$. Then $L(\RR)\cap G(\RR)'\subset G(\RR)'$ is Levi in the sense of \cite[(4.77a)]{knappvogan} (recall that $\lieq'$ is germane by the paragraph below \cite[(4.75)]{knappvogan}).
			\item Let $K^1$ be the real algebraic subgroup of $K$ generated by $K^0$ and $L\cap K$. Then its complexification $K^1_\CC$ is generated by $K^0_\CC$ and $\bar{Q}'\cap K_\CC$.
			\item We have $L\cap K^1=L\cap K$ and $(L\cap K)^0=L\cap K^0$.
			\item The canonical map
			$K(\RR)/K^1(\RR)\to K(\CC)/K^1(\CC)$
			is a bijection.
			\item The canonical injection $(L\cap K)(\CC)/(L\cap K^0)(\CC)\hookrightarrow K^1(\CC)/K^0(\CC)$ is a bijection.
		\end{enumerate}
	\end{lemma}

	\begin{proof}
		In view of the proof of Proposition \ref{prop:stable_vs_theta_stable}, one can find a cocharacter $\mu\in X_\ast(T_\CC)$ such that $\bar{Q}'=P_{G_\CC}(\mu)$. In fact, notice that $H\cap K^0_\CC= Z_{K^0_\CC}(T_\CC)=T_\CC$ in order to see that the cocharacter $\mu$ by that construction is valued in $T_\CC$ (recall that $\mu$ factors through $K^0_\CC$ since $\GG_{m,\Spec\CC}$ is connected). Since $T(\RR)$ is compact, we have $Q'=P_{G_\CC}(-\mu)$. We thus get $L'=Q'\cap \bar{Q}'$ (see \cite[Lemma 2.1.5]{conradpseudo} and \cite[Example 4.1.9]{conrad2014}). This shows (2).
		
		For (3), recall that the Levi subgroup of $G(\RR)$ attached to $\lieq'$, i.e., $L(\RR)$ is the normalizer of $\lieq'$ in $G(\RR)$.
		Then (3) is formal:
		\[\begin{split}
			L(\RR)&=L'\cap G(\RR)\\
			&=(Q'\cap \bar{Q}')\cap G(\RR)\\
			&=(N_{G(\CC)}(\lieq')\cap N_{G(\CC)}(\bar{\lieq}'))\cap G(\RR)\\
			&=(N_{G(\CC)}(\lieq')\cap G(\RR))\cap (N_{G(\CC)}(\bar{\lieq}')\cap G(\RR))\\
			&=N_{G(\CC)}(\lieq')\cap G(\RR).
		\end{split}
		\]
		We caution that the third equality holds since we are working over a field of characteristic $0$. This implies $L(\RR)\cap G(\RR)'=N_{G(\CC)}(\lieq')\cap G(\RR)'$.
		
		For (1), define $U_{G_\CC}(\mu)$ and $U_{K_\CC}(\mu)$ as in \cite[Theorem 4.1.7]{conrad2014}. Then we have
		\begin{equation}
			\bar{Q}'\cap K_\CC=(L'\cap K_\CC)\ltimes U_{K_\CC}(\mu)	\label{eq:levidec}
		\end{equation}
		by \cite[Proposition 4.1.10 1 and Theorem 4.1.7 3]{conrad2014}. Observe that $L'\cap K_\CC$ is (possibly disconnected) reductive with $(L'\cap K_\CC)^0=L'\cap K^0_\CC$ by \cite[Proposition 4.1.7 2 and Example 4.1.9]{conrad2014}. In particular, the latter equality of (5) follows. Recall also that $U_{K_\CC}(\mu)$ is unipotent by \cite[Proposition 4.1.7 4]{conrad2014}. This shows that \eqref{eq:levidec} is a Levi decomposition, and the proof of (1) is completed.
		
		We next prove (4) and (7) simultaneously. Let $K'_\CC$ denote the complex algebraic subgroup generated by $K^0_\CC$ and $\bar{Q}'\cap K_\CC$ here for convenience. It is easy to show that $K^1_\CC$ is generated by $L'$ and $K^0_\CC$. This implies $K^1_\CC\subset K'_\CC$ and (7). To see the converse containment $K'_\CC\subset K^1_\CC$, it will suffice to show $\bar{Q}'\cap K_\CC\subset K^1_\CC$. Recall the Levi decomposition \eqref{eq:levidec}. Since $U_{K_\CC}(\mu)$ is connected (\cite[Proposition 4.1.10 2]{conrad2014}), $U_{K_\CC}(\mu)$ is contained in $K^0_\CC$. This verifies $K'_\CC\subset K^1_\CC$.
		
		The first equality of (5) is evident by definition of $K^1$.
		
		Part (6) is formal:
		$K(\CC)/K^1(\CC)\cong K(\RR)/(K(\RR)\cap K^1(\CC))=K(\RR)/K^1(\RR)$.
	\end{proof}
	
	\begin{remark}
		One can regard (2) of Lemma \ref{lem:realformofK^1_C} as an algebraic analog of \cite[Proposition 4.78]{knappvogan} in the $\theta$-stable case. Statements (5)--(7) are an algebraic refinement of \cite[Proposition 5.13]{knappvogan}.
	\end{remark}
	
	\begin{remark}
		Let $G(\RR)^1$ denote the open subgroup of $G(\RR)$ corresponding to $K^1$. Then in virtue of the equality $(K^1)^1=K^1$ and Lemma \ref{lem:realformofK^1_C} (5), (7), the Levi subgroup $L(\RR)\cap G(\RR)^1$ meets every component of $G(\RR)^1$.
	\end{remark}
	
	\begin{remark}
		In \cite[Appendix B]{hechtetal}, $K^1_\CC$ was defined as the complex algebraic subgroup generated by $K^0_\CC$ and $\bar{Q}'\cap K_\CC$ so that the quotient $K/K^1_\CC$ parameterizes the connected components of the $K_\CC$-orbit containing $\bar{Q}'$. For the algebraic context, the compact group $K^1(\RR)$ appeared in \cite{knappvogan} as $(L\cap K)K_0$ in their notation.
	\end{remark}
	
	\begin{example}
		In the case of the classical Lie groups at our hand, we have $K^1_\CC\in\{K^0_\CC,K_\CC\}$ since $K_\CC$ has only two connected components. For the explicit computation of $K^1_\CC$ (and therefore of $K^1$), see section \ref{sec:disconn}. This can easily fail in general (think of $\GL_2\times\GL_3$).
	\end{example}

	\begin{corollary}[Cf.~{\cite[(1.1)]{oshima2013}}]\label{cor:duality}
		Let $G,\theta,K,T,H,Q',L, K^1$ be as in Lemma \ref{lem:realformofK^1_C}. Set \[x'=rt(Q')\in\rtype(G,K)(\CC).\]
		Let $u$ be the dimension of the unipotent radical of $\bar{Q}'\cap K^0_{\CC}$.
		
		Let $\Aa'$ be a $G_{\CC}$-equivariant tdo on $\Pp_{G_{\CC},gt(\bar{x}')}$, and $\Mm'$ be a $K_{\CC}$-equivariant integrable left $i^\cdot_{\bar{x}'} \Aa'$-connection on $\Pp^{K_\CC\mathrm{-st}}_{G_{\CC},\bar{x}'}$. Let $M'$ be the fiber of $\Mm'$ at $\bar{Q}'\in \Pp^{K_\CC\mathrm{-st}}_{G_{\CC},\bar{x}'}(\CC)$. 
		\begin{enumerate}
			\item The vector space $M'$ is naturally equipped with the structure of a $(\bar{\lieq}',L'\cap K_\CC)$-module. Moreover, its restriction to the unipotent radical of $\bar{\lieq}'$ is trivial.
			\item There exist isomorphisms
			\[\begin{split}
				H^\bullet(\Pp_{G_{\CC},gt(\bar{x}')},(i_{\bar{x}'})_+\Mm')
				&\cong
				L_{u-\bullet}P^{\lieg_\CC,K_\CC}_{\bar{\lieq}',L'\cap K_\CC}(M'\otimes_\CC \wedge^{\mathrm{top}} \lieg_\CC/\bar{\lieq}')\\
				&\cong P^{\lieg_\CC,K_\CC}_{\lieg_\CC,K^1_\CC} L_{u-\bullet} P^{\lieg_\CC,K^1_\CC}_{\bar{\lieq}',L'\cap K_\CC}(M'\otimes_\CC \wedge^{\mathrm{top}} \lieg_\CC/\bar{\lieq}')
			\end{split}\]
			of $(\lieg_\CC,K_{\CC})$-modules.
			\item There exists an isomorphism
			\[\Ff^{\lieg_\CC,K^0_\CC}_{\lieg_\CC,K^1_\CC}L_{u-\bullet} P^{\lieg_\CC,K^1_\CC}_{\bar{\lieq}',L'\cap K_\CC}(
			M'\otimes_\CC \wedge^{\mathrm{top}} \lieg_\CC/\bar{\lieq}')
			\cong L_{u-\bullet} P^{\lieg_\CC,K^0_\CC}_{\bar{\lieq}',L'\cap K^0_\CC}(\Ff^{\bar{\lieq}',L'\cap K^0_\CC}_{\bar{\lieq}',L'\cap K_\CC} M'\otimes_\CC \wedge^{\mathrm{top}} \lieg_\CC/\bar{\lieq}')\]
			of $(\lieg_\CC,K^0_\CC)$-modules.
		\end{enumerate}
	\end{corollary}
	
	We avoid using the notation of usual (co)homological induction functor $\Ll_\bullet$ here in order to clarify which open subgroup of $G^\theta$ we are considering.
	
	\begin{proof}
		Part (1) is a restatement of Theorem \ref{thm:geometric_induction} (1).
		
		See \cite[(5.8)]{knappvogan} for the second isomorphism of (2). For the first isomorphism of (2), notice that $\theta$ is a Cartan involution for $L$ since $L$ is $\theta$-stable. In view of the Zuckerman duality (\cite[Corollary 3.7]{knappvogan}), it will suffice to show that the $(\lieg_\CC,L'\cap K_\CC)$-module
		$\wedge^{\mathrm{top}}\liek_\CC/(\liek_\CC\cap \liel')$
		of \cite[Lemma 3.4]{knappvogan} is trivial. It is evident by definition that this is trivial as a $\lieg_\CC$-module. Since $(L\cap K)(\RR)$ meets every component of $(L'\cap K_\CC)(\CC)$, we may only prove that it is trivial as an $(L\cap K)(\RR)$-module. This was verified in the last paragraph of \cite[Proof of Proposition 5.93]{knappvogan}. We remark that the one-to-one onto assumption in that statement is not needed for the last argument there.
		
		For (3), it will suffice to prove an equivalence \[\Ff^{\lieg_\CC,K^0_\CC}_{\lieg_\CC,K^1_\CC}\circ
		L P^{\lieg_\CC,K^1_\CC}_{\bar{\lieq}',L'\cap K_\CC}
		\simeq LP^{\lieg_\CC,K^0_\CC}_{\bar{\lieq}',L'\cap K^0_\CC}
		\circ \Ff^{\bar{\lieq}',L'\cap K^0_\CC}_{\bar{\lieq}',L'\cap K_\CC}\]
		of the derived functors. We are reduced to the level of the abelian categories by \cite[Proposition 2.77]{knappvogan}. The corresponding isomorphism follows as
		\[\begin{split}
			\Ff^{\lieg_\CC,K^0_\CC}_{\lieg_\CC,K^1_\CC}
			\circ P^{\lieg_\CC,K^1_\CC}_{\bar{\lieq}',L'\cap K_\CC}
			&\cong \Ff^{\lieg_\CC,K^0_\CC}_{\lieg_\CC,K^1_\CC} \circ 
			P^{\lieg_\CC,K^1_\CC}_{\lieg_\CC,L'\cap K_\CC} \circ \ind^{\lieg_\CC}_{\bar{\lieq}'}\\
			&\cong P^{\lieg_\CC,K^0_\CC}_{\lieg_\CC,L'\cap K^0_\CC}\circ 
			\Ff^{\lieg_\CC,L'\cap K^0_\CC}_{\lieg_\CC,L'\cap K^1_\CC}
			\circ \ind^{\lieg_\CC}_{\bar{\lieq}'}\\
			&\cong P^{\lieg_\CC,K^0_\CC}_{\lieg_\CC,L'\cap K^0_\CC}
			\circ \ind^{\lieg_\CC}_{\bar{\lieq}'}
			\circ \Ff^{\bar{\lieq}',L'\cap K^0_\CC}_{\bar{\lieq}',L'\cap K^1_\CC}\\
			&\cong P^{\lieg_\CC,K^0_\CC}_{\bar{\lieq}',L'\cap K^0_\CC}\circ
			\Ff^{\bar{\lieq}',L'\cap K^0_\CC}_{\bar{\lieq}',L'\cap K_\CC}
		\end{split}\]
		(see Lemma \ref{lem:realformofK^1_C} (7) and \cite[Proposition 5.14 (b)]{knappvogan} for the second isomorphism).
	\end{proof}
	
	\begin{remark}
		One can prove (3) directly from the geometric construction.
	\end{remark}
	
	\begin{remark}
		Although we have a similar induction formula of the second isomorphism of (2) in the setting of Theorem \ref{thm:geometric_induction}, we omit it in order to avoid conflict of notations of $K^1$.
	\end{remark}
	
	For the relation with the right derived functor modules, let us note:
	
	\begin{definition}\label{def:ess_uni}
		We say a $(\lieg_\CC,K(\RR))$-module $X$ is essentially unitarizable if there exists a one-dimensional $(\lieg_\CC,K(\RR))$-module $\tau$ such that $X\otimes_\CC\tau$ is unitarizable.
	\end{definition}
	
	\begin{example}\label{ex:split_central}
		We call a linear functional on $\lieg_\CC$ a split central character if it vanishes on $\liek_\CC$ and the derived subalgebra of $\lieg_\CC$. This is a $(\lieg_\CC,K(\RR))$-module for a character of the component group $\pi_0(K(\RR))$ of the compact Lie group $K(\RR)$.
	\end{example}
	
	\begin{corollary}\label{cor:irreducible_unitary}
		Consider the setting of Corollary \ref{cor:duality}. Suppose that $\lieh_\CC$ acts on $M'$ by $\lambda$ (see the proof of Theorem \ref{thm:geometric_induction} (1)). Choose a Borel subgroup $B'$ satisfying $H_\CC\subset B'\subset Q'$, and $\rho$ be the half sum of the corresponding positive roots. Let $\tau$ be a one-dimensional $(\lieg_\CC,K(\RR))$-module. Write $\nu$ for the corresponding character of $\lieg_\CC$. Moreover, assume the following conditions:
		\begin{enumerate}
			\renewcommand{\labelenumi}{(\roman{enumi})}
			\item No $H_\CC$-root $\alpha$ of the unipotent radical $\lieu'$ of $\lieq'$ satisfies
			\[\langle \alpha^\vee, \lambda+\rho+\nu|_{\lieh_\CC}\rangle\in\left(-\infty,-1\right]\cup\{0\}.\]
			Here the canonical pairing is taken in $\lieh_\CC$ and its dual by regarding coroots and roots as elements of $\lieh_\CC$ and $\lieh^\vee_\CC$ respectively.
			\item The character $\lambda+\nu|_{\lieh_\CC}$ is unitary, i.e., $\lambda|_{\lieh}+\nu|_{\lieh}$ is valued in $\RR\sqrt{-1}$.
			\item The $(i^\cdot_x\Aa,K_\CC)$-module $\Mm'$ is irreducible, equivalently, $M'$ is an irreducible $(L'\cap K_\CC)$-module.
		\end{enumerate}
		Then:
		\begin{enumerate}
			\item The cohomology $H^i(\Pp_{G_{\CC},gt(\bar{x}')},(i_{\bar{x}'})_+\Mm')$ vanishes for all $i>0$.
			\item We have an isomorphism
			\[\Gamma(\Pp_{G_{\CC},gt(\bar{x}')},(i_{\bar{x}'})_+\Mm')
			\cong R^{u}I^{\lieg_\CC,K_\CC}_{\lieq',L'\cap K_\CC}(M'\otimes_\CC \wedge^{\mathrm{top}} \lieg_\CC/\bar{\lieq}')\]
			of $(\lieg_\CC,K_\CC)$-modules. Moreover, $\Gamma(\Pp_{G_{\CC},gt(\bar{x}')},(i_{\bar{x}'})_+\Mm')$ is irreducible, admissible, and essentially unitarizable as a $(\lieg_\CC,K(\RR))$-module.
		\end{enumerate}
	\end{corollary}
	
	\begin{proof}
		In view of Corollary \ref{cor:duality}, we may replace $M'$ by $M'\otimes_\CC \pi$ to assume that $\pi$ is trivial. Then (1) and the isomorphism in (2) follow from \cite[Theorems 8.21]{knappvogan}. The rest follows from \cite[Theorems 5.64, 9.1, and Corollary 8.28]{knappvogan}. We remark that $M'$ is a unitarizable $(\liel_\CC,(L\cap K)(\RR))$-module by (ii).
	\end{proof}
	
	\begin{remark}
		In our context, we can numerically determine the unitary chacaracters $\lambda$ of $\lieh$ as follows: Firstly, compute the complex conjugate action on $X^\ast(H_\CC)$ (see \cite[section 4.1]{hayashilinebdl} for explicit formulas in the case of classical Lie groups). Secondly, extend it $\RR$-linearly. Set
		\[\begin{array}{cc}
			X^\ast(H_\CC)^\sigma_{\RR}=\{\lambda\in X^\ast(H_\CC)_\RR:~\bar{\lambda}=\lambda\},
			&X^\ast(H_\CC)^{-\sigma}_{\RR}=\{\lambda\in X^\ast(H_\CC)_\RR:~\bar{\lambda}=-\lambda\}.
		\end{array}\]
		Then the set of unitary characters of $\lieh$ is identified with
		$\sqrt{-1}X^\ast(H_\CC)^\sigma_{\RR}\oplus X^\ast(H_\CC)^{-\sigma}_{\RR}$ in $X^\ast(H_\CC)_{\CC}$. In fact, let us denote the differential of $\theta$ by the same symbol. Set
		\[\liea=\{x\in\lieh:~\theta(x)=-x\}.\]
		Then we have natural identifications
		\[\begin{array}{cc}
			X^\ast(H_\CC)^\sigma_{\RR}\cong\Hom_{\RR}(\liea,\RR),
			&X^\ast(H_\CC)^{-\sigma}_{\RR}\cong \sqrt{-1}\Hom_{\RR}(\liet,\RR).
		\end{array}\]
	\end{remark}
	
	\begin{example}
		Put $G=\GL_{2n-1},\GL_{2n}$ with $n\geq 1$. Under the notations of \cite[Examples 4.1.1 and 4.1.2]{hayashilinebdl}, the sets of unitary characters of $\lieh$ are identified with
		\[\sqrt{-1}e_{2n-1}\oplus\oplus_{i=1}^{n-1} \RR\sqrt{-1}(e_{2i-1}+e_{2i})\oplus \oplus_{i=1}^{n-1} \RR(e_{2i-1}-e_{2i})\]
		\[\oplus_{i=1}^n \RR\sqrt{-1}(e_{2i-1}+e_{2i})\oplus \oplus_{i=1}^n \RR(e_{2i-1}-e_{2i})\]
		respectively.
	\end{example}
	
	We are now able to deduce a general statement of how to obtain forms of (co)homogically induced modules in a geometric way: Suppose that we are given a commutative diagram
	\[\begin{tikzcd}
		\Spec\CC\ar[r]\ar[d, "s"']&\Spec\RR\ar[d]\\
		\tilde{S}\ar[r]&S
	\end{tikzcd}\]
	of $\ZZ[1/2]$-schemes with $s$ flat and $\tilde{S}$ locally Noetherian of finite Krull dimension. Let $G$ be a reductive group scheme, equipped with an involution $\theta$ with $\theta\times_S\Spec\RR$ Cartan. Let $K\subset G^\theta$ be an open and closed subgroup scheme with $\pi_0(K)$ finite \'etale. Let $x\in\rtype(G,K)(\tilde{S})$, and $\Aa$ be a $G\times_S \tilde{S}$-equivariant tdo on $\Pp_{G\times_S \tilde{S},gt(x)}$.

	\begin{corollary}\label{cor:sformofcohohomologicalinduction1}
		Choose a maximal torus $T\subset K\times_S \Spec\RR$ and a complex stable parabolic subgroup $Q'\subset G\times_S \Spec\CC$ containing $T\times_{\Spec\RR}\Spec\CC$ such that $rt(\bar{Q}')=x|_{\Spec\CC}$. Let $u$ be the dimension of the unipotent radical of $\lieq'\cap\liek_\CC$.
		
		Let $\Mm$ be a $K\times_S\tilde{S}$-equivariant integrable left $i^\cdot_{x}\Aa$-connection. Let $M'$ denote the geometric fiber of $\Mm$ at $\bar{Q}'\in \Pp^{K\times_S \tilde{S}\mathrm{-st}}_{G\times_S\tilde{S},x}(\CC)$. Then $R^\bullet (p_{G\times_S \tilde{S},x})_\ast(i_{x})_+\Mm$ exhibits an $\tilde{S}$-form of \[L_{u-\bullet}P^{\lieg_\CC,K_\CC}_{\bar{\lieq}',L'\cap K_\CC}(M'\otimes_\CC \wedge^{\mathrm{top}} \lieg_\CC/\bar{\lieq}').\]
	\end{corollary}

	This consequence asserts that if the objects $G$, $K$, $x$, $\Aa$, and $\Mm$ are defined over a smaller subring $k$ of $\CC$, we obtain a $k$-form of the (co)homologically induced module. By certain descent arguments, $k$ may be subrings of $\RR$. In particular, we can put $\tilde{S}=S$ in such case. For examples of $x$, see section \ref{sec:halfintegralformofrelativetype}. A typical example of $\Aa$ is the twist of $\D_{\Pp_{G\times_S \tilde{S},gt(x)}}$ by a $G\times_S \tilde{S}$-equivariant line bundle $\Ll$. In this case, one can put $\Mm=i^\ast_{\bar{x}}\Ll$. The construction and classification of $\Ll$ were partially established in \cite{hayashilinebdl} using Borel-Tits' obstruction class $\beta_\lambda$ (see \cite[Theorem B and section 4]{hayashilinebdl} for general statements and examples respectively). 
	
	We end this section with the following conclusion for a typical situation where $\tilde{S}=S$ after Galois descent:
	
	\begin{corollary}\label{cor:arithmeticmodels}
		Consider the Galois double covering
		\[S'=\Spec\ZZ\left[1/2,\sqrt{-1}\right]\to \Spec\ZZ\left[1/2\right]=S.\]
		Let $(G,K)$ be the symmetric pair of the standard $\ZZ\left[1/2\right]$-form of a classical connected Lie group and its maximal compact subgroup attached to the involution $\theta$ of $G$ in \cite[section 3.3]{hayashilinebdl}. Let $T$ be the maximal torus of $K$ in \cite[section 3.4]{hayashilinebdl}. Set $H=Z_G(T)$. Let $\Pi$ be the positive system of $(K\times_S S',T\times_S S')$ in section \ref{sec:halfintegralformofrelativetype}. Take a $\theta$-stable parabolic subset of $\Delta(G\times_S S',H\times_S S')$ such that $\Phi|_{T\times_S S'}$ contains $\Pi$. Let $w\in N_{K(S')}(T\times_S S')$ be the representative of the longest element of $W(K,T)(S')$ defined in \cite[section 4.1]{hayashilinebdl}. Let $Q'$ be the parabolic subgroup of $G\times_S S'$ corresponding to $\Phi$. Let $\lambda$ be a character of $H\times_S S'$. Suppose that the following conditions are satisfied:
		\begin{enumerate}
			\renewcommand{\labelenumi}{(\roman{enumi})}
			\item $w\Phi=-\Phi$.
			\item $\bar{\lambda}=w^{-1}\lambda$.
			\item $\langle \alpha^\vee,\lambda\rangle=0$ for $\alpha\in \Phi\cap (-\Phi)$.
			\item $\lambda(w\bar{w})=1$.
		\end{enumerate}
		Then:
		\begin{enumerate}
			\renewcommand{\labelenumi}{(\arabic{enumi})}
			\item The relative parabolic type $rt(Q')$ descends an element $x\in\rtype(G,K)(S)$.
			\item The associated bundle $(G\times_S S')\times^{\bar{Q}'} \lambda$ descends to a $G$-equivariant line bundle $\Ll$ on $\Pp_{G,gt(x)}$.
			\item $\Gamma(\Pp_{G,gt(x)},(i_x)_+i^\ast_x\Ll)$ exhibits a $\ZZ\left[1/2\right]$-form of $A_{\lieq'_\CC}(\lambda)$.
		\end{enumerate}
	\end{corollary}
	
	\subsection{Application to $\GL_n$ over totally real and CM number fields}\label{sec:gl_n}
	
	This section aims to give arithmetic models of cohomological irreducible unitary representations of $\GL_n$ over totally real and CM number fields towards applications to integral models of automorphic representations.
	
	\subsubsection{Potentially defined characters}\label{sec:potential_character}
	
	To obtain cohomologically induced modules, we want line bundles on partial flag varieties. This time, we hope to work with the Weil restriction of a partial flag variety of $\GL_n$ over a number field to $\QQ$ for simplicity. In fact, the Weil restriction does not respect disjoint unions, so that there are more partial flag varieties of the Weil restriction of $\GL_n$ after base changes.
	
	For broader treatment of characters, we aim to work with cohomological representations over the complex numbers with all coefficients and their arithmetic models. Since $\GL_n$ is split over $\ZZ$, it is enough to study fractionality of the characters of a certain splitting base change of the Weil restriction of a split torus (see \cite[Example 2.3.2]{hayashilinebdl}).
	
	Before we perform this, let us collect some basic facts from algebraic number theory. Let $E/F$ be an extension of number fields. Write $E^{\Gal}$ for a normal hull of $E$ over $F$ and let $\Gamma$ denote the Galois group of $E^{\Gal}/F$. Write $\Hom_F(E,E^{\Gal})$ for the set of $F$-algebra homomorphisms from $E$ to $E^{\Gal}$.
	
	Firstly, let us record a well-known fact for convenience:
	
	\begin{lemma}\label{lem:ramification}
		A nonzero prime ideal $\liep\subset\OO_F$ is ramified in $E/F$ if and only it is ramified in $E^{\Gal}/F$.
	\end{lemma}
	
	\begin{proof}
		It follows by definition that $E^{\Gal}$ is the composition of all embeddings $E\hookrightarrow E^{\Gal}$ over $F$ in $E^{\Gal}$. Then the equivalence follows from \cite[Chapter III, $\S$ 2, Theorem 22]{frohlichtaylor} and \cite[Theorem 5.12]{mollin2011}.
	\end{proof}
	
	\begin{lemma}\label{lem:fractional_splitting}
		Let $f$ be an element of the intersection of the prime ideals of $\OO_F$ ramified in $E/F$. Then the inclusion $\OO_F\left[1/f\right]\subset \OO_{E^{\Gal}}\left[1/f\right]$ is finite \'etale. Moreover, the well-defined map
		\[\phi:\OO_E\left[1/f\right]\otimes_{\OO_F\left[1/f\right]} \OO_{E^{\Gal}}\left[1/f\right]
		\to \prod_{\sigma\in \Hom_F(E,E^{\Gal})}\OO_{E^{\Gal}}\left[1/f\right];~
		a\otimes b\mapsto (\sigma(a)b)_\sigma
		\]
		is an isomorphism of $\OO_{E^{\Gal}}$-algebras. In particular, if $E/F$ is Galois then so is
		\[\OO_F\left[1/f\right]\subset \OO_{E}\left[1/f\right].\]
	\end{lemma}
	
	\begin{proof}
		It is evident by definitions of $\OO_E$ and $\OO_{E^{\Gal}}$ that each map $\sigma\in\Hom_F(E,E^{\Gal})$ restricts to a homomorphism $\OO_E\to \OO_{E^{\Gal}}$. Hence $\phi$ is well-defined.
		
		It is an elementary fact that $\OO_F\left[1/f\right]\subset \OO_{E^{\Gal}}\left[1/f\right]$ is finite and faithfully flat (see \cite[Chapter II, $\S$ 1, Theorem 5 and its corollary]{frohlichtaylor} if necessary). In view of the faithfully flat descent, the proof will be completed by showing that $\phi$ is an isomorphism.
		
		Since $\OO_{E}$ is a locally free of finite rank as an $\OO_F$-module, it is sufficient to treat such elements $f$ with the property that $\OO_{E}\left[1/f\right]$ is free over $\OO_F\left[1/f\right]$. Fix a basis $a_1,\dots,a_m\in\OO_{E}\left[1/f\right]$. Then $\phi$ is the composition of the isomorphism
		\[\OO_E\left[1/f\right]\otimes_{\OO_F\left[1/f\right]} \OO_{E^{\Gal}}\left[1/f\right]
		\cong \OO_E\left[1/f\right]^m;~\sum_{i=1}^m a_i\otimes b_i\mapsto \left(\begin{array}{c}
			b_1 \\
			b_2 \\
			\vdots \\
			b_m
		\end{array}\right)
		\]
		and the $\OO_{E^{\Gal}}\left[1/f\right]$-linear endomorphism of $\OO_{E^{\Gal}}\left[1/f\right]^m$ determined by the matrix \[A\coloneqq(\sigma(a_i))_{\sigma,i}.\]
		It remains to prove that $A$ is invertible.
		
		The following calculation is standard:
		\[(\det A)^2=
		\det\left(A\cdot A^T\right)
		=\det\left(\sum_{\sigma}\sigma\left(a_ia_j\right)\right)_{i,j}.\]
		This value is known as the discriminant $d\left(a_1,\dots,a_m\right)$ of the chosen basis $\OO_F\left[1/f\right]$-module $\OO_E\left[1/f\right]$ (see \cite[Chapter III, $\S$ 2 and Chapter I, (1.24)]{frohlichtaylor}). Lemma \ref{lem:ramification} implies that every prime ideal of $\OO_F\left[1/f\right]$ is unramified in $E^{\Gal}$. Hence in view of \cite[Chapter III, $\S$ 2, Theorem 22]{frohlichtaylor}, $d(a_1,\dots,a_m)$ and therefore $\det A$ are units of $\OO_F\left[1/f\right]$. This shows that $A$ is invertible as desired.
	\end{proof}
	
	Let us fix $f$ as in Lemma \ref{lem:fractional_splitting} above.
	
	\begin{example}\label{ex:minimal_f_killing_ramified_primes}
		One can find an element $f\in\OO_F$ which defines the principal open subset of $\Spec \OO_F$ consisting of the zero and unramified prime ideals in $E/F$ (see \cite[Chapter 2, Exercise 3.19]{liu2002}). In this case, inverting $f$ corresponds to tensoring with the coordinate ring of this principal open subset over $\OO_F$, hence is independent of the specific choice of $f$. If $F=\QQ$, we can take $f=d_E$ by \cite[Chapter III, $\S$ 2, Theorem 22]{frohlichtaylor} since $\ZZ$ is a PID.
	\end{example}
	
	\begin{lemma}
		In the setting of Lemma \ref{lem:fractional_splitting}, take an intermediate extension $E^{\Gal}/E'/F$. Write $\Gamma'\subset \Gamma$ for the corresponding subgroup. Then 
		$\OO_{E'}\left[1/f\right]\subset \OO_{E^{\Gal}}\left[1/f\right]$ is a Galois extension of Galois group $\Gamma'$.
	\end{lemma}
	
	\begin{proof}
		Recall that $E^{\Gal}/E'$ is a Galois extension. Hence in view of the former lemma, it suffices to see that $f$ is contained in all ramified prime ideals of $\OO_{E'}$ in $E^{\Gal}$. Take any ramified prime ideal $\mathfrak{P}\subset \OO_{E^{\Gal}}$ in $E^{\Gal}/E'$. Set $\liep=\mathfrak{P}\cap\OO_{E'}$. In view of the tower formula of the ramification indices 
		(\cite[Chapter III, (1.13.a)]{frohlichtaylor}), $\mathfrak{P}$ is ramified in $E^{\Gal}/F$. Lemma \ref{lem:ramification} thus implies that $\liep$ is a ramified prime in $E^{\Gal}/F$. Our hypothesis on $f$ implies $f\in \liep$. It follows by definition of $\liep$ that $f\in\mathfrak{P}$ as desired.
	\end{proof}
	
	\begin{theorem}\label{thm:minimal_form}
		In the situation of Lemma \ref{lem:fractional_splitting}, let $H_{\spl}$ be a split torus over $\OO_E\left[1/f\right]$. Take $\lambda\in X^\ast((\res_{E/F} H_{\spl})\otimes_F E^{\Gal})$. Set
		\[\begin{array}{ccc}
			\Gamma(\lambda)\coloneqq \{\sigma\in\Gamma:~{}^\sigma\lambda=\lambda\},
			&F(\lambda)\coloneqq (E^{\Gal})^{\Gamma(\lambda)},
			&\OO(\lambda)\coloneqq \OO_{F(\lambda)}.
		\end{array}\]
		Then there is a unique character of $(\res_{\OO_E\left[1/f\right]/F\left[1/f\right]} H_{\spl})\otimes_{\OO_F\left[1/f\right]}
		\OO_{E(\lambda)}\left[1/f\right]$ such that its base change to $E^{\Gal}$ is $\lambda$.
	\end{theorem}
	
	\begin{proof}
		In view of Lemma \ref{lem:fractional_splitting}, $\res_{\OO_E\left[1/f\right]/\OO_F\left[1/f\right]} H_{\spl}$ is a torus which is split over $\OO_{E^{\Gal}}\left[1/f\right]$. Therefore the base change to $E$ determines a $\Gamma$-equivariant isomorphism
		\begin{flalign*}
			&X^\ast((\res_{\OO_E\left[1/f\right]/\OO_F\left[1/f\right]} H)\otimes_{\OO_F\left[1/f\right]}
			\OO_{E(\lambda)}\left[1/f\right]
			\otimes_{\OO_{E(\lambda)}\left[1/f\right]} \OO_{E^{\Gal}}\left[1/f\right])\\
			&\cong X^\ast((\res_{\OO_E\left[1/f\right]/\OO_F\left[1/f\right]} H)\otimes_{\OO_F\left[1/f\right]} \OO_{E^{\Gal}}\left[1/f\right])\\
			&\cong X^\ast((\res_{E/F} H_{\spl})\otimes_F E^{\Gal}).
		\end{flalign*}
		Since $\Gamma(\lambda)$ fixes $\lambda$, we get the desired character of $(\res_{\OO_E\left[1/f\right]/F\left[1/f\right]} H_{\spl})\otimes_{\OO_F\left[1/f\right]}
		\OO_{E(\lambda)}\left[1/f\right]$. The uniqueness follows from the above bijection.
	\end{proof}
	
	In this way, we determined the rationality of $\lambda$, and proved that it has a unique lift to a character over $\OO_{E(\lambda)}\left[1/f\right]$. This lift is essentially independent of the choice of $f$. I.e., such a character is obtained by the base change of that in the case of $f$ in Example \ref{ex:minimal_f_killing_ramified_primes} by the naturality of base changes, and the lift for $f$ in Example \ref{ex:minimal_f_killing_ramified_primes} is independent of the choice of $f$ by the proof (recall also the third line of Example \ref{ex:minimal_f_killing_ramified_primes}).

	\subsubsection{Totally real case}\label{sec:totally_real}
	
	Let $E$ be a totally real number field. To save space, write
	\[\begin{array}{cc}
		k=\ZZ\left[1/2d_E\right],&k'=\OO_E\left[1/2d_E\right].
	\end{array}\]
	Recall that $\OO_E\left[1/d_E\right]$ is a finite \'etale $\ZZ\left[1/d_E\right]$-algebra. Hence we obtain a reductive group scheme $\res_{k'/k}\GL_n$, which is a model of $\res_{E/\QQ} \GL_n$. Similarly, we have a smooth model $\res_{k'/k}\Oo(n)$ of a maximal compact subgroup of $(\res_{E/\QQ}\GL_n)(\RR)$. Indeed, this is the fixed point subgroup scheme by the involution $\theta\coloneqq((-)^T)^{-1}$ on $\res_{k'/k}\GL_n$.
	
	We wish to construct line bundles on partial flag schemes (cf.~Corollary \ref{cor:irreducible_unitary}). For this, we would firstly like to think of characters of the restriction of split and fundamental Cartan subgroups. The classification was essentially done for the split Cartan subgroup in the previous section. Here we aim to compare characters of the two tori, based on the argument of \cite[Examples 3.4.1 and 3.4.2]{hayashilinebdl} for the case $E=\QQ$. I.e., let $T\subset\GL_n$ be the torus over $\ZZ\left[1/2\right]$ as in \cite[Examples 3.4.1 and 3.4.2]{hayashilinebdl}. Set $H=Z_G(T)$. We take the base change to regard it as a torus over $\QQ$. Let $H_{\spl}$ be the split maximal torus of $\GL_n$ of diagonal matrices over $\QQ$. Let $E^{\Gal}$ be a normal hull of the extension $E/\QQ$. We note that $E^{\Gal}$ is still totally real by construction since $E$ is totally real. In particular, $E^{\Gal}(\sqrt{-1})$ is a CM number field. Let $g_n\in \GL_n(\QQ(\sqrt{-1}))$ as in \cite[section 1.5]{hayashilinebdl}. Then we have $g_n H g^{-1}_{n}=H_{\spl}$ over $\QQ(\sqrt{-1})$. This leads us to a commutative diagram
	\[\begin{tikzcd}
		X^\ast((\res_{E/\QQ} H_E)_{E^{\Gal}(\sqrt{-1})})\ar[r]\ar[d, "g_n"']
		&\prod_{\sigma\in\Hom_\QQ(E,E^{\Gal}(\sqrt{-1}))}
		X^\ast(H\otimes_\QQ E^{\Gal}(\sqrt{-1}))\ar[d, "g_n"]\\
		X^\ast((\res_{E/\QQ} H_{\spl,E})_{E^{\Gal}(\sqrt{-1})})\ar[r]
		&\prod_{\sigma\in\Hom_\QQ(E,E^{\Gal}(\sqrt{-1}))}
		X^\ast(H_{\spl}\otimes_{\QQ} E^{\Gal}(\sqrt{-1}))\\
		X^\ast((\res_{E/\QQ} H_{\spl,E})_{E^{\Gal}})\ar[r]\ar[u]
		&\prod_{\sigma\in\Hom_\QQ(E,E^{\Gal})}
		X^\ast(H_{\spl}\otimes_{\QQ} E^{\Gal})\ar[u]
	\end{tikzcd}\]
	of isomorphisms of abelian groups. Here the upper and middle horizontal arrows are (resp.~the bottom arrow is) induced from the canonical isomorphism
	\[E\otimes_\QQ E^{\Gal}(\sqrt{-1}) \cong\prod_{\sigma\in\Hom_\QQ(E,E^{\Gal}(\sqrt{-1}))}
	E^{\Gal}(\sqrt{-1})
	\]
	\[(\mathrm{resp.~}E\otimes_\QQ E^{\Gal} \cong\prod_{\sigma\in\Hom_\QQ(E,E^{\Gal})}E^{\Gal});\]
	the upper vertical arrows are obtained from the conjugation by $g_n\in \GL_n(E^{\Gal}(\sqrt{-1}))$; the lower vertical arrows are obtained by the base change $-\otimes_{E^{\Gal}} E^{\Gal}(\sqrt{-1})$. Henceforth we will identify the upper two vertices and the four vertices of the lower square respectively. We also recall that one obtains a free basis of $X^\ast(H_{E^{\Gal}(\sqrt{-1})})$ by the pullback of the standard one in $X^\ast(H_{\spl,E^{\Gal}(\sqrt{-1})})$. Under these identifications, we will denote their elements as
	\[\lambda=\sum_\sigma\lambda_{\sigma}\sigma=\sum_\sigma(\lambda_{i,\sigma})\sigma\]
	in order to emphasize the summands, where $\sigma$ under $\sum$ runs through \[\Hom_\QQ(E,E^{\Gal})=\Hom_\QQ(E,E^{\Gal}(\sqrt{-1})).\]
	
	Let us keep on discussing characters of $(\res_{E/\QQ} H_{E})_{E^{\Gal}(\sqrt{-1})}$. The involution on
	\[\prod_{\sigma\in\Hom_\QQ(E,E^{\Gal}(\sqrt{-1}))}
	X^\ast(H_{E^{\Gal}(\sqrt{-1})})\]
	induced from $\theta$ is $\sigma$-componentwisely given by 
	\begin{flalign*}
		&(\lambda_{1,\sigma},\lambda_{2,\sigma},\ldots,\lambda_{n,\sigma})\sigma\\
		&\mapsto
		\begin{cases}
			(-\lambda_{2,\sigma},-\lambda_{1,\sigma},-\lambda_{4,\sigma},-\lambda_{3,\sigma},
			\ldots,-\lambda_{n,\sigma},-\lambda_{n-1,\sigma})\sigma&(n:~\mathrm{even})\\
			(-\lambda_{2,\sigma},-\lambda_{1,\sigma},-\lambda_{4,\sigma},-\lambda_{3,\sigma},
			\ldots,-\lambda_{n-1,\sigma},-\lambda_{n-2,\sigma},\lambda_{n,\sigma})\sigma
			&(n:~\mathrm{odd})
		\end{cases}
	\end{flalign*}
	(cf.~Strategy \ref{strategy}, \cite[Examples 4.1.1 and 4.1.2]{hayashilinebdl}).
	It follows by construction that
	\[\Delta^+\coloneqq \coprod_{\sigma}
	\Delta^+(\GL_n, H_{E^{\Gal}(\sqrt{-1})})\sigma
	\subset \prod_{\sigma}
	X^\ast(H_{E^{\Gal}(\sqrt{-1})})\]
	gives a $\theta_{E^{\Gal}(\sqrt{-1})}$-stable positive system, where
	$\Delta^+(\GL_n,H_{E^{\Gal}(\sqrt{-1})})$
	is the positive system in \cite[Examples 4.1.1 and 4.1.2]{hayashilinebdl}.
	
	To work with partial flag schemes, recall the sequence of bijections
	\[\begin{split}
		\rtype(\res_{k'/k}\GL_n,\res_{k'/k}\Oo(n))(k)
		&\cong (\res_{k'/k} \rtype(\GL_n,\Oo(n)))(k)\\
		&\cong \rtype(\GL_n,\Oo(n))(k')\\
		&\cong \SPS^\theta,
	\end{split}\]
	where $\SPS^\theta$ is as in section \ref{sec:stdthetastablepss}. Henceforth fix $\Phi\in\SPS^\theta$. Write $x\in\rtype(\GL_n,\Oo(n))(k')$ for the corresponding relative parabolic type. Then the attached partial flag scheme of $\res_{k'/k}\GL_n$ is expressed as
	$\res_{k'/k}\Pp_{\GL_n,gt(x)}$. Let $Q'_{k'\left[\sqrt{-1}\right]}\subset\GL_n$ be the parabolic subgroup over $k'\left[\sqrt{-1}\right]$ corresponding to $\Phi$. Its conjugate parabolic subgroup $\bar{Q}'_{k'\left[\sqrt{-1}\right]}$ with respect to the involution $\bar{\ }:\sqrt{-1}\mapsto -\sqrt{-1}$ corresponds to $-\Phi$. Note that they have the same parabolic type since $\GL_n$ is split over $k'$. We next observe that $g_n \bar{Q}'_{k'\left[\sqrt{-1}\right]}g^{-1}_n$ contains the diagonal subgroups by definition of $g_n$ (cf.~\cite[Examples 3.4.1 and 3.4.2]{hayashilinebdl}). This parabolic subgroup scheme admits a unique $k$-form which we will denote by $Q_k$. It is evident by definition that its base change $Q_{k'\left[\sqrt{-1}\right]}$ to $k'\left[\sqrt{-1}\right]$ belongs to $(\res_{k'/k} \Pp_{\GL_n,gt(x)})(k)$.
	
	Take a character
	$\lambda\in X^\ast((\res_{E/\QQ} H)_{E^{\Gal}(\sqrt{-1})})$. We will write $\lambda=(\lambda_\sigma)_{\sigma}=((a_{i,\sigma}))$ through the identification
	\[X^\ast((\res_{E/\QQ} H)_{E^{\Gal}(\sqrt{-1})})
	\cong \prod_{\sigma\in\Hom_\QQ(E,E^{\Gal}(\sqrt{-1}))} X^\ast(H_{E^{\Gal}(\sqrt{-1})})
	\cong \prod_{\sigma\in\Hom_\QQ(E,E^{\Gal}(\sqrt{-1}))} \ZZ^n.\]
	Apply the argument of the former section to
	$g_n\lambda\in X^\ast((\res_{E/\QQ} H_{\spl})_{E^{\Gal}})$
	to lift $g_n\lambda$ to a character over $\OO(g_n\lambda)\left[1/2d_E\right]$ which we will denote by the same symbol. Apply $g^{-1}_n$ to obtain a character over $\OO(g_n\lambda)\left[1/2d_E,\sqrt{-1}\right]$. To relate the numerical description of $\lambda$ with the real setting, we fix an embedding $E^{\Gal}\hookrightarrow\RR$, rather than $E^{\Gal}\supset\QQ(g_n\lambda)\hookrightarrow\RR$. We then extend it to $E^{\Gal}(\sqrt{-1})\hookrightarrow \CC$ by $\sqrt{-1}\mapsto \sqrt{-1}$. We will identify $X^\ast((\res_{E/\QQ} H)_{E^{\Gal}(\sqrt{-1})})$ with $X^\ast((\res_{E/\QQ} H)_\CC)$ through the base change. This is needed particularly to parameterize the complex representations $V$ in Definitions \ref{def:A(V)_totally_real_even} and \ref{def:A(V)_totally_real_odd}.
	
	\begin{theorem}\label{thm:gl_n_totally_real}
		Assume that $\langle \alpha^\vee,\lambda_{\sigma}\rangle=0$ for $\alpha\in \Phi\cap (-\Phi)$ and $\sigma\in\Hom_\QQ(E,E^{\Gal}(\sqrt{-1}))$.
		\begin{enumerate}
			\item One can extend $g_n\lambda$ to a character of $(\res_{k'/k} Q)_{\OO(g_n\lambda)\left[1/2d_E\right]}$,
			which we will denote by the same symbol.
			\item Define a $(\res_{k'/k} \GL_n)_{\OO(g_n\lambda)\left[1/2d_E\right]}$-equivariant line bundle on
			\[(\res_{k'/k}\Pp_{\GL_n,gt(x)})_{\OO(g_n\lambda)\left[1/2d_E\right]}\]
			by
			\[\Ll\coloneqq (\res_{k'/k} \GL_n)_{\OO(g_n\lambda)\left[1/2d_E\right]}
			\times^{(\res_{k'/k} Q)_{\OO(g_n\lambda)\left[1/2d_E\right]}} g_n\lambda.\]
			Then
			\[\Gamma((\res_{k'/k}\GL_n/Q)_{\OO(g_n\lambda)\left[1/2d_E\right]},((\res_{k'/k} i_x)_{\OO(g_n\lambda)\left[1/2d_E\right]})_+((\res_{k'/k} i_x)_{\OO(g_n\lambda)\left[1/2d_E\right]})^\ast\Ll)\]
			exhibits a projective $\OO(g_n\lambda)\left[1/2d_E\right]$-form of the cohomologically induced module in Corollary \ref{cor:duality} (2) with $G=(\res_{E/\QQ}\GL_n)_{\RR}$, $M'=\lambda$, and $\bullet=0$. The parabolic subgroup in Corollary \ref{cor:duality} is
			$(\res_{k'\left[\sqrt{-1}\right]/k\left[\sqrt{-1}\right]} Q'_{k'\left[\sqrt{-1}\right]})_\CC$.
		\end{enumerate}
	\end{theorem}
	
	Let us introduce the following notation:
	
	\begin{definition}\label{def:A(V)_totally_real_even}
		Consider $\GL_{2n}$ with $n\geq 1$. For each $\sigma\in \Hom_\QQ(E,E^{\Gal}(\sqrt{-1}))$, assume
		\begin{enumerate}
			\renewcommand{\labelenumi}{(\roman{enumi})}
			\item $a_{1,\sigma}>a_{3,\sigma}>\cdots>a_{2n-1,\sigma}$,
			\item $a_{1,\sigma}+a_{2,\sigma}=a_{3,\sigma}+a_{4,\sigma}=\cdots=a_{2n-1,\sigma}+a_{2n,\sigma}$, and
			\item $a_{2n-1,\sigma}>a_{2n,\sigma}$.
		\end{enumerate}
		Let $V$ be the complex irreducible representation of $(\res_{E/\QQ} \GL_{2n})_\CC$
		with lowest weight $-\lambda$ with respect to $\Delta^+$. Write $\OO(V)\coloneqq \OO(-g_n\lambda)$\footnote{Since $\res_{E/\QQ}\GL_{2n}$ is quasi-split, $\QQ(-g_n\lambda)$ agrees with the field of rationality of $V$.}. Then we denote the $\OO(V)\left[1/2d_E\right]$-module in Theorem \ref{thm:gl_n_totally_real} by $A(V)_{\OO(V)\left[1/2d_E\right]}$.
	\end{definition}
	
	\begin{definition}\label{def:A(V)_totally_real_odd}
		Consider $\GL_{2n+1}$ with $n\geq 0$. For each $\sigma\in \Hom_\QQ(E,E^{\Gal}(\sqrt{-1}))$, assume
		\begin{enumerate}
			\renewcommand{\labelenumi}{(\roman{enumi})}
			\item $a_{1,\sigma}>a_{3,\sigma}>\cdots>a_{2n-1,\sigma}$,
			\item $a_{1,\sigma}+a_{2,\sigma}=a_{3,\sigma}+a_{4,\sigma}=\cdots=a_{2n-1,\sigma}+a_{2n,\sigma}=2a_{n+1,\sigma}$, and
			\item $a_{2n-1,\sigma}>a_{2n,\sigma}$.
		\end{enumerate}
		Let $V$ be the complex irreducible representation of $(\res_{E/\QQ}\GL_{2n+1})_\CC$
		with lowest weight $-\lambda$ with respect to $\Delta^+$. Write $\OO(V)\coloneqq \OO(-g_n\lambda)$. Write $\OO(V)\coloneqq \OO(-g_n\lambda)$. Then we denote the $\OO(V)\left[1/2d_E\right]$-module in Theorem \ref{thm:gl_n_totally_real} by $A(V)_{\OO(V)\left[1/2d_E\right]}$.
	\end{definition}
	
	In view of Corollary \ref{cor:unique}, $A(V)_{\OO(V)\left[1/2d_E\right]}$ is an $\OO(V)\left[1/2d_E\right]$-form of the complex Harish-Chandra module $A(V)$ introduced there in the sense that they only depend on $V$ by their construction.

	\subsubsection{CM case}
	
	Let $F$ be a totally real number field, and $E$ be a totally imaginary quadratic extension of $F$. In particular, $E$ is a CM number field. We wish to consider the arithmetic model $\res_{\OO_E\left[1/d_E\right]/\ZZ\left[1/d_E\right]} \GL_n$ of the Lie group $(\res_{E/\QQ}\GL_n)(\RR)$. For an arithmetic model of a maximal compact subgroup of this Lie group, let us introduce the $n$-th standard unitary group scheme. Henceforth we write
	\[\begin{array}{ccc}
		k=\ZZ\left[1/d_E\right],&k'=\OO_F\left[1/d_E\right],&k''=\OO_E\left[1/d_E\right].
	\end{array}\]
	to save space. We remark that $d_E$ kills all ramified prime ideals of $\OO_{F}$ in $E/F$ by \cite[Chapter III, $\S$ 2, Theorem 22]{frohlichtaylor}. In fact, every ramified prime ideal of $\OO_{E}$ in $E/F$ is ramified in $E/\QQ$ from the tower formula of the ramification indices (\cite[Chapter III, $\S$ 1, (1.13.a)]{frohlichtaylor}). Define an anti-involution $(-)^\ast$ on $\res_{k''/k'}\GL_n$ by the composition of the transpose map $(-)^T$ and the involution $\bar{\ }$ induced from the nontrivial element of the Galois group of the quadratic Galois extension $k'\subset k''$ (recall Lemma \ref{lem:fractional_splitting}). We now define the $n$-th standard unitary group $\U(n)$ over $k'$ in a similar way to \cite[Example 3.3.2]{hayashilinebdl}. Namely, the standard unitary group $\U(n)$ is given by
	\[\U(n,R)=\{g\in \GL_n(R\otimes_{k'} k''):~g^{-1}=g^\ast\}\]
	for commutative $k'$-algebras $R$. If we denote $\theta=((-)^\ast)^{-1}$, $\U(n)$ is the fixed point subgroup scheme of $\res_{k''/k'}\GL_n$ by $\theta$.
	
	\begin{proposition}\label{prop:U(n)}
		The counit $(\res_{k''/k'} \GL_n)\otimes_{k'} k''\to \GL_n$ restricts to an isomorphism
		\[\U(n)\otimes_{k'} k''\cong\GL_n.\]
		In particular, $\U(n)$ is a reductive group scheme.
	\end{proposition}
	
	\begin{proof}
		Let $A$ be a commutative $k''$-algebra, and $g=(g_{ij})\in\GL_n(A)$.
		Apply $A\otimes_{k''}-$ to the splitting isomorphism $k''\otimes_{k'} k''\cong (k'')^2$
		to obtain
		\[\phi_A:A\otimes_{k'} k''\cong A^2;~a\otimes c\mapsto (ac,a\bar{c}).\]
		Write $\sigma$ for the involution of $A\otimes_{k'} k''$ defined by $A\otimes_{k'} \bar{}$. Let $C$ be the involution of $A^2$ defined by switching the factors. Then it is evident by definition that $\phi_A\circ\sigma=C\circ \phi_A$. Take the inverses to get $\sigma\circ\phi^{-1}_A=\phi^{-1}_A\circ C$. One can prove by using this equality that $(\phi^{-1}_A(g_{ij},\bar{g}_{ij}))$ belongs to $\U(n,A)$. Moreover, this gives the inverse map.
	\end{proof}
	
	The reductive group scheme $\res_{k'/k}\U(n)$ is a model of a maximal compact subgroup of $(\res_{E/\QQ}\GL_n)(\RR)$. Let $T$ (resp.~$H_{\spl}$) be the subgroup scheme of diagonal matrices of $\U(n)$ (resp.~$\GL_n$). We remark that $H_{\spl}$ is identified with $T\otimes_{k'} k''$ through the isomorphism
	\[\U(n)\otimes_{k'} k''\cong \GL_n.\]
	In particular, $T$ is a maximal torus of $\U(n)$, which is split over $k''$. Let us also note for applications that $\res_{k''\left[1/2\right]/k\left[1/2\right]} H_{\spl}$ is the fundamental Cartan subgroup attached to
	\[\res_{k'\left[1/2\right]/k\left[1/2\right]} T.\]
	
	We again study characters as in the last section. There is a canonical identification
	\[X^\ast((\res_{E/\QQ} H_{\spl})\otimes_\QQ E^{\Gal}(\sqrt{-1}))\cong
	\prod_{\sigma\in\Hom_\QQ(E,E^{\Gal})} X^\ast(H_{\spl}\otimes_E E^{\Gal}),\]
	where $E^{\Gal}$ is a normal hull of $E/\QQ$. For $\lambda \in X^\ast((\res_{E/\QQ} H_{\spl})\otimes_\QQ E^{\Gal}(\sqrt{-1}))$, using the standard basis of $X^\ast(H_{\spl,E^{\Gal}})$, we will denote the corresponding element of \[\prod_{\sigma\in\Hom_\QQ(E,E^{\Gal})} X^\ast(H_{\spl}\otimes_E E^{\Gal})\]
	as
	$\sum_\sigma\lambda_{\sigma}\sigma=\sum_\sigma(\lambda_{i,\sigma})\sigma$.
	The composition with the nontrivial element of the Galois group of the Galois extension $E/F$ induces a free involution $\bar{\ }$ on $\Hom_\QQ(E,E^{\Gal})$. Unwinding the definitions, one can find that the involution on
	\[\prod_{\sigma\in\Hom_\QQ(E,E^{\Gal})}
	X^\ast(H_{\spl}\otimes_E E^{\Gal})\cong \prod_{\sigma\in\Hom_\QQ(E,E^{\Gal})}
	X^\ast(H_{\spl})\]
	induced from $\theta$ is given by 
	$\sum_\sigma\lambda_{\sigma}\sigma\mapsto -\sum_{\sigma}\lambda_{\sigma}\bar{\sigma}$. Choose a complete system
	$I_E\subset\Hom_\QQ(E,E^{\Gal})$
	of representatives of the quotient by the involution $\bar{\ }$. Take
	\[\Delta^+\coloneqq \left(\coprod_{\sigma\in I_E}
	\Delta^+(\GL_n, H_{\spl})\sigma\right)\coprod \left(\coprod_{\sigma\in I_E}-
	\Delta^+(\GL_n, H_{\spl})\bar{\sigma}\right)
	\subset \prod_{\sigma}
	X^\ast(H_{\spl})\]
	for a $\theta\otimes_\QQ E^{\Gal}$-stable positive system, where $\Delta^+(\GL_n,H_{\spl})$ is the standard positive system in \cite[Appendix C]{knapp2002}. Let us identify the Weyl group of
	\[((\res_{E/\QQ}\GL_n)_{E^{\Gal}},(\res_{E/\QQ} H_{\spl})_{E^{\Gal}})\]
	with the product of copies of the Weyl group of $(\GL_n, H_{\spl})$ indexed by $\Hom_\QQ(E,E^{\Gal})$. Define its element $\vec{w}=(w_\sigma)$ by
	\[w_\sigma=\begin{cases}
		w&(\sigma\in I_E)\\
		e&(\sigma\in\Hom_\QQ(E,E^{\Gal})\setminus I_E),
	\end{cases}\]
	where $e$ is the unit, and $w$ is the longest element of the Weyl group of $(\GL_n, H_{\spl})$ with respect to the positive system in \cite[Appendix C]{knapp2002}.
	
	To get models of orbits, recall that 
	\[\begin{split}
		\rtype(\res{k''/k}\GL_n,\res_{k'/k}\U(n))(k)
		&\cong (\res_{k'/k} \rtype(\res_{k''/k'}\GL_n,\U(n)))(k)\\
		&\cong (\res_{k'/k} \type \U(n))(k)\\
		&\cong (\type \U(n))(k')
	\end{split}
	\]
	(Corollary \ref{cor:restrictionsetting}). To compute $\type \U(n)$, we wish to know the Dynkin scheme of $\U(n)$. This is done in a similar way to \cite[Example 4.1.3]{hayashilinebdl} since the conjugate action on $X^\ast(T\otimes_{k'} k'')$ is equal to $-1$. Pick $x\in\rtype(\res_{k''/k'} \GL_n,\U(n))(k')$. Let $Q$ be the parabolic subgroup of $\GL_n$ over $k''$ containing the upper triangular matrices such that $Q$ has parabolic type $kt(x)|_{\Spec k''}$. Then $i_x$ can be expressed as $\Pp_{\U(n),kt(x)}\hookrightarrow \res_{k''/k'} \GL_n/Q$. 
	
	Let $\Phi\subset X^\ast(H_{\spl})$ denote the parabolic subset corresponding to $Q$. We remark that the set of roots of $(\res_{E/\QQ} Q)_{E^{\Gal}}$ is given by $\coprod_{\sigma\in\Hom_\QQ(E,E^{\Gal})} \Phi\sigma$. We denote the Lie algebra of the parabolic subgroup over $E^{\Gal}$ corresponding to the parabolic subset
	\[\left(\coprod_{\sigma\in I_E}
	\Phi\sigma\right)\coprod \left(\coprod_{\sigma\in I_E} w\Phi\bar{\sigma}\right)\]
	by $\lieq'$. 
	
	Take $\lambda\in X^\ast((\res_{E/\QQ} H_{\spl})_{E^{\Gal}})$. We will write $\lambda=(\lambda_\sigma)_{\sigma}=((\lambda_{\sigma,i})_i)_{\sigma}$ through the identification
	\[X^\ast((\res_{E/\QQ} H_{\spl})_{E^{\Gal}})
	\cong \prod_{\sigma\in\Hom_\QQ(E,E^{\Gal})} X^\ast(H_{E^{\Gal}})
	\cong \prod_{\sigma\in\Hom_\QQ(E,E^{\Gal})} \ZZ^n.\]
	One can lift $\vec{w}\lambda$ to a character over $\OO(\vec{w}\lambda)\left[1/d_E\right]$ which we will denote by the same symbol. 
	
	Fix $E^{\Gal}\hookrightarrow\CC$. It is easy to show that $E^{\Gal}$ is closed under the complex conjugation. Moreover, it restricts to the involution of $E$ corresponding to the nontrivial element of the Galois group of $E/F$.

	\begin{theorem}\label{thm:AqlambdaGLnCM}
		Assume the equalities
		$\langle\alpha^\vee,w\lambda_{\sigma}\rangle
		=\langle \alpha^\vee,\lambda_{\bar{\sigma}}\rangle=0$
		for all $\alpha\in \Phi\cap (-\Phi)$ and $\sigma\in I_E$.
		\begin{enumerate}
			\item One can extend $\vec{w}\lambda$ to a character of $(\res_{k''/k} Q)_{\OO(\lambda)\left[1/d_E\right]}$,
			which we will denote by the same symbol.
			\item Set
			\[\Ll\coloneqq (\res_{k''/k} \GL_n)_{\OO(\lambda)\left[1/d_E\right]}
			\times^{(\res_{k''/k} Q)_{\OO(\lambda)\left[1/d_E\right]}} \vec{w}\lambda.\]
			Then
			\[\Gamma((\res_{k''/k} \GL_n/Q)_{\OO(\lambda)\left[1/d_E\right]},((\res_{k'/k} i_x)_{\OO(\lambda)\left[1/d_E\right]})_+((\res_{k'/k} i_x)_{\OO(\lambda)\left[1/d_E\right]})^\ast\Ll)\]
			exhibits a projective $\OO(\vec{w}\lambda)\left[1/d_E\right]$-form of
			$A_{\lieq'_\CC}(\lambda)$.
		\end{enumerate}
	\end{theorem}
	
	\begin{definition}\label{def:A(V)_CM}
		Suppose $kt(x)=\emptyset$. For each $\sigma\in I_E$, we also assume
		\begin{enumerate}
			\renewcommand{\labelenumi}{(\roman{enumi})}
			\item $\lambda_{\sigma,1}>\lambda_{\sigma,2}>\cdots>\lambda_{\sigma,n}$ and
			\item $\lambda_{\sigma,1}+\lambda_{\sigma,n+1}=\lambda_{\sigma,2}+\lambda_{\sigma,n+2}=\cdots=\lambda_{\sigma,n}+\lambda_{\sigma,2n}$.
		\end{enumerate}
		Let $V$ be the complex irreducible representation of $(\res_{E/\QQ}\GL_{n})_\CC$ with lowest weight $-\lambda$. Write $\OO(V)=\OO(\vec{w}\lambda)$. Then we denote the $\OO(V)\left[1/d_E\right]$-module in Theorem \ref{thm:AqlambdaGLnCM} by $A(V)_{\OO(V)\left[1/d_E\right]}$.
	\end{definition}
	
	Again, this gives an $\OO(V)\left[1/d_E\right]$-form of the corresponding complex Harish-Chandra module introduced in Corollary \ref{cor:unique}.
	
	\subsection{Relative Lie algebra cohomology}\label{sec:(g,K)-cohomology}
	
	According to \cite{voganzuckerman1984}, the significant feature of $A_{\lieq}(\lambda)$-modules is the nonvanishing of the $(\lieg,K)$-cohomology. In this section, we introduce arithmetic models of the relative Lie algebra cohomology over the complex numbers, and establish finiteness properties for our models of cohomologically induced modules constructed above. We will work with an arithmetic analog of $(\lieg,S(\RR)_0 K(\RR))$-cohomology towards the reductive setting as in \cite{harderraghuram}. In fact, with a view towards applications to automorphic representations, we take the action of the center into consideration.
	
	Let $k$ be a commutative ring, $G$ be a smooth affine group scheme over $k$, and $K\subset G$ be a smooth closed subgroup. Let $Z_\lieg$ be the center of $\lieg$, and $\lies$ be a $k$-submodule of $Z_\lieg$. For a $(\lieg,K)$-module $X$ over $k$, consider the complex
	\begin{equation}
		\Hom_{\lies,K}\left(\bigwedge^\bullet\lieg/\left(\lies+\liek\right),X\right)
		\label{eq:gSKcohomology}
	\end{equation}
	with the usual differential
	\[\begin{split}
		(df)(x_1\wedge \cdots\wedge x_{q+1})
		&=
		\sum_{i=1}^{q+1} (-1)^{i+1} x_i f(x_1,x_2,\ldots,\hat{x}_i,\ldots,x_{q+1})\\
		& +\sum_{i<j} (-1)^{i+j} f([x_i,x_j],x_1,\ldots,\hat{x}_i,\ldots,\hat{x}_j,\ldots,x_{q+1}).
	\end{split}\]
	We define $(\lieg,\lies,K)$-cohomology of $X$ as the cohomology group $H^\bullet\left(\lieg,\lies, K; X\right)$ of the complex \eqref{eq:gSKcohomology}. This notion of $(\lieg,\lies, K)$-cohomology commutes with arbitrary flat base change of commutative Noetherian rings by \cite[Corollary 3.2.10]{hayashi2018}. If $\lies=\{0\}$ then we simply denote it by $H^\bullet\left(\lieg, K; X\right)$ as usual.
	
	\begin{example}
		Let $G$ be a connected real reductive algebraic group, $K\subset G$ be the fixed point subgroup by a Cartan involution of $G$, $S$ be a split central torus in $G$, and $X$ be a $(\lieg_\CC,K_\CC)$-module. Since the trivial group is the only compact subgroup of $S(\RR)_0$, \eqref{eq:gSKcohomology} recovers the classical relative Lie algebra cohomology:
		\[
		H^\bullet\left(\lieg_\CC,\lies_\CC, K_\CC; X\right)
		\cong
		H^\bullet\left(\lieg_\CC,S(\RR)_0 K(\RR); X\right).
		\]
	\end{example}
	
	\begin{example}\label{ex:relativeLiealgebracohomology}
		Suppose that $k=F$ is a field of characteristic zero. Then \eqref{eq:gSKcohomology} recovers the classical relative Lie algebra cohomology, i.e.,
		\[
		H^\bullet\left(\lieg,\lies, K; X\right)
		=
		H^\bullet\left(\lieg,\lies+\liek; X\right)^{\pi_0(K)}.
		\]
		Assume $\liek$ is reductive in $\lieg$ in the sense of \cite[Chapter I, section 2.4]{borelwallach}, i.e., $\lieg$ is semisimple as a $\liek$-module for the adjoint representation. If $H^\bullet\left(\lieg,\lies, K; X\right)\neq 0$ and $Z_\lieg$ acts on $X$ by a scalar then the action of $Z_{\lieg}$ on $X$ is trivial by \cite[Chapter I, section 5.3, Theorem (ii)]{borelwallach}.
	\end{example}
	
	In applications, we can reduce the computation of the cohomology to the case $\lies=\{0\}$:
	
	\begin{proposition}\label{prop:relativeliealgcoh}
		Let $\lieg^{\mathrm{ss}}$ be the derived subalgebra of $\lieg$. Assume the following conditions:
		\begin{enumerate}
			\renewcommand{\labelenumi}{(\roman{enumi})}
			\item $k=F$ is a field of characteristic zero,
			\item $\lieg$ is reductive,
			\item $\liek =(Z_\lieg\cap \liek)\oplus (\lieg^{\mathrm{ss}}\cap\liek)$\footnote{This equality holds if $\liek$ is a symmetric subalgebra of $\lieg$.},
			\item $\lies\cap \liek=\{0\}$,
			\item $K$ centralizes $\lies$, and
			\item $Z_{\lieg}$ acts trivially on $X$.
		\end{enumerate} 
		Let $\lies'$ be a subspace of $Z_{\lieg}$ containing $\lies$. Choose a complementary subspace $\liec'$ to $\lies'\cap \liek$ in $Z_{\lieg}\cap \liek$. We next choose a complementary subspace $\liec$ to $\lies'$ in $Z_{\lieg}$ containing $\liec'$. Set $\lieg'\coloneqq \liec\oplus \lieg^{\mathrm{ss}}$.
		\begin{enumerate}
			\item We have an isomorphism
			\[H^\bullet\left(\lieg,\lies, K; X\right)\cong 
			H^\bullet\left(\lieg',\liek\cap\lieg'; X\right)^{\pi_0(K)}
			\otimes_F \wedge^\bullet (\lies'/(\lies+(\liek\cap \lies')))^\vee 
			\]
			of graded vector spaces. In particular, if $H^\bullet\left(\lieg,\lies, K; X\right)\neq \{0\}$ then we have
			\[H^\bullet\left(\lieg',\liek\cap\lieg'; X\right)\neq \{0\}.\]
			\item If $\liek\subset\lieg'$, we have
			\[H^\bullet\left(\lieg,\lies, K; X\right)\cong 
			H^\bullet\left(\lieg',K; X\right)
			\otimes_F \wedge^\bullet (\lies'/(\lies+(\liek\cap \lies')))^\vee .
			\]
		\end{enumerate}
	\end{proposition}
	
	We remark that $\liec$ in the statement always exists since $\liec'\cap \lies'=\{0\}$.
	
	\begin{proof}
		We have decompositions 
		\[\begin{array}{cc}
			\lieg=\lieg'\oplus \lies',&\liek=(\liek\cap \lieg')\oplus (\liek\cap \lies')
		\end{array}\]
		in virtue of the hypotheses. The K\"unneth formula of \cite[Chapter I, section 1.3]{borelwallach} implies an isomorphism
		\[\begin{split}
			H^\bullet\left(\lieg,\lies+\liek; X\right)
			&\cong H^\bullet\left(\lieg',\liek\cap \lieg'; X\right)
			\otimes_F H^\bullet\left(\lies',\lies+(\liek\cap \lies');F\right)\\
			&\cong H^\bullet\left(\lieg',\liek\cap \lieg'; X\right)
			\otimes_F \wedge^\bullet (\lies'/(\lies+(\liek\cap \lies')))^\vee 
		\end{split}
		\]
		of graded vector spaces (recall $\lies\cap \liek=\{0\}$).
		In virtue of (v), take $\pi_0(K)$ to obtain
		\[H^\bullet\left(\lieg,\lies, K; X\right)\cong 
		H^\bullet\left(\lieg',\liek\cap \lieg'; X\right)^{\pi_0(K)}\otimes_F
		\wedge^\bullet(\lies'/(\lies+(\liek\cap \lies')))^\vee 
		\]
		of graded vector spaces. See Example \ref{ex:relativeLiealgebracohomology} for (2).
	\end{proof}
	
	\begin{example}\label{ex:cohomology_decomposition_1}
		Put $\lies'=\lies$ to obtain
		\[H^\bullet\left(\lieg,\lies, K; X\right)\cong H^\bullet\left(\lieg',K; X\right).\]
		We note that $\liec'=Z_{\lieg}\cap\liek$ in this case. Therefore $\liec$ is a complementary subspace to $\lies$ in $Z_\lieg$.
		
		A similar argument to Proposition \ref{prop:relativeliealgcoh} implies
		\[H^\bullet\left(\lieg,K; X\right)
		\cong
		H^\bullet\left(\lieg',K; X\right)\otimes_F H^\bullet(\lies;F)
		\cong H^\bullet\left(\lieg',K; X\right)\otimes_F (\wedge^\bullet\lies)^\vee.
		\]
		In particular, $H^\bullet\left(\lieg,K; X\right)$ is nonzero if and only if $H^\bullet\left(\lieg,\lies, K; X\right)$ is nonzero.
	\end{example}
	
	\begin{example}\label{ex:cohomology_decomposition_2}
		Put $\lies'=Z_{\lieg}$. In this case, we have $\liec'=\liec=\{0\}$ and $\lieg'=\lieg^{\mathrm{ss}}$. We thus obtain
		\[H^\bullet\left(\lieg,\lies, K; X\right)\cong 
		H^\bullet(\lieg^{\mathrm{ss}},\liek\cap\lieg^{\mathrm{ss}};X)^{\pi_0(K)}\otimes_F
		\wedge^\bullet(Z_\lieg/(\lies+(\liek\cap Z_\lieg)))^\vee.\]
	\end{example}

	For the condition above on the action of $Z_{\lieg}$, and more generally, the existence of infinitesimal characters, let us note that a relative analog of Dixmier's variant of Schur's lemma holds without the hypothesis (Q) in \cite{januszewskirationality}:
	
	\begin{proposition}\label{prop:schur}
		Let $(\lieg,K)$ be a Harish-Chandra pair over a field $F$ with $\dim_F \lieg<\infty$, and $X$ be an irreducible $(\lieg,K)$-module.
		\begin{enumerate}
			\item Every element $\varphi\in\End_{\lieg,K}(X)$ is algebraic over $F$.
			\item As a $U(\lieg)$-module, $X$ is finitely generated.
			\item Suppose that $K$ centralizes the center of $U(\lieg)$. The $(\lieg,K)$-module $X$ has an infinitesimal character if $F'\otimes_F X$ is irreducible for an algebraic closed field $F$' containing $F$.
		\end{enumerate}
	\end{proposition}
	
	\begin{proof}
		Part (1) follows by a minor modification of \cite{quillen1969}. Indeed, define a Harish-Chandra pair
		\[(F\left[\varphi\right]\otimes_F U(\lieg),K)\]
		by putting the trivial action of $K$ on $F\left[\varphi\right]$ and setting the zero map from $\liek$ to $F\left[\varphi\right]$ (cf.~\cite[Example 2.1.5 (8)]{hayashi2019}). We take a nonzero finite dimensional $K$-submodule of $X$ instead of a nonzero element. 
		
		Part (2) is standard (cf.~the proof of \cite[Proposition 3.2.4]{hayashi2018}). Part (3) follows from (1), (2), and \cite[Theorem 3.1.6]{hayashi2018}. 
	\end{proof}

	We know:
	
	\begin{theorem}[{\cite[Theorem 5.5]{voganzuckerman1984}}]\label{thm:a_q->cohomological}
		Let $G,\theta,K,T,H,Q',L$ be as in Lemma \ref{lem:realformofK^1_C}. Let $U'$ be the unipotent radical of $Q'$, and $\Delta(U',H_\CC)$ be the set of $H_\CC$-roots in $U'$. Write $u=\dim U'\cap K^0_\CC$ and $v=\dim U'-u$. Take a split central torus $S\subset G$. Let $\liep\subset\lieg$ be the $-1$-eigenspace of the differential of $\theta$. Let $\pi$ and $\tau$ be one-dimensional modules over $(\liel_\CC,(L\cap K)(\RR))$ and $(\lieg_\CC,K(\RR))$ such that $\pi\otimes_\CC\tau$ is unitarizable. Write $\lambda$ and $\nu$ for the linear functionals of $\lieh_\CC$ obtained by restriction of their Lie algebra actions\footnote{In this notation, $\pi\otimes_\CC\tau$ is unitarizable if and only if $(\lambda+\nu)|_{\lieh}$ is valued in $\RR\sqrt{-1}$.}. Choose a Borel subgroup $B'\subset G_\CC$ satisfying $H_\CC\subset B'\subset Q'$ to fix a positive system. We denote the half sum of the positive roots by $\rho$. Let $V$ be a $(\lieg_\CC,K(\RR))$-module which is irreducible of lowest weight $\mu$ as a $\lieg_\CC$-module. Assume that no $\alpha\in \Delta(U',H_\CC)$ satisfies $\langle \alpha^\vee,\lambda+\nu+\rho\rangle\in\{-1,-2,-3\ldots\}$. Write
		\[X\coloneqq R^uI^{\lieg_\CC,K_\CC}_{\lieq',(L\cap K)_\CC}(\pi\otimes_\CC \wedge^{\mathrm{top}} \lieg_\CC/\bar{\lieq}').\]
		Let $X_{\min}$ be the subspace spanned by minimal $K_\CC$-types of $X$.
		\begin{enumerate}
			\renewcommand{\labelenumi}{(\arabic{enumi})}
			\item The cohomology $H^\bullet(\lieg_\CC,\lies_\CC,K_\CC;X\otimes_\CC V)$ vanishes unless $\mu=-\lambda$.
			\item Assume that $\mu=-\lambda$.
			Then we have $X_{\min}\cong I^{\liek_\CC,K_\CC}_{\bar{\lieq}'\cap\liek_\CC,L'\cap K_\CC}(\pi\otimes_\CC \wedge^{\mathrm{top}}(\lieu'\cap\liep_\CC))$ and
			\[\begin{split}
				H^\bullet(\lieg_\CC,\lies_\CC,K_\CC;X\otimes_\CC V)
				&\cong \Hom_{L'\cap K_\CC}(\wedge^\bullet\liep_\CC/\lies_\CC,
				X_{\min}\otimes_\CC V)\\
				&\cong\Hom_{L'\cap K_\CC}(\wedge^{\bullet-v} (\liel'\cap\liep_\CC)/\lies_\CC,\pi\otimes_\CC H^0(\bar{\lieu},V)).
			\end{split}\]
			(recall Corollary \ref{cor:centraltorus} to see that the quotients in the right hand side make sense).
		\end{enumerate}
		In particular, the dimension of the cohomology of the bottom degree is one if $\pi^\vee\cong H^0(\bar{\lieu},V)$ as a representation of $L'\cap K_\CC$.
	\end{theorem}
	
	\begin{remark}
		Since $\lambda$ extends to a character of $\liel_\CC$, $\dim_\CC H^0(\bar{\lieu}',V)=1$. Indeed, $H^0(\bar{\lieu}',V)$ is the lowest weight space of $V$.
	\end{remark}
	
	\begin{remark}
		If $\mu=-\lambda|_{\lieh_\CC}$ then $\lambda$ is dominant. In particular, we have \[\langle\alpha^\vee,\lambda+\nu+\rho\rangle=\langle \alpha^\vee,\lambda+\rho\rangle>0\]
		for all $\alpha\in \Delta(U',H_\CC)$ in this case.
	\end{remark}
	
	\begin{remark}[{\cite[Theorems 8.21]{knappvogan}}]\label{rem:duality}
		We have $X\cong L_u P^{\lieg_\CC,K_\CC}_{\bar{\lieq}',(L\cap K)_\CC}(\pi\otimes_\CC \wedge^{\mathrm{top}} \lieg_\CC/\bar{\lieq}')$
		under the hypothesis of Theorem \ref{thm:a_q->cohomological}.
	\end{remark}

	The (non-)vanishing in the case $\mu=-\lambda$ is delicate:
	
	\begin{example}\label{ex:gl3}
		Put $G=\GL_3$, $K=\Oo(3)$. Let $S$ be the center of $G$. Assume $Q'=B'$ to be Borel. Let $V$ be an irredicible representation of $\GL_3$ as a complex algebraic group with unitary lowest weight $-\lambda$. Set $\pi=\lambda$. In this case, we have
		\[H^\bullet(\lieg_\CC,\lies_\CC,K_\CC;X\otimes_\CC V)
		\cong\Hom_{H\cap \Oo(3)}(\wedge^{\bullet-2}\lieh_\CC/(\lies_\CC+\lieso(2,\CC)),\CC)
		\neq 0.\]
		Let $\tau$ be the nontrivial character of $\pi_0(G(\RR))$. Regard it as a one-dimensional $(\lieg_\CC,K_\CC)$-module. Then $V\otimes_\CC\tau$ still has lowest weight $-\lambda$, but
		\[H^\bullet(\lieg_\CC,\lies_\CC,K_\CC;X\otimes_\CC V\otimes_\CC\tau)
		\cong \Hom_{H\cap \Oo(3)}(\wedge^{\bullet-2}\lieh_\CC/(\lies_\CC+\lieso(2,\CC)),\tau)
		=0\]
		since $\pi_0(H\cap \Oo(3))\cong\pi_0(\Oo(3))$ acts trivially on $\lieh_\CC/(\lies_\CC+\lieso(2,\CC))$. This kind of vanishing does not happen for $G=\GL_2$ with $Q'=B'$ Borel since $H\cap \Oo(2)=\SO(2)$ (cf.~$K^1=K^0$ in this case).
	\end{example}
	
	\begin{example}\label{ex:non-vanishing}
		Put $G=\GL_2$ and $K=\Oo(2)$. Define a unitary character $\chi$ by the quotient map
		\[\GL_2(\RR)\to \GL_2(\RR)/\GL_2(\RR)_0\cong\{\pm 1\}.\]
		We denote the corresponding Harish-Chandra module by the same symbol. Then it is easy to show
		\[\dim_\CC H^p(\lieg_\CC,K(\RR);\chi)=\begin{cases}
			1&(p=2,3)\\
			0&(\mathrm{otherwise}).
		\end{cases}\]
		Similarly, let $S$ be the center of $G$. Then we have
		\[\dim_\CC H^p(\lieg_\CC,\lies_\CC,K(\RR);\chi)=\begin{cases}
			1&(p=2)\\
			0&(\mathrm{otherwise}).
		\end{cases}\]
	\end{example}
	
	The second example tells us that there may be two or more $\pi$ with \[H^\bullet(\lieg_\CC,\lies_\CC,K_\CC;R^uI^{\lieg_\CC,K_\CC}_{\lieq',(L\cap K)_\CC}(\pi\otimes_\CC \wedge^{\mathrm{top}} \lieg_\CC/\bar{\lieq}')\otimes_\CC V)\neq 0.\]
	This issue does not occur when $Q'$ is Borel, i.e., for a given $V$, there is up to isomorphism at most a single $\pi$ with non-zero cohomology since $L$ is commutative in this case:
	
	\begin{example}\label{ex:abelian_case}
		Let $H$ be a real torus with a Cartan involution $\theta$, $T$ be an open subgroup of $H^\theta$, and $\chi$ be a character of $T_\CC$ with trivial differential. We regard $\chi$ as a one-dimensional $(\lieh_\CC,T_\CC)$-module for the trivial action of $\lieh_\CC$ (recall Example \ref{ex:split_central}). Let $\lies'$ be any vector subspace of $\lieh_\CC$. Then $H^\bullet(\lieh_\CC,\lies',T_\CC;\chi)= 0$ unless $\chi$ is trivial since $T_\CC$ and $\lies'$ act trivially on $\lieh_\CC/(\lies'+\liet_\CC)$.
	\end{example}

	\begin{theorem}[{\cite[Theorem 5.6]{voganzuckerman1984}}]\label{thm:cohomological->aq}
		Let $G,S,K,\theta$ be as in Theorem \ref{thm:a_q->cohomological}. Let $X$ be an irreducible essentially unitarizable $(\lieg_\CC,K(\RR))$-module. Assume that there exists a finite dimensional $(\lieg_\CC,K(\RR))$-module $V$ with the following properties:
		\begin{enumerate}
			\renewcommand{\labelenumi}{(\roman{enumi})}
			\item $V$ is irreducible as a $\lieg_\CC$-module, and
			\item $H^\bullet(\lieg_\CC,\lies_\CC,K_\CC;X\otimes_\CC V)$ is nonzero.
		\end{enumerate}
		Then there exist a $\theta$-stable parabolic subgroup $Q'$ with unipotent radical $U'$ containing the complexification of a maximal torus of $K$, a one-dimensional $(\liel',L'\cap K_\CC)$-module $\pi$, and a one-dimensional $(\lieg_\CC,K(\RR))$-module $\tau$ such that
		\begin{enumerate}
			\item $X\cong R^uI^{\lieg_\CC,K_\CC}_{\lieq',L'\cap K_\CC}
			(\pi\otimes_\CC \wedge^{\mathrm{top}} \lieg_\CC/\bar{\lieq}')$
			as a $(\lieg_\CC,K_\CC)$-module, and
			\item $\pi\otimes_\CC\tau$ is unitary,
		\end{enumerate}
		where $L'=Q'\cap \bar{Q}'$ and $u=\dim (U'\cap K^0_\CC)$.
	\end{theorem}
	
	We outline how to verify this in the disconnected algebraic setting. Let $\tau$ be a one dimensional $(\lieg_\CC,K(\RR))$-module such that $X\otimes_\CC \tau$ is unitarizable. Since $\bar{\lieu}'$ lies in the semisimple part of $\lieg_\CC$, $\bar{\lieu}'$ acts trivially on $\tau$. We may therefore replace $(X,V)$ with $(X\otimes_\CC \tau,\tau^\vee\otimes_\CC V)$ to assume that $X$ is unitarizable. This is a cohomologically induced module for certain $Q'$. To see that $X$ is induced from a one-dimensional representation $\pi$, use the fact that the component group of $L'\cap K_\CC$ is abelian (\cite[14.5. Corollaire]{boreltits1965}). We find that $\pi$ is unitary by seeing the Lie algebra action (recall \cite[Theorem 5.6]{voganzuckerman1984}).
	
	\begin{remark}\label{rem:uniqueness}
		Put the natural structure of a $(\lieq',L'\cap K_\CC)$-module on $\dim H^0(\bar{\lieu}',V)$. Suppose that $Q'=B'$ is Borel (e.g.~assume a highest weight of $V$ to be regular). Then in view of Theorem \ref{thm:a_q->cohomological} and Example \ref{ex:abelian_case}, we have $\pi\cong H^0(\bar{\lieu}',V)^\vee$ and $\dim H^v(\lieg_\CC,\lies_\CC,K_\CC;X\otimes_\CC V)=1$, where $v=\dim U'-\dim (U'\cap K^0_\CC)$.
		
		More generally, we have $\pi\cong H^0(\bar{\lieu}',V)^\vee$ if $\dim H^v(\lieg_\CC,\lies_\CC,K_\CC;X\otimes_\CC V)=1$ (use Theorem \ref{thm:a_q->cohomological} (2)). This does not hold in general without replacing $V$ (Example \ref{ex:non-vanishing}, $v=1$ in this case). The replacement of $V$ is possible if $K=K^1$ since $\pi\otimes_\CC H^0(\bar{\lieu}',V)$ is a character of $\pi_0(L'\cap K_\CC)\cong\pi_0(K_\CC)$. 
	\end{remark}

	For applications to cohomological cuspidal automorphic representations of $\GL_n$, let us also record a uniqueness property under strong conditions. The remark above and the following result are inspired by \cite[Proof of Proposition 4.2]{lischwermer}.
	
	\begin{corollary}\label{cor:unique}
		Consider the setting of Theorem \ref{thm:cohomological->aq}. Assume the following conditions:
		\begin{enumerate}
			\renewcommand{\labelenumi}{(\roman{enumi})}
			\item The flag variety of $G_\CC$ has a unique closed $K_\CC$-orbit.
			\item Fix a maximal torus $T$ of $K^0$. Set $H=Z_G(T)$. Choose a $\theta_\CC$-stable Borel subgroup of $G_\CC$ containing $H_\CC$. Let $\lambda$ denote the lowest weight of $V$ with respect to the positive system attached to $(H_\CC,B')$. We regard $\lambda$ as a character of $\lieh_\CC$. Then there exists a one-dimensional $(\lieg_\CC,K(\RR))$-module $\tau$ such that $\lambda+\nu$ is regular and unitary, where $\nu$ is the functional of $\lieh_\CC$ obtained by restriction of $\tau$. Here we say $\lambda+\nu$ is unitary if $\lambda|_{\lieh}+\nu|_{\lieh}$ is valued in $\RR\sqrt{-1}$.
		\end{enumerate}
		Then there is a unique irreducible essentially unitarizable $(\lieg_\CC,K(\RR))$-module $X$ up to isomorphism such that $H^\bullet(\lieg_\CC,\lies_\CC,K_\CC;X\otimes_\CC V)\neq 0$. We will denote it by $A(V)$.
	\end{corollary}

	\begin{proof}
		For existence, see Corollary \ref{cor:irreducible_unitary} and Theorem \ref{thm:a_q->cohomological}. To prove the uniqueness, apply Theorem \ref{thm:cohomological->aq} and Remark \ref{rem:uniqueness} to $(X,V)$ to express
		\[X\cong R^uI^{\lieg_\CC,K_\CC}_{\lieq',L'\cap K_\CC}
		(H^0(\bar{\lieu}',V)^\vee\otimes_\CC \wedge^{\mathrm{top}} \lieg_\CC/\bar{\lieq}')\]
		with $Q'$ Borel. The isomorphism class of the representation expressed in this form is independent of the choice of $Q'$ by (i). This completes the proof.
	\end{proof}
	
	\begin{remark}
		Condition (ii) is independent of the choice of $T$ and $B'$.
	\end{remark}
	
	\begin{remark}
		The $(\lieg_\CC,K(\RR))$-module $A(V)$ is unitarizable if $\lambda|_{\lieh}$ is valued in $\RR\sqrt{-1}$.
	\end{remark}
	
	We end this section with the finiteness of the cohomology group:
	
	\begin{proposition}\label{prop:relativeLiealgebracohomology}
		Assume the following conditions:
		\begin{enumerate}
			\renewcommand{\labelenumi}{(\roman{enumi})}
			\item $k$ is a Dedekind domain;
			\item $X$ has an exhaustive filtration $X=\cup_{i\geq 0} F_i X$ of $K$-submodules with the property that $\Gr^iX$ is finitely generated and torsion-free as a $k$-module for every integer $i$;
			\item $X$ is generically admissible as a $K$-module, i.e., for any finite dimensional representation $W$ of $K\otimes_k \Frac(k)$, $\dim_{\Frac (k)}\Hom_{K\otimes_k \Frac(k)}(W,X\otimes_k \Frac(k))<\infty$;
			\item $V$ is finitely generated as a $k$-module.
		\end{enumerate}
		Then the cohomology $H^i\left(\lieg,\lies, K; X\otimes_k V\right)$ is finitely generated as a $k$-module for every integer $i$. In particular, its quotient by torsion is finitely generated and projective as a $k$-module.
	\end{proposition}
	
	For the validity of (ii) in our applications, see Theorem \ref{thm:filtration}.
	
	\begin{lemma}\label{lem:Hom_finite}
		Let $k$ be a Noetherian domain, $K$ be a flat affine group scheme over $k$, and $X$ be a generically admissible $K$-module. Assume the condition (ii) of Proposition \ref{prop:relativeLiealgebracohomology}. Then $\Hom_K(V,X)$ is finitely generated and torsion-free as a $k$-module for every finitely generated $K$-module $V$.
	\end{lemma}
	
	\begin{proof}
		Take the filtration $F_\bullet X$ of (ii) in Proposition \ref{prop:relativeLiealgebracohomology}. Since $X$ is generically admissible, one can and do choose a nonnegative integer $i_0$ such that the canonical map
		\[\Hom_{K\otimes_k \Frac(k)}(V\otimes_k \Frac(k),F_{i_0}X\otimes_k \Frac(k))\to 
		\Hom_{K\otimes_k \Frac(k)}(V\otimes_k \Frac(k),F_{i}X\otimes_k \Frac(k))\]
		is an isomorphism for every integer $i\geq i_0$. For each $i> i_0$, we have a canonical map 
		\[\begin{tikzcd}
			0\ar[d]&0\ar[d]\\
			\Hom_K(V,F_{i-1}X)\ar[r]\ar[d]&\Hom_{K\otimes_k \Frac(k)}(V\otimes_k \Frac(k),F_{i-1}X\otimes_k \Frac(k))\ar[d, "\sim"{sloped, above}]\\
			\Hom_K(V,F_iX)\ar[r]\ar[d]
			&\Hom_{K\otimes_k \Frac(k)}(V\otimes_k \Frac(k),F_iX\otimes_k \Frac(k))\ar[d]\\
			\Hom_K(V,\Gr^iX)\ar[r]&\Hom_{K\otimes_k \Frac(k)}
			(V\otimes_k \Frac(k),\Gr^i X\otimes_k \Frac(k))
		\end{tikzcd}\]
		of exact sequences. Since the right middle vertical arrow is a bijection, the right bottom vertical arrow is zero. Since the bottom horizontal arrow is injective by the condition on the filtration $F_\bullet X$, the left bottom map is also zero. This implies that the left middle map is a bijection. Run though all $i$ to get
		\[\Hom_K(V,F_{i_0} X)\cong \varinjlim_{i} \Hom_K(V,F_iX)\cong \Hom_K(V,X)\]
		(cf.~\cite[Proposition 1.3.3]{hovey2004}). Observe that $F_{i_0}X$ is finitely generated and torsion-free by a routine inductive argument. The assertion now follows since $V$ is finitely generated.
	\end{proof}

	\begin{proof}[Proof of Proposition \ref{prop:relativeLiealgebracohomology}]
		It will suffice to show that $\Hom_{\lies,K}\left(\bigwedge^i\lieg/\left(\lies+\liek\right),X\otimes_k V\right)$ is finitely generated for every integer $i$. Since $V$ is finitely generated, it is identified with
		\[\Hom_{\lies,K}\left(V^\vee\otimes_k\bigwedge^i\lieg/\left(\lies+\liek\right),X\right).\]
		Then the assertion follows from Lemma \ref{lem:Hom_finite}.
	\end{proof}
	
	In our applications, we obtain a concrete description of the cohomology by thinking of forms of the situation of Theorem \ref{thm:a_q->cohomological}. The precise statement is as follows:
	
	\begin{variant}\label{var:gsKcohomology}
		Suppose that we are given a diagram
		\[\begin{tikzcd}
			\RR\ar[r,hook]&\CC\\
			k\ar[r]\ar[u]&\tilde{k}\ar[u,hook]
		\end{tikzcd}\]
		of commutative $\ZZ\left[1/2\right]$-algebras with the following properties:
		\begin{enumerate}
			\renewcommand{\labelenumi}{(\roman{enumi})}
			\item $\tilde{k}$ is a Dedekind domain, and
			\item the right vertical arrow is injective.
		\end{enumerate}
		Let $G$ be a reductive group scheme over $k$, equipped with an involution $\theta$ such that $\theta_{\RR}$ is Cartan. Let $S$ be a split central torus in $G$. Let $K$ be an open and closed subgroup scheme of $G^\theta$ with $\pi_0(K)$ finite \'etale. Let $x\in\rtype(G,K)(\tilde{k})$, $\Aa$ be a $G_{\tilde{k}}$-equivariant tdo on $\Pp_{G_{\tilde{k}},gt(x)}$, and $\Ll$ be a line bundle on $\Pp^{K_{\tilde{k}}\mathrm{-st}}_{G_{\tilde{k}},x}$, equipped with the structure of a $K_{\tilde{k}}$-equivariant integrable left $i^\cdot_x\Aa$-connection. Let $F_\bullet \Gamma(\Pp_{G_{\tilde{k}},gt(x)},(i_x)_+\Ll)$ be the filtration on $\Gamma(\Pp_{G_{\tilde{k}},gt(x)},(i_x)_+\Ll)$ in Theorem \ref{thm:filtration}.
		Let $V$ be a $(\lieg_{\tilde{k}},K_{\tilde{k}})$-module which is finitely generated and projective as a $\tilde{k}$-module such that $V_\CC$ is irreducible as a $\lieg_\CC$-module.	
		
		Choose a maximal torus $T\subset K_{\RR}$ and a $\theta_\CC$-stable parabolic subgroup $Q'\subset G_\CC$ containing $T_\CC$ such that $rt(\bar{Q}')=x|_{\Spec\CC}$. Let $u$ be the dimension of the unipotent radical of $Q'\cap K_\CC$. Write $L'=Q'\cap \bar{Q}$ and $H=Z_{G_{\RR}}(T)$. We fix a Borel subgroup $B'\subset G_\CC$ satisfying $H_\CC\subset B'\subset Q'$. Let $\rho$ denote the half sum of roots of $B'$ with respect to $H_\CC$. Let $\mu$ denote the lowest weight of $V_{\CC}$.
		
		Let $\pi$ be the geometric fiber of $\Ll$ at $\bar{Q}'\in \Pp^{K\times_S \tilde{S}\mathrm{-st}}_{G\times_S\tilde{S},x}(\CC)$. Let $\lambda$ be the character of $\lieh_\CC$ for which $\lieh_\CC$ acts on $\pi$.
		Suppose that there exists a one-dimensional $(\lieg_\CC,K(\RR))$-module such that $\pi\otimes_\CC \tau$ is a unitary $(\liel',L'\cap K_\CC)$-module. Moreover, assume that $\langle \alpha^\vee,\lambda+\rho\rangle\not\in(-\infty,-1]$ for any root $\alpha$ of the unipotent radical of $Q'$ with respect to $H_\CC$. Then we have
		\begin{flalign*}
			&H^\bullet(\lieg_{\tilde{k}},\lies_{\tilde{k}}, K_{\tilde{k}},\Gamma(\Pp_{G_{\tilde{k}},gt(x)},(i_x)_+\Ll)\otimes_{\tilde{k}} V)\\
			&\cong 
			\Hom_{K_{\tilde{k}}} (\wedge^{\bullet}\lieg_{\tilde{k}}/(\liek_{\tilde{k}}+\lies_{\tilde{k}}),
			F_0\Gamma(\Pp_{G_{\tilde{k}},gt(x)},(i_x)_+\Ll)\otimes_{\tilde{k}} V).
		\end{flalign*}
		Moreover, it vanishes unless $\lambda=-\mu$.
		In particular, the cohomology has no torsion.
	\end{variant}

	\begin{proof}
		
		The complex
		\[\Hom_{\lies,K_{\tilde{k}}} (\wedge^{\bullet}\lieg_{\tilde{k}}/(\liek_{\tilde{k}}+\lies_{\tilde{k}}),\Gamma(\Pp_{G_{\tilde{k}},gt(x)},(i_x)_+\Ll)\otimes_{\tilde{k}} V)\]
		is embedded into
		\[\Hom_{\lies_\CC,K_\CC} (\wedge^{\bullet}\lieg_\CC/(\liek_\CC+\lies_\CC),\Gamma(\Pp_{G,gt(x)},(i_x)_+\Ll)_\CC\otimes_{\CC} V_\CC).\]
		Notice also that the $(\lieg_\CC,K(\RR))$-module $\Gamma(\Pp_{G,gt(x)},(i_x)_+\Ll)_\CC$ is unitarizable by \cite[Theorem 9.1]{knappvogan}.
		Since the latter complex has trivial differential by \cite[Chapter II, section 3.1]{borelwallach}, so does the first one. The assertion is then a consequence of Corollary \ref{cor:duality} (2), Remark \ref{rem:duality}, Theorem \ref{thm:a_q->cohomological}, and the proof of Lemma \ref{lem:Hom_finite}.
	\end{proof}
	
	\begin{remark}
		We assumed a stronger dominant condition of $\lambda$ in Variant \ref{var:gsKcohomology} than Theorem \ref{thm:a_q->cohomological} for the unitarizability. However, it does not affect the case $\lambda=-\mu$ since $\langle \alpha^\vee,\lambda+\rho\rangle$ is an integer for any $\alpha$.
	\end{remark}
	
	\begin{remark}
		One can determine the (non-)vanishing and the local rank of the cohomology from the complexification through the preceding arguments.
	\end{remark}

	\appendix
	\section{Elementary definitions and facts on sheaves}\label{appendix}
	\subsection{Filtered sheaves}\label{sec:filt}
	\begin{definition}\label{defn:filtsheaf}
		Let $X$ be a topological space.
		\begin{enumerate}
			\renewcommand{\labelenumi}{(\arabic{enumi})}
			\item A sheaf of filtered abelian groups on $X$ is a sheaf $\Mm$ of abelian groups on $X$, equipped with a sequence of subsheaves
			$F_0\Mm\subset F_1\Mm\subset F_2\Mm\subset
			\cdots\subset \Mm$.
			\item A sheaf of filtered rings on $X$ is a sheaf $\Rr$ of rings on $X$, equipped with a sequence of sheaves of additive subgroups
			$F_0\Rr\subset F_1\Rr\subset F_2\Rr\subset
			\cdots \subset\Rr$ such that for every pair $(n,m)$ of nonnegative integers and local sections $a\in F_n\Rr$, $b\in F_m\Rr$, we have $ab\in F_{n+m}\Rr$.
			\item Let $\Rr=F_\bullet\Rr$ be a sheaf of filtered rings on $X$. Then a filtered left $\Rr$-module is a left $\Rr$-module $\Mm$, equipped with a filtration $F_\bullet \Mm$ such that for every pair $(n,m)$ of nonnegative integers and local sections $a\in F_n\Rr$, $v\in F_m\Mm$, we have $av\in F_{n+m}\Mm$.
		\end{enumerate}
	\end{definition}
	
	In particular, for a sheaf of filtered rings on $X$,
	multiplication induces a morphism
	$$
	F_n\Rr\otimes_{\ZZ_X} F_m\Rr\to F_{n+m}\Rr.
	$$

	\begin{prop-defn}
		For a topological space $X$ and a sheaf $\Rr$ of filtered rings on $X$, the following conditions are equivalent:
		\begin{enumerate}
			\renewcommand{\labelenumi}{(\alph{enumi})}
			\item For nonnegative integers $n,m$, an open subset $U\subset X$, and $a\in F_n\Rr(U)$, $b\in F_m\Rr(U)$, we have $ab-ba\in F_{n+m-1}\Rr(U)$.
			\item The associated graded sheaf $\Gr \Rr$ is commutative.
		\end{enumerate}
		We say that $\Rr$ is almost commutative if these equivalent conditions hold.
	\end{prop-defn}

	\begin{definition}\label{defn:exhaustive}
		A filtration $F_\bullet \Mm$ of a sheaf $\Mm$ on a topological space $X$ is called exhaustive if the canonical map $\varinjlim_n F_n \Mm\to \Mm$ is an isomorphism.
	\end{definition}
	
	\begin{remark}[{\cite[\href{https://stacks.math.columbia.edu/tag/0738}{Tag 0738} (3)]{stacks-project}}]
		Even if a filtration $F_\bullet \Mm$ is exaustive, $\varinjlim\Gamma(X,F_n \Mm)\to \Gamma(X,\Mm)$ is not an isomorphism in general. The map is an isomorphism if $X$ is quasi-compact.
	\end{remark}
	
	\begin{remark}
		Let $(X,\OO_X)$ be a scheme. Then the above notions for quasi-coherent sheaves on schemes $X$ can be defined inside $\Mod(\OO_X)$ since $\Mod_{\qc}(\OO_X)$ is closed under formation of colimits and finite limits of $\Mod(\OO_X)$ from \cite[\href{https://stacks.math.columbia.edu/tag/01IC}{Tag 01IC}, \href{https://stacks.math.columbia.edu/tag/01ID}{Tag 01ID}]{stacks-project}. In fact, notice that the restriction functors to open subsets respect small colimits since they are left adjoint.
	\end{remark}
	
	\subsection{Torsors}\label{sec:torsor}
	
	In this section, we give a quick summary on torsors attached to complexes, based on \cite{beilinsonbernstein1993}. We work in the settings of rings and sheaves for our purpose in section \ref{sec:tdo_and_pic_alg}.
	
	\begin{definition}
		For a group $G$, a $G$-torsor is a nonempty simply transitive $G$-set. We denote the category of $G$-torsors and $G$-equivariant maps by $\Tor(G)$.
	\end{definition}
	
	\begin{definition}
		Let $M^\bullet$ be a complex of additive groups with differential $d$. Suppose that $M^\bullet$ is concentrated in positive degrees.
		\begin{enumerate}
			\item An $M^\bullet$-torsor is a pair $(C,c)$ of an $M^1$-torsor $C$ and a map $c:C\to M^2$ such that $d\circ c=0$ and $c(a+\varphi)=da+c(\varphi)$ for $a\in M^1$ and $\varphi\in C$.
			\item A homomorphism $(C,c)\to (C',c')$ of $M^\bullet$-torsors is a morphism $f:C\to C'$ of $M^1$-torsors such that $c'\circ f=c$.
		\end{enumerate}
		We denote the category of $M^\bullet$-torsors by $\Tor(M^\bullet)$.
	\end{definition}
	
	\begin{example}
		Let $M$ be an additive group. Let $M[-1]$ be the complex concentrated in degree one with $M[-1]^1=M$. Then we have $\Tor(M[-1])=\Tor(M)$.
	\end{example}
	
	\begin{definition}[{\cite[(16.5.15)]{ega44}}]
		Let $\Gg$ be a sheaf of groups on $X$.
		\begin{enumerate}
			\item A $\Gg$-torsor is a $\Gg$-sheaf $\Cc$ with the following properties:
			\begin{enumerate}
				\item[(i)] The action is simply transitive on every open subset;
				\item[(ii)] There is an open covering $X=\cup_\lambda U_\lambda$ such that $\Cc(U_\lambda)$ is nonempty for every index $\lambda$.
			\end{enumerate}
			\item A homomorphism $\Cc\to \Cc'$ of $\Gg$-torsors is a $\Gg$-equivariant map of sheaves.
		\end{enumerate}
		We denote the groupoid of $\Mm$-torsors by $\Tor(\Mm)$.
	\end{definition}
	
	For a complex of additive groups on a topological space, we define its torsors in a similar way. We use similar notations for their category.
	
	\begin{proposition}\label{prop:tor_affine}
		Let $X$ be an affine scheme, and $\Mm^\bullet$ be a complex of sheaves of additive groups on $X$ concentrated in positive degrees. Suppose that $\Mm^1$ is endowed with the structure of a quasi-coherent $\OO_X$-module. Then $\Gamma(X,-)$ gives rise to a categorical equivalence
		\[\Tor(\Mm^\bullet)\simeq \Tor(\Gamma(X,\Mm^\bullet)).\]
	\end{proposition}
	
	This will be used in section \ref{sec:picalg} to see the equivalence of Picard algebras and Picard algebroids for affine schemes.
	
	\begin{proof}
		This is an easy consequence of \cite[Proposition (16.5.16)]{ega44}.
	\end{proof}

	Let us record basic operations on torsors.
	
	\begin{construction}
		Let $f:G\to H$ be a homomorphism of groups, and $C$ be a $G$-set. Define a left action of $G$ on $H\times C$ by $g(h,c)=(hf(g)^{-1},gc)$. We denote its quotient by $H\times^G C$.
	\end{construction}
	
	\begin{definition}[Change of groups]
		\begin{enumerate}
			\item Let $G\to H$ be a homomorphism of groups. Then we define a functor $\Tor(G)\to \Tor(H)$ by $C\mapsto H\times^G C$ as above with the action of $H$ for the multiplication of $H$ from the left side.
			\item Let $f:M^\bullet\to N^\bullet$ be a homomorphism of complexes of additive groups concentrated in positive degrees. Then define $\Tor(M^\bullet)\to \Tor(N^\bullet)$ by $(C,c)\mapsto (N^1\times^{M^1} C,c')$, where
			\[c'(\varphi)=dn+f(c(\varphi)).\]
		\end{enumerate}
	\end{definition}
	
	We define similar operations to homomorphisms of sheaves of (additive) groups and of their complexes on a topological space.
	
	\begin{definition}[Pullback]
		Let $f:X\to Y$ be a continuous map of topological spaces, $\Mm^\bullet$ be a complex of sheaves of additive groups on $Y$. Then the sheaf-theoretic pullback $f^{-1}$ determines a functor $f^{-1}:\Tor(\Mm^\bullet)\to\Tor(f^{-1}\Mm^\bullet)$ in a natural way.
	\end{definition}

	\subsection{Derivations and Differentials}\label{sec:derivation}
	For a homomorphism $\varphi:\OO_1\to\OO_2$ of sheaves of rings on a topological space $X$ and an $\OO_2$-module $\Ff$, let $\Der_{\OO_1}(\OO_2,\Ff)$ be the sheaf of $\varphi$-derivations into $\Ff$. If $\Ff=\OO_2$, we refer to $\Der_{\OO_1}(\OO_2,\Ff)$ as $\Theta_{\OO_2/\OO_1}$ or $\Der_{\OO_1}(\OO_2)$. If we are given a morphism $f:X\to S$ of ringed spaces, we always regard $\OO_X$ as an $f^{-1}\OO_S$-algebra sheaf by $f^\sharp$ (unless specified otherwise). If we are given an $f^{-1}\OO_S$-algebra sheaf $\OO$, write
	\[\begin{array}{cc}
		\Der_S(\OO,-)=\Der_{f^{-1}\OO_S}(\OO,-),
		&\Theta_{X/S}=\Theta_{\OO_X/f^{-1}\OO_S}.
	\end{array}\]
	
	\begin{lemma}[{\cite[\href{https://stacks.math.columbia.edu/tag/08RM}{Tag 08RM}]{stacks-project}}]
		Let $\OO_1\to\OO_2$ be a homomorphism of sheaves of rings on a topological space $X$, there is an $\OO_2$-module $\Omega^1_{\OO_2/\OO_1}$, equipped with an $\OO_1$-module homomorphism $d:\OO_2\to\Omega^1_{\OO_2/\OO_1}$ such that $-\circ d:\iHom_{\OO_2}(\Omega^1_{\OO_2/\OO_1},\Ff)\to\Der_{\OO_1}(\OO_2,\Ff)$ is an isomorphism for every $\OO_2$-module $\Ff$. In particular, the copresheaf of derivations is co-represented by $\Omega^1_{\OO_2/\OO_1}$.
	\end{lemma}
	
	Let us record a basic fact for convenience:
	
	\begin{lemma}[{\cite[\href{https://stacks.math.columbia.edu/tag/08RR}{Tag 08RR}]{stacks-project}}]
		Let $f:X\to Y$ be a continuous map of topological spaces, and $\OO_1\to\OO_2$ be a homomorphism of sheaves of rings on $Y$. Then we have a canonical isomorphism $f^{-1}\Omega^1_{\OO_2/\OO_1}\cong \Omega^1_{f^{-1}\OO_2/f^{-1}\OO_1}$.
	\end{lemma}
	
	For a morphism $f:X\to S$ of schemes, we will denote $\Omega^1_{\OO_X/f^{-1}\OO_S}$ by $\Omega^1_{X/S}$ as usual.
	
	\begin{lemma}[{\cite[\href{https://stacks.math.columbia.edu/tag/08S2}{Tag 08S2}, \href{https://stacks.math.columbia.edu/tag/01V3}{Tag 01V3}]{stacks-project}}]\label{lem:omega^1schemecase}
		Let $f:X\to S$ be a morphism of schemes.
		\begin{enumerate}
			\renewcommand{\labelenumi}{(\arabic{enumi})}
			\item The $\OO_X$-module $\Omega^1_{X/S}$ is quasi-coherent.
			\item Suppose that $f$ is locally of finite presentation. Then $\Omega^1_{X/S}$ is locally of finite presentation. In particular, $\Theta_{X/S}$ is a quasi-coherent $\OO_X$-module.
		\end{enumerate}
	\end{lemma}

	For the functorial property of $\Omega^1_{X/S}$ in terms of derivations, the following result is helpful:
	
	\begin{lemma}\label{lem:derivationvspullbackofdifferential}
		Suppose that we are given a diagram of schemes depicted as
		\[\begin{tikzcd}
			X\ar[rr, "f"]\ar[rd, "h"']&&Y\ar[ld, "g"]\\
			&S.
		\end{tikzcd}\]
		Then we have a canonical isomorphism
		$\Der_S (f^{-1}\OO_Y,-)\cong
		\iHom_{\OO_X}(f^\ast\Omega^1_{Y/S},-)$.
	\end{lemma}

	\begin{proof}
		This is straightforward:
		\[\begin{split}
			\iHom_{\OO_X}(f^\ast\Omega^1_{Y/S},-)
			&\cong
			\iHom_{f^{-1}\OO_Y}(f^{-1}\Omega^1_{Y/S},-)\\
			&\cong\iHom_{f^{-1}\OO_Y}
			(\Omega^1_{f^{-1}\OO_Y/h^{-1}\OO_S},-)\\
			&\cong\Der_S (f^{-1}\OO_Y,-).
		\end{split}\]
	\end{proof}
	
	Use the isomorphism $\Theta_{X/S}\cong(\Omega^1_{X/S})^\vee$ to get a pairing $\Omega^1_{X/S}\otimes_{\OO_X}\Theta_{X/S}\to\OO_X$ which we denote by $\langle-,-\rangle$. As usual, we put $\Omega^p_{X/S}:=\wedge^p\Omega^1_{X/S}$. We can extend the pairing of $\Omega^1_{X/S}$ and $\Theta_{X/S}$ to
	\[\langle-,-\rangle:\Omega^p_{X/S}\otimes_{\OO_X}\wedge^p\Theta_{X/S}\to\OO_X.\]
	\begin{proposition}
		Let $(X,\OO_X)$ be a ringed space, and $\Ff$, $\Gg$ be $\OO_X$-modules. Suppose that we are given an $\OO_X$-bilinear form $B:\Ff\otimes_{\OO_X}\Gg\to\OO_X$. For a nonnegative integer $p$, there is an $\OO_X$-bilinear form $B:\wedge^p\Ff\otimes_{\OO_X}\wedge^p\Gg\to\OO_X$ such that 
		\begin{equation}
			B(a_1\wedge a_2\wedge\cdots a_p\otimes b_1\wedge b_2\wedge\cdots b_p)=\sum_{\sigma\in\mathfrak{S}_p} \mathrm{sgn}(\sigma)\prod_{i=1}^p B(a_i\otimes b_{\sigma(i)})=\det(B(a_i,b_j))
			\label{eq:wedgepairing}
		\end{equation}
		for every open subset $U\subset X$ and sections $a_1,a_2,\ldots,a_p\in\Ff(U)$, $b_1,b_2,\ldots,b_p\in\Gg(U)$.
	\end{proposition}
	\begin{proof}
		Identify $\iEnd_{\OO_X}(\OO_X^p)$ with the sheaf of square matrices of size $p$ with $\OO_X$-entries. Define a map of sheaves
		$
		B^{p\times p}:\Ff^{\oplus p}\oplus\Gg^{\oplus p}\to\iEnd_{\OO_X}(\OO_X^p)
		$
		by
		\[B^{p\times p}(a_1, a_2,\ldots, a_p, b_1, b_2,\ldots, b_p)
		=(B(a_i,b_j)).\]
		The composite map $\det\circ B^{p\times p}:\Ff^{\oplus p}\oplus\Gg^{\oplus p}\to\OO_X$ gives rise to an $\OO_X$-bilinear form
		\[B:\wedge^p\Ff\otimes_{\OO_X}\wedge^p\Gg\to\OO_X\]
		since $\det\circ B^{p\times p}$ is multilinear in all variables and alternating in first and last $p$ variables. This satisfies the equality \eqref{eq:wedgepairing} by definition.
	\end{proof}

	We put
	\[\begin{array}{ccc}
		\Ff=\Omega^1_{X/S},&\Gg=\Theta_{X/S},&B=\langle-,-\rangle
	\end{array}\]
	to get a pairing $\langle-,-\rangle:\Omega^p_{X/S}\otimes_{\OO_X}\wedge^p\Theta_{X/S}\to\OO_X$.
	
	\begin{corollary}
		Suppose that the morphism $X\to S$ is smooth. Then the induced map $\Omega^p_{X/S}\to(\wedge^p\Theta_{X/S})^\vee$ is an isomorphism of $\OO_X$-modules.
	\end{corollary}
	
	\begin{proof}
		Recall that $\Omega^1_{X/S}$ is locally free of finite rank (see \cite[Proposition 2.5.4 (i)]{fu2011} for example).
		Since all the constructions are local, we may assume that $\Omega^1_{X/S}$ is free. The assertion is clear if $\Omega^1_{X/S}=0$. Henceforth assume that $\Omega^1_{X/S}\neq 0$. Choose a free basis $\theta_1,\theta_2,\ldots,\theta_n$. Let $\theta_1^\vee,\theta_2^\vee,\ldots,\theta_n^\vee\in\Hom_{\OO_X}(\Omega^1_{X/S},\OO_X)$ be its dual basis. For an integer $1\leq i\leq n$, write $\partial_i=\theta_i^\vee\circ d$. Then $\partial_1,\partial_2,\ldots,\partial_n$ form the dual basis of $\theta_1,\theta_2,\ldots,\theta_n$ with respect to the pairing $\langle-,-\rangle:\Omega^1_{X/S}\otimes_{\OO_X}\Theta_{X/S}\to\OO_X$.
		
		Let $\Lambda$ be the set of sequences $I=(i_1,i_2,\ldots,i_p)$ of integers with $1\leq i_1<i_2<\cdots<i_p\leq n$. For $I=(i_1,i_2,\ldots,i_p)\in\Lambda$, we denote $\theta_I=\theta_{i_1}\wedge \theta_{i_2}\wedge\cdots\theta_{i_p}$ (resp.\ $\partial_I=\partial_{i_1}\wedge \partial_{i_2}\wedge\cdots\partial_{i_p}$). Then $\theta_I$ (resp.\ $\partial_I$) with $I\in\Lambda$ form a free basis of $\Omega^p_{X/S}$ (resp.\ $\wedge^p\Theta_{X/S}$). Let $\partial_I^\vee\in(\wedge^p\Theta_{X/S})^\vee$ be the dual basis of $\partial_I$.
		
		The proof is completed by showing that the map $\Omega^p_{X/S}\to(\wedge^p\Theta_{X/S})^\vee$ sends $\theta_I$ to $\partial_I^\vee$ for every $I=(i_1,i_2,\ldots,i_p)\in\Lambda$. For this, it will suffice to show that for $J=(j_1,j_2,\ldots,j_p)\in \Lambda$, the following conditions are equivalent:
		\begin{enumerate}
			\renewcommand{\labelenumi}{(\alph{enumi})}
			\item $I=J$;
			\item $\langle \theta_I,\partial_J\rangle=1$;
			\item $\langle \theta_I,\partial_J\rangle\neq 0$.
		\end{enumerate}
		The assertion (b)$\Rightarrow$(c) follow since $\Omega^1_{X/S}$ is nonzero and free. Recall that $\langle\theta_I,\partial_J\rangle=\det(\delta_{i_l j_m})$. The implication (a)$\Rightarrow$(b) is now clear. Finally, suppose that (c) is satisfied. Since $\det$ is multilinear, for each $1\leq l\leq p$, there exists an integer $1\leq m_l\leq p$ such that $i_l=j_{m_l}$; otherwise, $\langle\theta_I,\partial_J\rangle=0$. We fix $m_l$ for each $l$. Then we can regard $i$ as a function
		\[\{1,2,\ldots,p\}\ni l\mapsto j_{m_l}\in\{j_1,j_2,\ldots,j_p\}.\] 
		Since $i_1<i_2<\cdots<i_p$, we have
		$j_{m_1}<j_{m_2}<\cdots<j_{m_p}$,
		i.e., $i$ is strictly monotone increasing. In particular, $i$ is injective. Since the domain and the taget of $i$ have the same cardinality $p$, it is a bijection. We now conclude $j_{m_l}=j_l$ for every $l$ again by the fact that $i$ is strictly monotone increasing.
	\end{proof}
	
	\begin{construction}[interior product]
		Let $X\to S$ be a smooth morphism of schemes. We then define an $\OO_X$-linear map $i:\Theta_{X/S}\otimes_{\OO_X}\Omega^\bullet_{X/S}\to\Omega^{\bullet-1}_{X/S}$ as follows:
		\begin{enumerate}
			\renewcommand{\labelenumi}{(\roman{enumi})}
			\item Consider the $\OO_X$-linear map
			\[\begin{split}
				\Theta_{X/S}\otimes_{\OO_X}\Omega^\bullet_{X/S}
				\otimes_{\OO_X}\wedge^{\bullet-1} \Theta_{X/S}
				&\xcong{C_{\Theta_{X/S},\Omega^\bullet_{X/S}}}
				\Omega^\bullet_{X/S}
				\otimes_{\OO_X}\Theta_{X/S}\otimes_{\OO_X}
				\wedge^{\bullet-1} \Theta_{X/S}\\
				&\to \Omega^\bullet_{X/S}
				\otimes_{\OO_X}
				\wedge^{\bullet} \Theta_{X/S}\\
				&\overset{\langle-,-\rangle}{\to}
				\OO_X,	
			\end{split}\]
			where the arrow in the second row is given by the multiplication.
			\item Pass to the adjoint to get
			$\Theta_{X/S}\otimes_{\OO_X}\Omega^\bullet_{X/S}\to(\wedge^{\bullet-1} \Theta_{X/S})^\vee$.
			\item Identify
			$(\wedge^{\bullet-1} \Theta_{X/S})^\vee$ with $\Omega^{\bullet-1}_{X/S}$.
		\end{enumerate}
	\end{construction}
	
	For a morphism $X\to S$ of schemes (resp.\ a homomorphism $k\to A$ of commutative rings), let $(\Omega^p_{X/S},d)$ (resp.\ $(\Omega^p_{A/k},d)$) denote the algebraic de Rham complex (\cite[section 16.6]{ega44}). We denote the sheaf $\Ker (d:\Omega^p_{X/S}\to\Omega^{p+1}_{X/S})$ (resp.\ the $k$-module $\Ker (d:\Omega^p_{A/k}\to\Omega^{p+1}_{A/k})$) of closed $p$-forms by $\Omega^{p,\cl}_{X/S}$ (resp.\ $\Omega^{p,\cl}_{A/k}$).
	
	\begin{example}[log differential]\label{ex:logdiff}
		We define $d\log:\OO_X^\times\to \Omega^1_{X/S}$ by $a\mapsto \frac{da}{a}$ for an open subset $U$ and $a\in \OO_X^\times(U)$. Then $d\log$ is a group homomorphism into $\Omega^{1,\cl}_{X/S}$. In particular, $d\log$ determines a homomorphism $\OO_X^\times\left[-1\right]\to\Omega^\bullet_{X/S}$ of complexes. We call $d\log$ the log differential.
	\end{example}
	
	\begin{construction}[Functoriality of de Rham complexes, {\cite[\href{https://stacks.math.columbia.edu/tag/08RS}{Tag 08RS}, \href{https://stacks.math.columbia.edu/tag/0FKL}{Tag 0FKL}]{stacks-project}}]\label{cons:deRham}
		Suppose that we are given a diagram of schemes
		\[\begin{tikzcd}
			X\ar[r, "f"]\ar[d]&Y\ar[d, "g"]\\
			S\ar[r]&T.
		\end{tikzcd}\]
		Then we define a canonical homomorphism $f^\ast:\Omega^\bullet_{Y/T}\to f_\ast\Omega^\bullet_{X/S}$ in the following way: The composition of $f^\sharp$ and the differential $d:\OO_X\to\Omega^1_{X/S}$ determines an $f^{-1}\OO_Y$-linear map $\Omega^1_{f^{-1}\OO_Y/(g\circ f)^{-1}\OO_T}\to \Omega^1_{X/S}$ through the isomorphism
		\[\Der_T(f^{-1}\OO_Y,\Omega^1_{X/S})
		\cong\iHom_{f^{-1}\OO_Y}(\Omega^1_{f^{-1}\OO_Y/(g\circ f)^{-1}\OO_T},\Omega^1_{X/S})\]
		(Lemma \ref{lem:derivationvspullbackofdifferential}).
		The desired map is obtained from
		\[f^{-1}\Omega^\bullet_{Y/T}
		\cong \wedge^\bullet f^{-1}\Omega^1_{Y/T}
		\cong \wedge^\bullet \Omega^1_{f^{-1}\OO_Y/(g\circ f)^{-1}\OO_T}
		\to \Omega^\bullet_{X/S}.\]
		Similarly, for a commutative ring $k$ and a homomorphism $f:A\to B$ of commutative $k$-algebras, we can define a canonical map $f_\ast:\Omega^\bullet_{A/k}\to\Omega^\bullet_{B/k}$ of the de Rham complexes.
	\end{construction}

	A routine computation in the basic theory of manifolds shows the following result:
	
	\begin{lemma}\label{lem:pairingwithdifferential}
		Let $k\to A$ be a homomorphism of commutative rings. Then we have
		\[\begin{split}
			\langle d\eta,\partial_1\wedge\partial_2 \cdots\wedge \partial_{p+1}\rangle
			&=\sum_{i=1}^{p+1}(-1)^{i+1} \partial_i(\langle \eta,\partial_1\wedge\partial_2\wedge \cdots\wedge \overset{\vee}{\partial_i}\wedge \cdots\wedge \partial_{p+1}\rangle)\\
			&+\sum_{1\leq i<j\leq p+1} (-1)^{i+j}
			\langle\eta,\left[\partial_i,\partial_j\right]\wedge
			\partial_1\wedge\partial_2\wedge \cdots\wedge \overset{\vee}{\partial_i}\wedge \cdots\wedge\overset{\vee}{\partial_j}\wedge\cdots\wedge \partial_{p+1}\rangle
		\end{split}\]
		for a nonnegative integer $p$, $\eta\in\Omega^p_{A/k}$ and $\partial_1,\ldots,\partial_{p+1}\in\Theta_{A/k}$. A similar equality also holds in the sheaf setting.
	\end{lemma}
	
	\begin{definition}[Lie derivative]
		Let $f:X\to S$ be a smooth morphism of schemes. Define $L:\Theta_{X/S}\otimes_{f^{-1}\OO_S}\Omega^\bullet_{X/S}\to\Omega^{\bullet}_{X/S}$ by $d\circ i+i\circ (\id_{\Theta_{X/S}}\otimes_{f^{-1}\OO_S} d)$.
	\end{definition}
	
	Regard $\Theta_{X/S}$ as a Lie algebra sheaf over $f^{-1}\OO_S$ for the commutator. The standard argument for manifolds shows the following result:
	
	\begin{proposition}\label{prop:liederivativeproperties}
		Let $f:X\to S$ be a smooth morphism of schemes, and $p\geq 0$ be an integer.
		\begin{enumerate}
			\renewcommand{\labelenumi}{(\arabic{enumi})}
			\item For local sections $\eta\in\Omega^p_{X/S}$, $\partial,\partial_1,\ldots,\partial_p\in\Theta_{X/S}$, we have
			\begin{flalign*}
				&\langle L_\partial\eta,\partial_1\wedge\cdots\wedge\partial_p\rangle\\
				&=\partial(\langle\eta,\partial_1\wedge\cdots\wedge\partial_p\rangle)
				-\sum_{i=1}^p \langle \eta,\partial_1\wedge\cdots\wedge
				\left[\partial,\partial_i\right] \wedge\cdots\partial_p\rangle.
			\end{flalign*}
			\item The Lie derivative $L$ defines the structure of a $\Theta_{X/S}$-module on $\Omega^p_{X/S}$.
			\item For an open subset $U\subset X$ and local sections $a\in\OO_X(U)$, $\partial\in\Theta_{X/S}(U)$, $\eta\in\Omega^p_{X/S}(U)$, we have
			$L_\partial (a\eta)=a L_\partial\eta +\partial(a)\eta$.
			\item We have
			$da\wedge i_\partial(\eta)+i_\partial(da\wedge \eta)
			=\partial(a)\eta$
			for local sections $a\in\OO_X$, $\eta\in\Omega^p_{X/S}$, $\partial\in\Theta_{X/S}$.
		\end{enumerate}
	\end{proposition}
	
	\begin{definition}[Canonical bundle]\label{defn:canonicalbundle}
		Let $f:X\to S$ be a smooth morphism of schemes. Then we define a line bundle $\omega$ on $X$ as follows: for $n\geq 0$, let $U_n=\{x\in X:~\rank\Omega^1_{X/S,x}=n\}$. Then $U_n$ form a disjoint open covering of $X$ since $\Omega^1_{X/S}$ is locally free of finite rank. We define $\omega_{X/S}$ by $\omega_{X/S}|_{U_n}=\Omega^n_{X/S}|_{U_n}$.
	\end{definition}
	
	\begin{corollary}\label{cor:omegaisrightDmod}
		The $\OO_X$-module $\omega_{X/S}$ is a left $\Theta_{X/S}$-module for the Lie derivative. Moreover, for an open subset $U\subset X$ and local sections $a\in\OO_X(U)$, $\partial\in\Theta_{X/S}(U)$, $\omega\in\omega_{X/S}(U)$, the following conditions are satisfied:
		\begin{enumerate}
			\renewcommand{\labelenumi}{(\roman{enumi})}
			\item $L_\partial(a\omega)=a L_\partial \omega
			+\partial(a)\omega$;
			\item $L_{a\partial}\omega=L_{\partial}(a\omega)$.
		\end{enumerate}
	\end{corollary}
	
	This implies that $\omega_{X/S}$ is a right $\D_{X/S}$-module for the minus of the Lie derivative.
	
	Recall that in the classical theory of $\D$-modules, relative local coordinates were used on the course of proofs of Kashiwara's equivalence and the base change theorem. Let us see the existence of relative local coordinates over general bases:
	
	\begin{lemma}[Local coordinate]\label{lem:coordinate}
		Let $S$ be a scheme, and
		\[j:U\to \bfA^n_S=S\times_{\Spec\ZZ} \Spec\ZZ\left[y_1,y_2,\ldots,y_n\right]\]
		be an \'etale morphism of smooth $S$-schemes. Set $x_p=j^\sharp(y_p)$ for $1\leq p\leq n$. Then
		\begin{enumerate}
			\item $\{dx_p\in\Omega^1_{U/S}(U):~1\leq p\leq n\}$ is a free $\OO_U$-basis of $\Omega^1_{U/S}$.
			\item The dual basis $\{\partial_p\coloneqq (dx_p)^\vee\in\Theta_{U/S}(U)\}$ satisfies
			\begin{enumerate}
				\item[(i)] $\left[\partial_p,x_q\right]=
				\left\{\begin{array}{cc}
					1&(p=q)\\
					0&(p\neq q),
				\end{array}\right.$
				\item[(ii)] $\left[\partial_p,\partial_q\right]=0$ for all $i,j$.
			\end{enumerate}
		\end{enumerate}
	\end{lemma}
	
	\begin{proof}
		The first claim follows from the canonical isomorphism $j^\ast\Omega^1_{\bfA^n_S/S}\cong\Omega^1_{U/S}$.
		
		Part (i) in (2) is immediate from the definitions:
		$\left[\partial_p,x_q\right]
		=\langle dx_q,\partial_p\rangle$.
		For (ii), it will suffice to show $\langle \theta,\left[\partial_p,\partial_q\right]\rangle=0$ for all $\theta\in\Omega^1_{U/S}$. We may assume $\theta=dx_r$ for some $r$ by (1). Then the equality follows from
		$\langle dx_r,\left[\partial_p,\partial_q\right]\rangle
		=\langle d(\langle \partial_q,x_r\rangle),\partial_p\rangle
		-\langle d(\langle \partial_p,x_r\rangle),\partial_q\rangle
		=0$.
		This completes the proof.
	\end{proof}

	\begin{theorem}[Relative local coordinate]\label{thm:relativecoordinate}
		Let $i:Y\hookrightarrow X$ be an immersion of smooth $S$-schemes, and $y\in Y$. Then there exist an open neighborhood $U\ni x=j(y)$ in $X$, a set
		\[\{x_p\in\OO_X(U):~1\leq p\leq n\}\]
		of local sections, and $r\in\{0,1,\ldots,n-1\}$ with the following properties:
		\begin{enumerate}
			\renewcommand{\labelenumi}{(\roman{enumi})}
			\item $\{dx_p\in\Omega^1_{X/S}(U):~1\leq p\leq n\}$ is a free basis of $\Omega^1_{X/S}$,
			\item $\{(dx_p)|_{Y\cap U}\in\Omega^1_{Y/S}(Y\cap U):~
			1\leq p\leq r\}$ is a free basis of
			$\Omega^1_{(Y\cap U)/S}$,
			\item $(dx_p)|_{Y\cap U}=0$ for $r+1\leq p\leq n$, 
			\item $\left[\partial_p,x_q\right]=\left\{
			\begin{array}{cc}
				1&(p=q)\\
				0&(p\neq q),
			\end{array}\right.$
			\item $\left[\partial_p,\partial_q\right]=0$ for all $1\leq p,q\leq n$.
		\end{enumerate}
		Here $(\partial_i)$ is the dual basis of $(dx_i)$.
	\end{theorem}
	
	\begin{proof}
		We may replace $X$ by an open neighborhood of $x$ to assume that $i$ is a closed immersion. Then in view of \cite[Corollaire 17.12.2]{ega44}, one can find an \'etale morphism $j:U\to \bfA^n_S$ with $x\in j(U)$ and a Cartesian diagram
		\[\begin{tikzcd}
			Y\cap U\ar[r, "i|_{Y\cap U}"]\ar[d]
			&U\ar[d, "j"]\\
			\bfA_S^r\ar[r, hook]&\bfA^n_S.
		\end{tikzcd}\]
		If we write $\bfA^r_S=S\times_{\Spec\ZZ} \Spec\ZZ\left[y_1,y_2,\ldots,y_r\right]$ and $\bfA^n_S=S\times_{\Spec\ZZ} \Spec\ZZ\left[y_1,y_2,\ldots,y_n\right]$ then
		the bottom horizontal arrow is given by
		\[y_p\mapsto \left\{\begin{array}{cc}
			y_p&(1\leq p\leq r)\\
			0&(r+1\leq p\leq n).
		\end{array}\right.\]
		Then $x_p=j^\sharp(y_p)$ satisfy the conditions by Lemma \ref{lem:coordinate}.
	\end{proof}
	
	\subsection{Equivariant structure on the (co)tangent sheaf}\label{sec:equivsheaf}
	Let $x:X\to S$ be a morphism of schemes, and $G$ be a group scheme over $S$ with an action $a:G\times_S X\to X$. In this section, we recall the definition of the canonical action of $G$ on the (co)tangent sheaf on $X$.
	
	Define morphisms $\pr_1,\pr_2,\pr_{23},m_G,i$ by
	\[\begin{array}{cc}
		\pr_1:G\times_S X\to G;~(g,x)\mapsto g,
		&\pr_2:G\times_S X\to G;~(g,x)\mapsto x,
	\end{array}\]
	\[\pr_{23}:G\times_S G\times_S X\to G\times_S X;~
	(g_1,g_2,x)\mapsto (g_2,x),\]
	\[\begin{array}{cc}
		m_G:G\times_S G\to G:~(g_1,g_2)\mapsto g_1g_2,&i:X\to G\times_S X;~x\mapsto (e,x).
	\end{array}\]
	Recall that a quasi-coherent $\OO_X$-module $\Ff$, equipped with an isomorphism $I:a^\ast \Ff\cong \pr_2^\ast\Ff$ is called $G$-equivariant if the equalities
	\[\begin{array}{cc}
		\pr_{23}^\ast I\circ(\id_G\times a)^\ast I=
		(m_G\times \id_X)^\ast I,
		&i^\ast I=\id_{\Ff}
	\end{array}\]
	hold.
	
	Notice that the canonical homomorphisms
	\[\begin{array}{cc}
		\pr_1^\ast\Omega^1_{G/S}\to\Omega^1_{G\times_S X/S},
		&\pr_2^\ast\Omega^1_{X/S}\to\Omega^1_{G\times_S X/S}
	\end{array}\]
	(recall Construction \ref{cons:deRham}) sum up to an isomorphism
	\begin{equation}
		\pr_1^\ast\Omega^1_{G/S}\oplus \pr_2^\ast\Omega^1_{X/S}
		\cong\Omega^1_{G\times_S X/S}\label{eq:cot_product}.
	\end{equation}
	The composite map $a^\ast \Omega^1_{X/S}\to \pr_2^\ast\Omega^1_{X/S}$ of the canonical homomorphism
	$a^\ast\Omega^1_{X/S}\to \Omega^1_{G\times_S X/S}$
	with the projection
	\[\Omega^1_{G\times_S X/S}\overset{\eqref{eq:cot_product}}{\cong}\pr_1^\ast\Omega^1_{G/S}\oplus \pr_2^\ast\Omega^1_{X/S}\overset{\pr_2}{\to}\pr_2^\ast\Omega^1_{X/S}\]
	is an isomorphism.
	Moreover, this is the structure of a $G$-equivariant sheaf on $\Omega^1_{X/S}$.
	
	Finally, suppose that $x$ is smooth. Then $\Omega^1_{X/S}$ is locally free of finite rank. Therefore take the dual to get the structure of a $G$-equivariant sheaf on $\Theta_{X/S}$.
	
	\subsection{Homological algebra of $\Aa$-modules}\label{sec:homologicalalgebra}
	In this section, we collect basic facts on unbounded derived categories of modules over sheaves of (possibly noncommutative) rings. In particular, we aim to introduce the derived functors $-\otimes^L_{\Aa}-$ and $Rf_\ast$. In fact, we can readily generalize \cite[section 2]{lipman2009} to the setting of noncommutative rings by careful discrimination of right and left modules.
	
	Let $Y$ be a topological space, and $\Aa$ be a sheaf of rings on $Y$.
	
	We first review the derived tensor product. From perspectives of homotopy theory, we temporally adopt the homological grading. To construct the derived tensor product, let us recall a reasonable class of resolutions:
	
	\begin{definition}[{\cite[Definition 18.5.1. (iii)]{kashiwaraschapira2}}]
		A right (resp.~left) $\Aa$-module $\Ff$ is flat if the functor $-\otimes_{\Aa}\Ff:\Mod(\Aa)\to \Mod(\ZZ_Y)$ (resp.~$\Ff\otimes_{\Aa}-:\Mod_{\mathrm{r}}(\Aa)\to \Mod(\ZZ_Y)$) is exact.
	\end{definition}
	
	\begin{definition}[{\cite[Definition (2.5.1)]{lipman2009}}]
		A complex $\Pp_\bullet$ of left (resp.\ right) $\Aa$-modules is called q-flat if for every quasi-isomorphism $\Ff_\bullet\to\Gg_\bullet$ of right (resp.\ left) $\Aa$-modules, the induced map $\Ff_\bullet\otimes_{\Aa}\Pp_\bullet\to \Gg_\bullet\otimes_{\Aa}\Pp_\bullet$ (resp.\ $\Pp_\bullet\otimes_{\Aa}\Ff_\bullet\to \Pp_\bullet\otimes_{\Aa}\Gg_\bullet$) is a quasi-isomorphism.
	\end{definition}
	
	\begin{remark}[{\cite[Example (2.5.2)]{lipman2009}}]\label{rem:q-flat_stalk}
		A complex of $\Aa$-modules is q-flat if and only if so are its stalks.
	\end{remark}
	
	\begin{remark}[{\cite[Examples (2.5.3), (2.5.4)]{lipman2009}}]\label{rem:q-flat_triangulatedsubcategory}
		The q-flat complexes form a triangulated full subcategory of $K(\Aa)$. Moreover, it is closed under small filtered colimits and therefore under small coproducts.
	\end{remark}
	
	\begin{remark}[{\cite[Example (2.5.4)]{lipman2009}}]\label{rem:q-flat_boundedabove}
		A bounded-below chain complex $\mathcal{P}_\bullet$ of left $\Aa$-modules is q-flat if and only if it is degree-wisely flat, i.e., $\mathcal{P}_n$ is flat as a left $\Aa$-module for every integer $n$.
	\end{remark}	
	
	Since we work with bimodules in this paper, suppose that we are given more algebras: Let $\OO_Y$ be a sheaf of commutative rings\footnote{In this paper, the main example is $y^{-1}\OO_S$ for $y:Y\to S$ a morphism of schemes, not the structure sheaf of $Y$.}. Let $\Aa$ and $\Bb$ be $\OO_Y$-algebras. We remark that $\OO_Y$ is central in $\Aa$ and $\Bb$. Notice that we can identify the category of $(\Aa,\Bb)$-bimodules over $\OO_Y$ (i.e., the $(\Aa,\Bb)$-bimodule objects in $\Mod(\OO_Y)$) with $\Mod(\Aa\otimes_{\OO_Y}\Bb^{\op})$.
	
	\begin{proposition}\label{prop:projectiveclass}
		\begin{enumerate}
			\item Let $\Pp\subset \Mod_{\mathrm{r}}(\Aa)$ be the full subcategory consisting of flat right $\Aa$-modules. Then $(\Pp,\Mod(\Aa\otimes_{\OO_Y}\Bb^{\op}))$ and $(\Pp,\Mod(\Aa))$ are $-\otimes_{\Aa}-$-projective in the dual sense of \cite[Definition 13.4.2]{kashiwaraschapira2}.
			\item Let $\Pp\subset K_{\mathrm{r}}(\Aa)$ be the full subcategory consisting of q-flat complexes of right $\Aa$-modules. Then $(\Pp,K(\Aa))$ is $-\otimes_{\Aa}-$-projective in the sense of \cite[Definition 10.3.9]{kashiwaraschapira2}. Similarly, $(\Pp,K(\Aa\otimes_{\OO_Y}\Bb))$ is $-\otimes_{\Aa}-$-projective.
			\item Let $\Pp\subset K(\Aa\otimes_{\OO_Y}\Bb)$ be the full subcategory consisting of complexes of $(\Aa,\Bb)$-modules over $\OO_Y$ which are q-flat as complexes of left $\Aa$-modules. Suppose that $\Bb$ is flat over $\OO_Y$. Then $(K_{\mathrm{r}}(\Aa),\Pp)$ is $-\otimes_{\Aa}-$-projective.
		\end{enumerate}
	\end{proposition}
	
	\begin{proof}
		Part (1) follows from the definition of flat right $\Aa$-modules and \cite[Propositions 18.5.4 and 13.4.4]{kashiwaraschapira2}.
		
		We next prove (2). Every complex $\Ff_\bullet$ of right or left $\Aa$-modules admits a q-flat resolution by (1), \cite[Lemma 14.4.1]{kashiwaraschapira2} and Remarks \ref{rem:q-flat_triangulatedsubcategory}, \ref{rem:q-flat_boundedabove} (see \cite[Definition 11.2.5 and Proposition 11.2.8]{kashiwaraschapira2} if necessary to recall the definition of the structure of the triangulated category on $K(\Aa)$). Suppose that we are given a q-flat complex $\Ff_\bullet$ of right $\Aa$-modules and a complex $\Mm_\bullet$ of left $\Aa$-modules or of $(\Aa,\Bb)$-bimodules. If $\Mm_\bullet$ is acyclic then so is $\Ff_\bullet\otimes_{\Aa}\Mm_\bullet$ by definition of q-flat complexes. If $\Ff_\bullet$ is acyclic then so is $\Ff_\bullet\otimes_{\Aa}\Mm_\bullet$ by the argument in the first paragraph of \cite[(2.5.7)]{lipman2009}. This completes the proof of (2).
		
		For (3), we see that every complex $\Gg_\bullet$ in $K(\Aa\otimes_{\OO_Y}\Bb^{\op})$ is quasi-isomorphic to a member of $\Pp$. In fact, choose a q-flat resolution $\Ff_\bullet$ of $\Gg_\bullet$ as a complex of left $\Aa\otimes_{\OO_Y}\Bb^{\op}$-modules (use (2)). Since $\Bb$ is flat over $\OO_Y$, $\Ff_\bullet$ is q-flat as a complex of left $\Aa$-modules. The rest goes the same line as (2).
	\end{proof}
	
	\cite[Proposition 10.3.10]{kashiwaraschapira2} immediately implies:
	
	\begin{corollary}\label{cor:derivedtensorproduct}
		\begin{enumerate}
			\item The functor $-\otimes_{\Aa}-:\Mod_{\mathrm{r}}(\Aa)\times\Mod(\Aa)\to\Mod(\ZZ_Y)$ admits a left derived functor $-\otimes^L_{\Aa}-:D_{\mathrm{r}}(\Aa)\times D(\Aa)\to D(\ZZ_Y)$.
			\item The functor $-\otimes_{\Aa}-:\Mod_{\mathrm{r}}(\Aa)\times\Mod(\Aa\otimes_{\OO_Y}\Bb^{\op})\to\Mod_{\mathrm{r}}(\Bb)$ admits a left derived functor $-\otimes^L_{\Aa}-:D_{\mathrm{r}}(\Aa)\times D(\Aa\otimes_{\OO_Y}\Bb^{\op})\to D_{\mathrm{r}}(\Bb)$.
			\item Then the derived functors of (1) and (2) are compatible in the sense that the diagram
			\[\begin{tikzcd}
				D_{\mathrm{r}}(\Aa)\times D(\Aa\otimes_{\OO_Y}\Bb^{\op})
				\ar[r, "-\otimes^L_{\Aa}-"]\ar[d]
				&D_{\mathrm{r}}(\Bb)\ar[d]\\
				D_{\mathrm{r}}(\Aa)\times D(\Aa)\ar[r, "-\otimes^L_{\Aa}-"]
				&\Mod(\ZZ_Y)
			\end{tikzcd}\]
			is 2-commutative, where the vertical arrows are defined by restriction.
		\end{enumerate}
	\end{corollary}
	
	We next aim to introduce the noncommutative analog of the adjunction $(f^\ast,f_\ast)$ of functors between the categories of modules on ringed spaces. Namely, let $f:X\to Y$ be a continuous map of topological spaces, $\Aa$ be a sheaf of rings on $Y$, and $f^{-1}\Aa\to \Bb$ be a homomorphism of sheaves of rings on $X$. Then we have an adjunction
	\begin{equation}
		\Bb\otimes_{f^{-1}\Aa} f^{-1}(-):\Mod(\Aa)\rightleftarrows\Mod(\Bb):f_\ast.
		\label{eq:adjunction}
	\end{equation}
	
	We have two possible ways to define $D(\Aa)\to D(\Bb)$ from $\Bb\otimes_{f^{-1}\Aa} f^{-1}(-)$:
	\begin{itemize}
		\item Compose the functor $f^{-1}:D(\Aa)\to D(f^{-1}\Aa)$ induced from the exact functor
		\[f^{-1}:\Mod(\Aa)\to \Mod(f^{-1}\Aa)\]
		with $\Bb\otimes^L_{f^{-1}\Aa}-:D(f^{-1}\Aa)\to D(\Bb)$. Here we use q-flat resolutions in the second variable to define $\Bb\otimes^L_{f^{-1}\Aa}-$.
		\item Take the left derived functor of $\Bb\otimes_{f^{-1}\Aa} f^{-1}(-):\Mod(\Aa)\to\Mod(\Bb)$;
	\end{itemize}
	
	\begin{proposition}\label{prop:leftderivedfunctor}
		\begin{enumerate}
			\item The flat left $\Aa$-modules form a $\Bb\otimes_{f^{-1}\Aa} f^{-1}(-)$-projective full subcategory of $\Mod(\Aa)$.
			\item The functor
			$\Bb\otimes_{f^{-1}\Aa} f^{-1}(-):\Mod(\Aa)\to\Mod(\Bb)$
			admits a left derived functor
			\[D(\Aa)\to D(\Bb).\]
			Moreover, the left derived functor is computed by q-flat resolutions.
		\end{enumerate}
		
	\end{proposition}
	
	\begin{proof}
		We first prove (1). For this, we check the dual conditions of \cite[Corollary 13.3.8]{kashiwaraschapira2}. This is immediate from \cite[Proposition 18.5.4 (i)]{kashiwaraschapira2} since $f^{-1}:\Mod(\Aa)\to\Mod(f^{-1}\Aa)$ is exact and respects flat modules.
		
		We next prove (2). In view of \cite[Proposition 10.3.3 (b)]{kashiwaraschapira2}, it will suffice to show that q-flat complexes of left $\Aa$-modules are $\Bb\otimes_{f^{-1}\Aa} f^{-1}(-)$-projective in the sense of \cite[Definition 10.3.2]{kashiwaraschapira2}. We have seen in Proposition \ref{prop:projectiveclass} that every complex of left $\Aa$-modules admits a q-flat resolution. Suppose that $\Ff_\bullet$ is an acyclic q-flat complex of left $\Aa$-modules. Then $f^{-1}\Ff_\bullet$ is an acyclic q-flat complex of left $f^{-1}\Aa$-modules by Remark \ref{rem:q-flat_stalk}. Then \cite[(2.5.7)]{lipman2009} implies that $\Bb\otimes_{f^{-1}\Aa} f^{-1}\Ff_\bullet$ is acyclic. This completes the proof of (2).
	\end{proof}
	
	The two functors from $D(\Aa)$ to $D(\Bb)$ are equivalent since $f^{-1}$ respects q-flat complexes. 
	
	To work with right derived functors, recall
	
	\begin{definition}[{\cite[Definition (2.3.1)]{lipman2009}}]
		A complex $\Ii\in K(\Aa)$ is called q-injective\footnote{This terminology is taken from \cite{lipman2009}. Such a complex is called K-injective in other literature.} if for every quasi-isomorphism $\Gg^\bullet_1\to\Gg^\bullet_2$, the induced map $\Hom_{K(\Aa)}(\Gg^\bullet_2,\Ii^\bullet)\to\Hom_{K(\Aa)}(\Gg^\bullet_1,\Ii^\bullet)$ is surjective.
	\end{definition}
	
	\begin{lemma}\label{lem:q-inj}
		For a complex $\Ii\in K(\Aa)$, the following conditions are equivalent:
		\begin{enumerate}
			\renewcommand{\labelenumi}{(\alph{enumi})}
			\item $\Ii$ is q-injective;
			\item $\Ii$ is homotopically injective in the sense of \cite[Definition 14.1.4 (i)]{kashiwaraschapira2}, i.e., for any acyclic complex $\Gg^\bullet\in K(\Aa)$, $\Hom_{K(\Aa)}(\Gg^\bullet,\Ii^\bullet)\cong\{0\}$;
			\item for every quasi-isomorphism $\Gg^\bullet_1\to\Gg^\bullet_2$, the induced map \[\Hom_{K(\Aa)}(\Gg^\bullet_2,\Ii^\bullet)\to\Hom_{K(\Aa)}(\Gg^\bullet_1,\Ii^\bullet)\]
			is bijective.
		\end{enumerate}
	\end{lemma}
	
	\begin{proof}
		Condition (c) clearly implies (a). Suppose that (a) is satisfied. Let $\Gg^\bullet$ be an acyclic complex. Then apply (a) to $\Gg^\bullet_1=\Gg^\bullet$ and $\Gg^\bullet_2=0$ to deduce $\Hom_{K(\Aa)}(\Gg^\bullet,\Ii^\bullet)\cong\{0\}$. Hence (b) follows. For the equivalence of (b) and (c), see \cite[Theorem 14.3.1 (i)]{kashiwaraschapira2}.
	\end{proof}

	\begin{proposition}\label{prop:rightderivedfunctor}
		\begin{enumerate}
			\item The functor $f_\ast:\Mod(\Bb)\to\Mod(\Aa)$ admits a right derived functor $D(\Bb)\to D(\Aa)$. Moreover, the right derived functor is computed by q-injective resolutions.
			\item The right derived functor of (1) is right adjoint to the left derived functor of Proposition \ref{prop:leftderivedfunctor} (2).
			\item For a q-flat complex $\Ff^\bullet$\footnote{We use the cohomological grading here since a q-injective complex appears here.} of left $\Aa$-modules and a q-injective complex $\Ii^\bullet$ of left $\Bb$-modules, the canonical map
			\[\Hom_{K(\Aa)}(\Ff^\bullet,f_\ast \Ii^\bullet)
			\to \Hom_{D(\Aa)}(\Ff^\bullet,f_\ast \Ii^\bullet)\]
			is an isomorphism.
		\end{enumerate}
	\end{proposition}
	
	\begin{proof}
		Part (1) follows from \cite[Theorem 14.3.1]{kashiwaraschapira2} and \cite[Corollary (2.3.2.3)]{lipman2009}. Part (2) follows from Proposition \ref{prop:leftderivedfunctor} (1) and \cite[Theorem 14.4.5]{kashiwaraschapira2}. One can verify (3) and (4) in a similar way to \cite[Corollary (3.2.2)]{lipman2009}.
	\end{proof}
	
	In particular, we obtain a right derived functor $Rf_\ast:D(f^{-1}\Aa)\to D(\Aa)$ by putting $\Bb=f^{-1}\Aa$. Therefore we have two ways to define $D(\Bb)\to D(\Aa)$ again:
	
	\begin{itemize}
		\item Compose the restriction functor $D(\Bb)\to D(f^{-1}\Aa)$ with the right derived functor $Rf_\ast:D(f^{-1}\Aa)\to D(\Aa)$.
		\item Take the right derived functor of $f_\ast:\Mod(\Bb)\to\Mod(\Aa)$;
	\end{itemize}
	
	These are equivalent by \cite[Proposition 14.4.7]{kashiwaraschapira2}

	One can also prove by \cite[Proposition 14.4.7]{kashiwaraschapira2} that the functor $Rf_\ast$ is essentially independent of $\Aa$ and $\Bb$ in the following sense:
	
	\begin{proposition}\label{prop:independent}
		Suppose that we are given a homomorphism $p:\Aa\to\Aa'$ of sheaves of rings on $X$, a homomorphism $q:\Bb\to\Bb'$ of sheaves of rings on $Y$, and homomorphisms $u:\Bb\to f^{-1}\Aa$, $v:\Bb'\to f^{-1}\Aa'$ such that $u\circ f^{-1}p=q\circ v$;
		\[\begin{tikzcd}
			\Bb\ar[r, "u"]\ar[d, "q"']&f^{-1}\Aa\ar[d, "f^{-1}p"]\\
			\Bb'\ar[r, "v"]&f^{-1}\Aa'.
		\end{tikzcd}\]
		Then the diagram
		\[\begin{tikzcd}
			D(\Aa')\ar[r, "Rf_\ast"]\ar[d]&D(\Bb')\ar[d]\\
			D(\Aa)\ar[r, "Rf_\ast"]&D(\Bb),
		\end{tikzcd}\]
		is 2-commutative, where the vertical arrows are obtained by restriction of the actions. In particular, for $\Mm^\bullet\in D(f^{-1}\Aa)$, $Rf_\ast\Mm^\bullet$ can be computed by the underlying structure of the complex of sheaves of abelian groups.
	\end{proposition}
	
	As an application of this generalized formalism, one can define the cohomology groups for complexes of left $\Aa$-modules. Write $A=\Gamma(Y,\Aa)$. Let $Z$ be the single point space. We denote the unique map $Y\to Z$ by $p$. Regard $A$ as a sheaf of rings on $Z$. Then one can identify the functor
	$\Gamma(Y,-):\Mod(\Aa)\to\Mod(A)$
	with $p_\ast:\Mod(\Aa)\to\Mod(A)$ attached to the counit $p^{-1}A\to \Aa$. We thus obtain a right derived functor $R\Gamma(Y,-):D(\Aa)\to D(A)$. For each integer $i$, set $H^i(Y,-)=H^i\circ R \Gamma(Y,-):D(\Aa)\to \Mod(A)$.
	For each left $\Aa$-module $\Mm$, the underlying structure of an abelian group on $H^i(Y,\Mm)$ is independent of $\Aa$. That is, forget the action of $\Aa$ to regard $\Mm$ just as a sheaf of abelian groups on $Y$. Then the attached $i$-th cohomology group is isomorphic to $H^i(Y,\Mm)$ as an abelian group since both groups can be computed by an $\Aa$-linear flasque resolution.
	
	Finally, let us give a basic observation:
	
	\begin{construction}[Projection map]\label{cons:projmap}
		Let $f:X\to Y$ be a morphism of ringed spaces, and $\Aa$, $\Bb$ be $\OO_Y$-algebras. Then define a natural transformation 
		\begin{equation}
			Rf_\ast(-)\otimes^L_{\Aa}-\to Rf_\ast(-\otimes^L_{f^{-1}\Aa}f^{-1}(-))
			\label{eq:projmap}
		\end{equation}
		of functors from $D_{\mathrm{r}}(f^{-1}\Aa)\times D(\Aa\otimes_{\OO_Y}\Bb^{\op})$ to $D_{\mathrm{r}}(\Bb)$ as follows: Apply the unit of the adjunction $(f^{-1},Rf_\ast)$ to the domain of \eqref{eq:projmap} to get
		\[Rf_\ast(-)\otimes^L_{\Aa}-
		\to Rf_\ast(f^{-1}(Rf_\ast(-)\otimes^L_{\Aa}-)).\]
		Since $f^{-1}$ respects the tensor product, we get
		\[Rf_\ast(f^{-1}(Rf_\ast(-)\otimes^L_{\Aa}-))
		\simeq R f_\ast (f^{-1}Rf_\ast(-)\otimes^L_{f^{-1}\Aa}f^{-1}(-)).\]
		Finally, apply the counit of the adjunction $(f^{-1},Rf_\ast)$ to get
		\[R f_\ast (f^{-1}Rf_\ast(-)\otimes^L_{f^{-1}\Aa}f^{-1}(-))
		\to Rf_\ast (-\otimes^L_{f^{-1}\Aa}f^{-1}(-)).\]
		Compose the above homomorphisms.
	\end{construction}
	
	\begin{theorem}[Projection formula]\label{thm:projformula}
		Let $i:Y\to X$ be a closed immersion of ringed spaces, $\Aa$, $\Bb$ be $\OO_X$-algebras. Then the projection map
		$i_\ast(-)\otimes^L_{\Aa}-\to i_\ast(-\otimes^L_{i^{-1}\Aa}i^{-1}(-))$
		of functors from $D_{\mathrm{r}}(i^{-1}\Aa)\times D(\Aa\otimes_{\OO_X}\Bb^{\op})$ is an equivalence.
	\end{theorem}
	
	\begin{proof}
		Let $\Ff_\bullet$ be a q-flat complex of right $i^{-1}\Aa$-modules, and $\Mm_\bullet$ be a complex of $(\Aa,\Bb)$-bimodules. Since q-flatness is determined stalk-wisely, $i_\ast(\Ff_\bullet)$ is a q-flat complex of right $\Aa$-modules. Hence the assertion is reduced to showing that the canonical map
		\[i_\ast(\Ff_\bullet)\otimes_{\Aa}\Mm_\bullet\to i_\ast(\Ff_\bullet\otimes_{i^{-1}\Aa}i^{-1}\Mm_\bullet)\]
		is a (quasi-)isomorphism. This follows by seeing the stalks.
	\end{proof}
	
	This is used as basic techniques to study the derived direct image functor $f_+$ of twisted $\D$-modules.

	\subsection{Base change theorem}\label{sec:bc}
	
	First, we review the general construction of the base change map:
	
	\begin{construction}[Base change map]\label{cons:basechangemap}
		Consider a diagram
		\[\begin{tikzcd}
			\Cc_1\ar[r, "U'"]\ar[d, "V'"']&\Cc_2\ar[d, "V"]\\
			\Cc_3\ar[r, "U"]&\Cc_4,
		\end{tikzcd}\]
		of functors between categories. Suppose that this is lax 2-commutative, i.e., a natural transformation $\eta:V\circ U\to U\circ V'$ is given. Assume also that $U$ and $U'$ admit left adjoint functors $F$ and $F'$ respectively. Then we define a natural transformation $F\circ V\to V'\circ F'$ called the base change map as follows: Apply $V$ to the unit of $(F',U')$ to get $V\to V\circ U'\circ F'$. Apply $\eta$ to this natural transformation to get $V\to U\circ V'\circ F'$. Finally, pass to the adjuncton of $(F,U)$ to get $F\circ V\to V'\circ F'$. 
	\end{construction}
	
	\begin{example}[{Flat base change theorem, \cite[Proposition (3.9.5)]{lipman2009}}]\label{ex:flatbasechange}
		Let $x:X\to S$ and $s:S'\to S$ be morphisms of schemes. Apply Construction \ref{cons:basechangemap} to
		\[\begin{array}{cccc}
			\Cc_1=D(\OO_{X'}), & \Cc_2=D(\OO_{X}),
			& \Cc_3=D(\OO_{S'}), & \Cc_4=D(\OO_{S}),\\
			U=Rs_\ast, & V=Rx_\ast, & U'=R(s_X)_\ast, & V'=Rx'_\ast
		\end{array}\]
		to get
		\begin{equation}
			Ls^*Rx_*\to
			Rx'_* Ls^*_X\mathcal,
			\label{eq:bcmapofderivedfunctor}
		\end{equation}
		where $x':X'=X\times_S S'\to S'$ and $s_X:X'\to X$ are the canonical projection. Suppose that $x$ is concentrated (i.e., quasi-compact and quasi-separated), and that $s$ is flat. Then the natural transformation \eqref{eq:bcmapofderivedfunctor} is an equivalence on $D_{\qc}(\OO_X)$, i.e., for $\Mm^\bullet\in D_{\qc}(\OO_X)$, we have a canonical equivalence
		$s^*Rx_*\Mm^\bullet\simeq
		Rx'_* s^*_X \Mm^\bullet$.
	\end{example}
	
	\begin{example}\label{ex:noflatbcmap}
		Suppose that we are given a Cartesian diagram
		\[\begin{tikzcd}
			X'\ar[r, "s_X"]\ar[d, "x'"']&X\ar[d, "x"]\\
			S'\ar[r, "s"]&S
		\end{tikzcd}\]
		of schemes. Apply Construction \ref{cons:basechangemap} to
		\[\begin{array}{cccc}
			\Cc_1=\Mod(\OO_{X'}), & \Cc_2=\Mod(\OO_{X}),
			& \Cc_3=\Mod(\OO_{S'}), & \Cc_4=Mod(\OO_{S}),\\
			U=s_\ast, & V=x_\ast, & U'=(s_X)_\ast, & V'=x'_\ast
		\end{array}\]
		to get a canonical map
		\begin{equation}
			s^\ast x_\ast\Mm\to x'_\ast s^\ast_X\Mm
			\label{eq:bcmapofsheaf}
		\end{equation}
	\end{example}
	
	\begin{example}\label{ex:bcthmforopimm}
		Suppose that we are given a Cartesian diagram
		\[\begin{tikzcd}
			V\ar[r, "j_X"]\ar[d, "f_U"']&X\ar[d, "f"]\\
			U\ar[r, "j"]&Y
		\end{tikzcd}\]
		of topological spaces. Assume that $j$ is an open immersion. For a sheaf $\Aa$ of rings on $Y$, apply Construction \ref{cons:basechangemap} to
		\[\begin{array}{cccc}
			\Cc_1=D(f^{-1}_Uj^{-1}\Aa), & \Cc_2=D(f^{-1}\Aa),
			& \Cc_3=D(\Aa|_U), & \Cc_4=D(\Aa),\\
			U=Rj_\ast ,& V=Rf_\ast ,& U'=R(j_X)_\ast, & V'=R(f_U)_\ast
		\end{array}\]
		to get a canonical equivalence
		$Rf_\ast(-)|_U\simeq R(f_U)_\ast(-|_V)$.
		In fact, $j^{-1}$ and $j^{-1}_X$ respect q-injective complexes since they admit exact left adjoint functors $j_!$ and $(j_X)_!$ respectively (\cite[Chapter II, 6.3, Propositions 6.6, 6.9 ii)]{iversen1986}). Then the equivalence follows from $f_\ast(-)|_U=(f_U)_\ast(-|_V)$ on $\Mod(f^{-1}\Aa)$.
	\end{example}

	\subsection{Upper dimension of $Rf_\ast$}\label{sec:boundRf_ast}
	
	Let $f:X\to Y$ be a continuous map of topological spaces, and $\Rr$ be a sheaf of (possibly noncommutative) rings on $Y$. Let us recall the finiteness condition on $Rf_\ast$ in \cite[Definition (1.11.1)]{lipman2009}:
	
	\begin{definition}\label{defn:dim^+(f,R)}
		Set $\dim^+(f,\Rr)\coloneqq 
		\inf\{d\in\ZZ:~Rf_\ast(D^{\leq 0}_{\mathrm{r}}(f^{-1}\Rr))\subset D^{\leq d}_{\mathrm{r}}(\Rr)\}$.
	\end{definition}
	
	\begin{definition}\label{defn:finitedimension}
		We say $Rf_\ast:D_{\mathrm{r}}(f^{-1}\Rr)\to D_{\mathrm{r}}(\Rr)$ is bounded if $\dim^+(f,\Rr)<\infty$, i.e., there exists an integer $d$ such that $R f_\ast$ sends $D^{\leq 0}_{\mathrm{r}}(f^{-1}\Rr)$ to $D^{\leq d}_{\mathrm{r}}(\Rr)$.
	\end{definition}
	
	We will sometimes assume it in section \ref{sec:deriveddirectimage} as a working hypothesis.
	
	\begin{example}\label{ex:cloimm_dim0}
		Suppose that $f=i$ is a closed immersion. Then $R i_\ast=i_\ast$ sends
		$D^{\leq 0}_{\mathrm{r}}(i^{-1}\Rr)$ to $D^{\leq 0}_{\mathrm{r}}(\Rr)$ since 
		$i_\ast:\Mod_{\mathrm{r}}(i^{-1}\Rr)\to \Mod_{\mathrm{r}}(\Rr)$ is exact. In particular, we have $\dim^+(i,\Rr)\leq 0<\infty$.
	\end{example}
	
	\begin{example}\label{ex:findim}
		The functor $Rf_\ast:D_{\mathrm{r}}(f^{-1}\Rr)\to D_{\mathrm{r}}(\Rr)$ is bounded if $X$ is Noetherian of finite Krull dimension. In fact, a similar argument to \cite[Proposition (3.9.2)]{lipman2009} works by \cite[Chapter III, Theorem 2.7]{hartshorneag} and the remark below Proposition \ref{prop:independent}. For $d$ in \cite[Proof of Proposition (3.9.2)]{lipman2009}, take the Krull dimension of $X$.
	\end{example}
	
	Definition \ref{defn:finitedimension} is local in the following weak sense:
	
	\begin{proposition}\label{prop:localityoffindim}
		\begin{enumerate}
			\item Let $U$ be an open subset of $Y$. If $Rf_\ast:D_{\mathrm{r}}(f^{-1}\Rr)\to D_{\mathrm{r}}(\Rr)$ is bounded, so is its restriction $R(f_U)_\ast:D_{\mathrm{r}}(f^{-1}_U\Rr|_U)\to D_{\mathrm{r}}(\Rr|_U)$.
			\item Suppose that we are given a finite open covering $Y=\cup_{i=1}^n U_i$. Then $\dim^+(f,\Rr)<\infty$ if and only if so are $\dim^+(f_{U_i},\Rr|_{U_i})<\infty$ for all $i$.
		\end{enumerate}
	\end{proposition}
	
	\begin{proof}
		Suppose $\dim^+(f,\Rr)<\infty$. Take an integer $d\geq \dim^+(f,\Rr)$. We remark that this inequality is equivalent to the containment $Rf_\ast(D^{\leq 0}_{\mathrm{r}}(f^{-1}\Rr))\subset D^{\leq d}_{\mathrm{r}}(\Rr)$ by definition. Let $U$ be an open subset of $Y$. The proof of (1) is completed by showing $d\geq \dim^+(f_U,\Rr|_U)$. Let $\Ff^\bullet\in K_{\mathrm{r}}(f^{-1}_U \Rr|_U)$ be a q-injective object whose cohomology sheaves are concentrated in nonpositive degrees. Let $j$ denote the embedding $f^{-1}(U)\hookrightarrow X$. Since the functor $j_!$ of extension by zero from $f^{-1}(U)$ to $X$ is exact, $j_!\Ff^\bullet$ has nonzero cohomology only in nonpositive degrees. Hence we have
		$R^i(f_U)_\ast(\Ff^\bullet)
		\cong R^i(f_U)_\ast((j_!\Ff^\bullet)|_{f^{-1}(U)})
		\cong R^if_\ast(j_!\Ff^\bullet)|_U=0$
		for $i>d$. For the second isomorphism, see Example \ref{ex:bcthmforopimm}. This shows (1).
		
		The ``only if'' direction of (2) follows from (1). Suppose that $\dim^+(f_{U_i},\Rr|_{U_i})<\infty$ for all $i$. Then it is easy to show $\dim(f,\Rr)\leq \max\{d_1,d_2,\ldots,d_n\}<\infty$ by using Example \ref{ex:bcthmforopimm}. This completes the proof.
	\end{proof}
	
	Let us also introduce a local analog of Definition \ref{defn:finitedimension}:
	
	\begin{definition}\label{defn:locallyoffindim}
		\begin{enumerate}
			\item We say $Rf_\ast:D_{\mathrm{r}}(f^{-1}\Rr)\to D_{\mathrm{r}}(\Rr)$ is locally bounded if there exists an open covering $Y=\cup_\lambda U_\lambda$ such that $\dim^+(f_{U_\lambda},\Rr|_{U_\lambda})<\infty$ for every index $\lambda$.
			\item We say $Rf_\ast:D(f^{-1}\Rr)\to D(\Rr)$ is locally bounded if $Rf_\ast:D_{\mathrm{r}}(f^{-1}\Rr^{\op})\to D_{\mathrm{r}}(\Rr^{\op})$ is so.
		\end{enumerate}
		
	\end{definition}
	
	We remark that this is local in the following stronger sense:
	
	\begin{proposition}\label{prop:localityoflocalfindim}
		\begin{enumerate}
			\item Let $U$ be an open subset of $Y$. If $Rf_\ast:D_{\mathrm{r}}(f^{-1}\Rr)\to D_{\mathrm{r}}(\Rr)$ is locally bounded, so is $R(f_U)_\ast:D_{\mathrm{r}}(f^{-1}_U\Rr|_U)\to D_{\mathrm{r}}(\Rr|_U)$.
			\item Suppose that we are given an open covering $Y=\cup_{\mu} V_\mu$. Then $Rf_\ast:D_{\mathrm{r}}(f^{-1}\Rr)\to D_{\mathrm{r}}(\Rr)$ is locally bounded if and only if so is $R(f_{V_\mu})_\ast:D_{\mathrm{r}}(f^{-1}_{V_\mu}\Rr|_{V_\mu})\to D_{\mathrm{r}}(\Rr|_{V_\mu})$ for every $\mu$.
			\item Suppose that $Y$ is quasi-compact. Then $Rf_\ast:D_{\mathrm{r}}(f^{-1}\Rr)\to D_{\mathrm{r}}(\Rr)$ is locally bounded if and only if it is bounded.
			\item Suppose that we are given a covering $Y=\cup_\mu V_\mu$ consisting of quasi-compact open subsets. Then $Rf_\ast:D_{\mathrm{r}}(f^{-1}\Rr)\to D_{\mathrm{r}}(\Rr)$ is locally bounded if and only if
			\[R(f_{V_\mu})_\ast:D_{\mathrm{r}}(f^{-1}_{V_\mu}\Rr|_{V_\mu})\to D_{\mathrm{r}}(\Rr|_{V_\mu})\]
			is bounded for every $\mu$.
		\end{enumerate}
	\end{proposition}
	
	A typical case of (4) in our paper is when $Y$ is a scheme and $V_\mu$ are the affine open subsets.
	
	\begin{proof}
		The first two assertions follow from the basic observation that $(f_U)_V=f_V$ for open subsets
		$V\subset U\subset Y$ (also use Proposition \ref{prop:localityoffindim} (1) for (1)). Part (3) is immediate from Proposition \ref{prop:localityoffindim} (2). Part (4) is combination of (2) and (3).
	\end{proof}
	
	\begin{corollary}\label{cor:loc_noe_fin_dim}
		Suppose that $X$ is locally Noetherian of finite Krull dimension, and that $f$ is quasi-compact. Then $Rf_\ast:D_{\mathrm{r}}(f^{-1}\Rr)\to D_{\mathrm{r}}(\Rr)$ is locally bounded.
	\end{corollary}
	
	As an easy application, notice that the proofs of \cite[Lemma (3.9.3.1) and Corollary (3.9.3.2)]{lipman2009} readily assert:
	
	\begin{lemma}\label{lem:partialomegaaccessibility}
		Let $f:X\to Y$ be a concentrated morphism of schemes, $\Aa$ be a sheaf of rings on $Y$, and $n\in\ZZ$. Suppose that the functor
		$Rf_\ast:D(f^{-1}\Aa)\to D(\Aa)$
		is locally bounded. Then for a filtered diagram $C^\bullet_\lambda$ of complexes of $f^{-1}\Aa$-modules, we have a canonical equivalence
		$\varinjlim_\alpha Rf_\ast( C^\bullet_\lambda)\simeq Rf_\ast(\varinjlim_\alpha C^\bullet_\lambda)$.
	\end{lemma}
	
	This will be used in section \ref{sec:deff+}.

	
	\bibliographystyle{plain}

\end{document}